\documentclass[11pt]{article}

\usepackage{graphicx}
\usepackage{caption}
\usepackage{subcaption}

\usepackage[margin = 1in]{geometry}
\usepackage{amsmath,amssymb,amsfonts,amsthm}
\usepackage{color}
\usepackage{hyperref} % For creating hyperlinks in cross references

\usepackage{tikz}
\usetikzlibrary{arrows}
\usetikzlibrary{calc}
\usetikzlibrary{matrix}

\def \p{\partial}

\def \div{\nabla \cdot}
\def \grad{\nabla}

\def \bM{\mathbb{M}}

\def \cTh{\mathcal{T}_h}
\def \cTc{\mathcal{T}_{c,h}}
\def \cTco{\hat{\mathcal{T}}_{c,h}}
\newcommand{\vect}[1]{\boldsymbol{#1}}
\newcommand{\trace}{\text{trace}}

\newcommand{\bey}{\begin{eqnarray}}
\newcommand{\eey}{\end{eqnarray}}

\newcommand{\beq}{\begin{equation}}
\newcommand{\eeq}{\end{equation}}
\theoremstyle{plain}% default

\theoremstyle{definition}

\theoremstyle{remark}
\newtheorem{exam}{\hspace{6mm}Example}[section]

\newenvironment{keywords}%
   {\begin{trivlist}\item[]{\bfseries\sffamily Key words:}~}
   {\end{trivlist}}
\newenvironment{AMS}%
   {\begin{trivlist}\item[]{\bfseries\sffamily AMS subject classifications:}~}
   {\end{trivlist}}

\allowdisplaybreaks

\title{%
A Study on Moving Mesh Finite Element Solution of\\ the Porous Medium Equation
}

\author{%
   Cuong~Ngo%
   \thanks{%
      Department of~Mathematics, University of~Kansas, Lawrence, KS~66045, U.S.A.
      (\href{mailto:cngo@ku.edu}{\nolinkurl{cngo@ku.edu}}).%
   }%
   \and
   Weizhang~Huang%
   \thanks{%
      Department of~Mathematics, University of~Kansas, Lawrence, KS~66045, U.S.A.
      (\href{mailto:whuang@ku.edu}{\nolinkurl{whuang@ku.edu}}).%
   }%
 }

\date{}

\begin{document}
\maketitle

%--- abstract --------------------------------------------------------
\begin{abstract}
An adaptive moving mesh finite element method is studied for the numerical solution of the porous medium
equation with and without variable exponents and absorption.
The method is based on the moving mesh partial differential equation approach
and employs its newly developed implementation.
%{\color{red}
The implementation has several improvements over the traditional one, including
its explicit, compact form of the mesh velocities, ease to program,
and less likelihood of producing singular meshes.
Three types of metric tensor that correspond to uniform and arclength-based and Hessian-based adaptive meshes
are considered.
%}
The method shows first-order convergence
for uniform and arclength-based adaptive meshes, and second-order convergence for Hessian-based
adaptive meshes. It is also shown that the method can be used for situations with complex free boundaries,
emerging and splitting of free boundaries, and the porous medium equation with variable exponents and
absorption. Two-dimensional numerical results are presented.
\end{abstract}

%--- keywords --------------------------------------------------------

%--- ams --------------------------------------------------------

%%%%%%%%%%%%%
%% Section: Introduction 
\section{Introduction}
\label{SEC:intro}

We consider the numerical solution of the initial-boundary value problem (IBVP)
of the porous medium equation (PME) in two dimensions, 
\begin{equation}
    \begin{cases}
        u_t = \div(|u|^m \grad u) , \quad & \text{in} \quad \Omega \times (t_0,T]  \\
        u(\vect{x},t_0) = u_0(\vect{x}) , \quad & \text{on} \quad \Omega \\
        u(\vect{x},t) = 0 , \quad & \text{on} \quad \p \Omega  \times (t_0,T]
    \end{cases}
  \label{PME-1}
\end{equation}
where $\Omega$ is a bounded polygonal domain, $u_0(\vect{x})$ is a given function, and $m\ge 1$ is a physical
parameter. PME is a nontrivial generalization of the heat equation. It is found in many areas of
the physical sciences, including gas flow in porous medium, incompressible fluid dynamics,
nonlinear heat transfer, and image processing; e.g., see \cite{Vazquez2007} and references therein.
In the case of gas flow in porous medium, $u$ represents the density of the gas,
$u^m$ the pressure, $u \grad (u^m)$ the flux, $\grad (u^m)$ the velocity, and $m$ is the isentropic coefficient.
%{\color{red}
In the case of radiation diffusion in plasmas, $m$ stands for the power of temperature appearing in the nonlinear diffusion
coefficient and can take values up to 5.5 \cite[Page 23]{Vazquez2007}.
%}
PME itself represents a nonlinear diffusion process.
One of its many interesting features is its degeneracy which induces a property of the finite propagation:
if $u_0(\vect{x})$ has a compact support, the solution
will have a compact support for any time $t > t_0$.
This in effect creates a free boundary that stays between the regions where $u$ is nonzero and
where $u$ vanishes and propagates at a finite speed for all time.
Contrary to the heat equation which smooths out its initial solutions, PME solution
can become nonsmooth even in cases where it has a smooth initial solution. Moreover,
for a certain type of initial solutions, the solution can exhibit a waiting-time phenomenon
for which the free boundary will not move until a finite amount of time has elapsed.
PME has been studied extensively in theory and there is a vast literature,
including the earlier work by Ole{\u\i}nik et al. \cite{Oleinik1958},
Kala{\v{s}}nikov \cite{Kalashnikov1967}, Aronson \cite{Aronson1969},
the more recent work by Shmarev \cite{Shmarev2005,Shmarev2003},
and the monograph by V{\'a}zquez \cite{Vazquez2007} and references therein.

The numerical solution of PME has also received considerable attention from researchers.
Particularly, error estimates have been obtained for various finite element approximations.
For example, using a regularization approach
%{\color{red}
(to avoid degenerate or negative diffusion, for instance, with the diffusion
coefficient $|u|^m$ being replaced by $\max (|u|, \epsilon/2)^m$ for some regularization parameter $\epsilon > 0$) 
and taking
%}
$\epsilon = \mathcal{O}(h^{\frac{2m+4}{m^2+4m+2}})$, Rose \cite{Rose-1983} shows that
the error for a P1 finite element (for space) -- backward Euler (for time)
approximation of PME is bounded by
\begin{equation}
\left (\sum_{n} \Delta t \| u_h^n - u \|_{L^{m+2}(\Omega)}^{m+2} \right )^{\frac{1}{m+2}}
\le C \left (\Delta t^{\frac{1}{m+1}} + \left (\ln \left (\frac{1}{h}\right )\right )^{\frac{1}{(m+1)(m+2)}} h^{\frac{2}{m+1}}\right ),
\label{rose-1}
\end{equation}
where $h$ is the maximum element diameter and $u_h^n$ is the numerical approximation of $u$ at $t= t_n$.
Nochetto and Verdi \cite{Nochetto-1988} consider a class of degenerate PDEs which includes PME
as a special example and improve the result of \cite{Rose-1983}. They show that the error
for a P1 finite element -- 1st-order semi-implicit approximation of PME is bounded by
\begin{align}
&\| u_h - u \|_{L^{\infty}(0,T; H^{-1}(\Omega))}
+ \| (u_{h})^{m+1} - u^{m+1} \|_{L^2(0,T; L^2(\Omega))}
+ \left \| \int_0^t ((u_{h})^{m+1} - u^{m+1})\right \|_{L^{\infty}(0,T; H^{1}(\Omega))}
\notag
\\
& \qquad \qquad 
\le C \left (\frac{h^2}{\epsilon} + \frac{h^4}{\epsilon^2 \Delta t} + \Delta t\right )^{\frac 1 2}
\label{nochetto-1}
\\
& \qquad \qquad = \mathcal{O} (h^{\frac{m+2}{2(m+1)}}),\quad \text{if}\quad 
\Delta t = \mathcal{O}(h^{\frac{m+2}{m+1}}),\quad \epsilon = \mathcal{O}(h^{\frac{m}{m+1}}).
\notag
\end{align}
For the P1 finite element -- backward Euler approximation, error bounds in various norm are obtained,
for instance,
\begin{equation}
\| u_h - u \|_{L^\infty(0,T; H^{-1}(\Omega))} \le C \left (\Delta t + \left (\ln \left (\frac{1}{h}\right )\right )^{\frac{2m+3}{2m+2}} h\right )
\label{rulla-1}
\end{equation}
by Rulla and Walkington \cite{Rulla1996}, 
\begin{equation}
\| u_h - u \|_{L^2(0,T; L^2(\Omega))} \le C h^{\frac{m^2+6m+8}{6m^2+14m+8}}
\quad \text{when}\quad \Delta t = \mathcal{O}(h^{\frac{5 m+4}{2m}})
\label{ebmeyer-1}
\end{equation}
by Ebmeyer \cite{Ebmeyer-1998}, and 
\begin{equation}
\| u_h - u \|_{L^{m+2}(0,T; L^{m+2}(\Omega))} \le
C \left ( \Delta t^{\frac 1 2} + h + h^{\frac{1}{m+1} \left (\frac{d m}{2m+4}+1\right )} \right )^{\frac{1}{m+2}}
\label{wei-1}
\end{equation}
by Wei and Lefton \cite{Wei1999},  where $d$ is the space dimension.
It is remarked that these estimates are obtained for quasi-uniform meshes. The convergence rate is first order at best
and decreases with $m$. Some of these estimates are shown to be optimal
in the corresponding norm in lieu of the known regularity of the solution of PME.
Moreover, it is worth mentioning that Ebmeyer and Liu \cite{Ebmeyer-2008} obtain error estimates
in quasi-norm and Emmrich and {\v{S}}i{\v{s}}ka \cite{Emmrich2012} prove that
a Galerkin finite element -- backward Euler approximation converges to the weak solution of PME
using the theory of monotone operators.
More recently,  Duque et al. \cite{Duque2013} establish $L^{1+\max (\gamma/2)}$ error bounds for
the approximation of a general order continuous Galerkin in space and a general order
discontinuous Galerkin in time for PME with a variable exponent $m = \gamma(\vect{x})$.
Zhang and Wu \cite{Zhang-2009} consider the numerical simulation of the one-dimensional PME
on a uniform mesh using a high-order local discontinuous
Galerkin finite element method. The method can effectively eliminate unwanted, nonphysical
oscillations in the computed solution near the free boundary and lead to a high-order
convergence rate within the solution support and away from the free boundary.

The low regularity (steep gradient and corner shape near the moving free boundary)
and evolving nature of the solution makes adaptive moving mesh methods
an attractive tool to improve accuracy in the numerical solution of PME.
A number of works exist in this direction. For example,
Budd et al. \cite{BCHR99} investigate the numerical simulation of
self-similarity solutions of one-dimensional PME using the Moving Mesh PDE (MMPDE)
moving mesh method \cite{HRR94a,HR11} and a specially designed monitor function
to preserve the scaling invariance of PME.
In a series of papers \cite{BHJ05,BHJ05a,BHJJ06} (also see the review paper \cite{Baines-2011}),
Baines and his co-workers study the numerical solution of PME in one and two dimensions using
a moving mesh finite element method that is based upon conserving a local proportion of the total mass
that is present in the projected initial data. Numerical results show that their method gives a second-order
convergence for $m=1$ but only a first-order convergence for $m=3$ when a uniform initial mesh
is used \cite{BHJ05a,Baines-2011}. For $m=3$, the second-order convergence can be recovered in one dimension
if an optimal initial mesh is used. Unfortunately, such an optimal mesh 
is significantly more expensive to compute in two dimensions than in one dimension.
Recently, Duque et al. \cite{Duque2014, Duque2015} present a moving mesh finite element method
based on an MMPDE for PME with variable exponents and with/without absorption.
The method shows a first-order convergence when tested for the Barenblatt-Pattle solution of PME.

The objective of this paper is to study an adaptive moving mesh finite element method for the numerical
solution of PME. The method is also based on an MMPDE but significantly different from
the method of \cite{Duque2014, Duque2015}.
The MMPDE we use is formulated by minimizing an energy (cf. (\ref{Ih}))
based on the equidistribution and alignment conditions and the mesh adaptation is controlled
through a matrix-valued function (i.e., a metric tensor) instead of a scalar function.
The advantage of using a metric tensor is that it provides information not only to control
the size of mesh elements but also their shape and orientation.
Generally speaking, a so generated mesh has better alignment with the geometry of the physical solution
than that with a scalar mesh adaptation function.
Moreover, a newly developed compact formulation of the method (cf. (\ref{mmpde-2}) and (\ref{mmpde-3}))
makes its implementation much easier and more efficient.
Mesh adaptation based on the gradient and Hessian of the solution will be considered.
%{\color{red}
The arclength metric tensor (a gradient-based metric tensor) has been widely used
in the context of moving mesh methods because it results in more stable mesh movement than
a Hessian-based metric tensor and works well for many problems. On the other hand,
there is no theoretical guarantee that the arclength metric tensor will lead to the optimal convergence order
for piecewise linear interpolation or finite element approximation since the error is determined
by the Hessian of the solution in these situations.
As a matter of fact, there are problems, although not very common, for which we have to use
a Hessian-based metric tensor in order to achieve the optimal convergence order for piecewise linear interpolation
or finite element approximation. Interestingly, PME is one of those problems. We shall show that the linear finite element approximation
of the Barenblatt-Pattle solution of PME shows a first-order convergence for arclength-based adaptive meshes
and a second-order convergence for Hessian-based adaptive meshes.
%}
%For Hessian-based adaptive meshes, the method shows a second-order convergence
%when applied to the Barenblatt-Pattle solution of PME.

Another feature of the method that is different from those of
Duque et al. \cite{Duque2014, Duque2015} and Baines et al.  \cite{BHJ05,BHJ05a,Baines-2011,BHJJ06}
is that PME is solved on a large domain that contains the free boundary
for the whole time period under consideration. In this way, there is no need to explicitly trace
the movement of the free boundary and thus the method can readily deal with more complicated structures
in the solution and in the differential equation. Numerical examples with simple free boundary
(such as the Barenblatt-Pattle solution) and more complex and even merging free boundaries
will be presented as well as those for PME with variable exponents and absorption.

%{\color{red}
In addition to the above mentioned MMPDE method, a number of other moving mesh methods
have been developed in the past; e.g., see
Hirt et al. \cite{HAC74} (ALE -- Arbitrary Lagrangian-Eulerian)
Miller and Miller \cite{MM81} (MFE -- moving finite element),
Liao and Anderson \cite{LA92} (deformation map),
Li et al. \cite{LTZ01} (mesh rezoning),
Cao et al. \cite{CHR02} (geometric conservation law),
Baines et al. \cite{BHJ05a} (conservation of fraction mass),
and
Budd and Williams \cite{BW06} (parabolic {M}onge-{A}mp\`ere equation).
The interested reader is also referred to the books/review articles
\cite{Bai94a,Baines-2011,Shashkov-2015,BHR09,HR11,Tan05} and references therein.
%}

The outline of the paper is as follows. Some properties of PME that are relevant to the numerical
simulation are described in Section~\ref{SEC:PME-theory}. Section~\ref{SEC:mmfem} is devoted
to the description of the moving mesh finite element method, including the linear finite element
discretion of PME on a moving mesh and the generation of an adaptive moving mesh using
the MMPDE approach. Numerical examples are presented in Sections~\ref{SEC:PME-numerics}
and \ref{SEC:PME-numerics-2} for PME and PME with variable exponents and absorption, respectively.
Finally, Section~\ref{SEC:conclusion} contains the conclusions and further comments. 

\vspace{10pt}
%%%%%%%%%%%%%%%%%%%%%%%
%% Section: Background
\section{Properties of the porous medium equation}
\label{SEC:PME-theory}

Before we describe the finite element approximation of IBVP (\ref{PME-1}), it is instructive to recall some
of its properties that are relevant to numerical simulation. First of all, 
it is known (e.g., see V{\'a}zquez \cite{Vazquez2007}) that when $u_0^{m+2} \in L^{1}(\Omega)$, 
IBVP \eqref{PME-1} has a weak solution $u$ satisfying
$u(\cdot, t)^{m+2} \in L^{1}(\Omega)$ for any $t  \in (t_0, T]$ and $u^{m+1} \in L^{2}(0,T; H^1_0(\Omega))$.
Moreover, a nonnegative weak solution exists
if $u_0^{m+2} \in L^{1}(\Omega)$, $u_0 \in L^1(\Omega)$, and $u_0 \ge 0$. The uniqueness
of the weak solution is guaranteed if it is further assumed that $u \in L^2(\Omega \times (0,T))$.

PME is degenerate whenever $u = 0$. Due to this degeneracy, PDE has the property of the finite speed
of propagation: if compact-supported initially,  its solution remains compact-supported at any finite
time, with the support monotonically expanding as time evolves. The boundary of the support forms a moving interface
$\Gamma(t)$ which is commonly referred to as a free boundary. The velocity of the free boundary is
given by Darcy's law (e.g., see Shmarev \cite{Shmarev2005}), i.e.,
\begin{equation}
    \Gamma'(t) = - \lim_{\vect{x} \to \Gamma(t)^-} \grad \left( \frac{u^m(\vect{x},t)}{m} \right),
    \label{Darcy-law}
\end{equation}
where the limit is taken from the interior of the support.
In addition, PME exhibits a waiting-time phenomenon: for a certain type of initial solutions
the free boundary does not move until a finite amount of time has elapsed.
Loosely speaking, from (\ref{Darcy-law}) we may expect to see this phenomenon for initial solutions having
vanishing $\grad (u^m)$ at the initial free boundary.

A few classes of special solutions to IBVP \eqref{PME-1} have been known, among which is the Barenblatt-Pattle
solution, viz.,
    \begin{equation}
        u(r,t) = 
        \begin{cases}
            \frac{1}{\lambda^d(t)} 
            \left( 
            1 - \left( \frac{r}{r_0 \lambda(t)} \right)^2 
            \right)^{\frac{1}{m}} , \quad & \text{for} \quad |r| \leq r_0 \lambda(t)  \\
                                0 , \quad & \text{for} \quad |r| > r_0 \lambda(t)  
        \end{cases}
        \label{BP-soln}
    \end{equation}
    where 
    \[ 
        r = | \vect{x} | , \quad 
        \lambda(t) = \left(\frac{t}{t_0} \right)^{\frac{1}{2+d m}} ,
        \quad 
        t_0 = \frac{r_0^2 m}{2 (2 + d m)} ,
    \]
and $r_0>0$ is a given parameter. 
It is radially symmetric, self-similar, and compact-supported for any finite time.
Moreover, $u^m$ is Lipschitz continuous in $\vect{x}$ and $t$
and $\nabla (u^m)$ is bounded in the support of $u(\cdot, t)$, $\text{supp}(u(\cdot, t))$.
Furthermore, the solution is H\"{o}lder continuous.
The slope of the solution at the free boundary is finite for $m = 1$
and becomes infinite when $m>1$, which causes challenges for the numerical solution of PME.
These regularity properties also hold for general compactly supported solutions of (\ref{PME-1}),
and their free boundaries can be shown to be at least Lipschitz continuous in both space and time;
e.g., see \cite{Aronson1969,Caffarelli1980,Daskalopoulos1998a,Shmarev2005}.

 %%%%%%%%%%%%%%%%%%%%%%%
%% Section 3: MMFEM 
\section{The moving mesh finite element method}
\label{SEC:mmfem}

In this section we describe the adaptive moving mesh finite element approximation of IBVP~\eqref{PME-1}.
To begin with, we note that there are roughly two approaches for solving the IBVP. 
The first, as used in \cite{Rose-1983,Zhang-2009}, is to solve PME in a large domain
containing the free boundary for the whole time period of the simulation.
With this approach, there is no need to explicitly treat the free boundary, which makes
the approach more amenable to problems with complex solution supports.
The main disadvantage of this approach is that the solution has a corner shape between
the regions of zero and nonzero solution values and thus its regularity on the whole domain
is at most $H^1$. An $H^1$ regularity often means at best a first-order convergence
in the numerical solution as the mesh is refined. 

The second approach is to solve the problem only in the region of compact support; e.g., see
\cite{BHJ05a,Duque2014,Duque2015}.
One of the advantages of this approach is that a smaller spatial domain is used
and thus fewer mesh points can be used to achieve the same computational accuracy.
Moreover, the regularity of the solution is better on the support than on a larger domain
since it does not have a corner shape. As a result,  the numerical solution
can have a higher convergence order than that with the first approach.
The main disadvantage is that the boundary movement has to be treated explicitly
using Darcy's law (\ref{Darcy-law}). 
%It is also worth mentioning that 
%Since the domain is changing, the interior mesh points
%also need to be moved accordingly.

We use the first approach in this work. We choose this approach due to its advantage of no need to
explicitly treat the free boundary and its potential to deal with problems
having complex solution supports. To better resolve the corner shape in the solution and 
improve the computational accuracy, we employ an MMPDE-based moving mesh
method \cite{HR11} to adaptively and dynamically concentrate the mesh points around the free boundary.
Linear finite elements and the fifth-order Radau IIA method (e.g., see  Hairer and Wanner \cite{HW96})
are used for the spatial and temporal discretization, respectively. As we will see in Section~\ref{SEC:PME-numerics},
a second-order convergence of the finite element approximation in space can be achieved when
a properly adapted mesh is used.

% section 3.1
\subsection{Finite element discretization}
\label{SEC:fem}

We now describe the finite element discretization. Denote the time instants by
\begin{equation}
    t_0 = 0 < t_1 < \ldots < t_{n_f} \equiv T .
    \label{time-grid}
\end{equation}
For the moment, we assume that the simplicial meshes $\cTh^n$, $n = 0,\ldots,n_f$
for the physical domain $\Omega$ at these time instants are known and
have the same connectivity and the same numbers of vertices and elements.
(Their generation will be discussed in the next subsection.)
Denote the coordinates of the vertices of $\cTh^n$
by $\vect{x}_j^n$, $j = 1, ..., N_v$, where $N_v$ is the number of all vertices.
The mesh $\cTh(t)$ between any two time instants $t_n$ and $t_{n+1}$
is defined through linear interpolation, i.e., 
\begin{align}
&    \vect{x}_j(t) =
    \frac{t-t_n}{t_{n+1}-t_n} \vect{x}_j^{n+1} +
    \frac{t_{n+1}-t}{t_{n+1}-t_n} \vect{x}_j^{n} , \quad \forall j = 1, ..., N_v 
\\
&\dot{\vect{x}}_j (t) = \frac{\vect{x}_j^{n+1}-\vect{x}_j^{n}}{t_{n+1}-t_n},\quad
j = 1, ..., N_v.
\end{align} 
Denote by $\phi_j(\vect{x}, t)$ the linear basis function associated with vertex $\vect{x}_j(t)$. 
For convenience, we assume that the vertices are arranged such that the first $N_{vi}$ vertices
are the interior vertices.
Let
\[
V_h(t) = \text{span}\{\phi_1(\cdot, t), ..., \phi_{N_{vi}}(\cdot, t)\} .
\]
Then, the linear finite element approximation to the solution of IBVP \eqref{PME-1} is defined
as $u_h(\cdot, t) \in V_h(t)$, $ t \in (t_0, T]$ such that
\begin{equation}
\begin{cases} 
	& \int_{\Omega} \frac{\p u_h}{\p t} v~d\vect{x} = - \int_{\Omega} |u_h|^m \grad u_h
	\cdot \grad v~d\vect{x} , 
	    \quad \forall v \in V_h(t),\quad t_0 < t \leq T  \\
	& \int_{\Omega} (u_h(\vect{x},0) - u^0(\vect{x})) v~d\vect{x}  = 0, \quad \forall v \in V_h(t) .
    \end{cases}
    \label{fem0}
\end{equation}
The above equation can be cast in matrix form. Indeed, expressing 
\[
    u_h(\vect{x},t) = \sum_{j=1}^{N_{v_i}} u_j(t) \phi_j(\vect{x},t)
\]
and differentiating it with respect to $t$, we have
\[
\frac{\p u_h}{\p t}  = \sum_{j=1}^{N_{v_i}} \frac{\p u_j}{\p t} \phi_j(\vect{x},t) + 
	\sum_{j=1}^{N_{v_i}} u_j(t) \frac{\p \phi_j}{\p t} . 
\]
It can be shown (e.g. see Jimack and Wathen \cite[Lemma 2.3]{Jimack-1991}) that
\[
\frac{\p \phi_j}{\p t}  =  - \nabla \phi_j \cdot \dot{\vect{X}},\quad a.e. \text{ in } \Omega
\]
where the mesh velocity $\dot{\vect{X}}$ is defined as
\[
\dot{\vect{X}} (\vect{x},t) = \sum_{j=1}^{N_v} \dot{\vect{x}}_j(t) \phi_j (\vect{x}, t) .
\]
Then, we get
\[
\frac{\p u_h}{\p t}  = \sum_{j=1}^{N_{v_i}} \frac{\p u_j}{\p t} \phi_j(\vect{x},t) - \nabla u_h\cdot  \dot{\vect{X}}  .
\]
From this and taking $v = \phi_i$ ($i = 1, ..., N_{vi}$) in (\ref{fem0}) successively, we obtain
\[
        \sum_{j=1}^{N_{vi}} \left(\int_{\Omega} 
	 \phi_j \phi_i~d\vect{x}\right) \frac{d u_j}{dt} 
        = \int_{\Omega} \grad u_h \cdot \left( \dot{\vect{X}} \phi_i
        - u_h^m \grad \phi_i \right) ~d\vect{x}, \quad i = 1, ..., N_{vi}, \quad t_0 < t \leq T
\]
which can be cast in the matrix form as
\begin{equation}
B(\vect{X}) \dot{\vect{U}} = F(\vect{U}, \vect{X}, \dot{\vect{X}}),
\label{fem-2}
\end{equation}
where $B$ is the mass matrix and $\vect{X}$ and $\vect{U}$ are the vectors representing
the mesh and solution, respectively.
This ODE system is integrated from $t_n$ to $t_{n+1}= t_n + \Delta t_n$ using the fifth-order Radau IIA method,
%{\color{red}
with $\Delta t_n$ being determined
by a standard time step size selection procedure (e.g., see Hairer et al. \cite[Section II.4]{HNW93}) and
%}
using a two-step error estimator of Gonzalez-Pinto et al.~\cite{Montijano2004}.
The relative and absolute tolerances $rtol = 10^{-6}$ and $atol = 10^{-8}$ are taken in the computation.

The whole computation alternates between the integration of PME and the generation of the mesh.
Starting with the current mesh $\cTh^n$ and a solution $u_h^n(\vect{x})\approx u(\vect{x}, t_n)$
defined thereon, a new mesh $\cTh^{n+1}$ is generated
using the moving mesh strategy to be described in the next subsection.
Then, the discrete PME (\ref{fem-2}) is integrated from $t_n$ to $t_{n+1}$
%{\color{red}
(as described above)
%}
to obtain the solution
approximation $u_h^{n+1}(\vect{x})$.
%with $\Delta t_n$ being determined
%by a standard time step size selection procedure for (\ref{fem-2}).

% section 3.2
\subsection{An MMPDE-based moving mesh strategy}
\label{SEC:mmpde} 

We now describe the generation of $\cTh^n$, $ n = 1, ..., n_f$. 
We assume that the mesh $\cTh^n$ and a computed solution
$u_h^n(\vect{x})$ are known at $t = t_n$. We also assume that a reference computational mesh
$\cTco = \{ \hat{\vect{\xi}}_j, j = 1, ..., N_v\}$ having the same
connectivity and the same numbers of vertices and elements as $\cTh^n$ has been chosen.
In our computation, it is taken as a uniform mesh
%{\color{red}
(in the Euclidean metric)
%}
defined on $\Omega$.
%{\color{red}
$\cTco$ stays fixed for the whole computation.
%}

\begin{figure}[t]
\centering
{\footnotesize
\begin{tikzpicture}[scale = 0.8]
\draw [thick,->] (0,0) -- (7,0);
\draw [right] (7,0) node {$\xi$};
\draw [thick,->] (0,0) -- (0,7);
\draw [above] (0,7) node {$x$};
\draw[line width=0.25ex] (0,0) -- (1.4, 1.2) -- (4.4,5) -- (6,6);
\draw (1.4,0) -- (1.4, 1.2) -- (0, 1.2);
\draw (4.4,0) -- (4.4, 5) -- (0, 5);
\draw (6,0) -- (6, 6) -- (0, 6);
\draw[fill] (0,0) circle (.5ex);
\draw[fill] (1.4,1.2) circle (.5ex);
\draw[fill] (4.4,5) circle (.5ex);
\draw[fill] (6,6) circle (.5ex);
\draw[fill] (1.4,0) circle (.5ex);
\draw[fill] (4.4,0) circle (.5ex);
\draw[fill] (6,0) circle (.5ex);
\draw[fill] (0,1.2) circle (.5ex);
\draw[fill] (0,5) circle (.5ex);
\draw[fill] (0,6) circle (.5ex);
\draw [below] (0,0) node {$\xi_1^{n+1}(\hat{\xi}_1)$};
\draw [below] (1.4,0) node {$\xi_2^{n+1}$};
\draw [below] (4.7,0) node {$\xi_3^{n+1}$};
\draw [below] (6,0) node {$\xi_4^{n+1}(\hat{\xi}_4)$};
\draw [left] (0,0.1) node {$x_1^{n}(x_1^{n+1})$};
\draw [left] (0,6) node {$x_4^{n}(x_4^{n+1})$};
\draw [left] (0,1.2) node {$x_2^{n}$};
\draw [left] (0,5) node {$x_3^{n}$};
\draw [left] (0,4.4933) node {$x_3^{n+1}$};
\draw [left] (0,1.96) node {$x_2^{n+1}$};

\draw[fill=red] (2,0) circle (.5ex);
\draw[fill=red] (4,0) circle (.5ex);
\draw [below] (2.1,0) node {$\hat{\xi}_2$};
\draw [below] (4,0) node {$\hat{\xi}_3$};
\draw[fill=red] (2,1.96) circle (.5ex);
\draw[fill=red] (0,1.96) circle (.5ex);
\draw[fill=red] (4,4.4933) circle (.5ex);
\draw[fill=red] (0,4.4933) circle (.5ex);
\draw [right] (6,6) node {$x=\Phi_h(\xi)$};

\draw[red] (2,0) -- (2, 1.96) -- (0, 1.96);
\draw[red,->] (2,0) -- (2, 1);
\draw[red,->] (2, 1.96) -- (1, 1.96);
\draw[red] (4,0) -- (4, 4.4933) -- (0, 4.4933);
\draw[red,->] (4,0) -- (4, 2.25);
\draw[red,->] (4, 4.4933) -- (2, 4.4933);
\end{tikzpicture}
}
\caption{A sketch of the relations among the meshes $\cTco = \{ \hat{\vect{\xi}}_j\}$,
$\cTc^{n+1} = \{\vect{\xi}_j^{n+1} \}$, $\cTh^n=\{\vect{x}_j^{n}\}$,
and $\cTh^{n+1}=\{\vect{x}_j^{n+1}\}$. The function $\vect{x} = \Phi_h(\vect{\xi})$ is determined
as the correspondence between $\cTc^{n+1}$ and $\cTh^{n}$; and $\cTh^{n+1}$ is computed
as $\cTh^{n+1} = \Phi_h(\cTco)$ using linear interpolation.}
\label{fig:mesh-relation}
\end{figure}
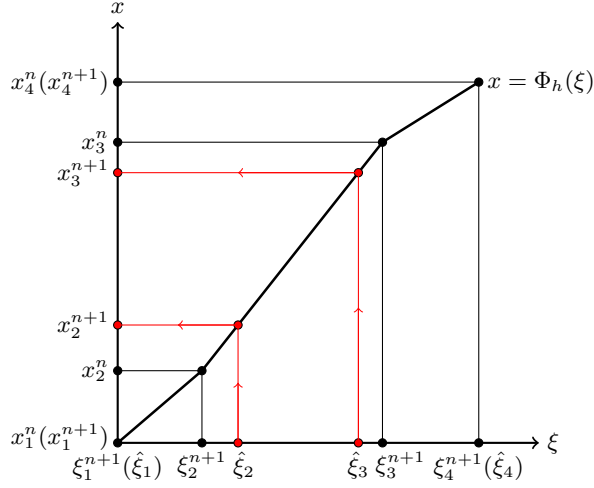

%{\color{red}
The generation of $\cTh^{n+1}$ is through the computational mesh
$\cTc = \{ \vect{\xi}_j, j = 1, ..., N_v\}$ which serves as an intermediate variable. 
(A sketch of the relations among the meshes $\cTco$,
$\cTc^{n+1}$, $\cTh^n$, and $\cTh^{n+1}$ is shown in Fig.~\ref{fig:mesh-relation}.)
First, an MMPDE-based mesh equation
(to be described below) for the velocities of the computational vertices is employed.
It takes the form (cf. (\ref{mmpde-2}))
\[
\begin{cases}
\frac{d \vect{\xi}_j}{d t} = \vect{v}_j(\bM, \cTh^n; \vect{\xi}_1, ..., \vect{\xi}_{N_v}),&\quad j = 1, ..., N_v, \quad t \in (t_n, t_{n+1}]
\\
\vect{\xi}_j(t_n) = \hat{\vect{\xi}}_j,& \quad j = 1, ..., N_v 
\end{cases}
\]
where $\vect{v}_j$ denotes the mesh velocity for the $j$-th node which depends on $\cTh^n$, the metric tensor
$\bM$ defined thereon, and $\cTc$. Here, the initial mesh is taken to be the reference computational mesh $\cTco$.
The system is integrated and the mesh $\cTc^{n+1} \approx \cTc(t_{n+1})$ is obtained. During the integration, both $\bM$
and $\cTh^n$ are kept fixed. Notice that $\cTc^{n+1}$ and $\cTh^n$
form a correspondence relation, say, $\vect{x}_j^n = \Phi_h( \vect{\xi}_j^{n+1}),\, j = 1, ..., N_v$ or $\vect{x} = \Phi_h( \vect{\xi} )$.
Then, the vertices of the new physical mesh $\cTh^{n+1}$ are defined as
$\vect{x}_j^{n+1} = \Phi_h(\hat{\vect{\xi}}_j), \, j = 1, ..., N_v$. Since $\Phi_h$ is defined only at
the vertices of $\cTc^{n+1}$, we need to compute $\Phi_h(\hat{\vect{\xi}}_j)$ using interpolation.
Linear interpolation is used since it is important to keep the nonsingularity of the mesh
while it is unnecessary to compute the mesh to high accuracy.
%}

% We will also use a computational mesh
%$\cTc = \{ \vect{\xi}_j, j = 1, ..., N_v\}$ as an intermediate variable.
%
%The procedure for generating $\cTh^{n+1}$ is as follows.
%An MMPDE-based mesh equation (to be described below)
%for the velocities of the computational vertices
%is integrated, starting with the reference computational mesh $\cTco$,
%from $t_n$ to $t_{n+1}$ for the new computational mesh $\cTc^{n+1}$.
%The equation is formulated 
%based on the current mesh $\cTh^n$ and a metric tensor $\bM$ defined thereon.
%This metric tensor is computed using the solution $u_h^n$. 
%The obtained computational mesh $\cTc^{n+1}$ forms a correspondence relation
%with $\cTh^n$. Denote the relation by $\cTh^n = \Phi_h( \cTc^{n+1} )$.
%The new physical mesh is defined as $\cTh^{n+1} = \Phi_h( \cTco )$, which is computed using linear interpolation.

The metric tensor is assumed to be symmetric and uniformly positive definite on $\Omega$.
It is used to control the size, shape, and orientation of the elements of the mesh to be generated.
We consider three types of mesh, one is uniform and the other two are arclength- and Hessian-based
adaptive meshes. The metric tensors associated with the adaptive meshes are defined as
\begin{align}
\bM & = (\mathbb{I} + \grad u_h^n (\grad u_h^n)^T)^{\frac{1}{2}} ,
\label{M-arclength}
\\
\bM & = \left[ \det \left( \mathbb{I} +  |H(u_h^n)| \right)\right]^{-\frac{1}{6}}
                \left( \mathbb{I} +  |H(u_h^n)|\right) ,
 \label{M-hessian}
 \end{align}
where $\mathbb{I}$ is the $d\times d$ identity matrix, $H(u_h^n)$ is a recovered Hessian
for the piecewise linear finite element solution $u_h^n$, and
$|H(u_h^n)| = Q \text{diag} (|\lambda_1|, ..., |\lambda_d|) Q^T$ with 
$Q \text{diag} (\lambda_1, ..., \lambda_d) Q^T$ being the eigen-decomposition of $H(u_h^n)$.
The tensor (\ref{M-arclength}) is the frequently used arclength monitor function which, loosely speaking,
places mesh points according to the uniformity in arclength. The tensor (\ref{M-hessian})
is optimal \cite{HS03} for the $L^2$ norm of linear interpolation error.
In our computation, we use a least squares fitting strategy for Hessian recovery
(e.g., see \cite{Kamens09PhD,KaHu2013}).

We now describe the formulation of the MMPDE-based mesh equation.
A key component of the formulation is the $\bM$-uniform mesh concept
with which any (nonuniform) adaptive mesh is viewed as a uniform
one in some metric. It is known \cite{Hua06,HR11} that such an $\bM$-uniform mesh $\cTh$
approximately satisfies the equidistribution and alignment conditions
\begin{align}
& |K| \sqrt{\det(\bM_K)} = \frac{\sigma_h |K_c|}{|\Omega_c|} , \quad \forall K \in \cTh
\label{eq-1}
\\
&\frac{1}{d} \trace 
      \left( 
	\left( F'_K \right)^{-1}
	  \bM_K^{-1} \left (F'_K \right )^{-T}
      \right)
      =
      \det \left( 
	\left( F'_K \right)^{-1}
	  \bM_K^{-1} \left (F'_K \right )^{-T}
      \right)^{\frac{1}{d}} , \quad \forall K \in \cTh
\label{ali-1}
\end{align}
where $|K|$ is the volume of $K$, $\bM_K$ is the average of $\bM$ over $K$, $\det(\cdot)$ and $\trace(\cdot)$
denote the determinant and trace of a matrix, respectively, $|K_c|$ is the volume of the element $K_c \in \cTc$
corresponding to $K$, $F'_K$ is the Jacobian matrix of the affine mapping
$F_K: K_c \to K$, and 
\[
\sigma_h = \sum_{K \in \cTh} |K| \sqrt{\det(\bM_K)}, \quad 
|\Omega_c | = \sum_{K_c \in \cTc} |K_c| .
\]
The equidistribution condition (\ref{eq-1}) requires that the volume of $K$ in the metric $\bM$ be proportional
to $|K_c|$ with constant proportionality while the alignment condition (\ref{ali-1}) requires that $K$ be similar to $K_c$.
The meshes that closely satisfy these conditions can be obtained by minimizing the energy function
\begin{align}
        I_h & = \theta \sum_{K \in \cTh} |K| \sqrt{\det(\bM_K)}
			\left( 
                        \trace({(F_K')}^{-1} {\bM}_{K}^{-1} {(F_K')}^{-T})
			\right)^{\frac{d p}{2}} \nonumber \\
			&  \quad \quad + (1-2\theta) d^{\frac{dp}{2}}
			\sum_{K \in \cTh} |K| \sqrt{\det (\bM_K)}
			\left( 
			\frac{|K_c|}{|K| \sqrt{\det(\bM_K)}} 
                        \right)^{p},
\label{Ih}
 \end{align}
which is a Riemann sum of a continuous functional developed in \cite{Hua01b}
based on equidistribution and alignment for variational mesh adaptation.
Here, $\theta \in (0, \frac{1}{2}]$ and $p > 1$ are non-dimensional parameters. We choose $\theta = 1/3$
and $p = 2$ in our computation.

Notice that $I_h$ is a function of the position of the computational vertices $\{ \vect{\xi}_j\}_{j=1}^{N_v}$ and
the physical coordinates $\{ \vect{x}_j\}_{j=1}^{N_v}$. For the current situation, we choose $\cTh$ to be $\cTh^n$
(the current physical mesh). Then, $I_h$ is the function of $\vect{\xi}_j$, $j = 1, ..., N_v$ only. 
Instead of taking direct minimization of $I_h$ with respect to these coordinates, we follow the MMPDE approach
\cite{HRR94a} and define the moving mesh equation as a gradient system of $I_h$,
\begin{equation}
\frac{d \vect{\xi}_j}{d t} = - \frac{P_j}{\tau} \left [ \frac{\partial I_h} {\partial \vect{\xi}_j} \right ]^T,
\quad j = 1, ..., N_v
\label{mmpde-1}
\end{equation}
where the derivative of $I_h$ with respect to $\vect{\xi}_j$, ${\partial I_h}/{\partial \vect{\xi}_j}$,
is considered as a row vector, $\tau > 0$ is a parameter used to control the response time of the mesh movement
to the change in the metric tensor, and $P_j = \det(\bM(\vect{x}_j))^{\frac{p-1}{2}} $ is chosen such that
(\ref{mmpde-1}) is invariant under the scaling transformation of $\bM$: $\bM \to c \bM$ for any positive
constant $c$. The derivative of $I_h$ with respect to $\vect{\xi}_j$ can be found
analytically using the notion of scalar-by-matrix differentiation; see \cite{HK2014}.
With these analytical formulas, we can rewrite (\ref{mmpde-1}) into 
\begin{equation}
\label{mmpde-2}
\frac{d \vect{\xi}_j} {d  t}= \frac{P_j}{\tau} \sum_{K \in \omega_j} |K| \vect{v}_{j_K}^K , \quad j = 1, \dotsc, N_v
\end{equation}
where $\omega_j$ is the element patch associated with the $j$-th vertex, $j_K$ is its local index of the vertex on $K$,
and $\vect{v}_{j_K}^K$ is the velocity contributed by the element $K$ to the vertex $j_K$.
The velocities contributed by $K$ to its vertices are given by
\begin{equation}
\label{mmpde-3}
\begin{bmatrix} {(\vect{v}_1^K)}^T \\ \vdots \\ {(\vect{v}_d^K)}^T \end{bmatrix}
= - E_K^{-1} \frac{\partial G}{\partial \mathbb{J}} - \frac{\partial G}{\partial \det(\mathbb{J})}
\frac{\det(\hat{E}_K)}{\det(E_K)} \hat{E}_K^{-1}, 
\quad
\vect{v}_0^K = - \sum_{i=1}^d \vect{v}_d^K ,
\end{equation}
where $E_K = [\vect{x}_1^K-\vect{x}_0^K, ..., \vect{x}_d^K-\vect{x}_0^K]$ is the edge matrix of $K$,
$\hat{E}_K$ is the edge matrix for $K_c$ which is defined similarly, the function $G$ is associated with
the energy (\ref{Ih}) and defined as
\[
G(\mathbb{J}, \det(\mathbb{J}), \bM)
= \theta \sqrt{\det(\bM)} \left ( \trace(\mathbb{J} \bM^{-1} \mathbb{J}^T) \right )^{\frac{dp}{2}}
+ (1-2\theta) d^{\frac{d p}{2}} \sqrt{\det(\bM)} \left (\frac{\det(\mathbb{J})}{\sqrt{\det(\bM)}} \right )^p,
\]
and its derivatives (evaluated at $(\mathbb{J}, \det(\mathbb{J}), \bM) = ((F_K')^{-1}, \det(F_K')^{-1}, \bM_K)$)
with respect to the first (the Jacobian matrix) and second arguments are given by
\begin{align*}
& \frac{\partial G}{\partial \mathbb{J}} = d p \theta \sqrt{\det(\bM)} \left ( \trace(\mathbb{J} \bM^{-1} \mathbb{J}^T)
\right )^{\frac{dp}{2}-1}\bM^{-1} \mathbb{J}^T,
\\
& \frac{\partial G}{\partial \det(\mathbb{J})} = p (1-2\theta) d^{\frac{d p}{2}} \det(\bM)^{\frac{1-p}{2}}
\det (\mathbb{J})^{p-1} .
\end{align*}
Notice that ${\partial G}/{\partial \mathbb{J}}$ is a $d$-by-$d$ matrix.

In practical computation, we can first compute the local velocities $\vect{v}_j^K$, $j = 0, ..., d$ for all elements
using (\ref{mmpde-3}) and then obtain the velocity for any mesh point
by summing the volume weighted contributions from its neighboring elements (cf. (\ref{mmpde-2})).
The mesh equation should be modified for boundary mesh points.
For fixed points (such as corners),
we can set the velocity to be zero. For those on a boundary edge or surface, the mesh velocities
should be modified such that they do not move out of the domain.

The mesh equation (\ref{mmpde-2}) (with proper modifications for boundary mesh points) is integrated
from $t=t_n$ to $t_{n+1}$ starting with $\cTco$ as the initial mesh. In our computation,
the Matlab ODE solver {\em ode15s} (an implicit scheme) is used to integrate (\ref{mmpde-2}).

Equation (\ref{mmpde-2}) is called the $\vect{\xi}$-formulation of the MMPDE moving mesh method
since it has been formulated in terms of the derivatives of $I_h$ with respect to $\vect{\xi}_j$ and
the velocities for the computational coordinates. We can obtain an $\vect{x}$-formulation by directly
differentiating $I_h$ with respect to $\vect{x}_j$ (with $\cTc$ being taken as $\cTco$ and fixed)
and the new physical mesh $\cTh^{n+1}$ by directly integrating this formulation. The main disadvantage
of this formulation is that its formula is more complicated than that of the $\vect{\xi}$-formulation
and the metric tensor, which is defined on $\cTh^{n}$, needs to be updated every time the physical
mesh is changed during the time integration of the mesh equation for $\cTh^{n+1}$.
It is analytically shown in \cite{HK2015} that the mesh governed by the $\vect{x}$-formulation will stay
nonsingular when it is nonsingular initially. Although such a theoretical result is not available for
the $\vect{\xi}$-formulation, our limited numerical experience shows that the $\vect{\xi}$-formulation
also produces nonsingular meshes.

%%%%%%%%%%%%%%%%%%%%%%%%%%%%%%%
%% Section: Numerical Results 
\section{Numerical Results for PME}
\label{SEC:PME-numerics}

\begin{figure}[ht]
    \begin{center}
    \begin{subfigure}[b]{0.35\linewidth}
       \includegraphics[scale=0.35]{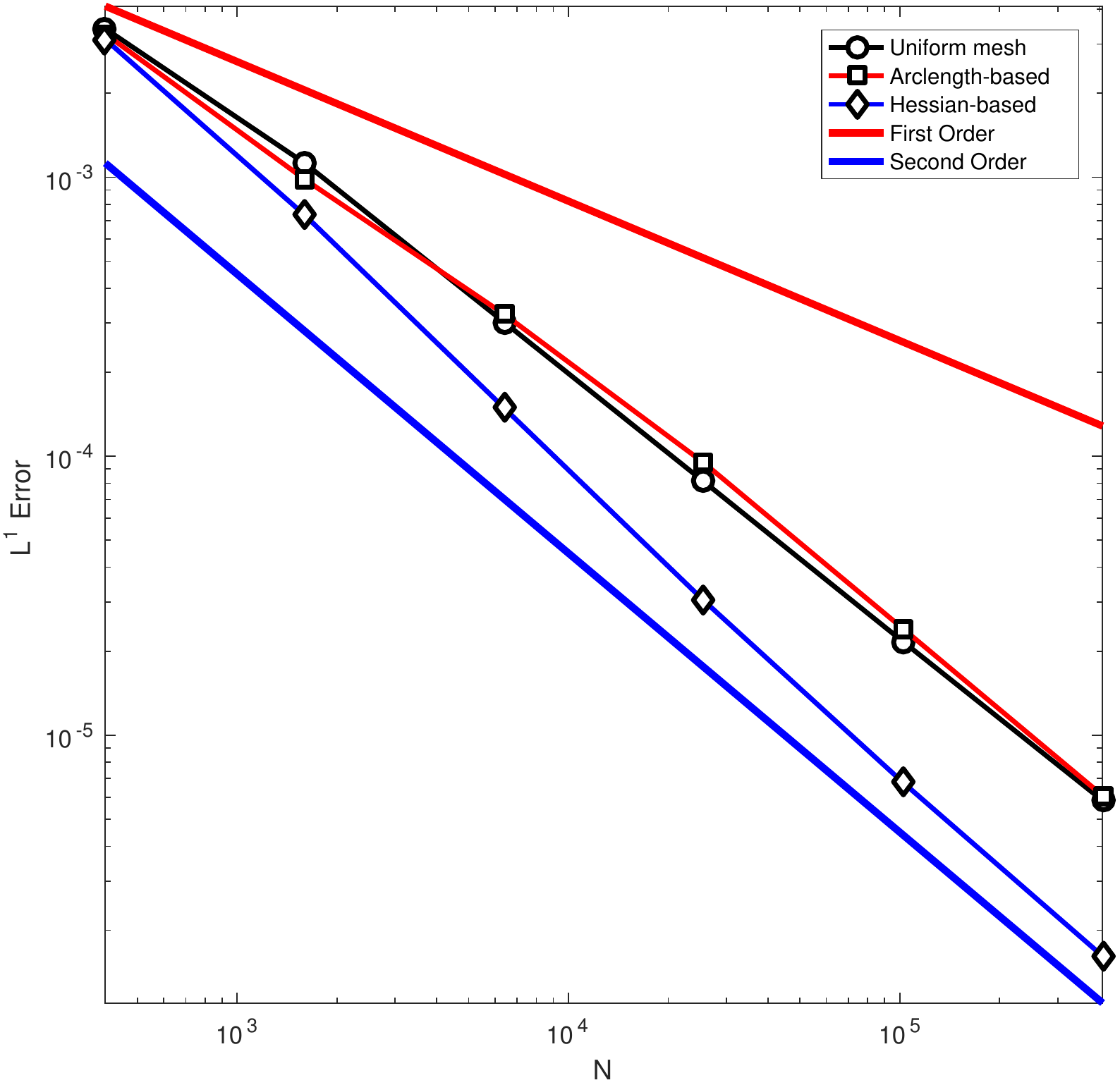}
       \caption{$m = 1$}
       \label{subfig:pme4-adaptivity-m1-L1}
    \end{subfigure}%
    \hspace{3mm}
    \begin{subfigure}[b]{0.35\linewidth}
       \includegraphics[scale=0.35]{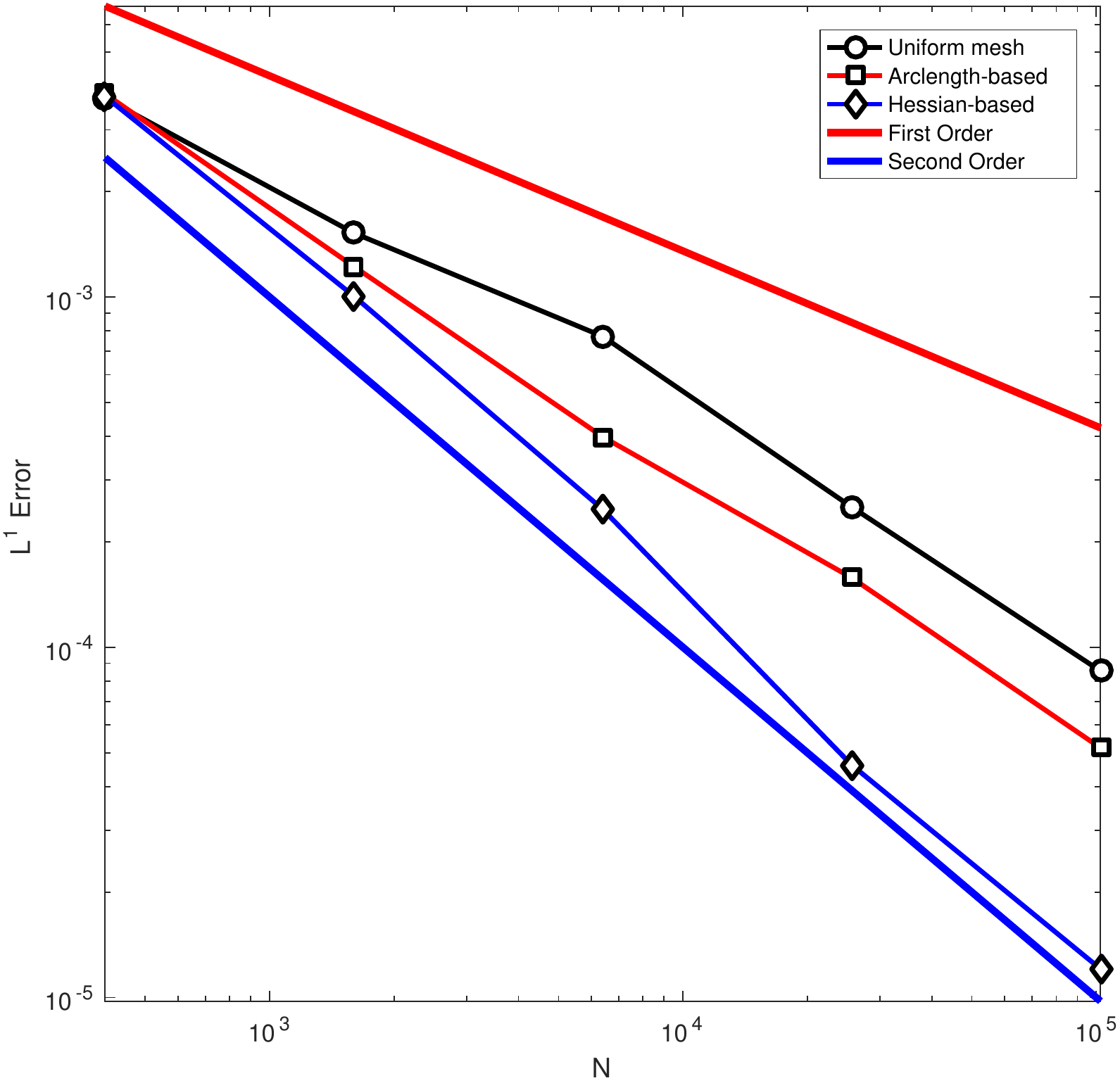}
       \caption{$m = 2$}
       \label{subfig:pme4-adaptivity-m2-L1}
    \end{subfigure}%
    \caption{Example~\ref{exam4.1}. Convergence history (in $L^1$ norm) for the three meshing strategies as $N$ increases.}
    \label{fig:pme4-adaptivity-L1}
    \end{center}
\end{figure}

In this section we present numerical results obtained with the moving mesh finite element method described
in the previous section for a number of PME examples. They include the Barenblatt-Pattle
solution and the generalizations of several one-dimensional examples studied by
Zhang and Wu \cite{Zhang-2009}. These examples are selected to demonstrate the accuracy of our method
as well as its ability to deal with solutions having complex support and the waiting-time phenomenon.
For the cases having an exact solution, the error in the computed
solution will be measured in the (global) $L^2$ norm, i.e.,
\[
\|e_h\|_{L^2(t_0,T; L^2(\Omega))} = \left ( \int_{t_0}^{T} \int_{\Omega} e_h^2(\vect{x}, t) d \vect{x}  d t \right )^{\frac{1}{2}}.
\]
We choose this norm because various error estimates are obtained in this norm, e.g., see (\ref{ebmeyer-1}).
%{\color{red}
(An exception is Fig.~\ref{fig:pme4-adaptivity-L1} where the convergence history in $L^1$ norm is plotted
for comparison purpose.)
%}
In our computation, we use $\tau = 10^{-4}$ (for the mesh movement), the maximal allowed time step size
$\Delta t_{max} = 10^{-3}$ (for integrating PME), and the Hessian-based metric tensor (\ref{M-hessian}),
unless stated otherwise.

\begin{figure}[ht]
    \begin{center} 
        \begin{subfigure}[b]{0.35\linewidth}\includegraphics[scale=0.29]{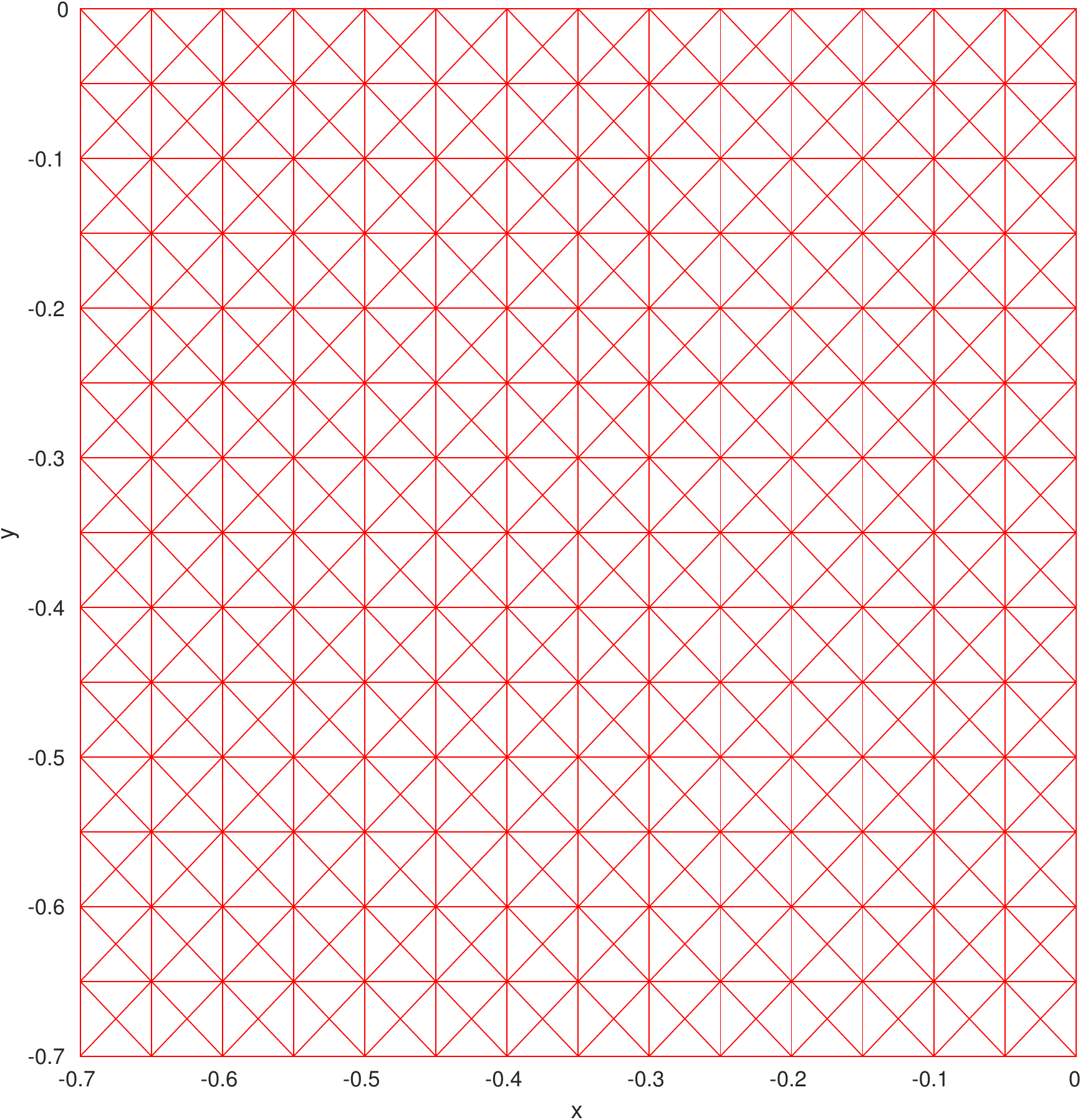}\caption{Uniform mesh}\end{subfigure}
        \begin{subfigure}[b]{0.35\linewidth}\includegraphics[scale=0.3]{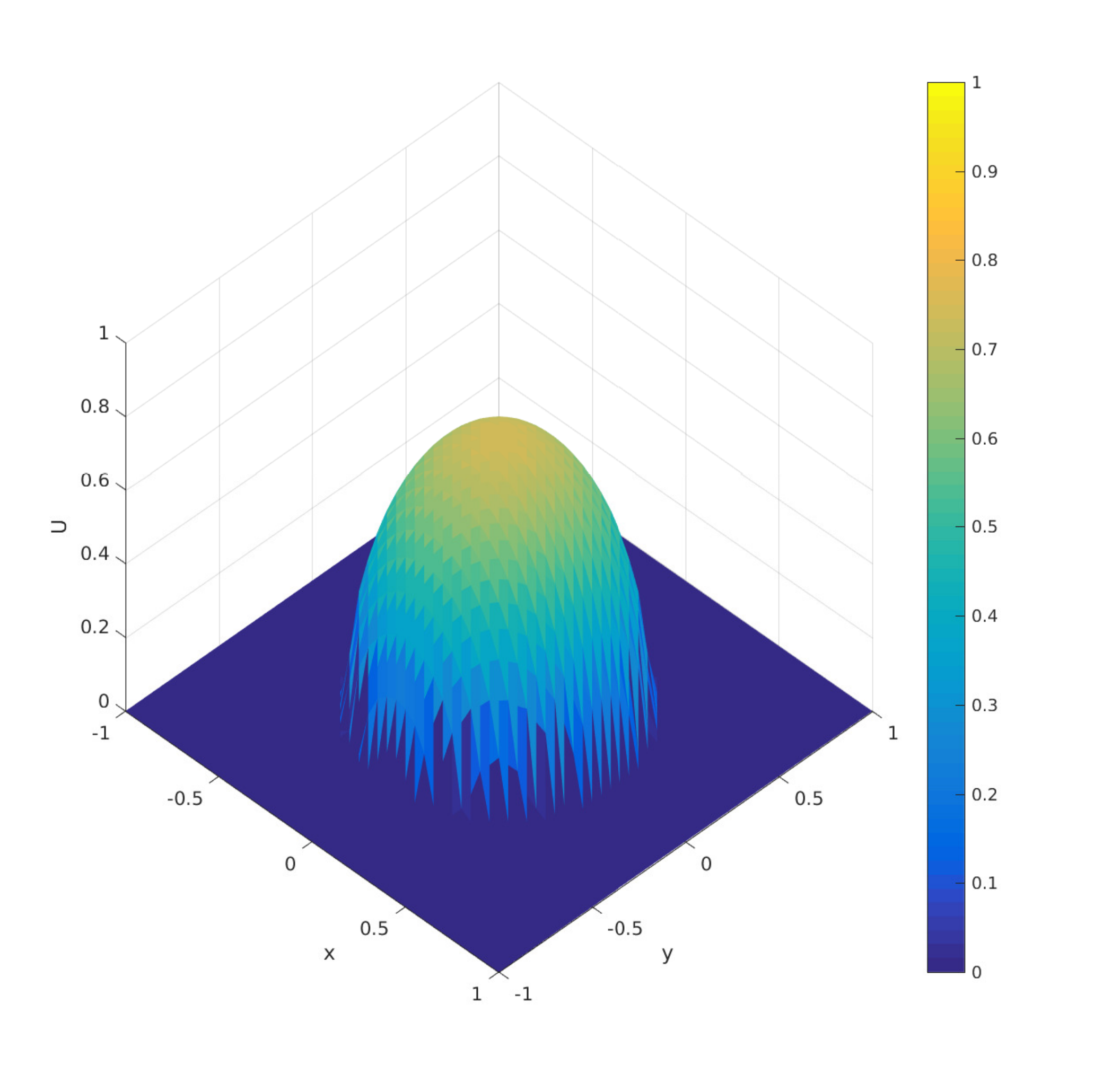}\caption{Uniform mesh}\end{subfigure}\\
        \begin{subfigure}[b]{0.35\linewidth}\includegraphics[scale=0.29]{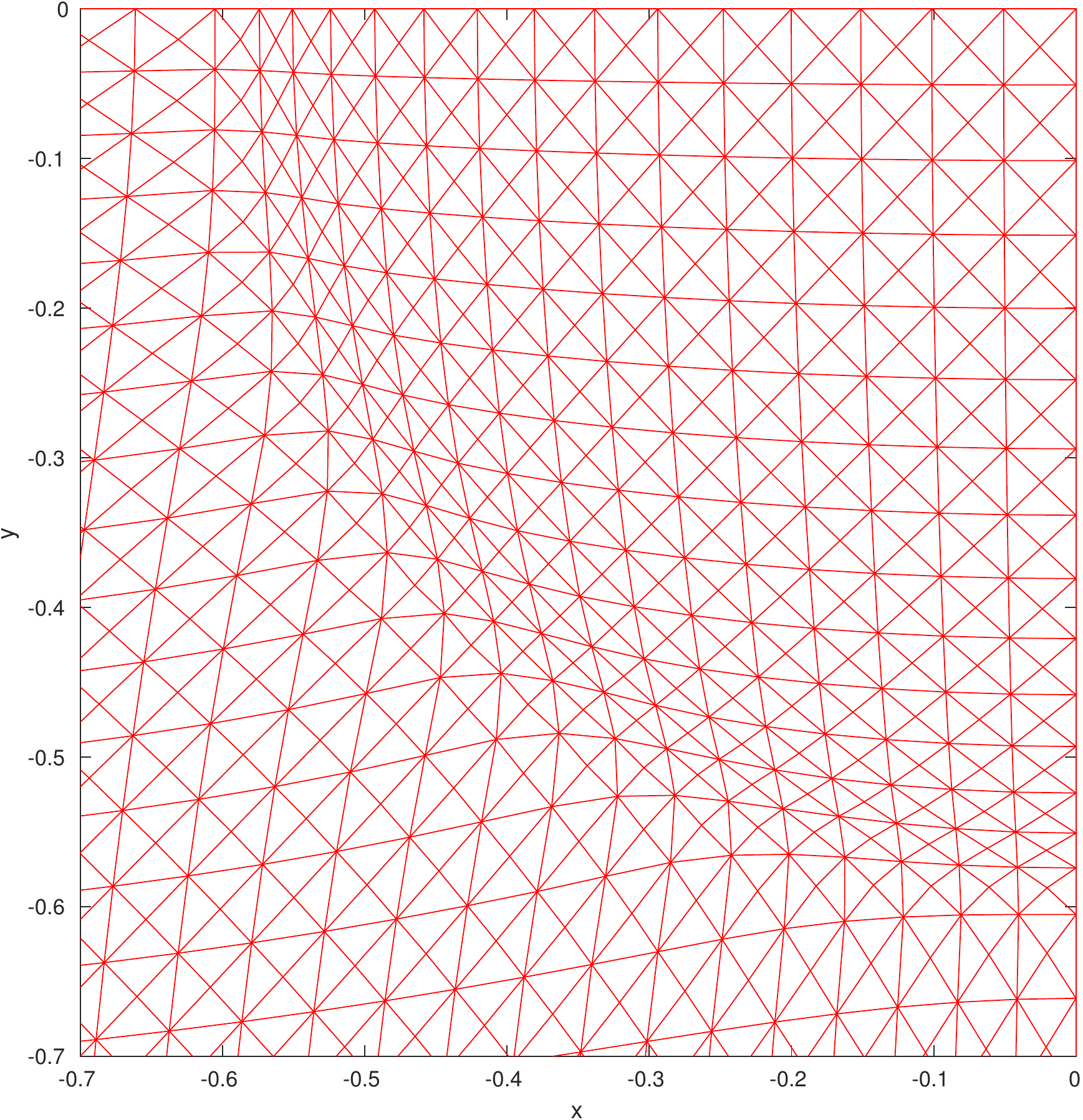}\caption{Arclength metric tensor}\end{subfigure}
        \begin{subfigure}[b]{0.35\linewidth}\includegraphics[scale=0.3]{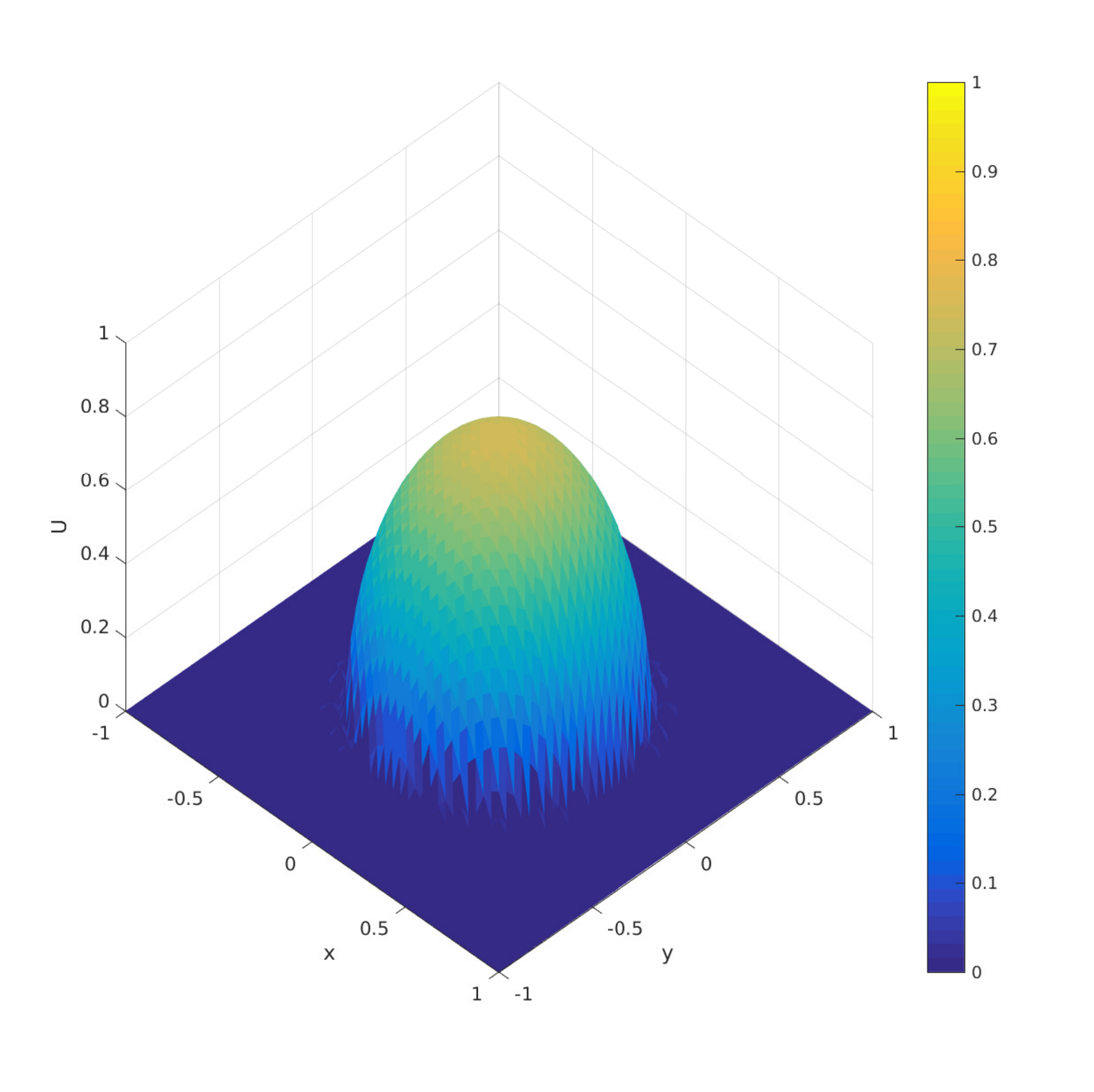}\caption{Arclength metric tensor}\end{subfigure}\\
        \begin{subfigure}[b]{0.35\linewidth}\includegraphics[scale=0.29]{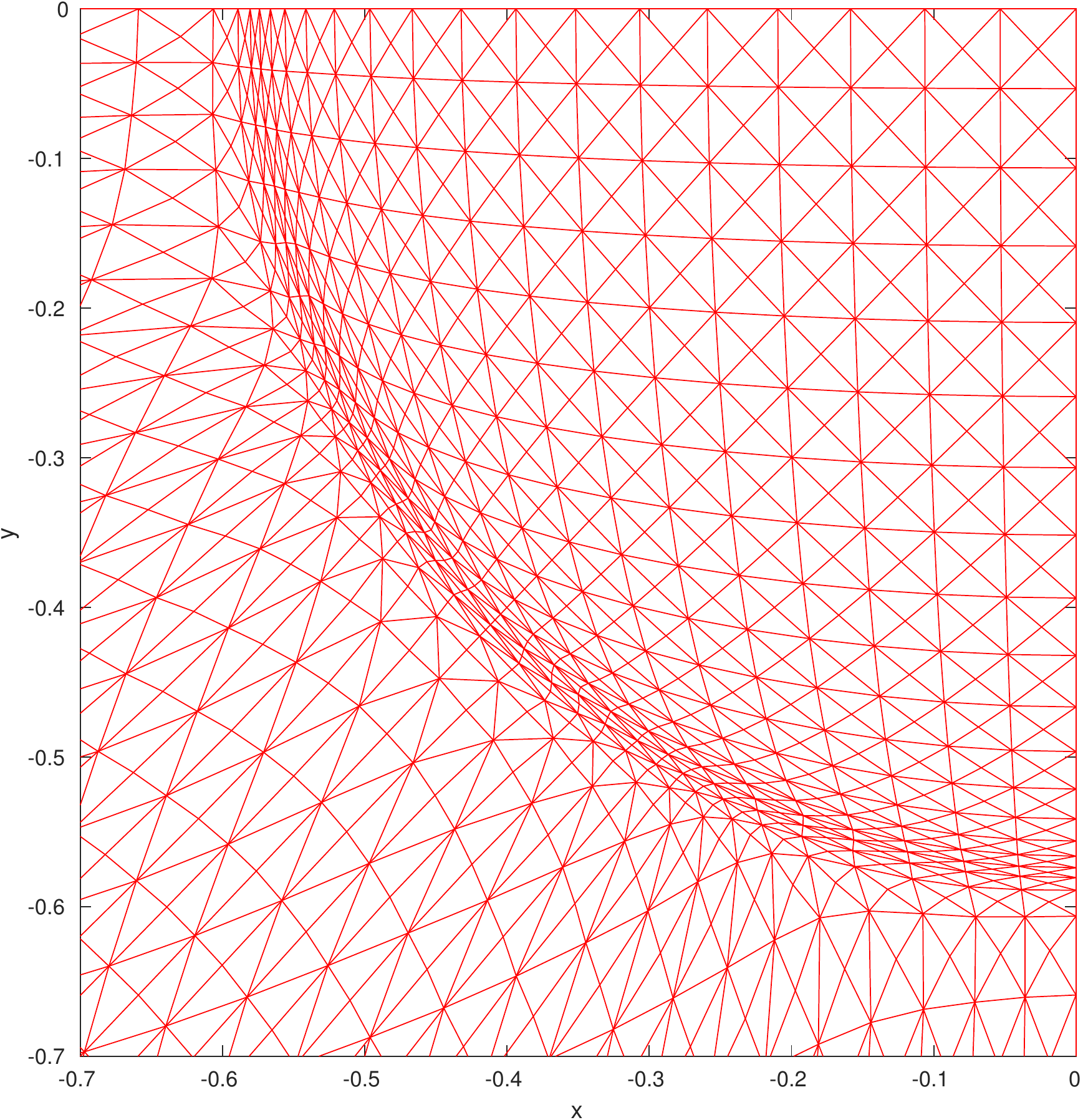}\caption{Hessian-based metric tensor}\end{subfigure}
        \begin{subfigure}[b]{0.35\linewidth}\includegraphics[scale=0.3]{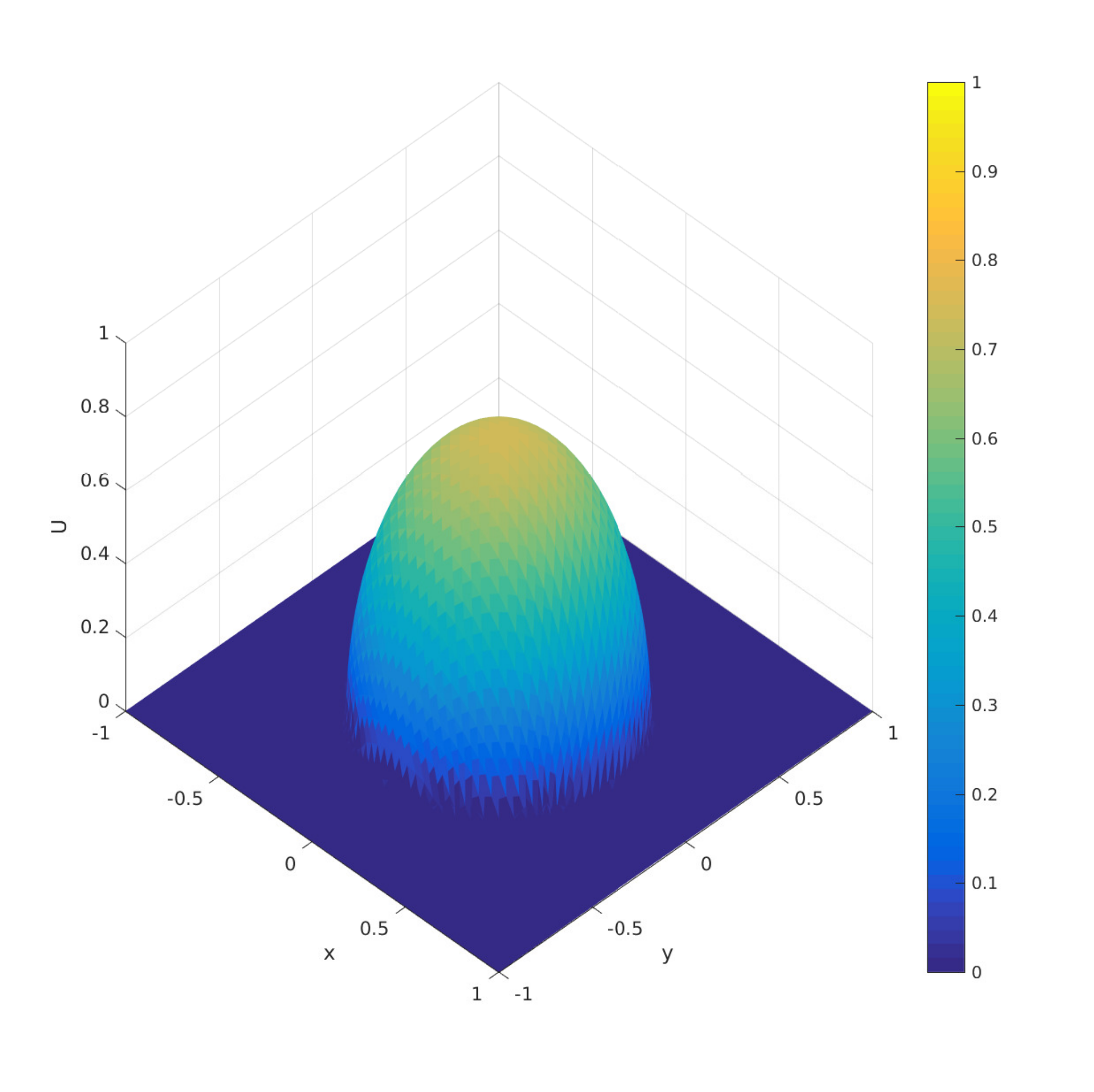}\caption{Hessian-based metric tensor}\end{subfigure}
    \caption{Example~\ref{exam4.1} with $m=2$.
    The meshes (closer view near (-0.35, -0.35)) and computed solutions at $t = T$ obtained with uniform and
    arclength- and Hessian-based adaptive meshes ($N = 25600$).}
\label{fig:mesh-soln-compare}
\end{center}
\end{figure}

\begin{figure}[ht]
    \begin{center}
    \begin{subfigure}[b]{0.35\linewidth}
       \includegraphics[scale=0.35]{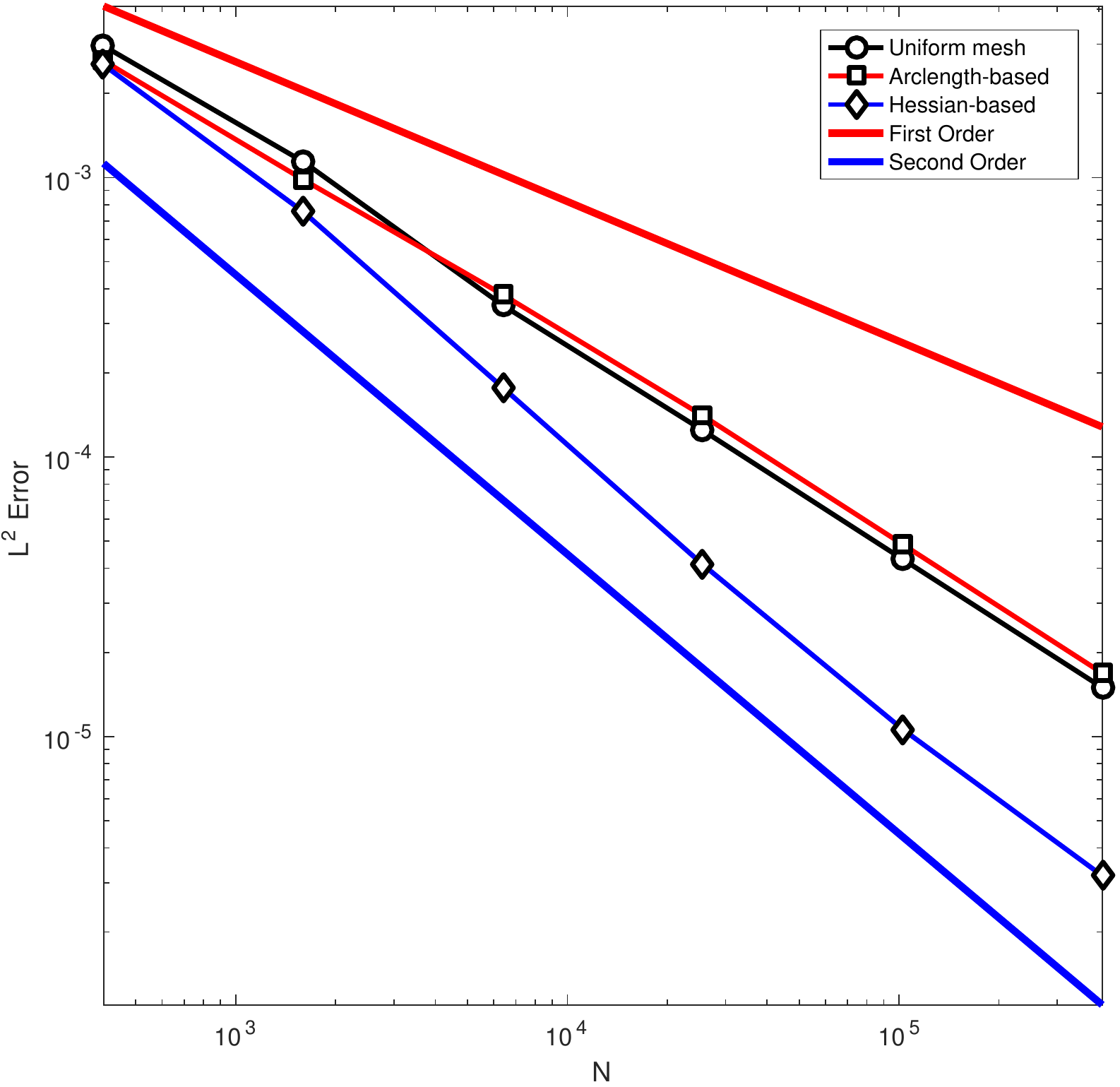}
       \caption{$m = 1$}
       \label{subfig:pme4-adaptivity-m1}
    \end{subfigure}%
    \hspace{3mm}
    \begin{subfigure}[b]{0.35\linewidth}
       \includegraphics[scale=0.35]{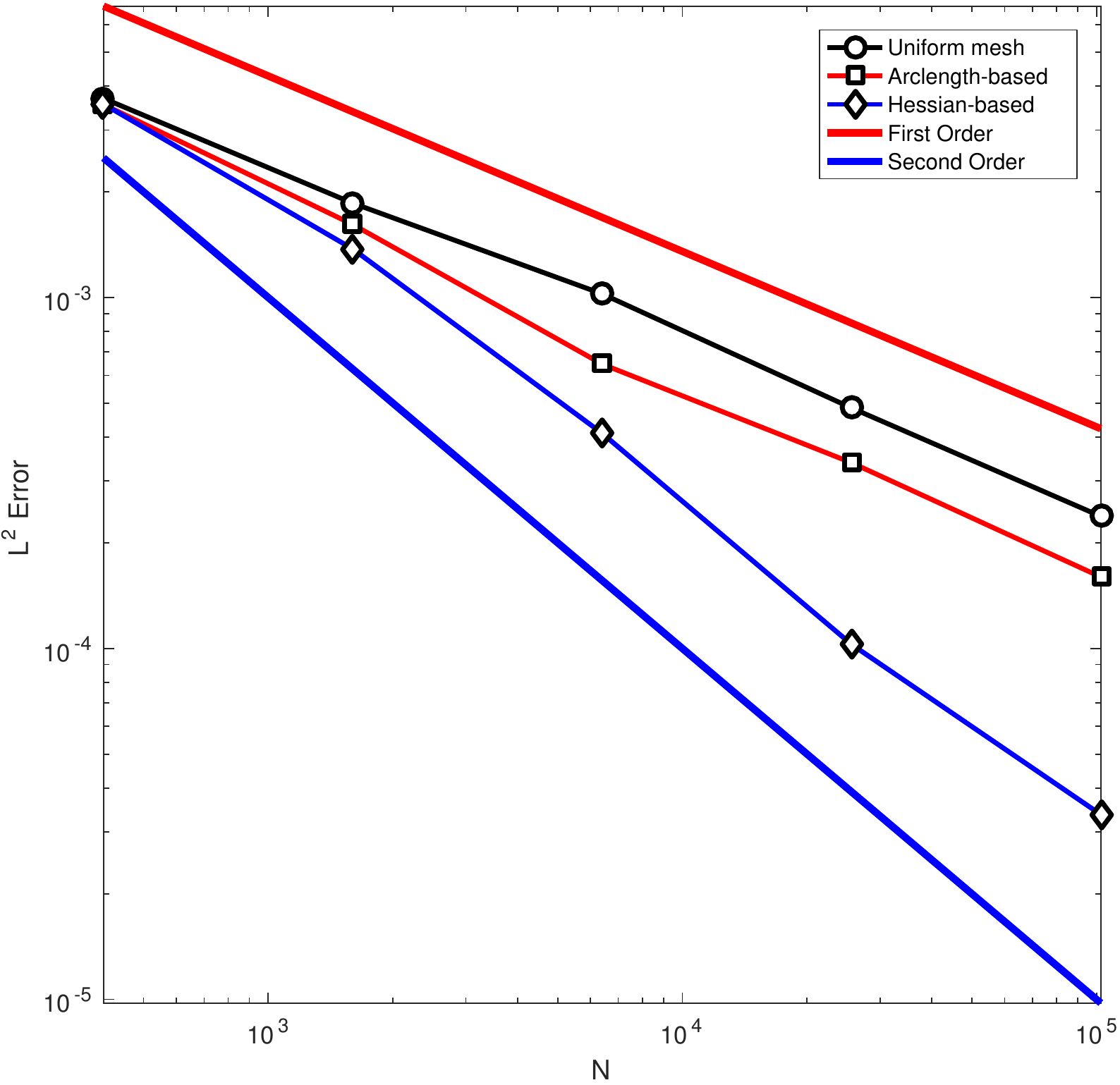}
       \caption{$m = 2$}
       \label{subfig:pme4-adaptivity-m2}
    \end{subfigure}%
    \caption{Example~\ref{exam4.1}. Convergence history (in $L^2$ norm) for the three meshing strategies as $N$ (the number of the elements) increases.}
    \label{fig:pme4-adaptivity}
    \end{center}
\end{figure}

\begin{figure}[ht]
   \begin{center}
       \includegraphics[scale=0.35]{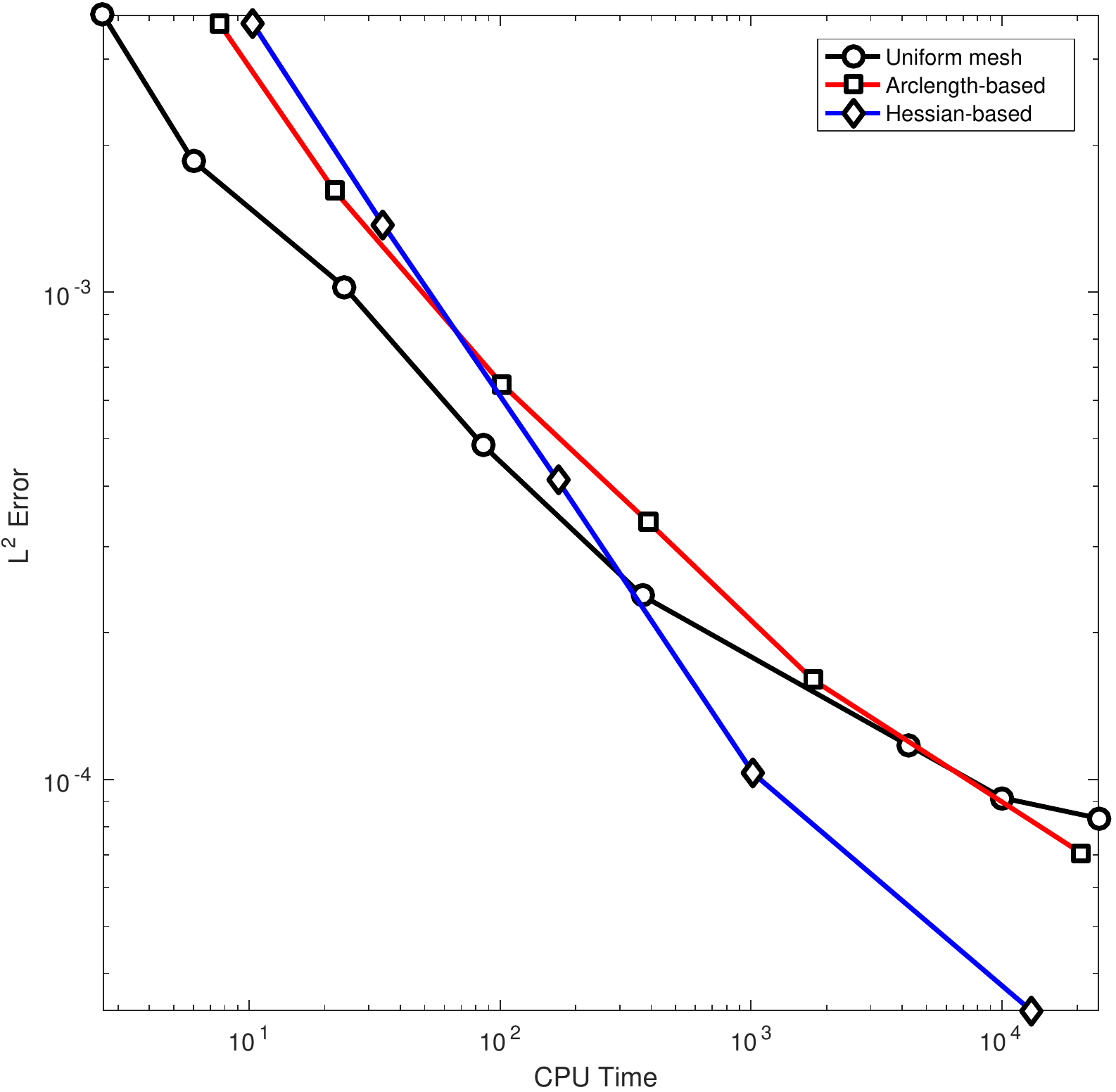}
    \caption{Example~\ref{exam4.1} with $m=2$. The $L^2$ norm of the error is plotted against the CPU time in seconds
    for the computation corresponding to Fig.~\ref{subfig:pme4-adaptivity-m2}.}
    \label{fig:pme4-cpu-m2}
    \end{center}
\end{figure}

\begin{figure}[ht]
    \begin{center}
    \begin{subfigure}[b]{0.35\linewidth}
       \includegraphics[scale=0.30]{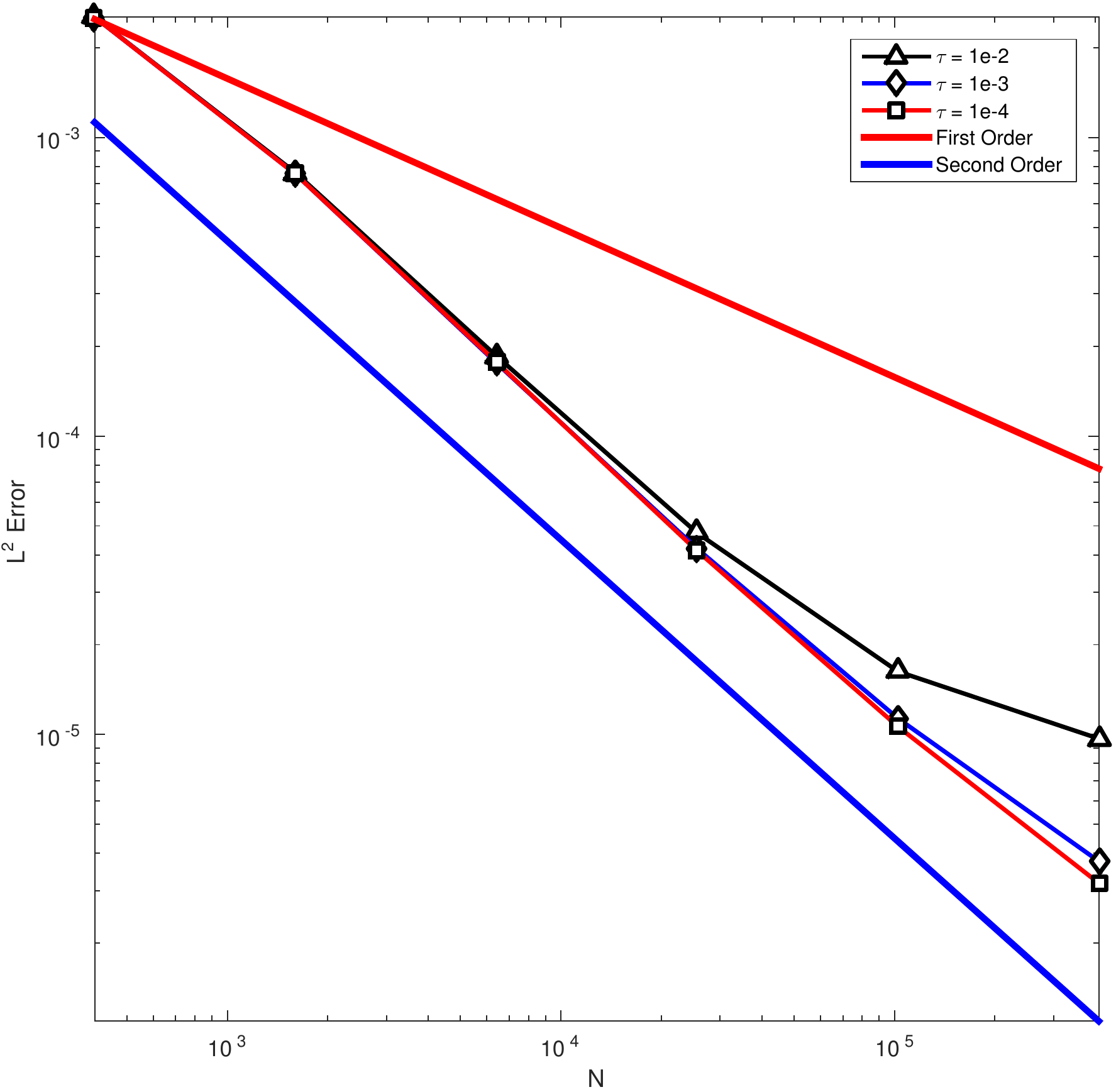}
       \caption{$m = 1$}
       \label{subfig:pme4-tauDelta-m1}
    \end{subfigure}%
    \begin{subfigure}[b]{0.35\linewidth}
       \includegraphics[scale=0.29]{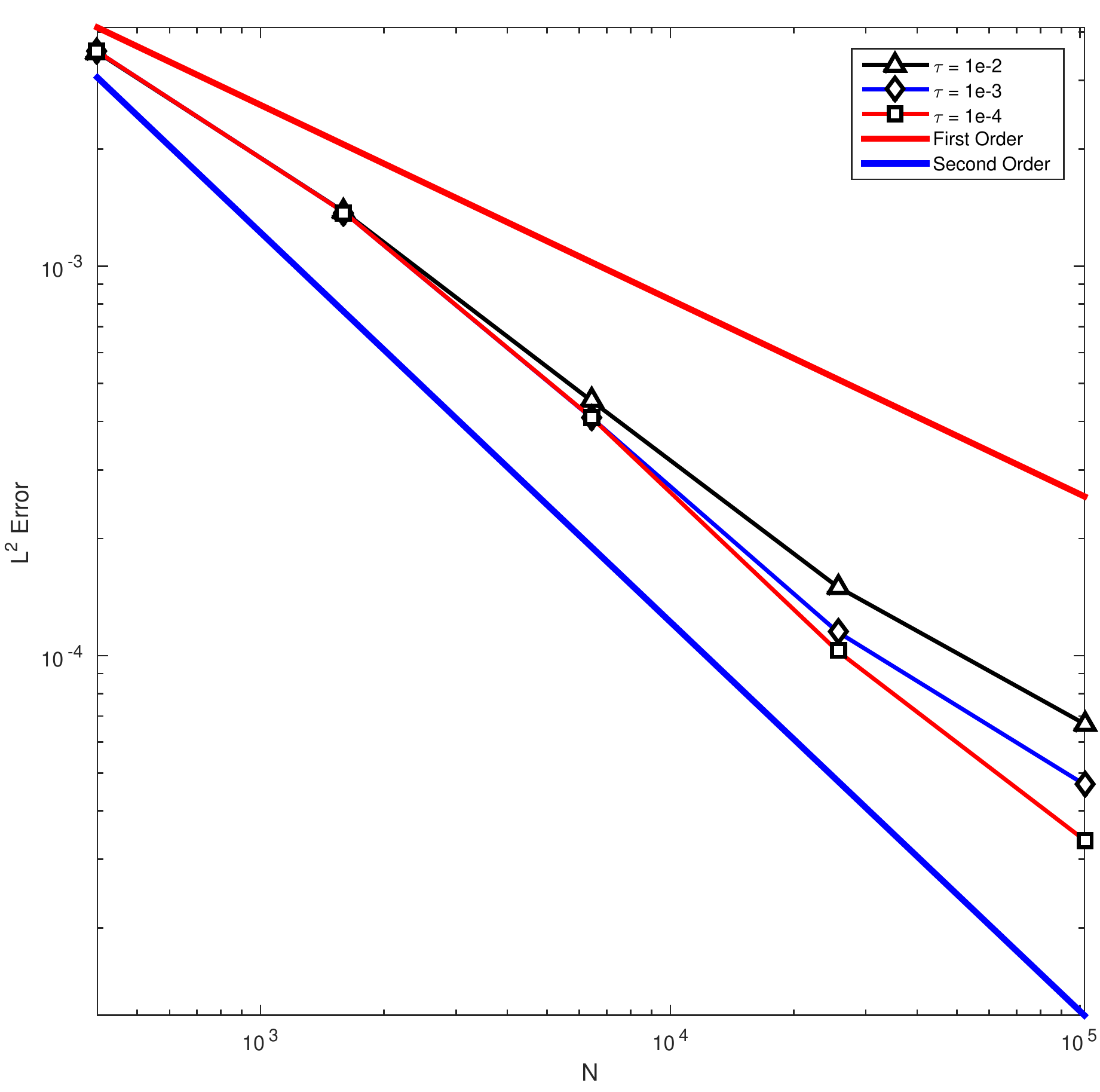}
       \caption{$m = 2$}
       \label{subfig:pme4-tauDelta-m2}
    \end{subfigure}%
    \caption{Example~\ref{exam4.1}. Convergence history for different values of $\tau$.}
    \label{fig:pme4-tauDelta}
    \end{center}
\end{figure}

\begin{figure}[ht]
    \begin{center}
    \begin{subfigure}[b]{0.35\linewidth}
       \includegraphics[scale=0.34]{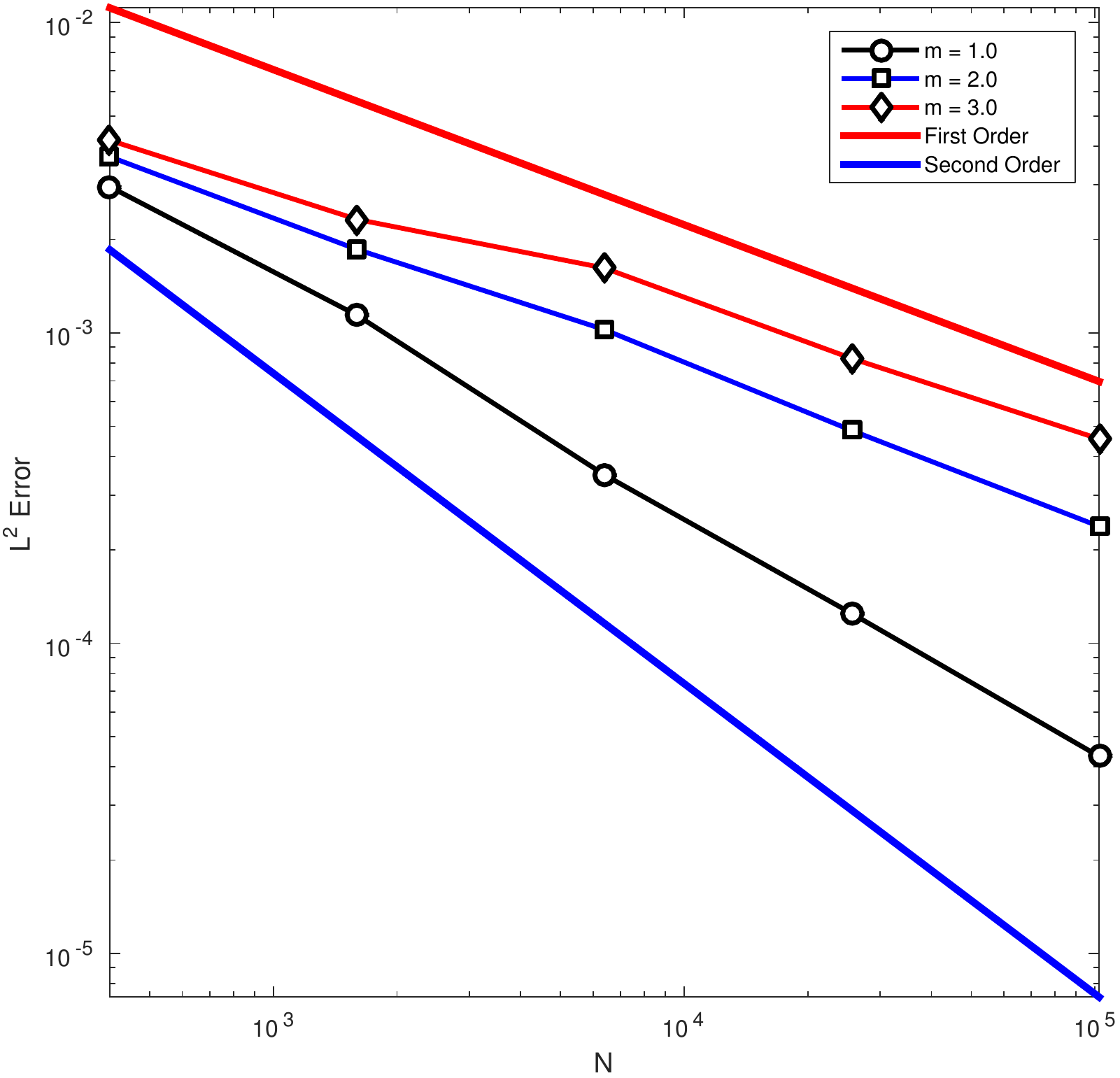}
       \caption{Uniform mesh}
       \label{fig:pme4-m-vary:uniform}
    \end{subfigure}%
    \hspace{1mm}
    \begin{subfigure}[b]{0.35\linewidth}
       \includegraphics[scale=0.34]{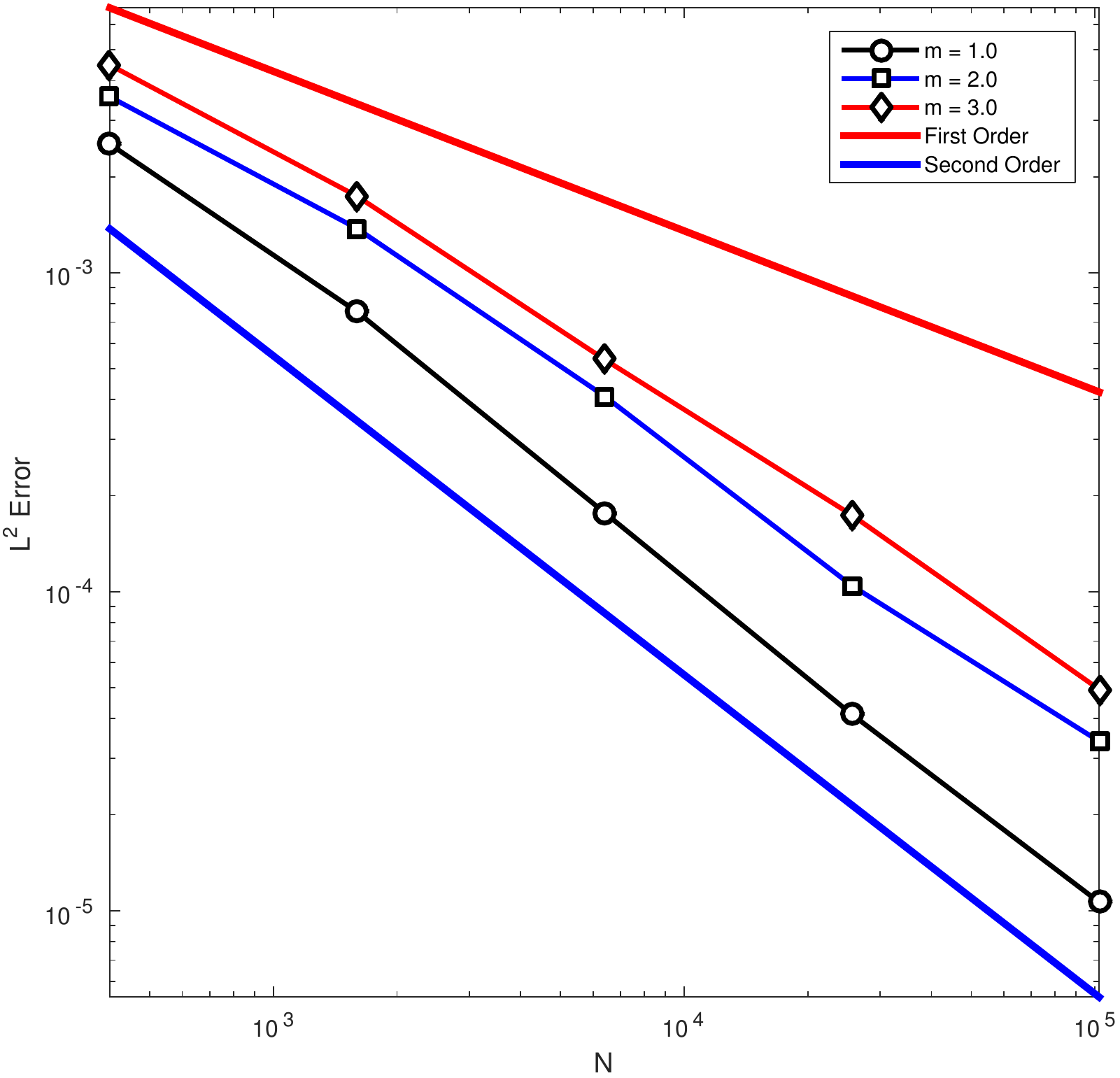}
       \caption{Hessian-based adaptive mesh}
       \label{fig:pme4-m-vary:adaptive}
    \end{subfigure}%
    \caption{Example~\ref{exam4.1}. Convergence history for different values of $m$.}
    \label{fig:pme4-m-vary}
    \end{center}
\end{figure}

\begin{exam}[Barenblatt-Pattle solution]
\label{exam4.1}

We first consider the Barenblatt-Pattle solution (\ref{BP-soln}) with $r_0 = 0.5$ and $T = (t_0+0.1)/2$.
We use it to verify the accuracy of the numerical method and the effects of the mesh adaptivity and
the physical parameter $m$ on the computational accuracy.

Typical meshes and computed solutions at the final time obtained with the uniform mesh and two adaptive mesh
strategies are shown in Fig.~\ref{fig:mesh-soln-compare} and the convergence history is shown in
Fig.~\ref{fig:pme4-adaptivity} for the cases $m = 1$ and $2$. We can see that for both uniform and
the arclength-based adaptive meshes, the convergence order is 
about 1.5 (i.e., $\mathcal{O}(N^{-\frac{1.5}{2}})$) for $m=1$ and 
1 for the case $m=2$, with an arclength-based adaptive mesh producing slightly more
accurate solutions for both cases. We notice that the exact solution (\ref{BP-soln}) is
in $H^1(\Omega)$ for $m = 1$ and $W^{1,\frac{m}{m-1}-\epsilon}(\Omega)$ for $m> 1$, where $\epsilon$ is
a small positive number. The observed convergence order is higher than what we can expect from
the solution regularity. (For example, the theoretical estimate (\ref{ebmeyer-1})
shows a convergence order of $15/28$ for $m=1$ and $6/15$ for $m=2$.)
 Even more surprisingly, Hessian-based adaptive meshes lead to a second-order
 convergence rate for both the $m=1$ and $2$ cases.
We do not have a rigorous explanation for this but would like to point out two relevant observations.
The first is that the mesh is denser near the free boundary with the Hessian-based metric tensor than
with the arclength metric tensor (e.g., see Fig.~\ref{fig:mesh-soln-compare}).
The other is that the exact solution has higher regularity in its support than on the whole domain $\Omega$.
Indeed, it can be directly verified that
\[ 
        \sqrt{\det(|H(u(\cdot, t))|)} \in L^{\frac{2m}{3m-2}-\epsilon}(\text{supp}(u(\cdot, t))),
\]
where $\epsilon$ is a small positive number and $H(u(\cdot, t))$ denotes the Hessian of $u$. It is known \cite{HS03} that for an $\mathbb{M}$-uniform
mesh associated with the metric tensor (\ref{M-hessian}), the linear interpolation error on a polygonal domain $D$
is bounded by   
\[
        \| u - \Pi_1 u \|_{L^2(D)} \leq C N^{-1} \| \sqrt{\det(|H(u)|)}\|_{L^{\frac{2}{3}}(D)} + h.o.t. ,
\]
where $h.o.t.$ stands for higher-order terms. From this we can expect a second-order convergence
if we consider linear interpolation only in the support of the solution with a Hessian-based adaptive mesh.
Although this analysis does not apply directly to our current situation with a larger domain than the support,
it may shed some light on why the scheme with Hessian-based adaptive meshes shows
a second-order convergence.

%{\color{red}
For comparison purpose, we plot the convergence history in the $L^1$ norm in Fig.~\ref{fig:pme4-adaptivity-L1}. It can be
seen that the the $L^1$ norm of the error behaves similarly as the $L^2$ norm.

We have seen from Figs.~\ref{fig:pme4-adaptivity} and \ref{fig:pme4-adaptivity-L1} that mesh adaptation, especially the
Hessan-based one, can significantly improve the accuracy. But this 
comes with additional cost. To show if the mesh adaptation can also improve the efficiency,
we plot the solution error against the required CPU time (in seconds) in Fig.~\ref{fig:pme4-cpu-m2}
for the computation corresponding to Fig.~\ref{subfig:pme4-adaptivity-m2}.
We can see that a uniform mesh is more efficient when low accuracy is desired while mesh adaptation shows
advantages for high accuracy. This is consistent with our limited experience with adaptive
moving mesh computation (also see \cite[Page 17]{HR11}). The location of the break-even point depends on
specific problems and specific mesh adaptation strategies. For the current situation,
we have $(N, \|e_h\|_{L^2}) \approx (300, 3\times 10^{-4})$
for Hessian-based adaptation and $(5000, 10^{-4})$ for arclength-based adaptation.
%}

We now examine the effects of the parameter $\tau$ on the accuracy. Recall that $\tau$ is used in
the moving mesh equation (\ref{mmpde-1}) to adjust the response time of the mesh movement to
the changes in the metric tensor.
%{\color{red}
The smaller $\tau$ is, the faster the response is. On the other hand, for smaller $\tau$,
the mesh equation (\ref{mmpde-2})  becomes stiffer and harder to integrate.
Fortunately, this only causes a slight increase in the cost when an implicit solver (Matlab solver {\em ode15s} in our computation)
is used for the mesh equation.
%}
The convergence history is shown in Fig.~\ref{fig:pme4-tauDelta}
for Hessian-based adaptive meshes for $\tau = 10^{-2}, 10^{-3}$, and $10^{-4}$.
We can see that for both cases with $m=1$ and $m=2$, the convergence with $\tau = 10^{-2}$ and
$10^{-3}$ slows down when the mesh is becoming finer  whereas that with $\tau=10^{-4}$ stays
second order at least for the same considered range of the number of mesh elements.
This indicates that the mesh concentration needs to follow up the movement of the free boundary
very closely or otherwise we may lose the accuracy improvements gained with mesh adaptation. 

Next, we consider the effects of the physical parameter $m$. As seen in Section~\ref{SEC:PME-theory},
the solution at the free boundary becomes steeper for larger $m$. It is not surprising that PME will also become
more difficult to solve numerically. Indeed, as we can see in Fig.~\ref{fig:pme4-m-vary:uniform},
the convergence rate for the uniform mesh decreases as $m$ increases.
This is qualitatively consistent with the theoretical analysis for various finite element
approximations for PME on quasi-uniform meshes
which also shows a decrease in convergence order with $m$; cf. (\ref{rose-1}), (\ref{nochetto-1}), (\ref{ebmeyer-1}),
and (\ref{wei-1}) and e.g., see \cite{Ebmeyer-1998,Ebmeyer-2008,Nochetto-1988,Rose1983}.
On the other hand, for Hessian-based adaptive meshes the convergence order is second
for $m=1$, 2, and 3, although the error is larger for larger $m$;
see Fig.~\ref{fig:pme4-m-vary:adaptive}. 
%[It should also be pointed out that for $m=2$ and $3$,
%a smaller $\tau = 10^{-5}$ was used in the computation
%(with additionally a smaller time step $\Delta t_{max} = 10^{-5}$ for the case $m = 3$), which is needed to show
%the second order convergence rate.]
The final mesh and computed solution obtained for $m=3$ with the Hessian-based
mesh adaptation are shown in Fig.~\ref{fig:pme4-m3}.

%{\color{red}
It is worth pointing out that there are small oscillations around the free boundary
in computed solutions; e.g., see Fig.~\ref{fig:oscillations-cross-section}.
This is due to the nature of the problem where the solution is steep and has a corner shape
near the free boundary and the loss of the maximum principle in the discretization.
A standard finite element discretization like the one we used here
typically leads to solutions with oscillations for this type of problems
(also see Zhang and Wu \cite{Zhang-2009} for the case with the one-dimensional PME).
The oscillations may be suppressed using, for instance, monotone schemes (e.g., see \cite{BaSo1991,NH2015,Oberman2006})
or structure-preserving schemes (e.g., see \cite{LSSV07,LiHu2013,Le05,YuSh2008,Zhang-2009,ZZS2013}). 
These schemes and their combination with adaptive mesh movement for PME are worth future investigations.
%}
\qed \end{exam}

\begin{figure}[ht]
    \begin{center}
    \begin{subfigure}[b]{0.35\linewidth}
       \includegraphics[scale=0.29]{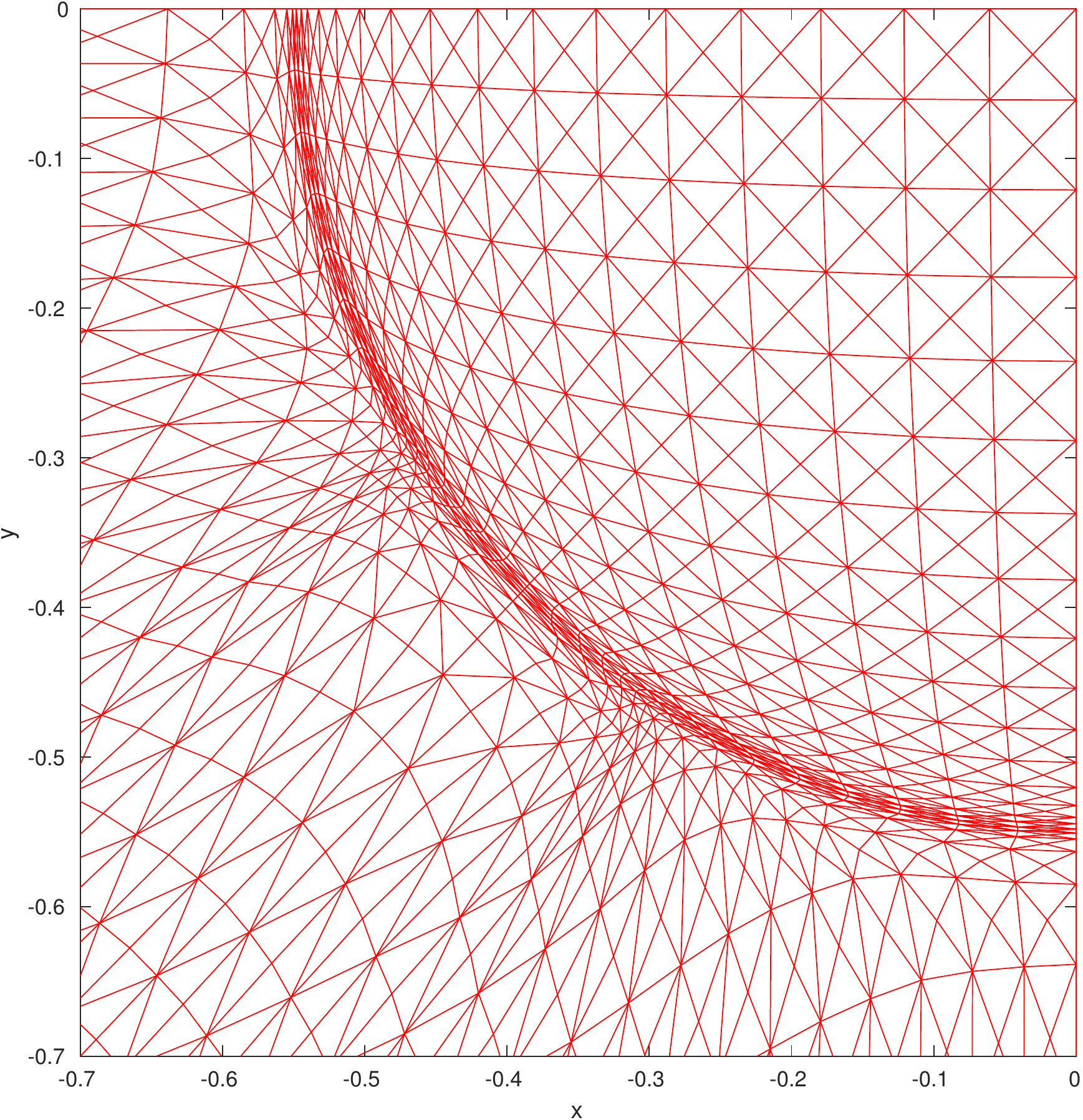}
       \caption{Mesh}
    \end{subfigure}%
    \begin{subfigure}[b]{0.35\linewidth}
       \includegraphics[scale=0.30]{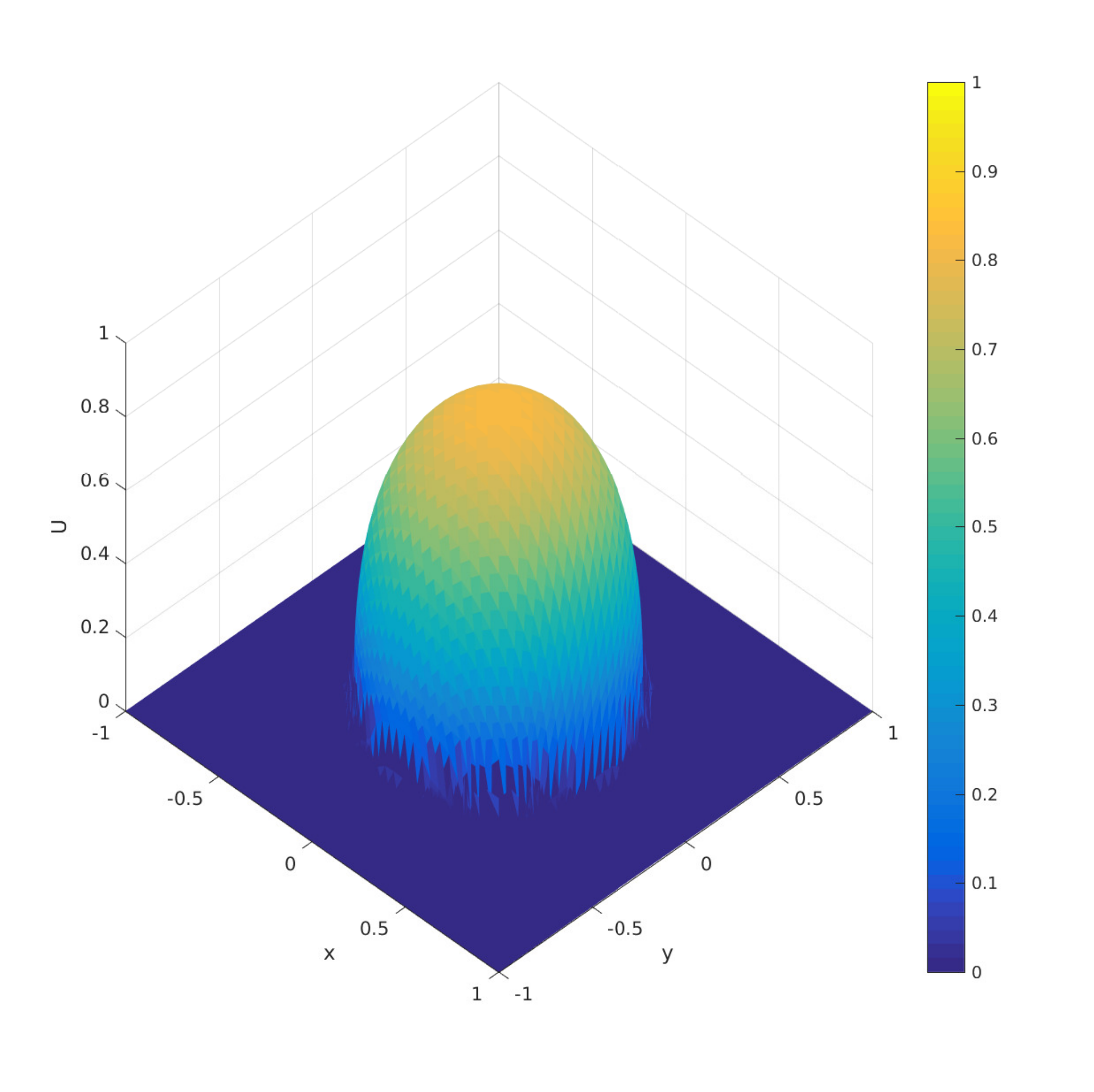}
       \caption{Computed solution}
    \end{subfigure}%
    \caption{Example~\ref{exam4.1}. The final mesh (close view near (-0.35, -0.35)) and computed solution for $m = 3 $
    		with the Hessian-based mesh adaptation ($N = 25600$).}
    \label{fig:pme4-m3}
    \end{center}
\end{figure}

 \begin{figure}[ht]
        \centering
        \begin{subfigure}[b]{0.25\linewidth}\includegraphics[scale=0.25]{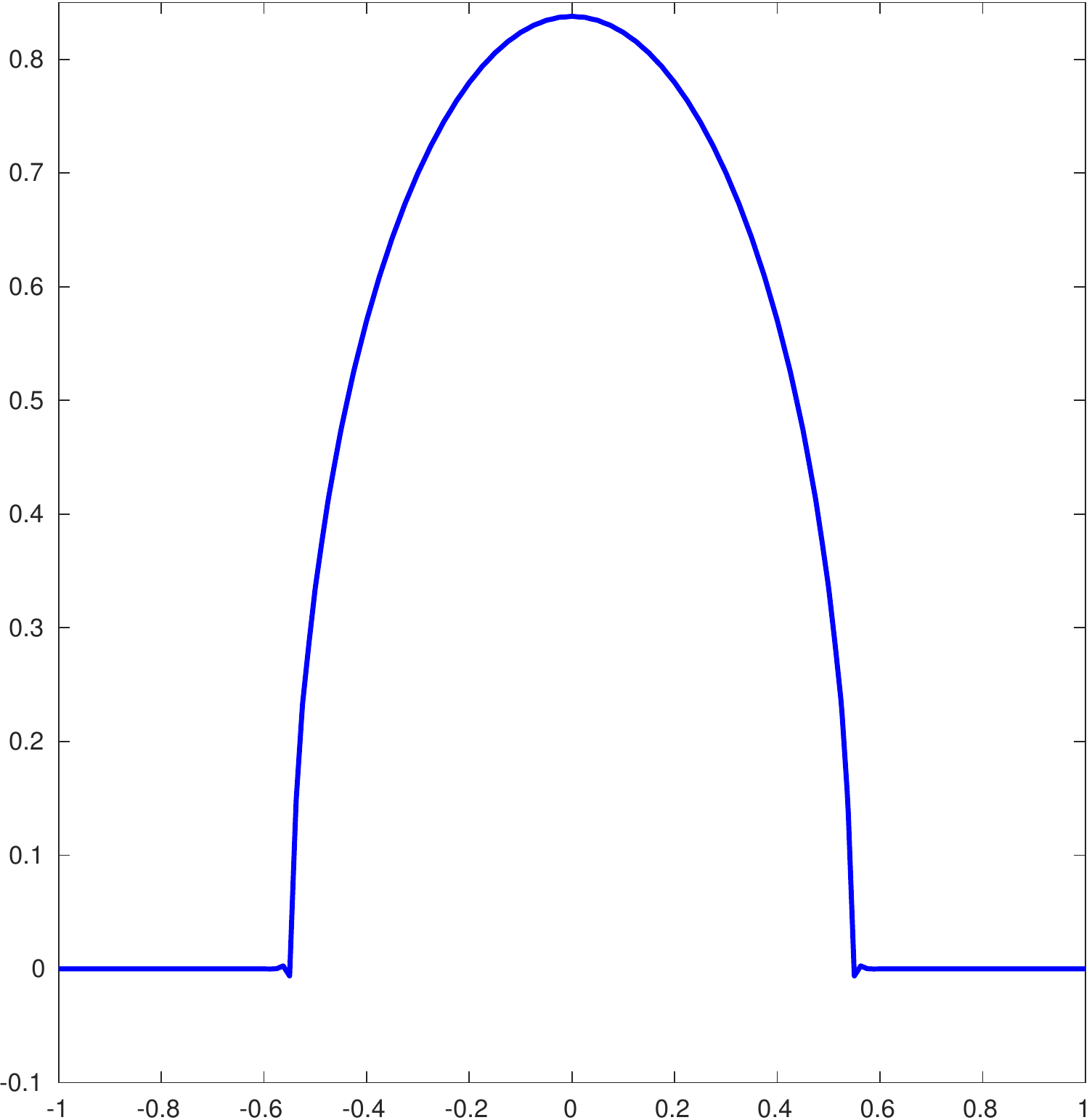}\caption{uniform mesh}\end{subfigure}%
        \hspace{3mm}
        \begin{subfigure}[b]{0.25\linewidth}\includegraphics[scale=0.25]{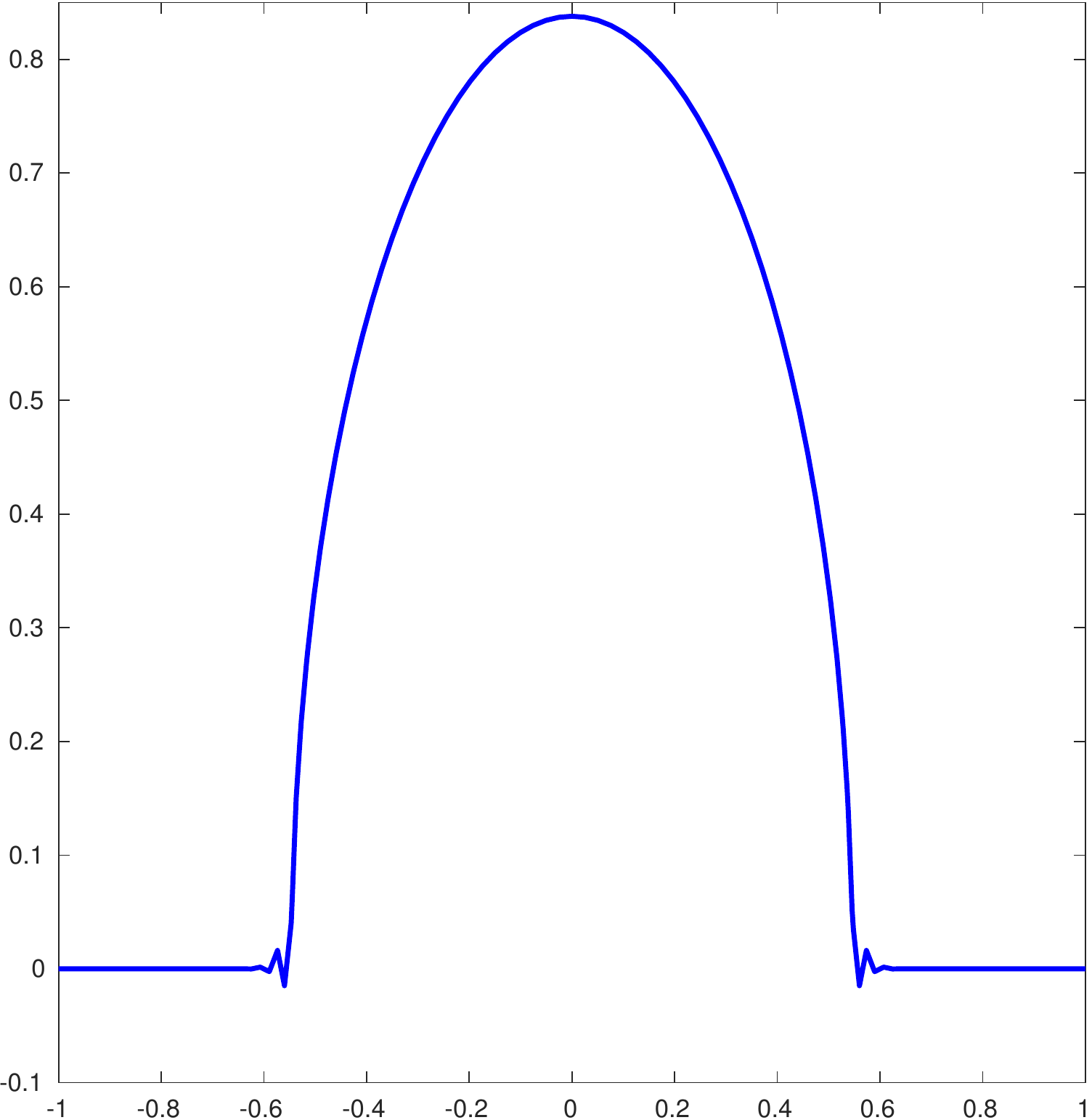}\caption{arclength-based mesh}\end{subfigure}%
         \hspace{3mm}
        \begin{subfigure}[b]{0.25\linewidth}\includegraphics[scale=0.25]{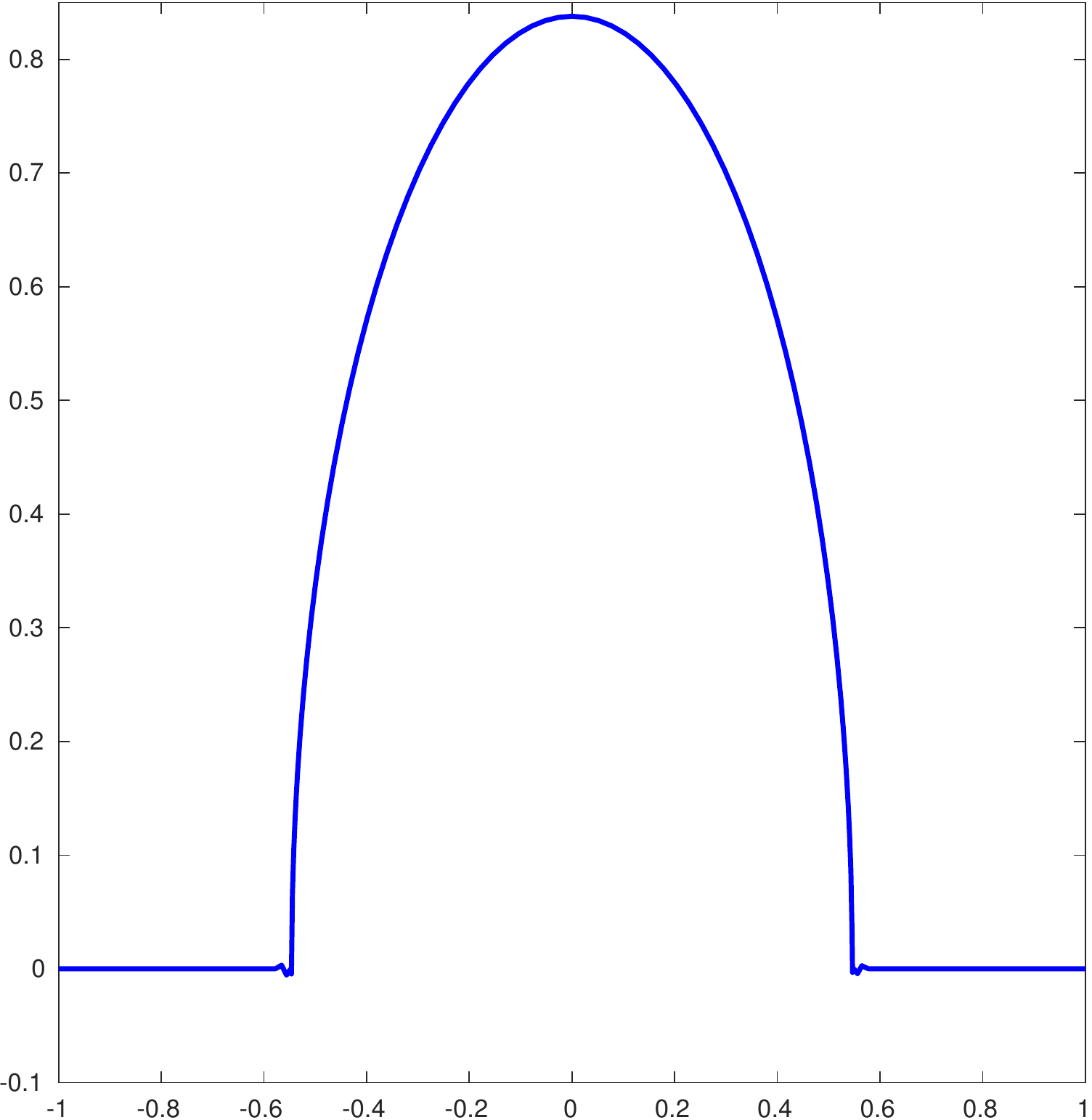}\caption{Hessian-based mesh}\end{subfigure}%
         \caption{Example~\ref{exam4.1} with $m=2$. The cross section at $y = 0$ of the computed solutions at $t=T$ obtained
         with a uniform mesh and arclength- and Hessian-based adaptive meshes  ($N = 102,400$).}
        \label{fig:oscillations-cross-section}
    \end{figure}

\begin{figure}[htp]
        \centering
        \begin{subfigure}[b]{0.35\linewidth}\includegraphics[scale=0.34]{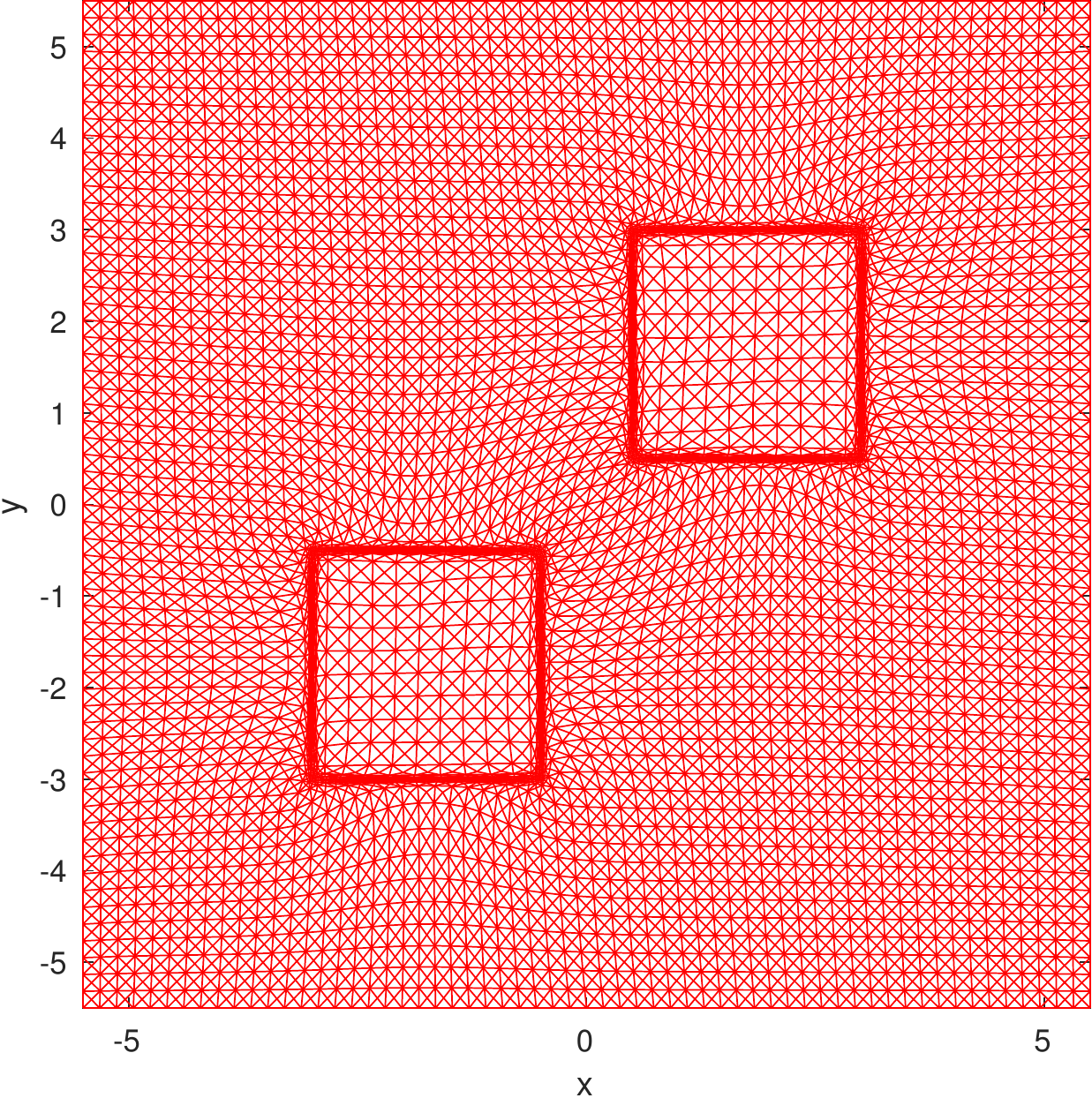}\caption{$t = 0$}\end{subfigure}\hspace{5mm}%
        \begin{subfigure}[b]{0.35\linewidth}\includegraphics[scale=0.32]{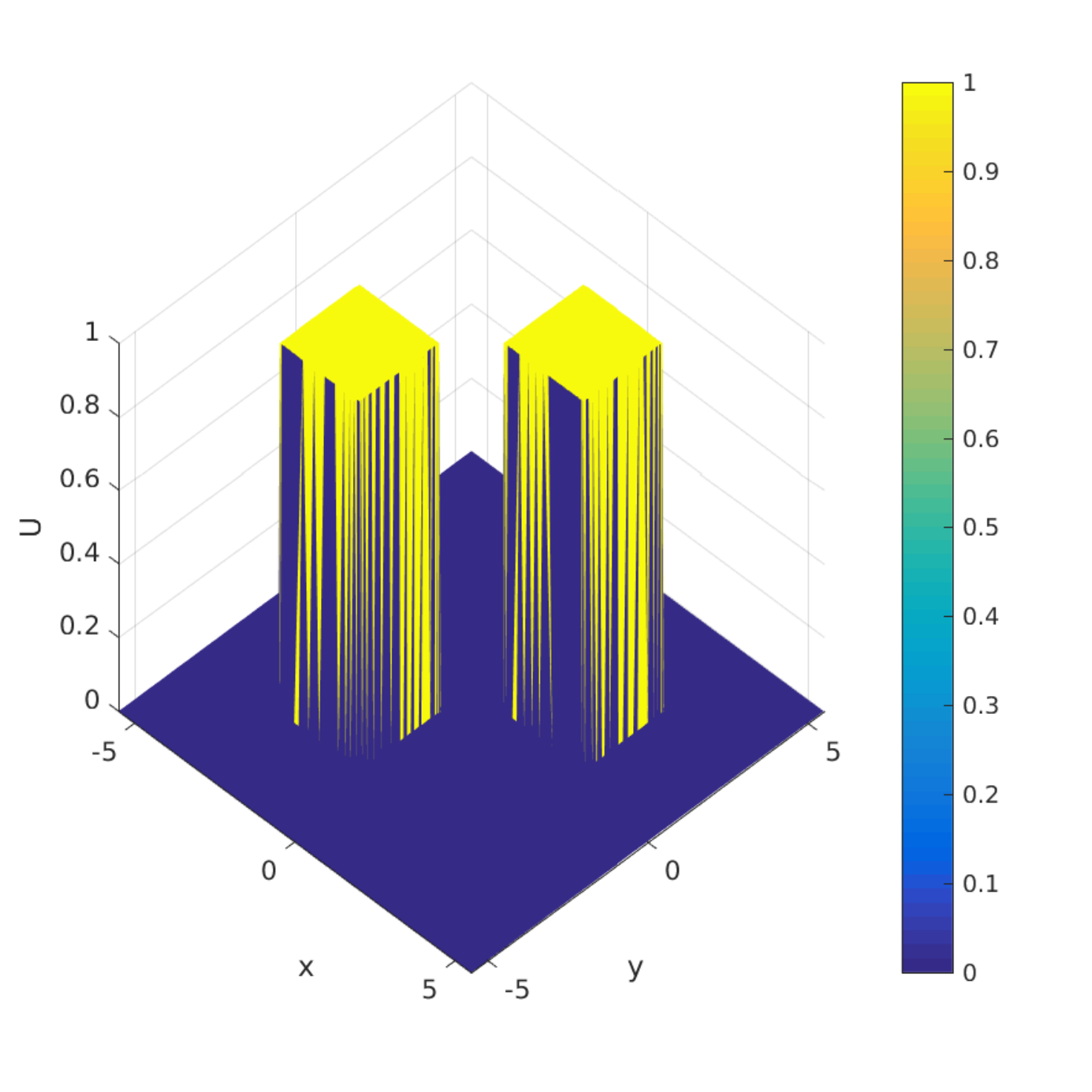}\caption{$t = 0$}\end{subfigure}\\%
        \vspace{-1mm}
        \begin{subfigure}[b]{0.35\linewidth}\includegraphics[scale=0.34]{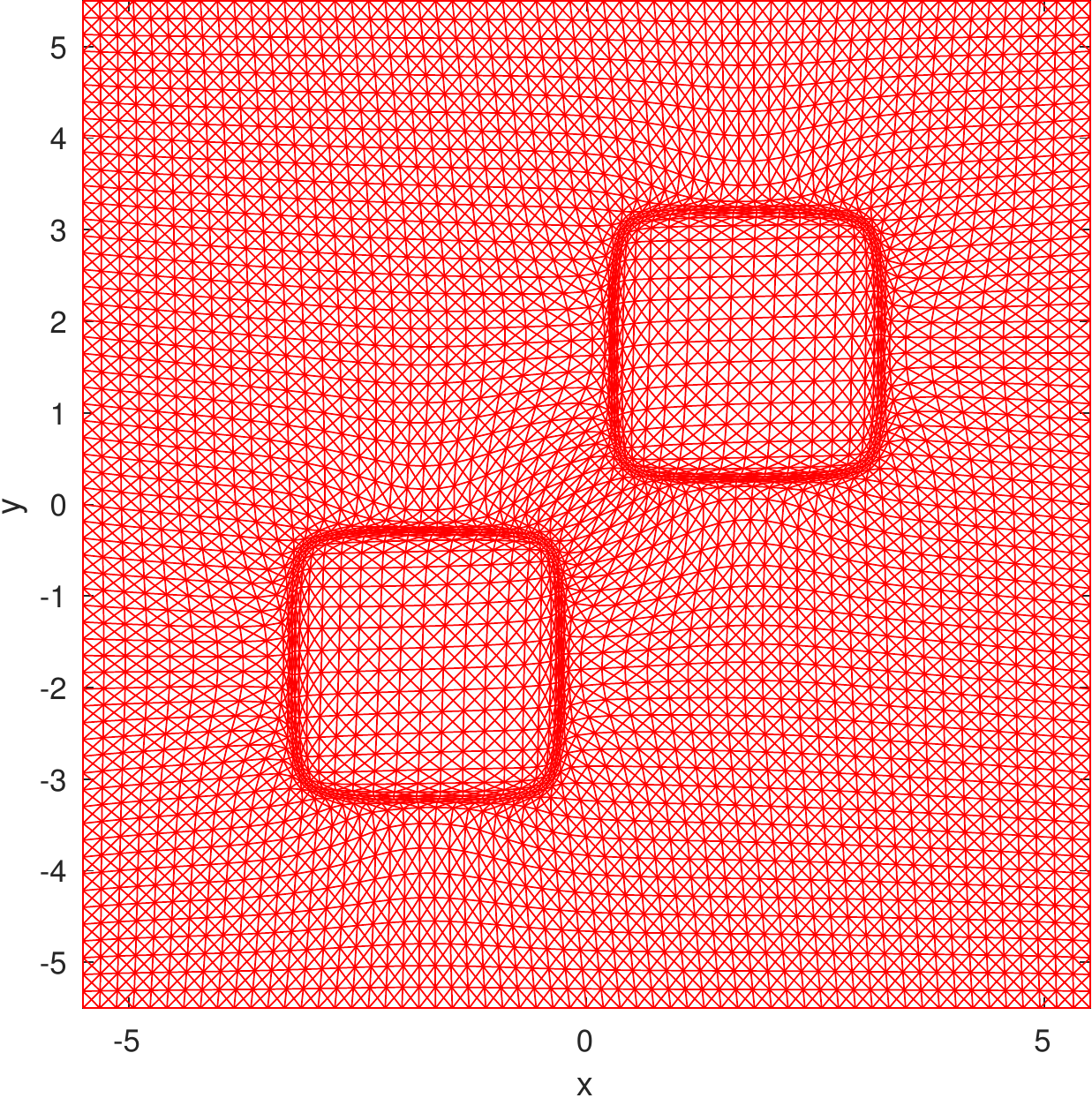}\caption{$t = 0.51$}\end{subfigure}\hspace{5mm}%
        \begin{subfigure}[b]{0.35\linewidth}\includegraphics[scale=0.32]{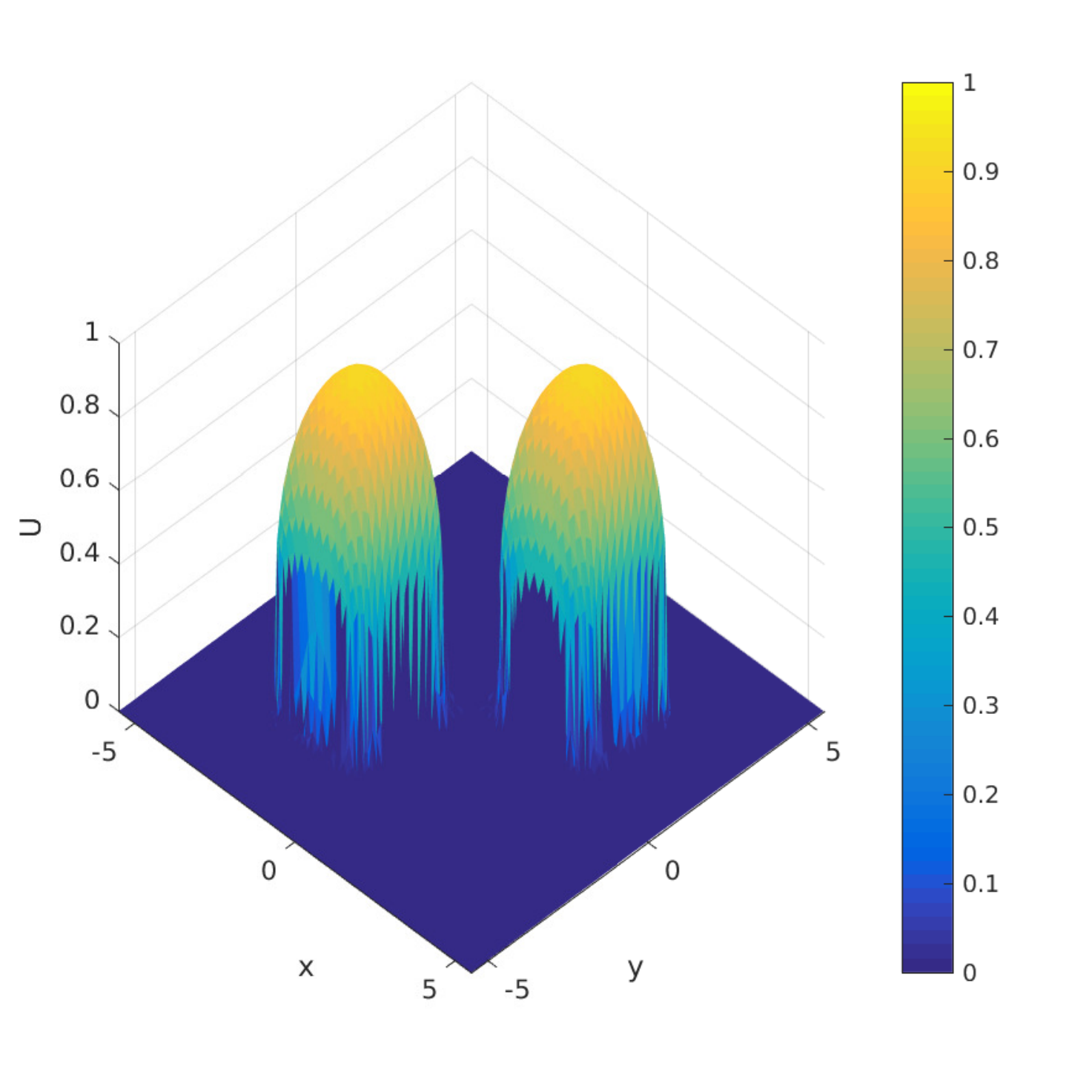}\caption{$t = 0.51$}\end{subfigure}\\%
        \vspace{-1mm}
        \begin{subfigure}[b]{0.35\linewidth}\includegraphics[scale=0.34]{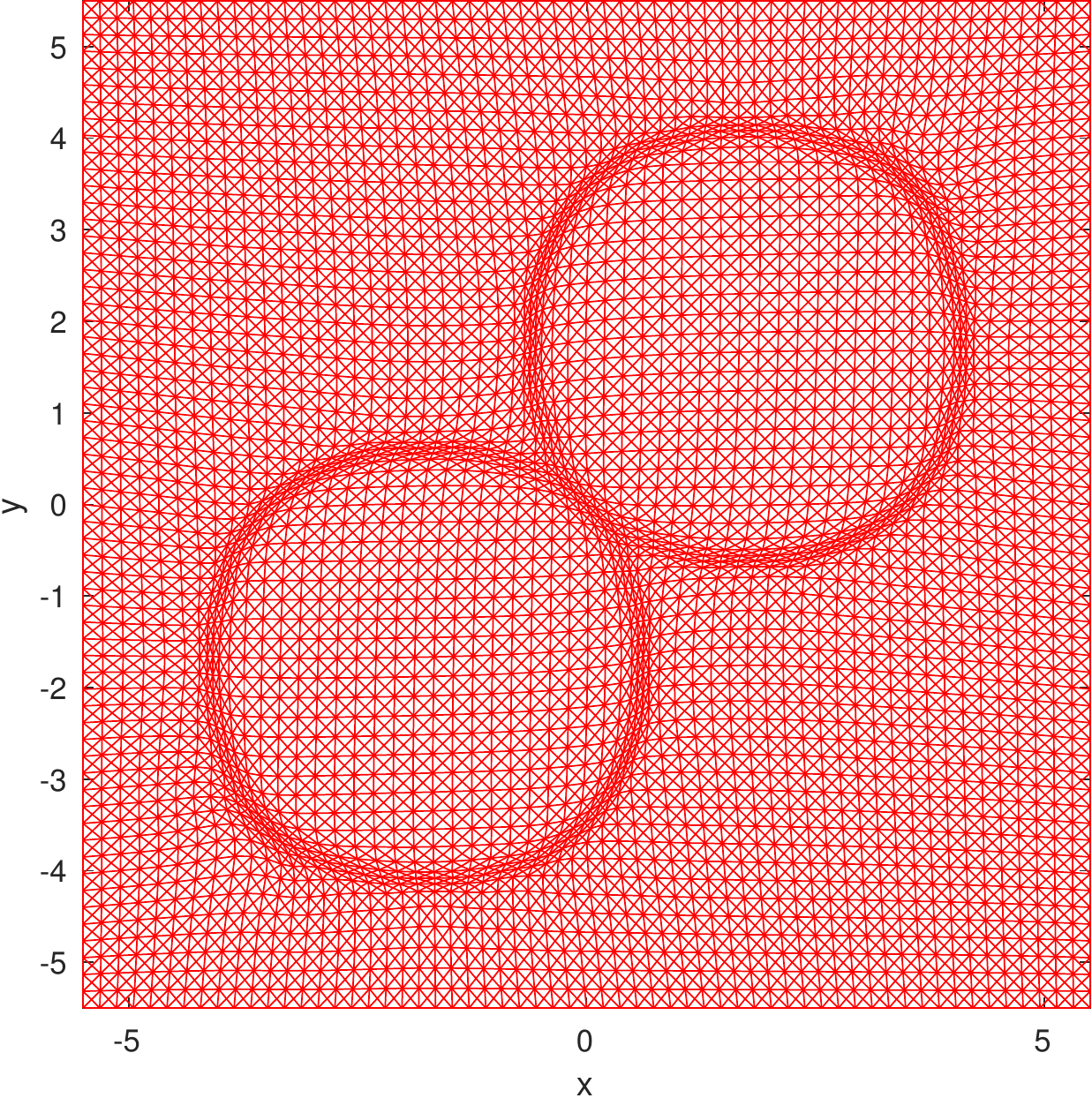}\caption{$t = 100.01$}\end{subfigure}\hspace{5mm}%
        \begin{subfigure}[b]{0.35\linewidth}\includegraphics[scale=0.32]{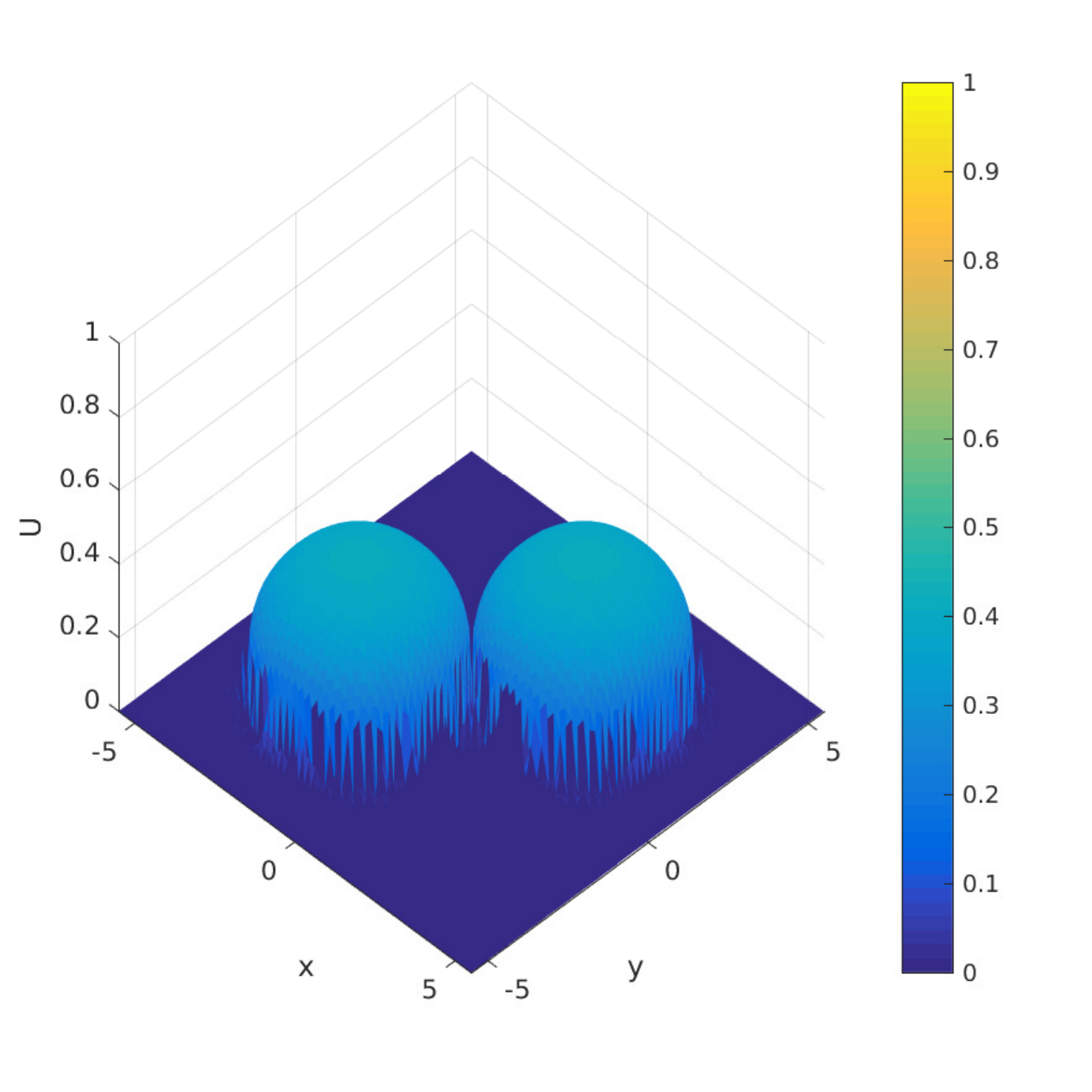}\caption{$t = 100.01$}\end{subfigure}\\%
        \vspace{-1mm}
        \begin{subfigure}[b]{0.35\linewidth}\includegraphics[scale=0.34]{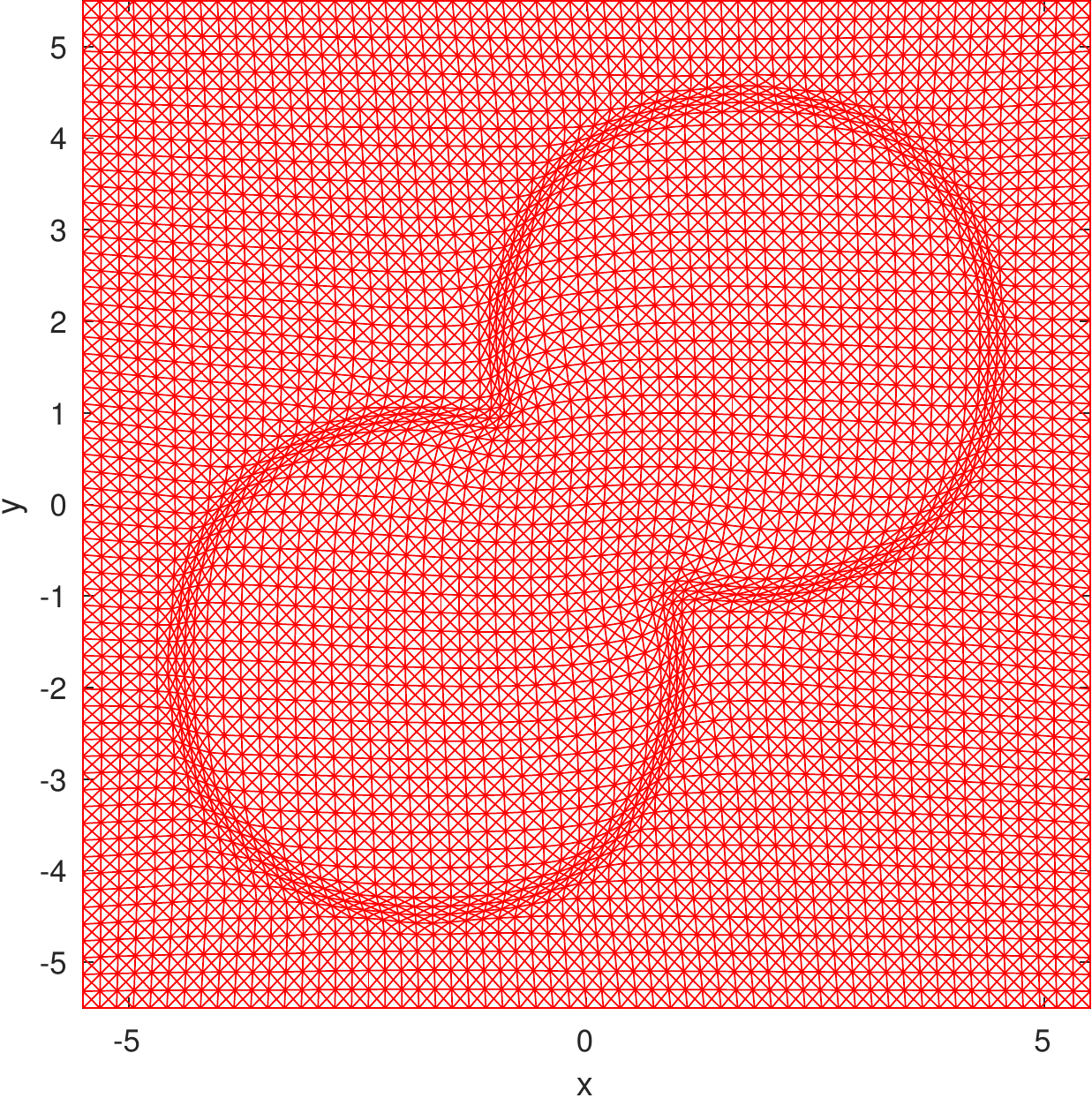}\caption{$t = 500$}\end{subfigure}\hspace{5mm}%
        \begin{subfigure}[b]{0.35\linewidth}\includegraphics[scale=0.32]{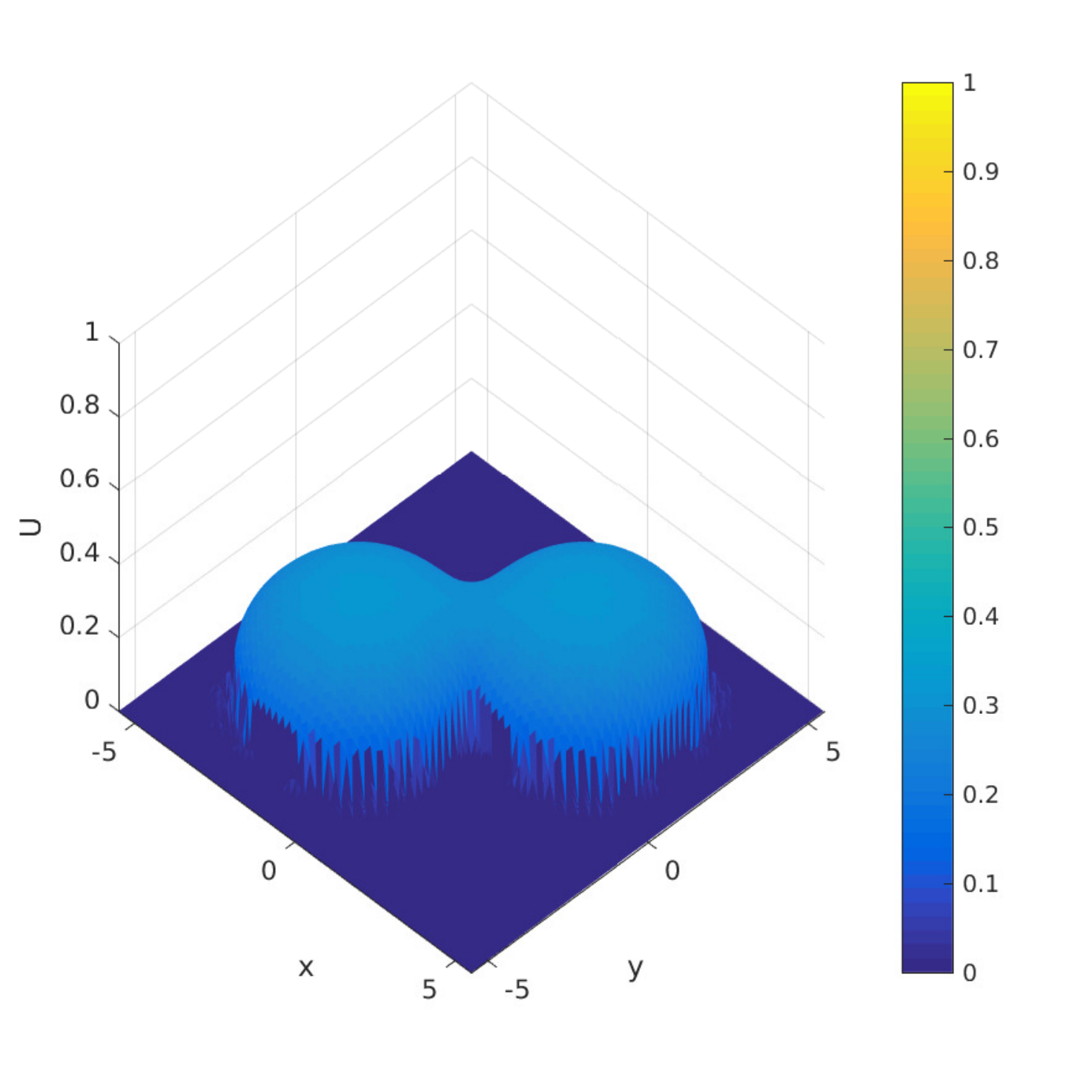}\caption{$t = 500$}\end{subfigure}%
        \caption{Example~\ref{exam4.2}. An adaptive mesh and the corresponding computed solution at various time instants ($N = 14400$).}
        \label{fig:Two-Box-I-soln-mesh}
    \end{figure}

\begin{exam}[Solution with complex support]
\label{exam4.2}

We now consider examples with complex solution support.
The first example models the movement and interaction of two columns of a substance, which start out with the same height.  It has
\[
        m = 5, \quad \Omega = (-5.5,5.5) \times (-5.5,5.5),
\]
\begin{equation}
        u_0(x,y) = 
        \begin{cases}
            1 , & \quad \text{for} \quad (x,y) \in (0.5,3) \times (0.5,3) \\
            1 , & \quad \text{for} \quad (x,y) \in (-3,-0.5) \times (-3,-0.5) \\
            0 , & \quad \text{otherwise} .
        \end{cases}
\label{two-box-1}
\end{equation}
A typical adaptive mesh and the corresponding computed solution obtained
with the Hessian-based mesh adaptation are shown in Fig.~\ref{fig:Two-Box-I-soln-mesh}.
It can be seen that as time evolves, the support of the solution expands from the two boxes,
and then merges into one big region. The mesh adaptation strategy works nicely for the current
example, with the mesh points moving to concentrate around the free boundary.
Particularly, the mesh stays concentrated and nonsingular even during the merging process of
the two separated support regions.
Moreover, the numerical results show that the support of the solution becomes smoother as time evolves,
consistent with the theoretical prediction (e.g., see \cite{Shmarev2005}).
\qed \end{exam}

\begin{figure}[ht]
        \centering
        \begin{subfigure}[b]{0.35\linewidth}\includegraphics[scale=0.34]{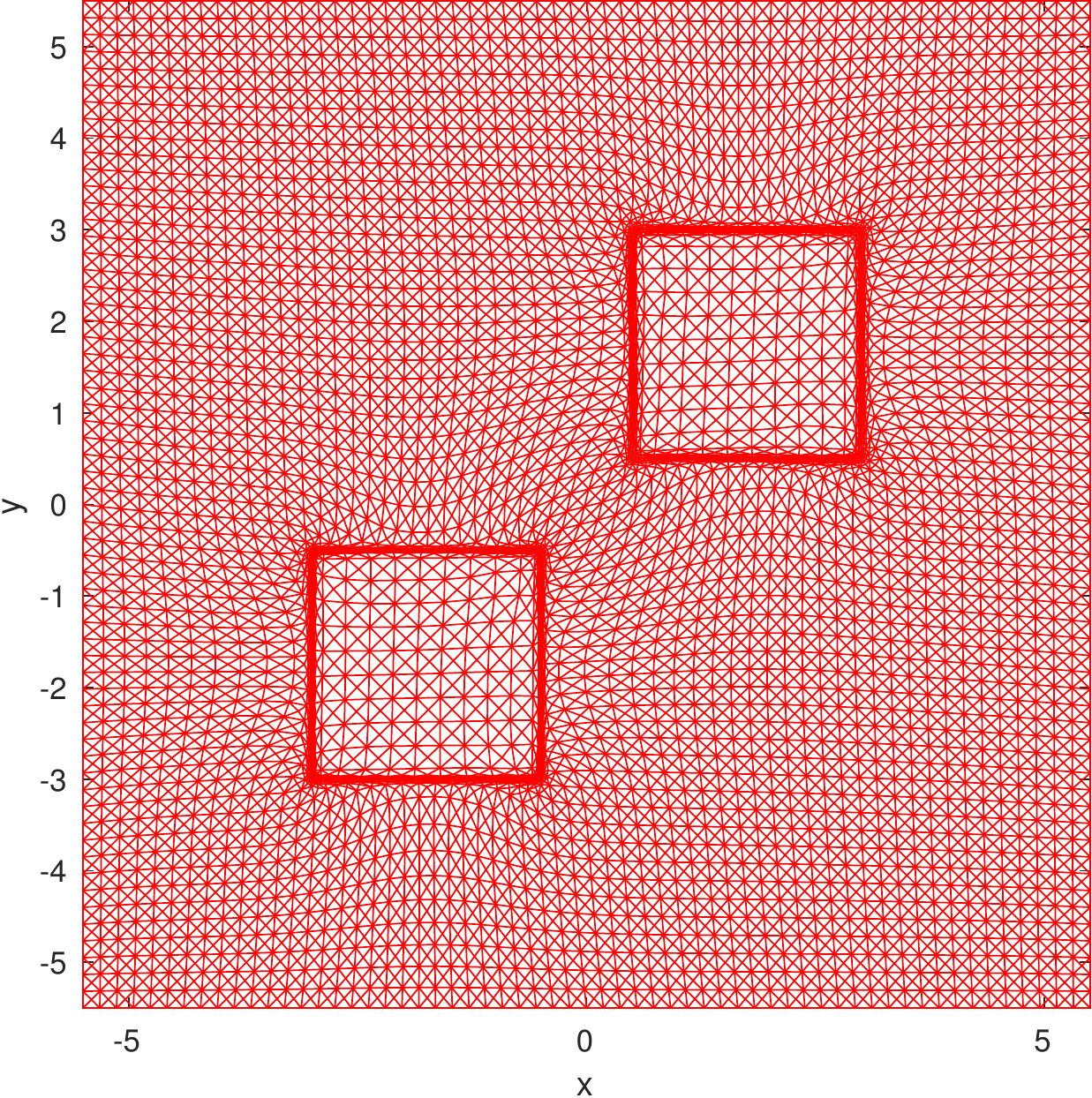}\caption{$t = 0$}\end{subfigure}\hspace{5mm}%
        \begin{subfigure}[b]{0.35\linewidth}\includegraphics[scale=0.32]{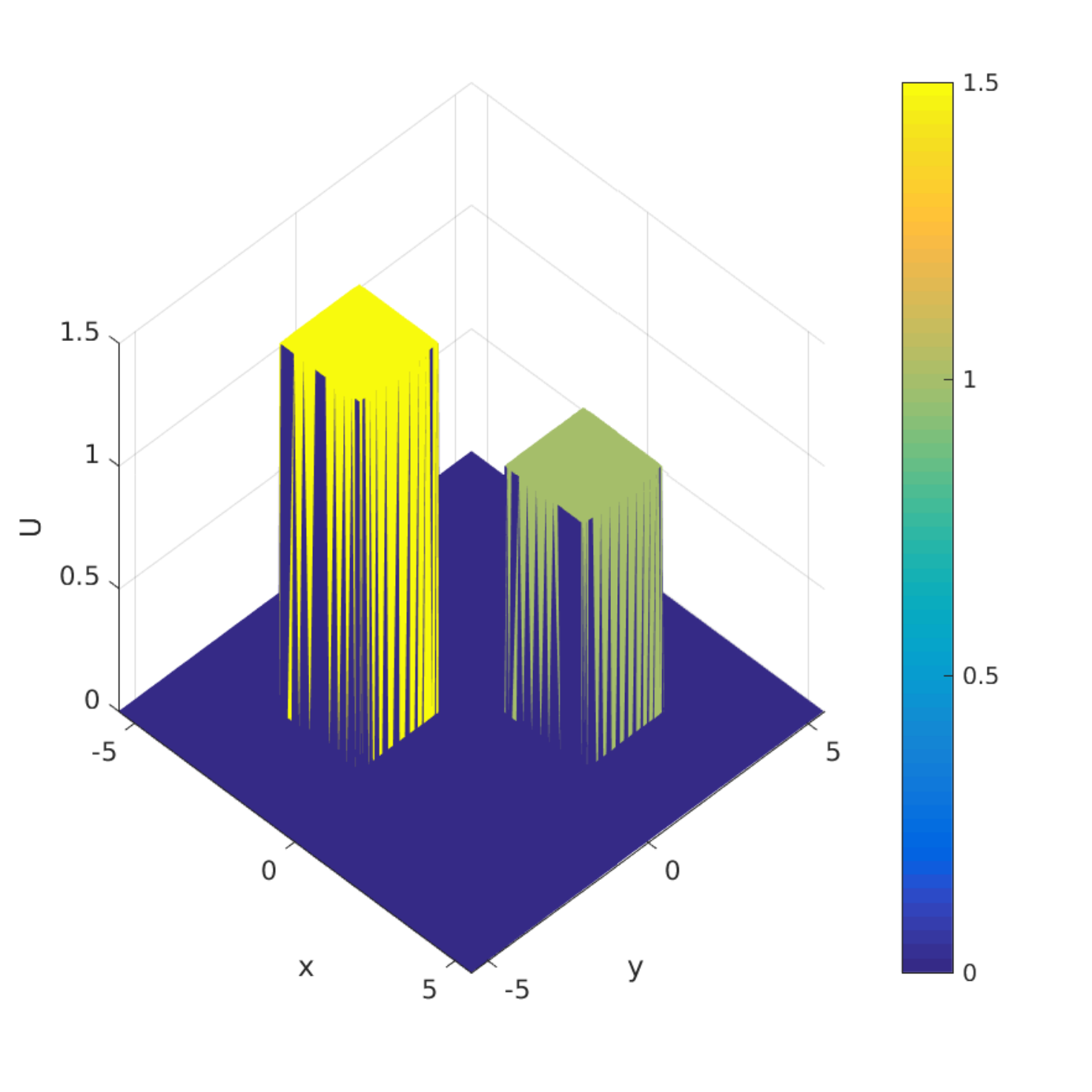}\caption{$t = 0$}\end{subfigure}\\%
        \vspace{-1mm}
        \begin{subfigure}[b]{0.35\linewidth}\includegraphics[scale=0.34]{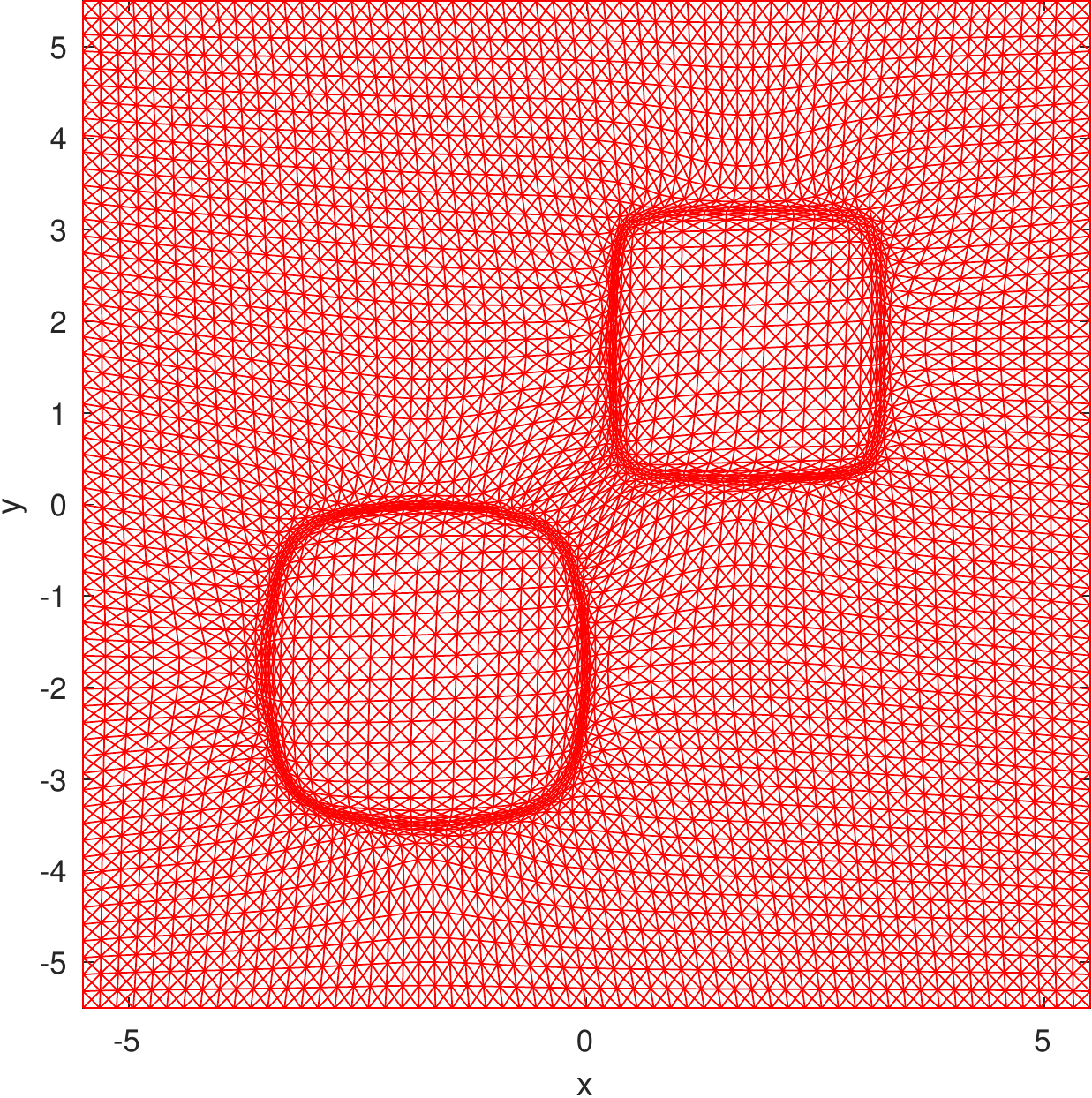}\caption{$t = 0.5$}\end{subfigure}\hspace{5mm}%
        \begin{subfigure}[b]{0.35\linewidth}\includegraphics[scale=0.32]{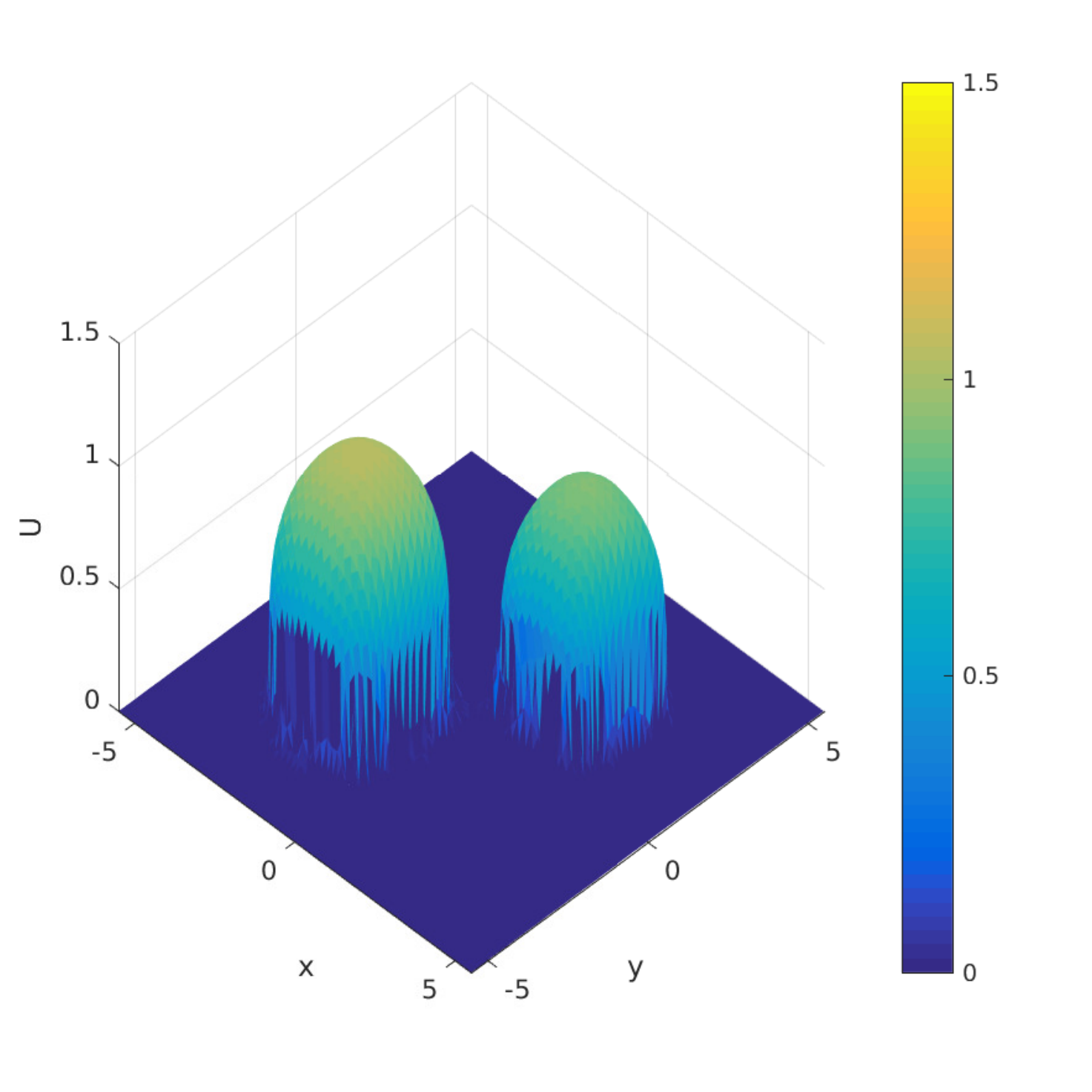}\caption{$t = 0.5$}\end{subfigure}\\%
        \vspace{-1mm}
        \begin{subfigure}[b]{0.35\linewidth}\includegraphics[scale=0.34]{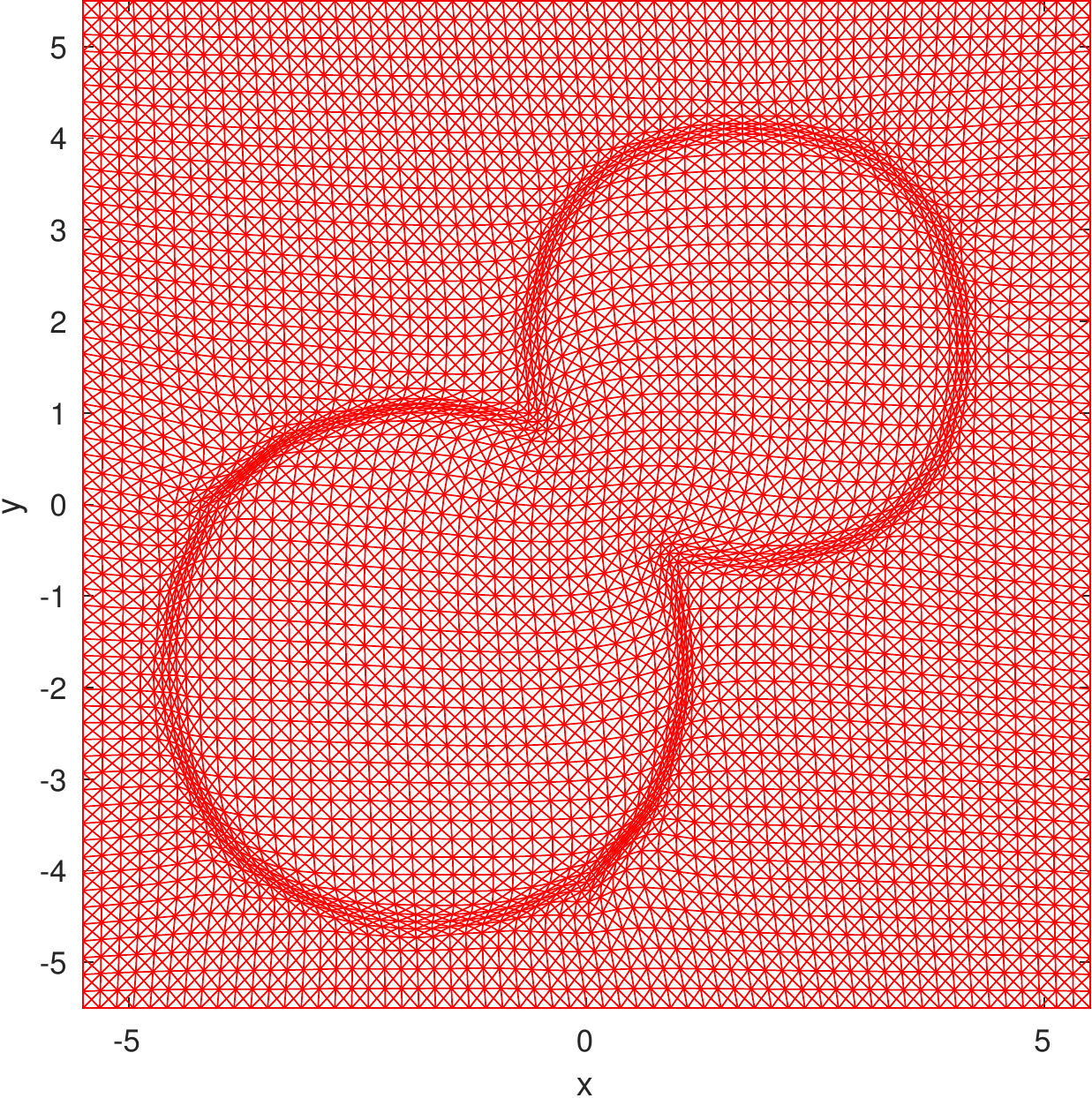}\caption{$t = 100$}\end{subfigure}\hspace{5mm}%
        \begin{subfigure}[b]{0.35\linewidth}\includegraphics[scale=0.32]{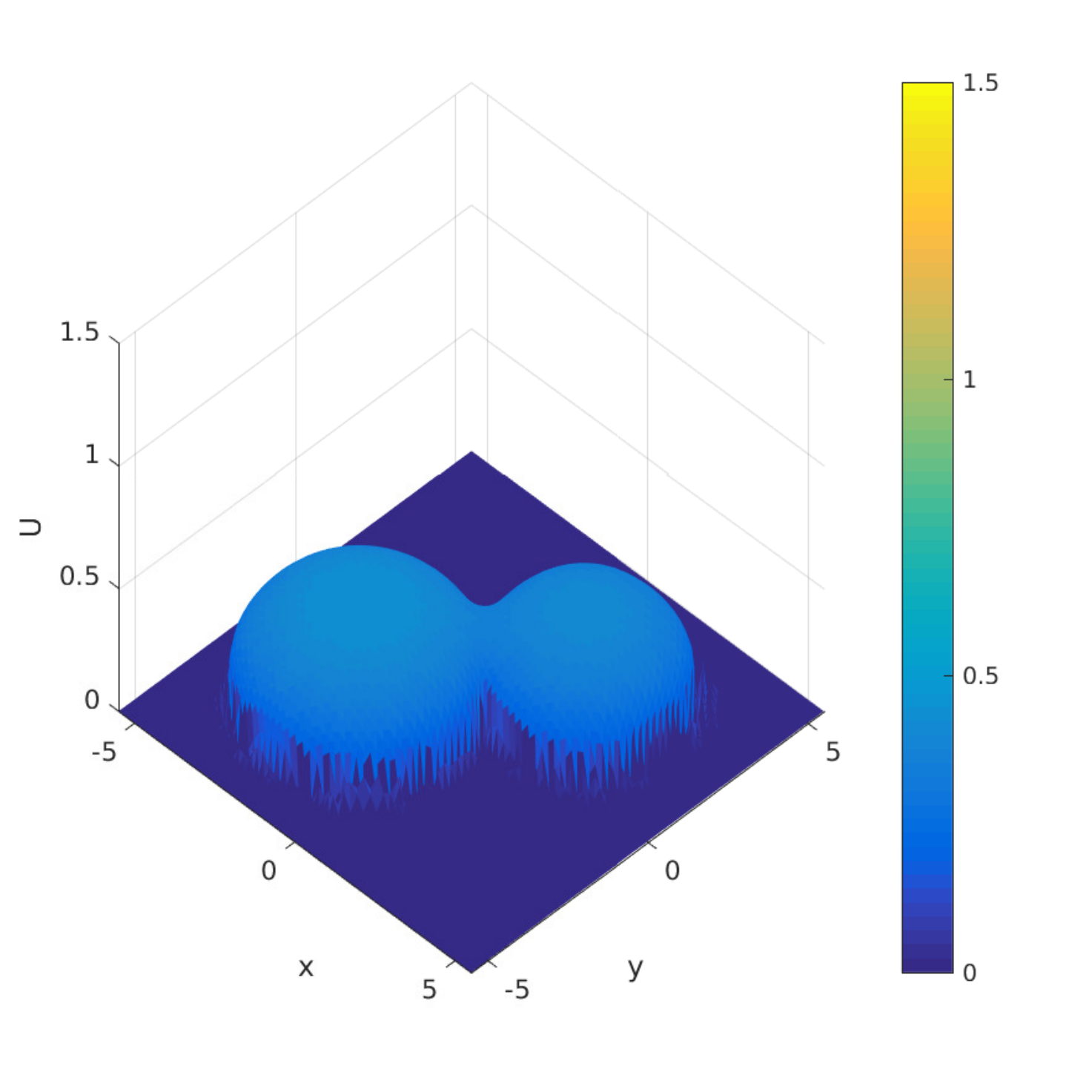}\caption{$t = 100$}\end{subfigure}\\%
        \vspace{-1mm}
        \begin{subfigure}[b]{0.35\linewidth}\includegraphics[scale=0.34]{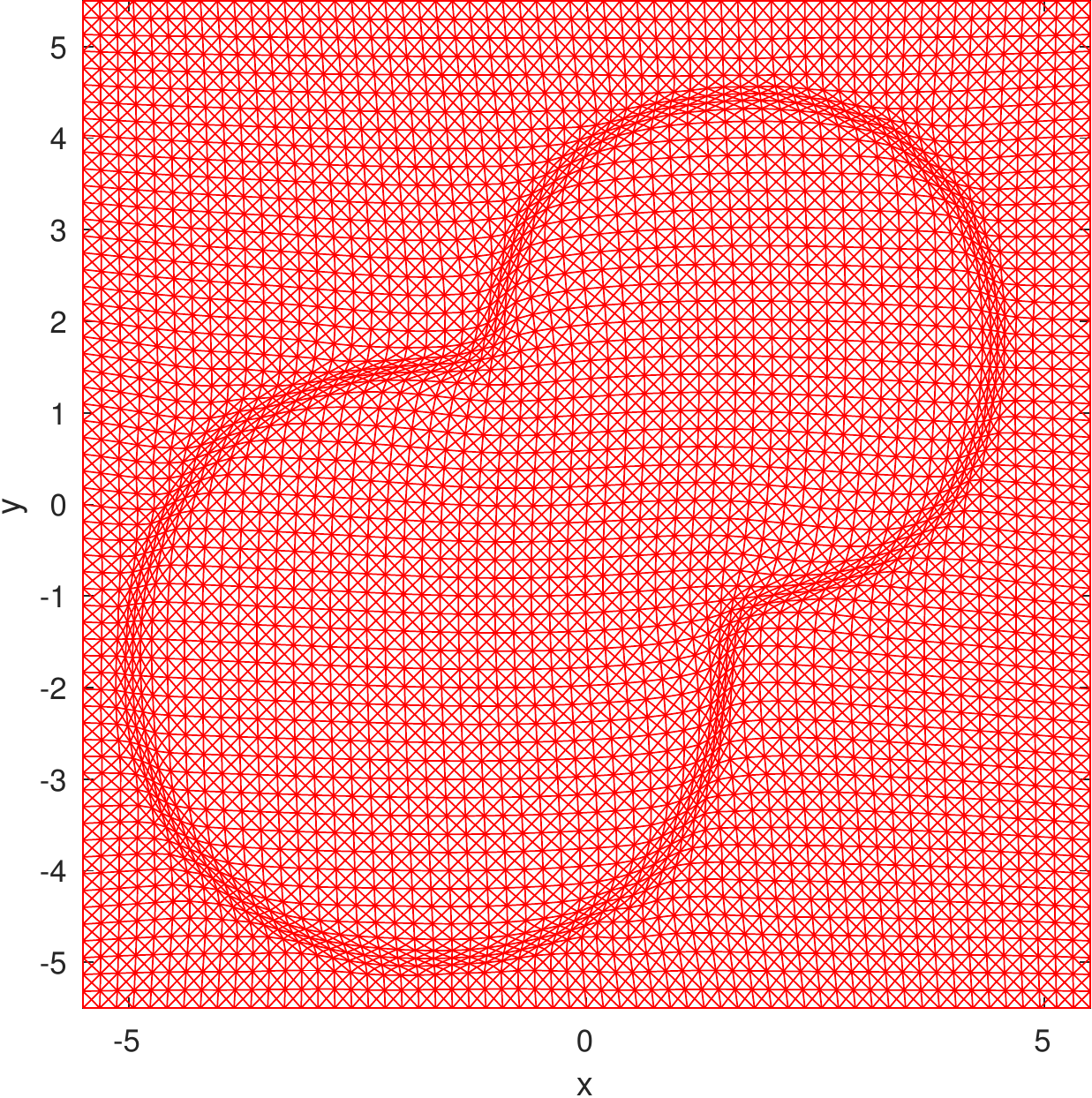}\caption{$t = 500$}\end{subfigure}\hspace{5mm}%
        \begin{subfigure}[b]{0.35\linewidth}\includegraphics[scale=0.32]{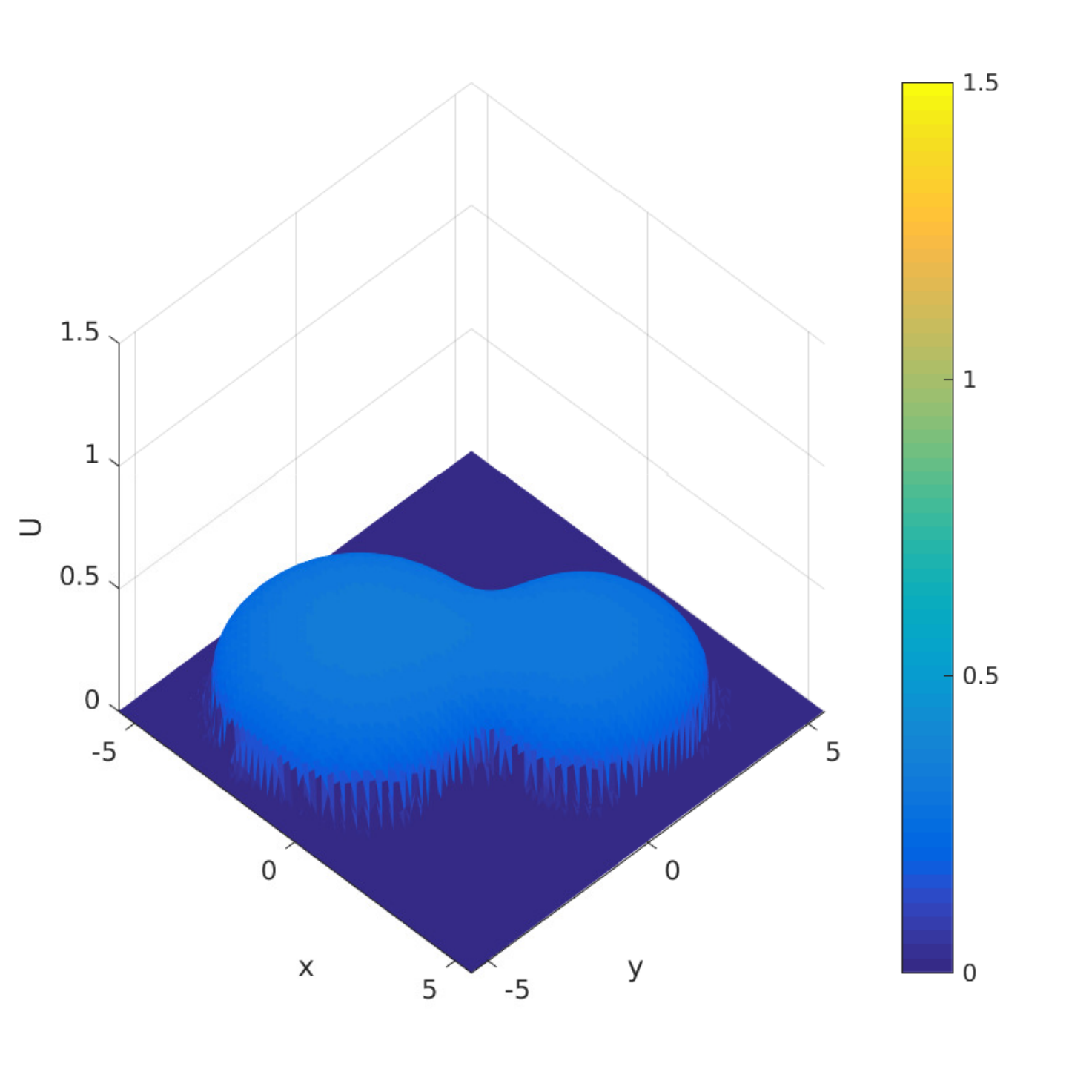}\caption{$t = 500$}\end{subfigure}%
        \caption{Example~\ref{exam4.3}. An adaptive mesh and the corresponding computed solution at various time instants ($N = 14400$).}
        \label{fig:Two-Box-II-soln-mesh}
    \end{figure}

\begin{exam}[Solution with complex support]
\label{exam4.3}

The next example is similar to the previous one except that the initial solution has different heights in the two boxes, 
    \begin{equation}
        u_0(x,y) = 
        \begin{cases}
            1 , & \quad \text{for} \quad (x,y) \in (0.5,3) \times (0.5,3) \\
            1.5 , & \quad \text{for} \quad (x,y) \in (-3,-0.5) \times (-3,-0.5) \\
            0 , & \quad \text{otherwise} .
        \end{cases}
        \label{two-box-2}
    \end{equation}
A typical adaptive mesh and the corresponding solution are shown in Fig.~\ref{fig:Two-Box-II-soln-mesh}.
Once again, the mesh is concentrated correctly around the free boundaries
as they evolve with time. Moreover, the region with larger initial solution values
expands faster than the region with smaller values. Overall, the support of the solution
for this example expands faster than that of the previous example. At $t=50$, the two boxes have already
merged into a single region of calabash shape.
\qed \end{exam}

\begin{figure}[ht]
        \centering
        \begin{subfigure}[b]{0.25\linewidth}\includegraphics[scale=0.20]{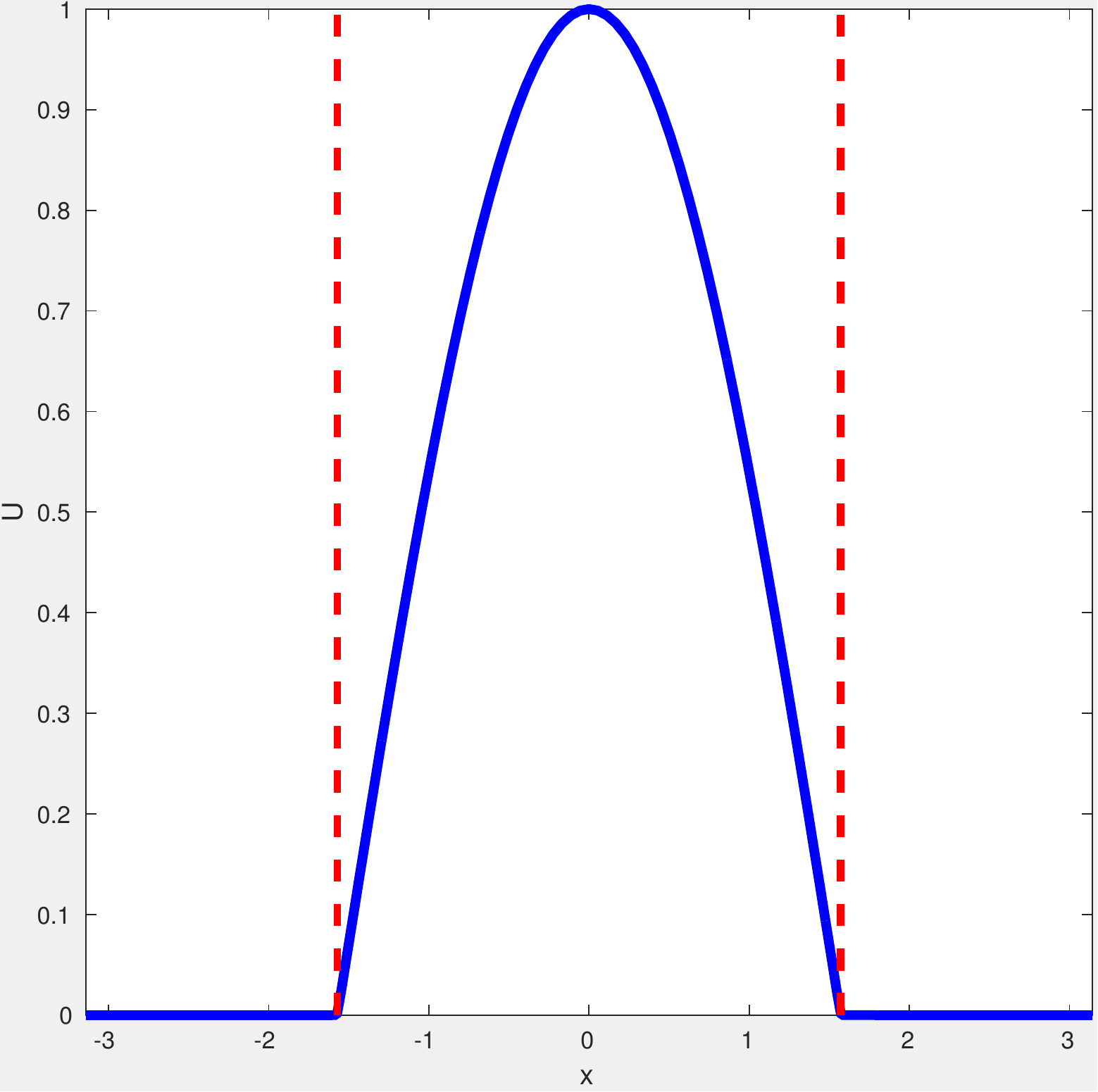}\caption{$t = 0$}\end{subfigure}%
        \begin{subfigure}[b]{0.25\linewidth}\includegraphics[scale=0.20]{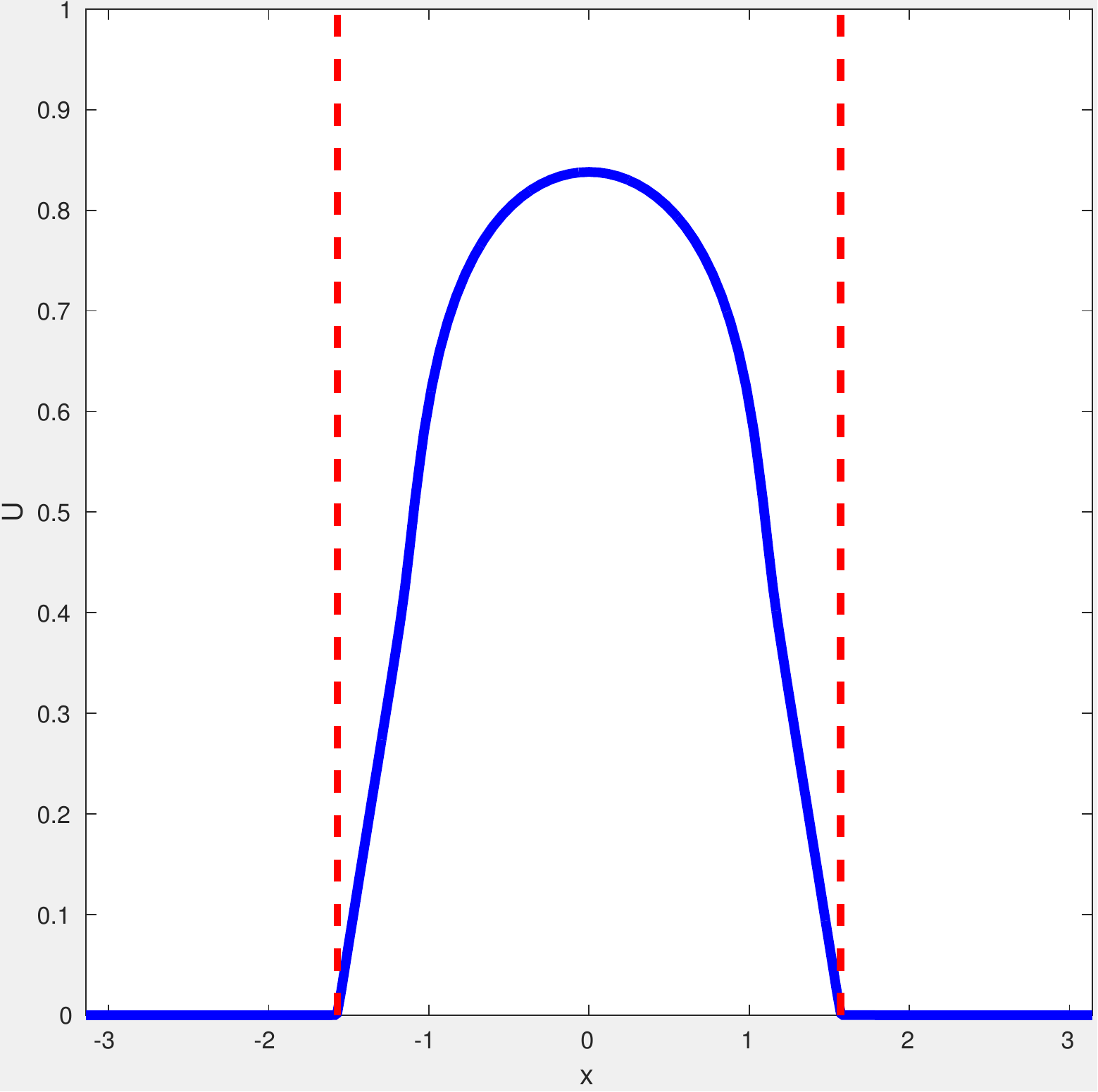}\caption{$t = 0.05$}\end{subfigure}%
        \begin{subfigure}[b]{0.25\linewidth}\includegraphics[scale=0.20]{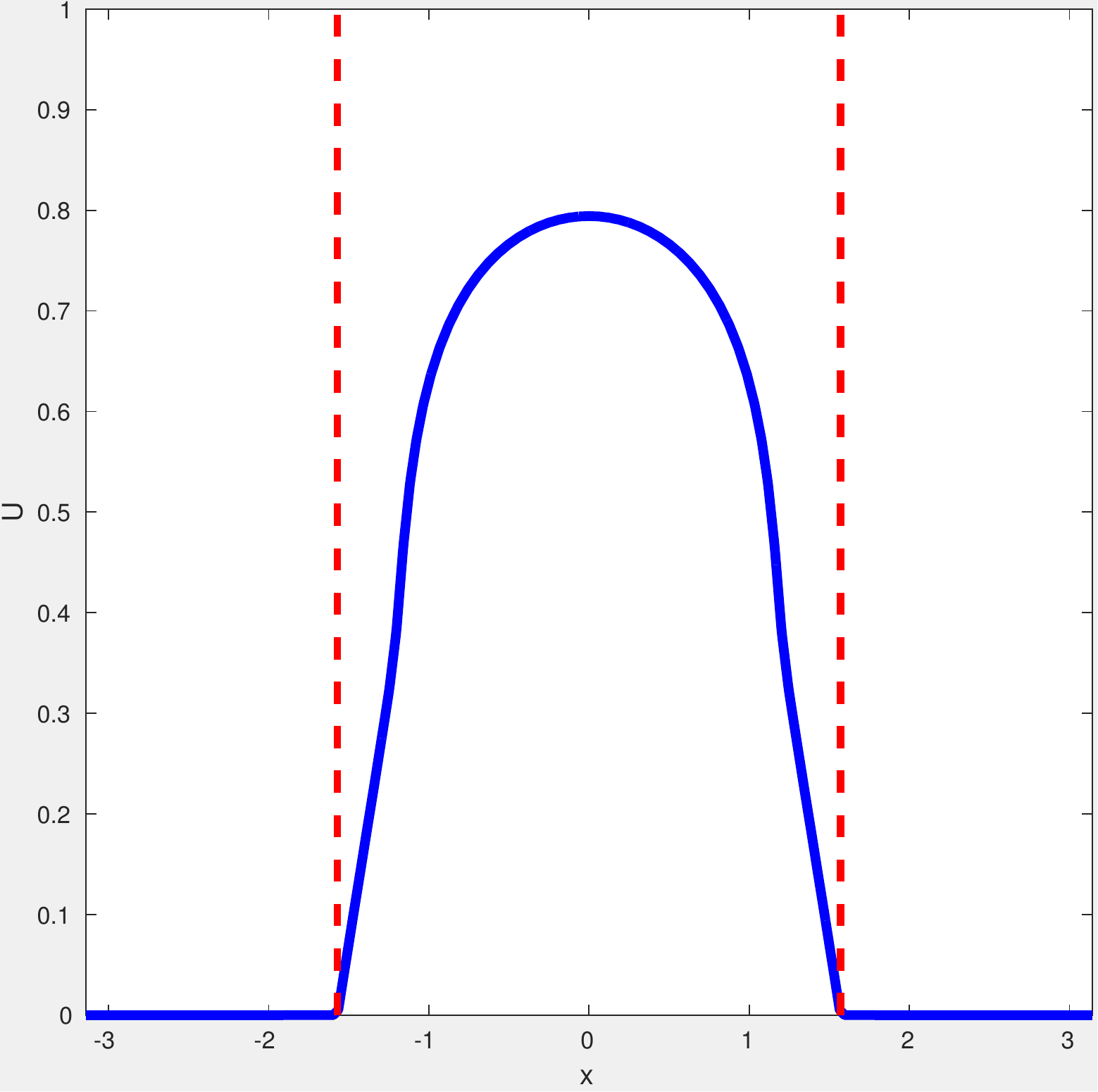}\caption{$t = 0.1$}\end{subfigure}%
        \begin{subfigure}[b]{0.25\linewidth}\includegraphics[scale=0.20]{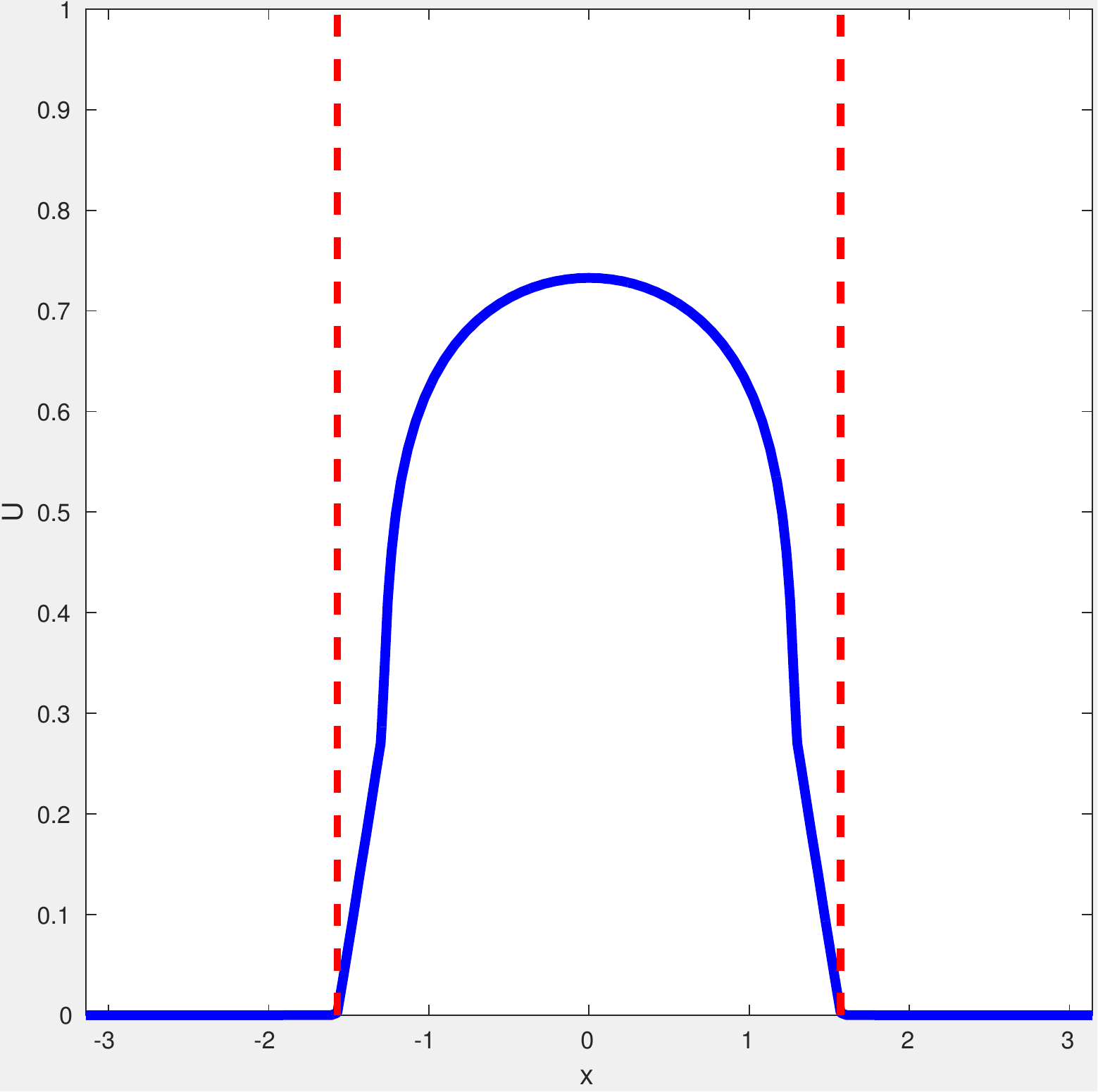}\caption{$t = 0.25$}\end{subfigure}\\%
        \begin{subfigure}[b]{0.25\linewidth}\includegraphics[scale=0.20]{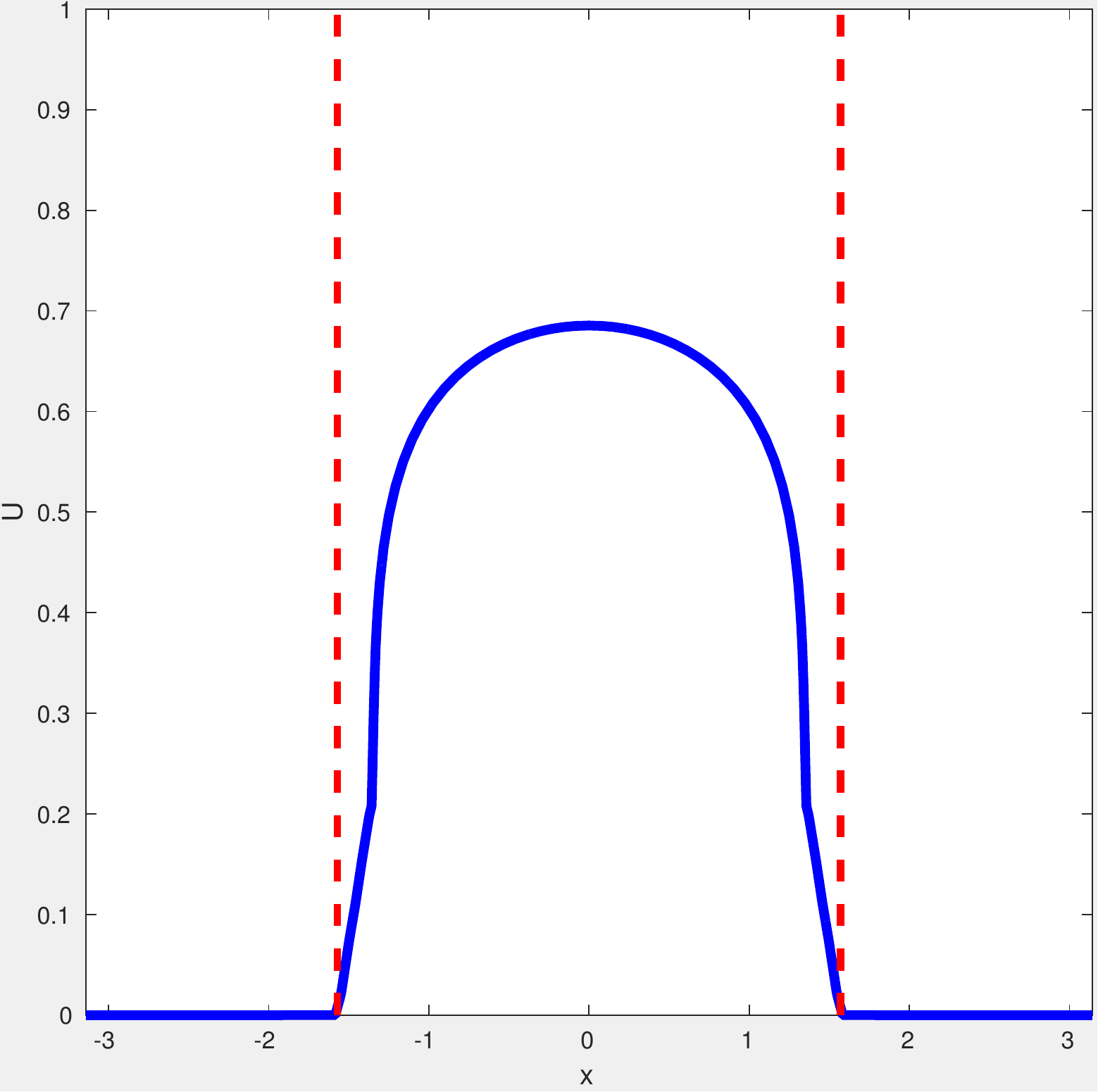}\caption{$t = 0.5$}\end{subfigure}%
        \begin{subfigure}[b]{0.25\linewidth}\includegraphics[scale=0.20]{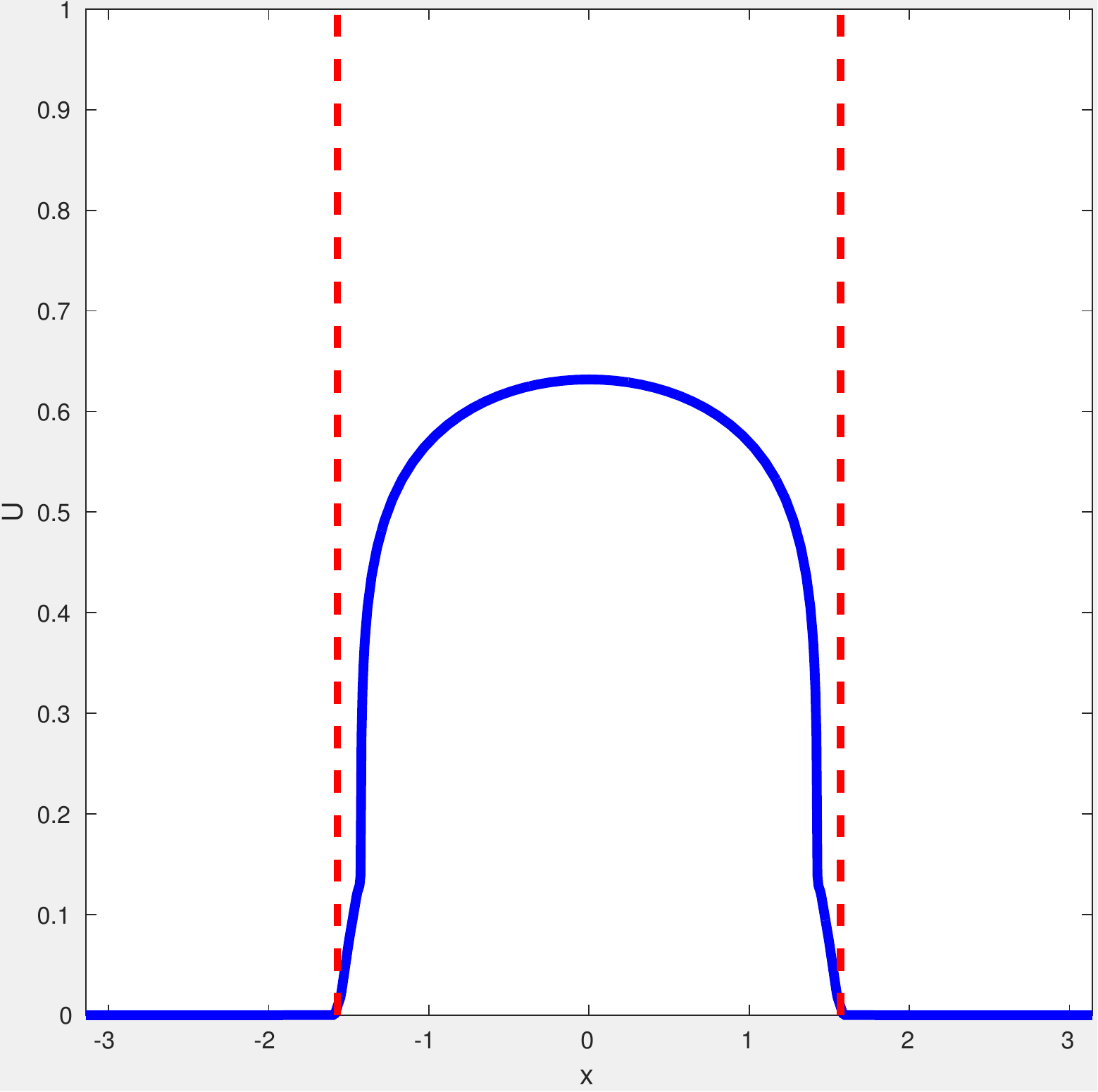}\caption{$t = 1.1$}\end{subfigure}%
        \begin{subfigure}[b]{0.25\linewidth}\includegraphics[scale=0.20]{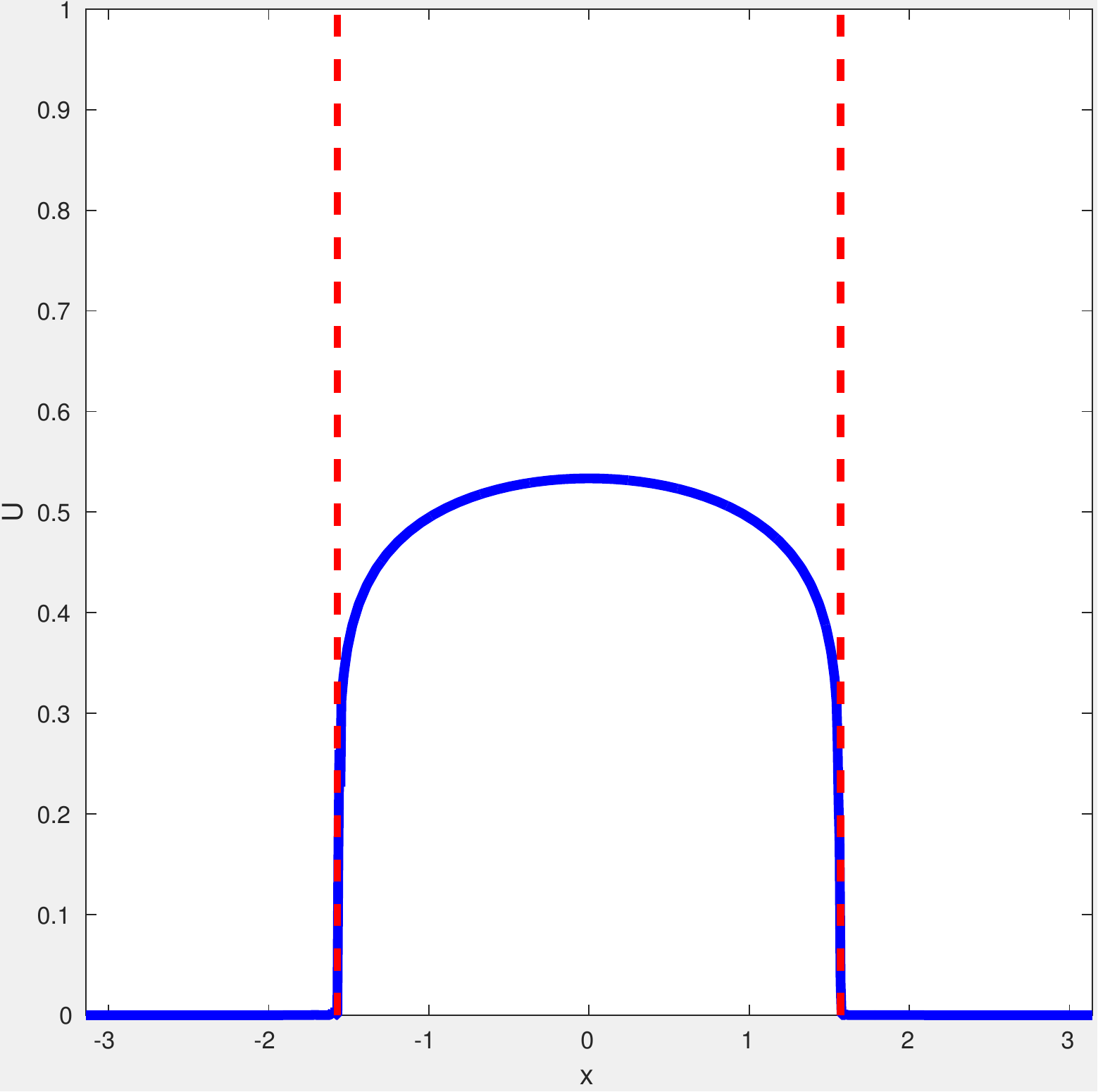}\caption{$t = 5$}\end{subfigure}%
        \begin{subfigure}[b]{0.25\linewidth}\includegraphics[scale=0.20]{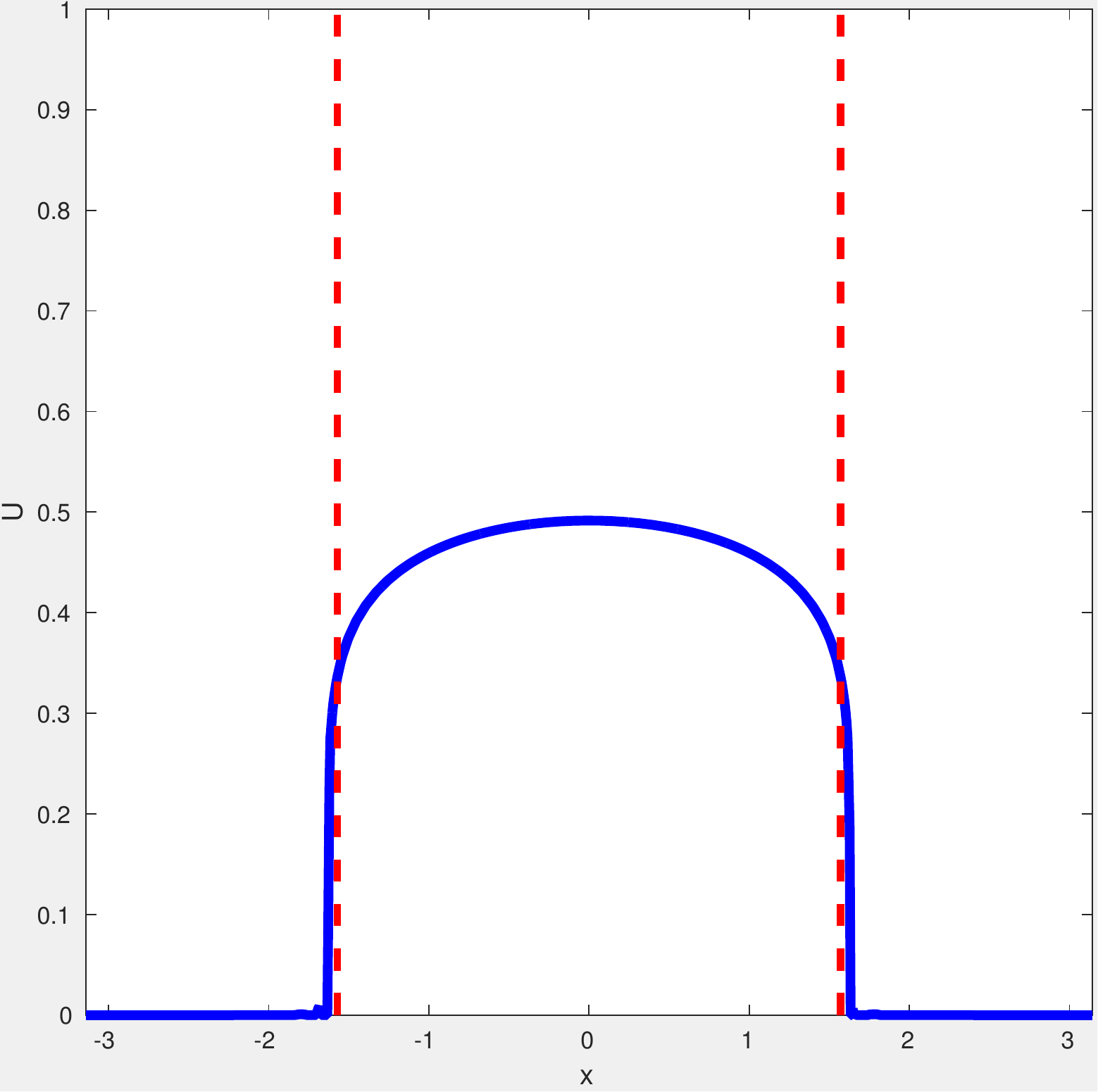}\caption{$t = 10$}\end{subfigure}\\%
        \begin{subfigure}[b]{0.25\linewidth}\includegraphics[scale=0.20]{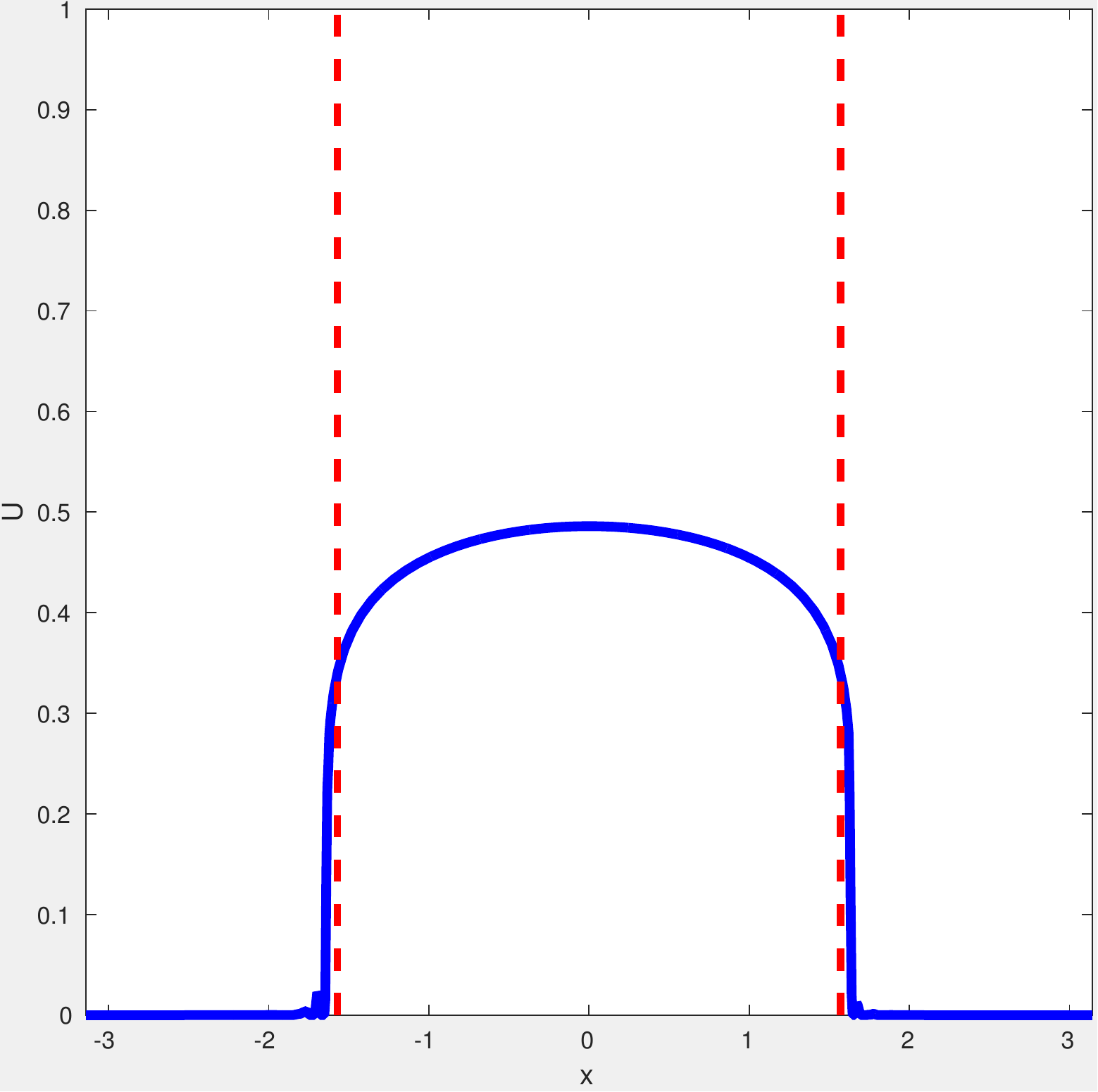}\caption{$t = 11$}\end{subfigure}%
        \begin{subfigure}[b]{0.25\linewidth}\includegraphics[scale=0.20]{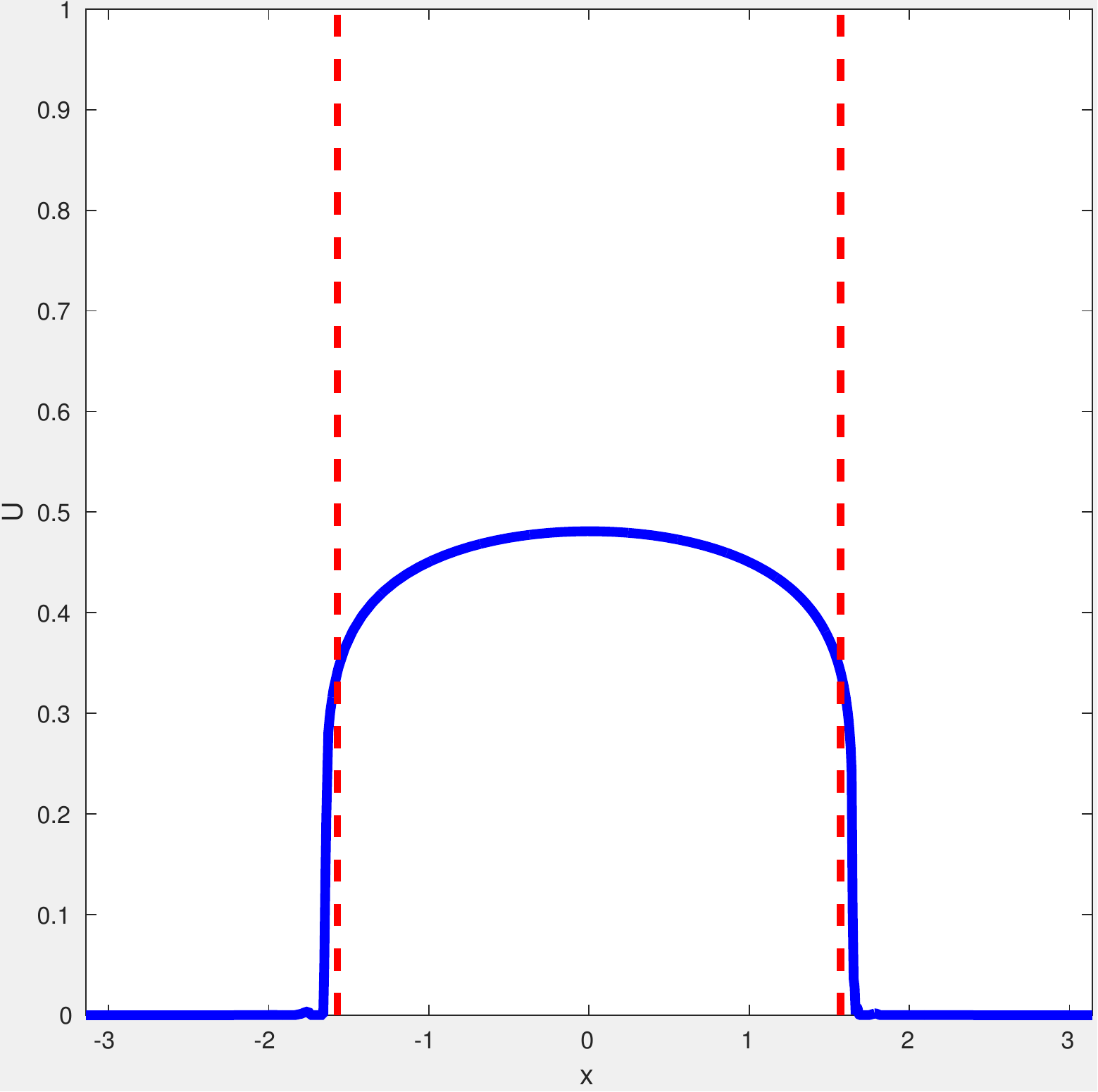}\caption{$t = 12$}\end{subfigure}%
        \begin{subfigure}[b]{0.25\linewidth}\includegraphics[scale=0.20]{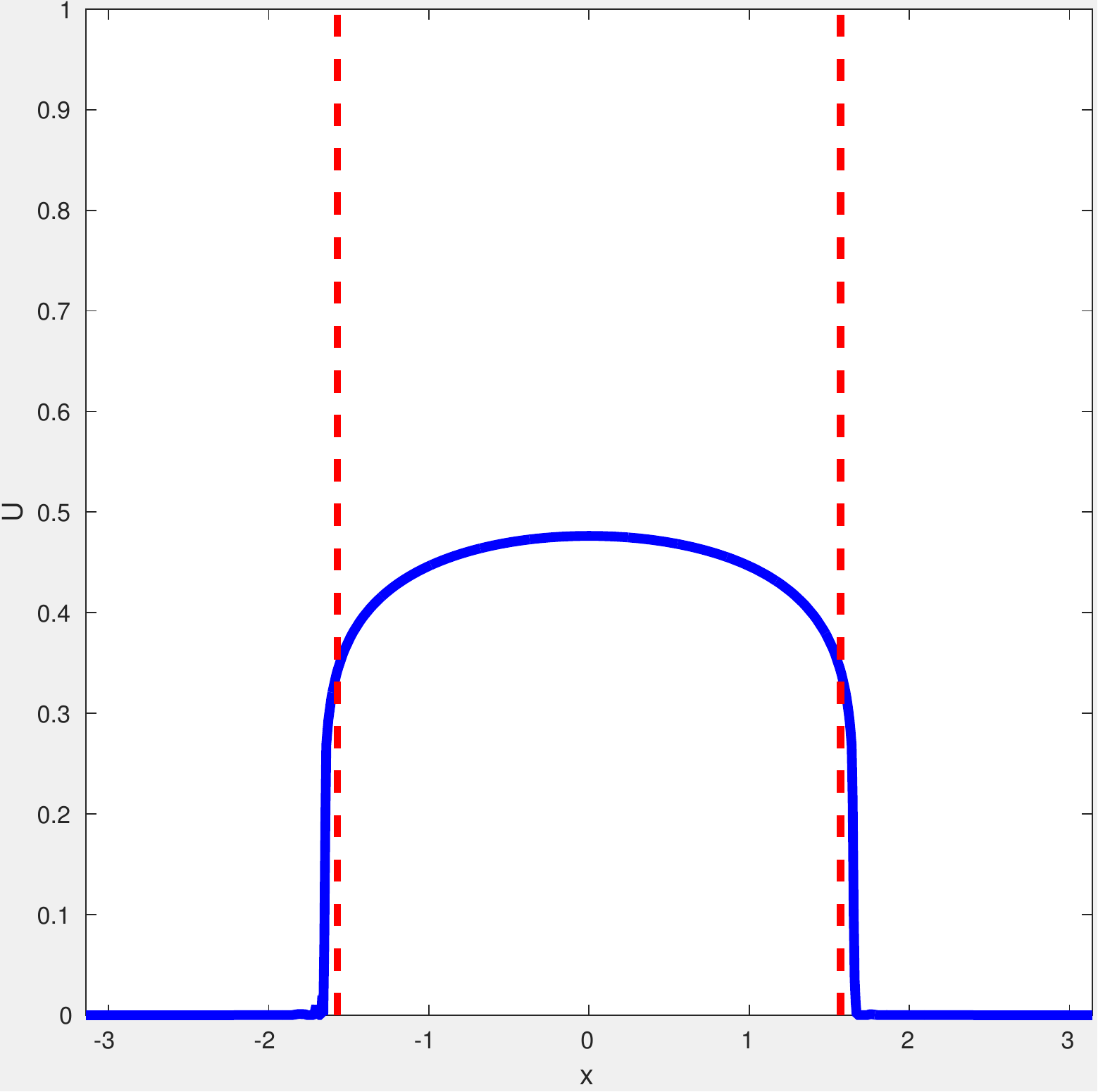}\caption{$t = 13$}\end{subfigure}%
        \begin{subfigure}[b]{0.25\linewidth}\includegraphics[scale=0.20]{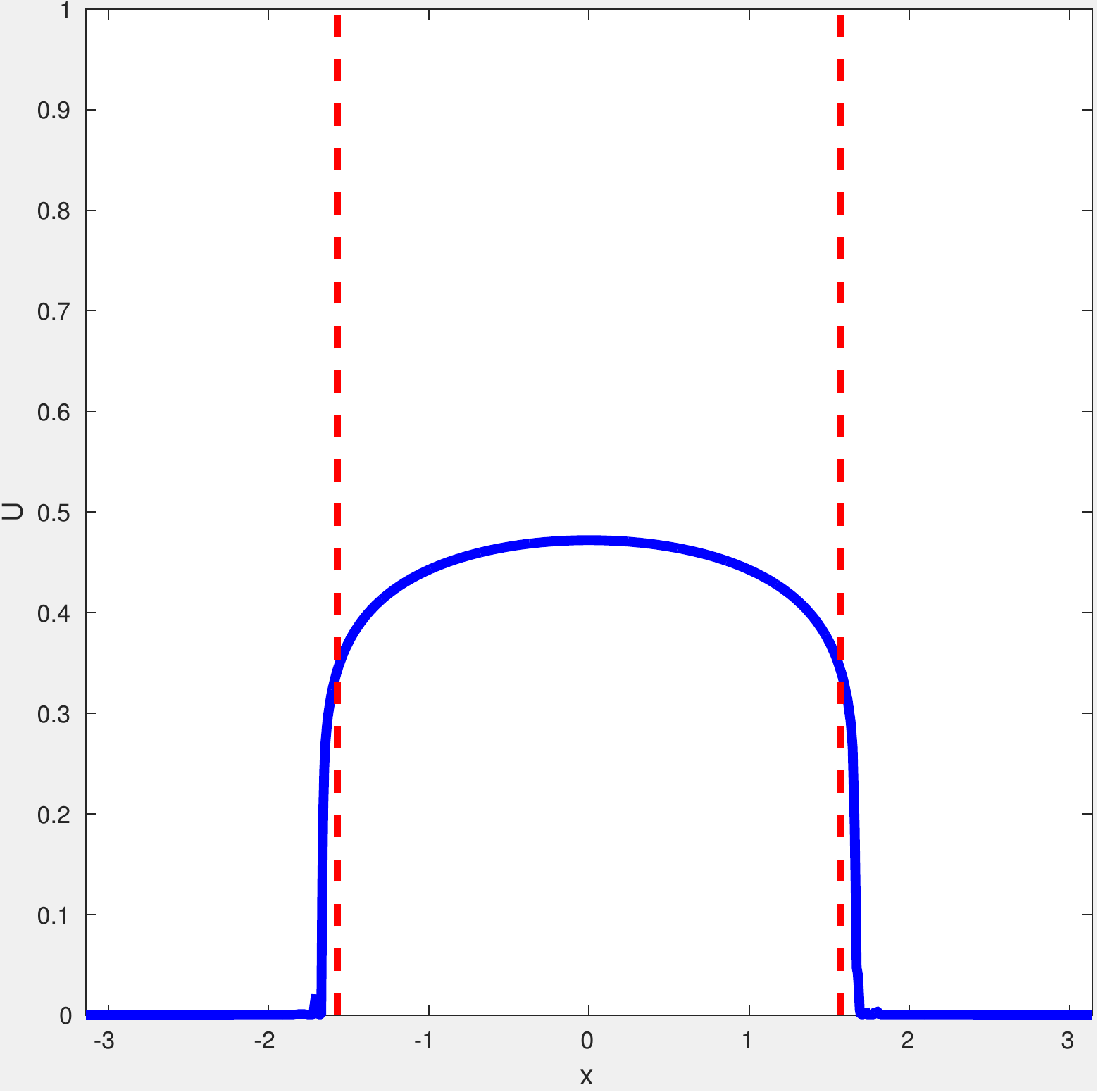}\caption{$t = 14$}\end{subfigure}\\%
        \begin{subfigure}[b]{0.25\linewidth}\includegraphics[scale=0.20]{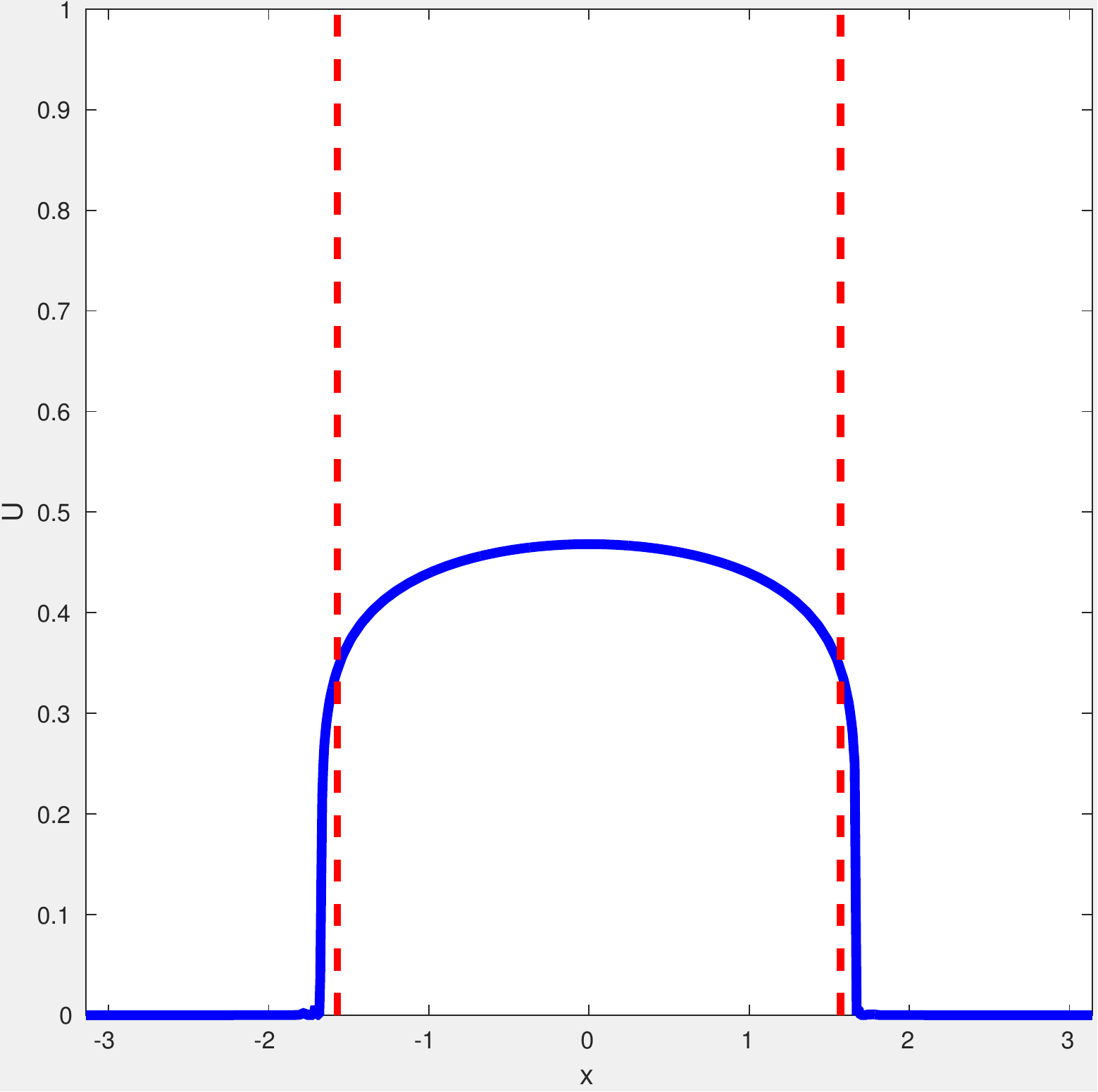}\caption{$t = 15$}\end{subfigure}%
        \begin{subfigure}[b]{0.25\linewidth}\includegraphics[scale=0.20]{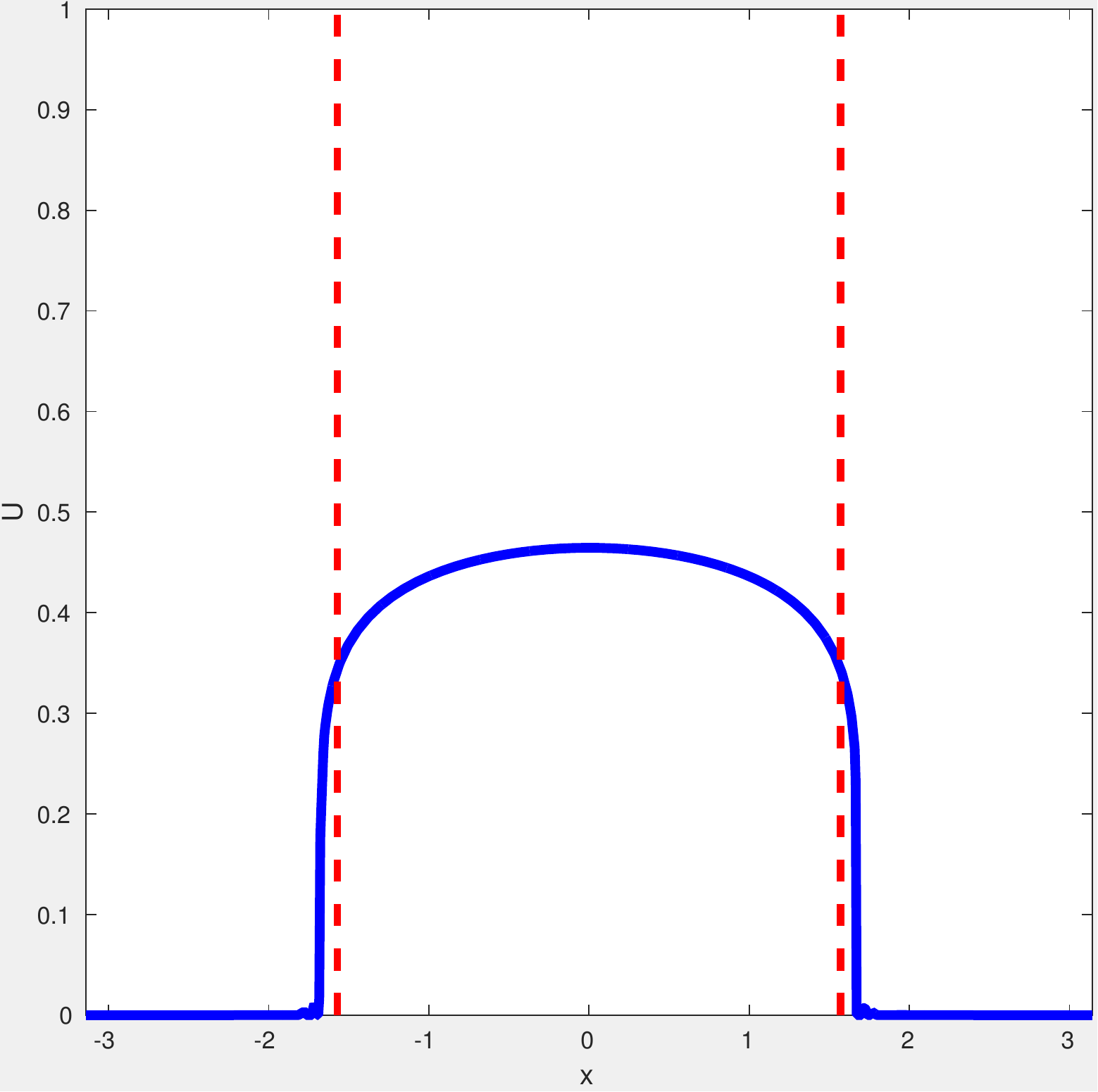}\caption{$t = 16$}\end{subfigure}%
        \begin{subfigure}[b]{0.25\linewidth}\includegraphics[scale=0.20]{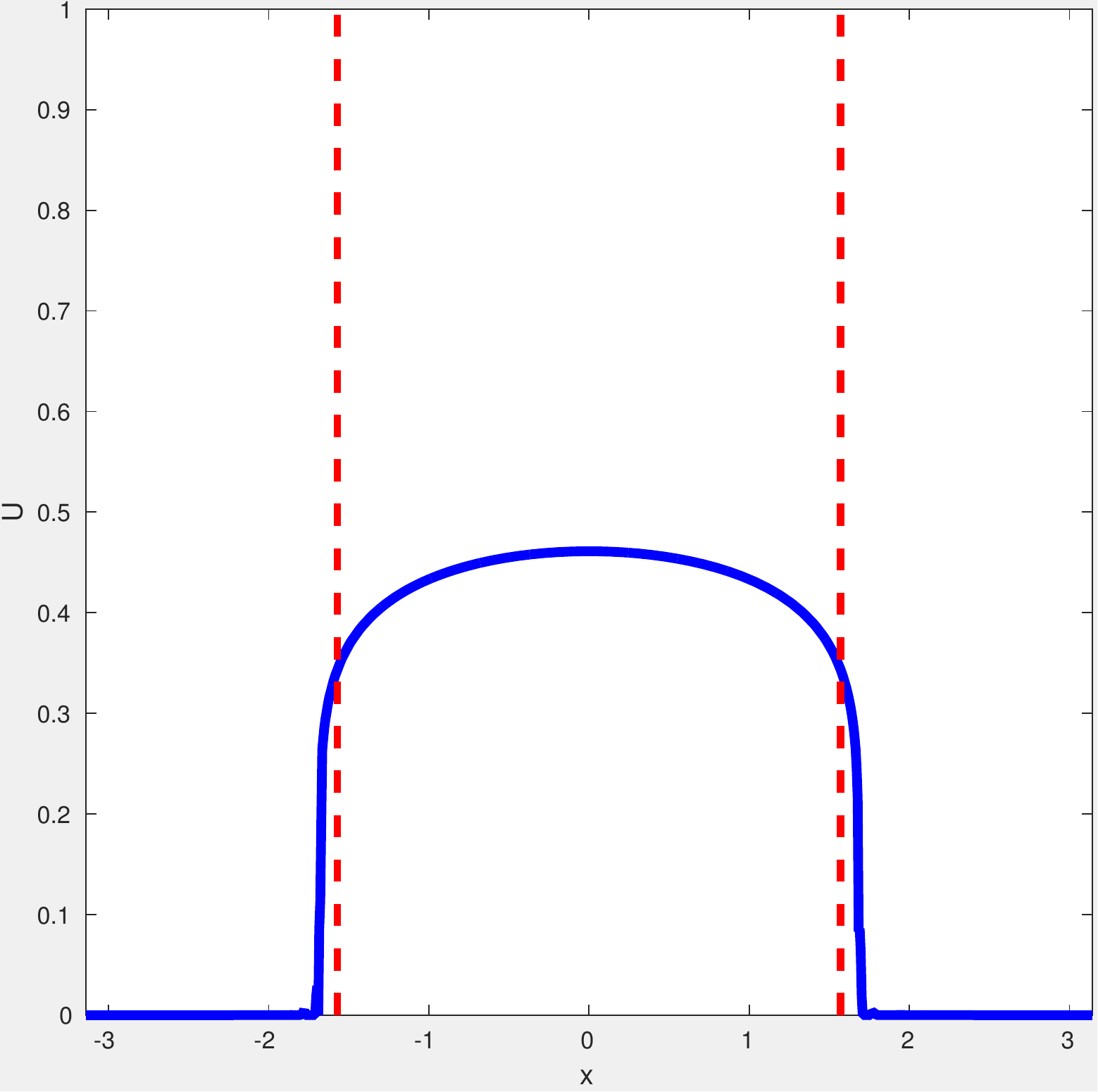}\caption{$t = 17$}\end{subfigure}%
        \begin{subfigure}[b]{0.25\linewidth}\includegraphics[scale=0.20]{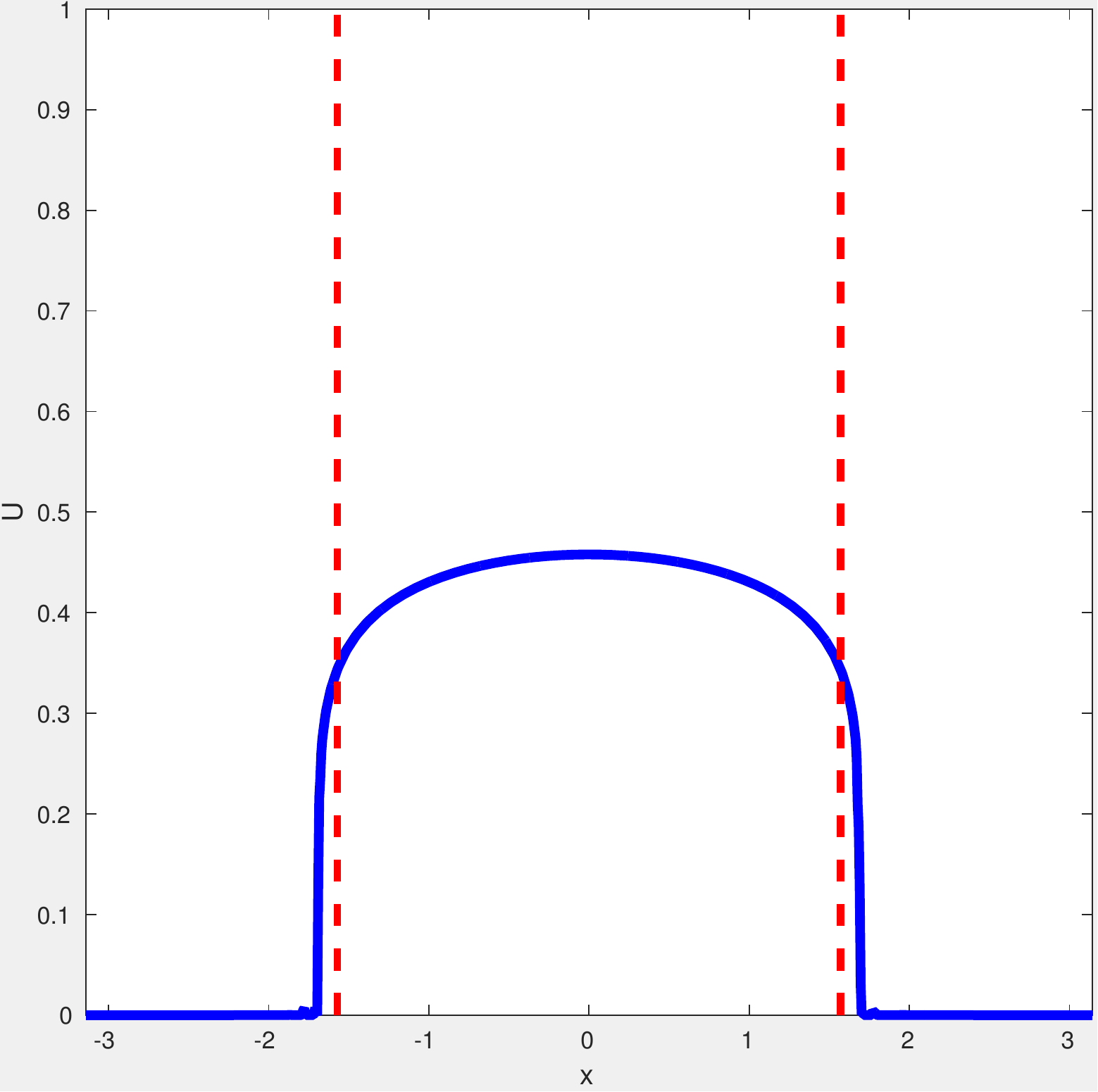}\caption{$t = 18$}\end{subfigure}\\%
        \caption{Example~\ref{exam4.4}. The cross section at $y = 0$ of a computed solution is shown at various time instants ($N = 40000$).}
        \label{fig:waiting-2D-cross-section}
    \end{figure}

\begin{figure}[ht]
        \centering
        \begin{subfigure}[b]{0.35\linewidth}\includegraphics[scale=0.34]{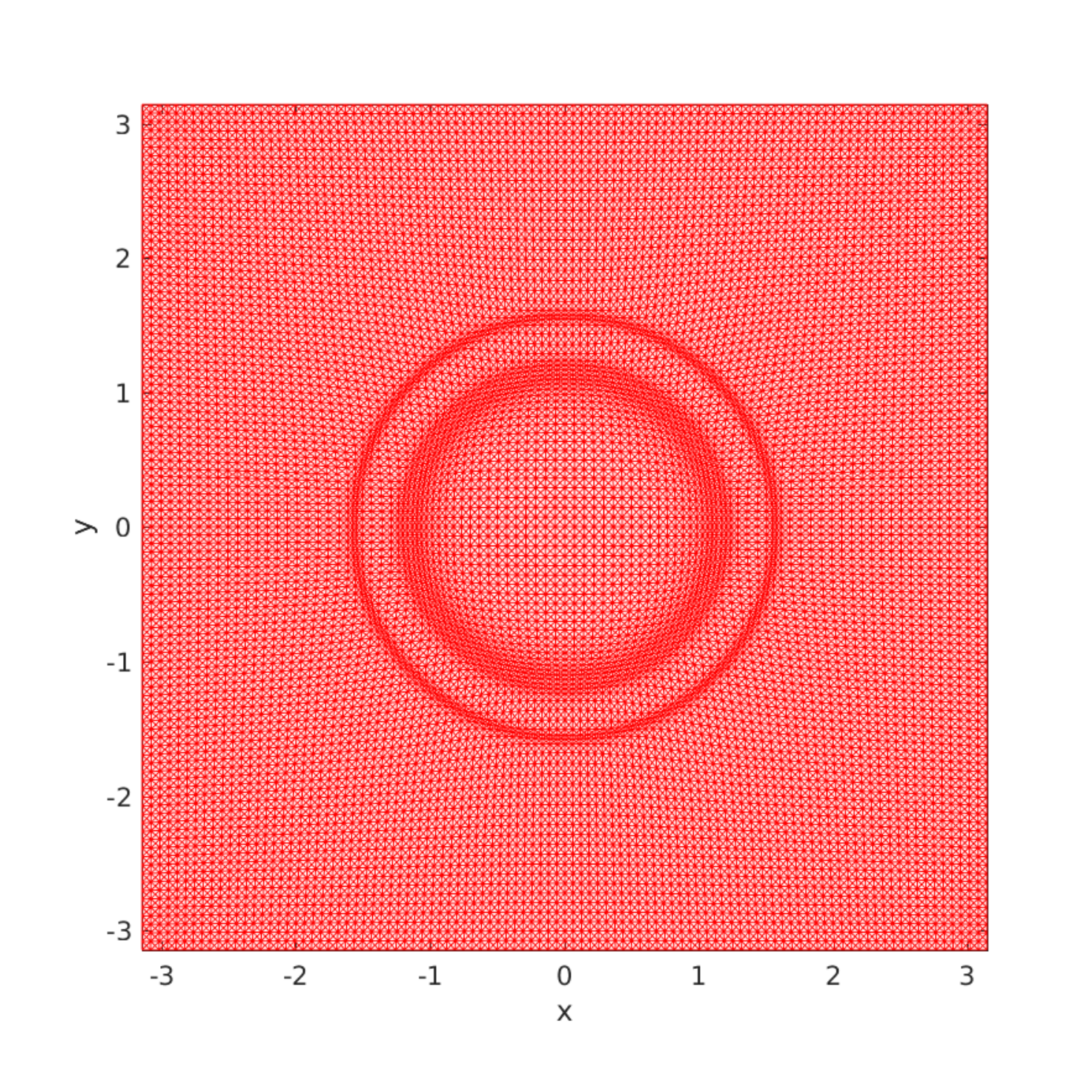}\caption{$t = 0.1$}\end{subfigure}\hspace{5mm}
        \begin{subfigure}[b]{0.35\linewidth}\includegraphics[scale=0.32]{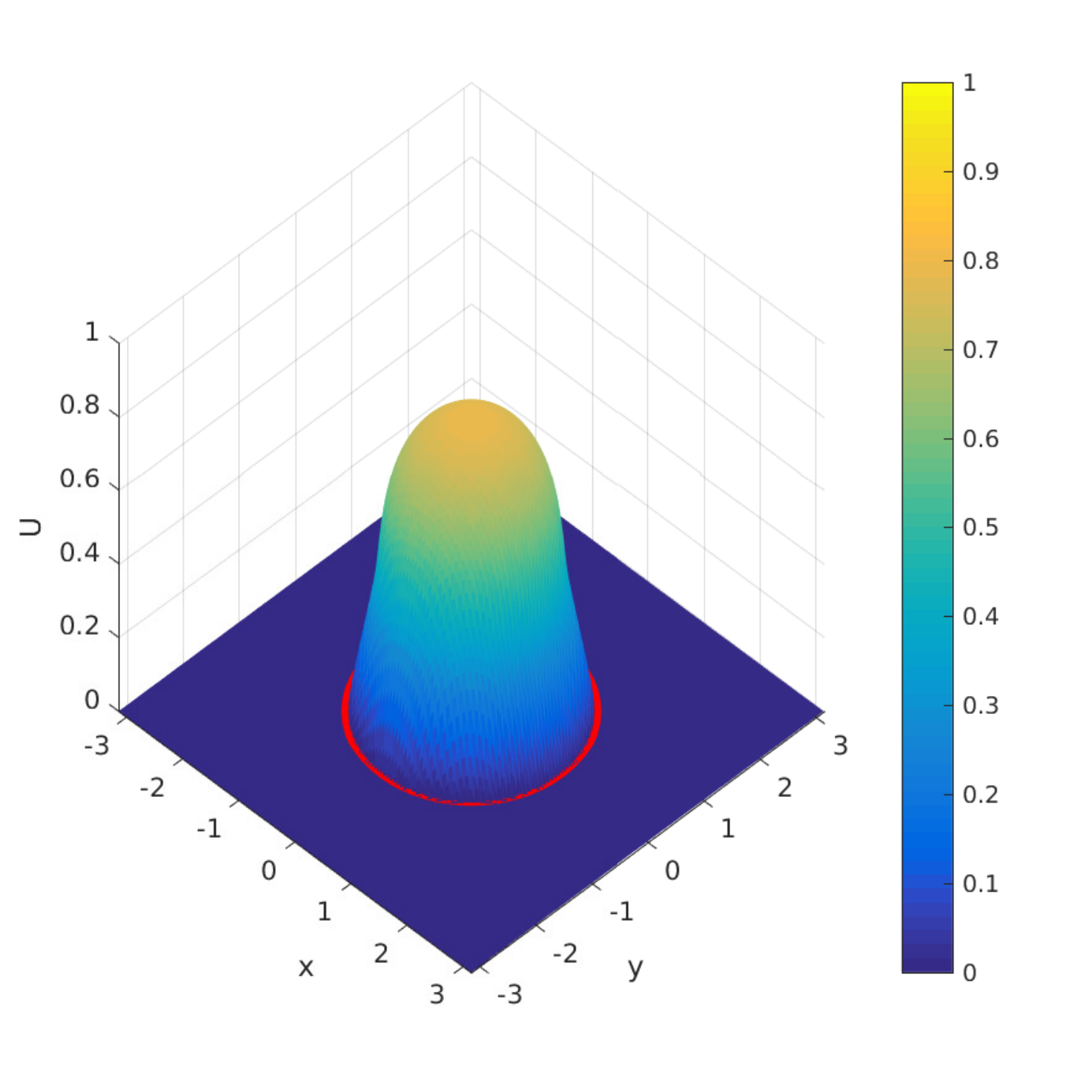}\caption{$t = 0.1$}\end{subfigure}\\%
         \vspace{-1mm}
        \begin{subfigure}[b]{0.35\linewidth}\includegraphics[scale=0.34]{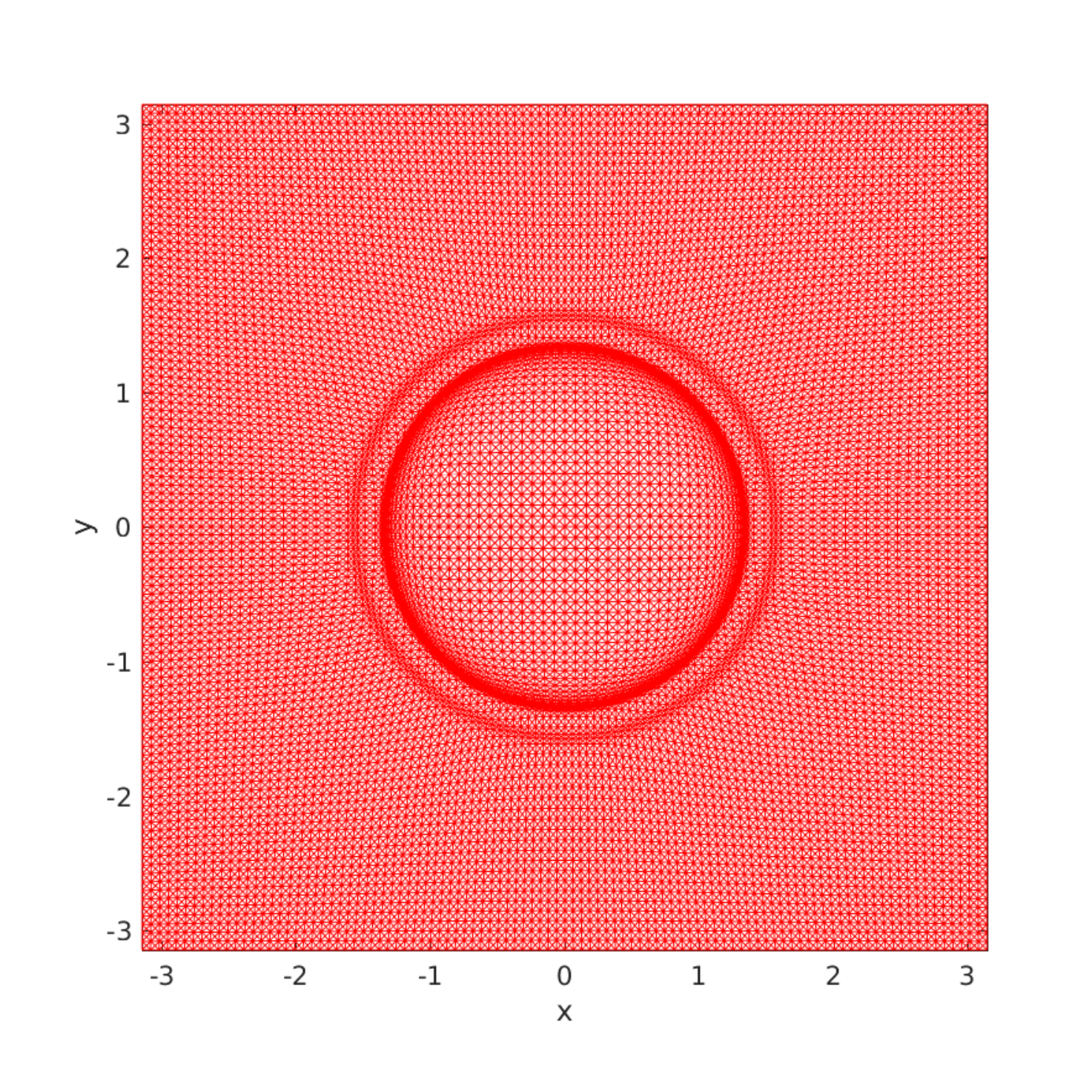}\caption{$t = 0.5$}\end{subfigure}\hspace{5mm}
        \begin{subfigure}[b]{0.35\linewidth}\includegraphics[scale=0.32]{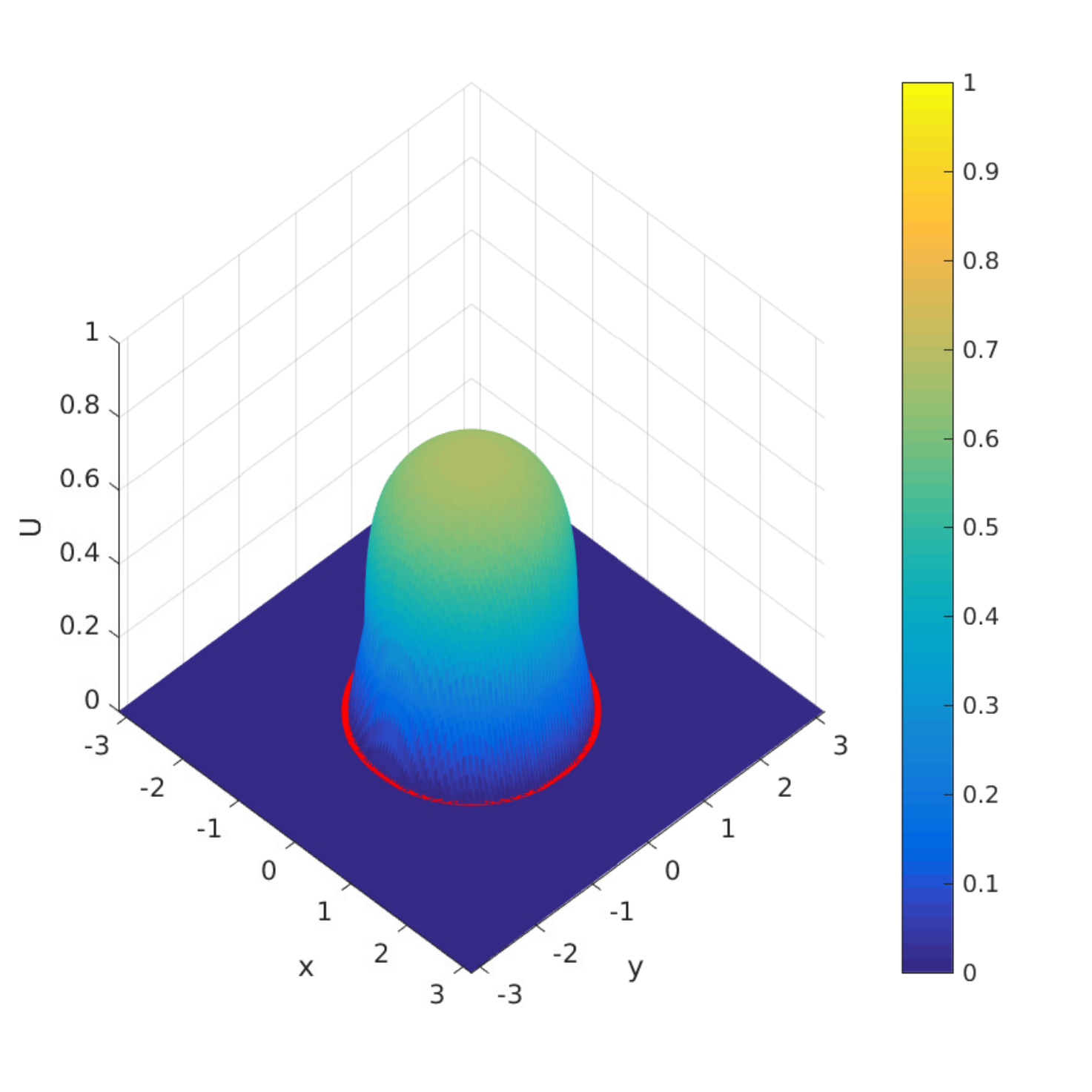}\caption{$t = 0.5$}\end{subfigure}\\%
         \vspace{-1mm}
        \begin{subfigure}[b]{0.35\linewidth}\includegraphics[scale=0.34]{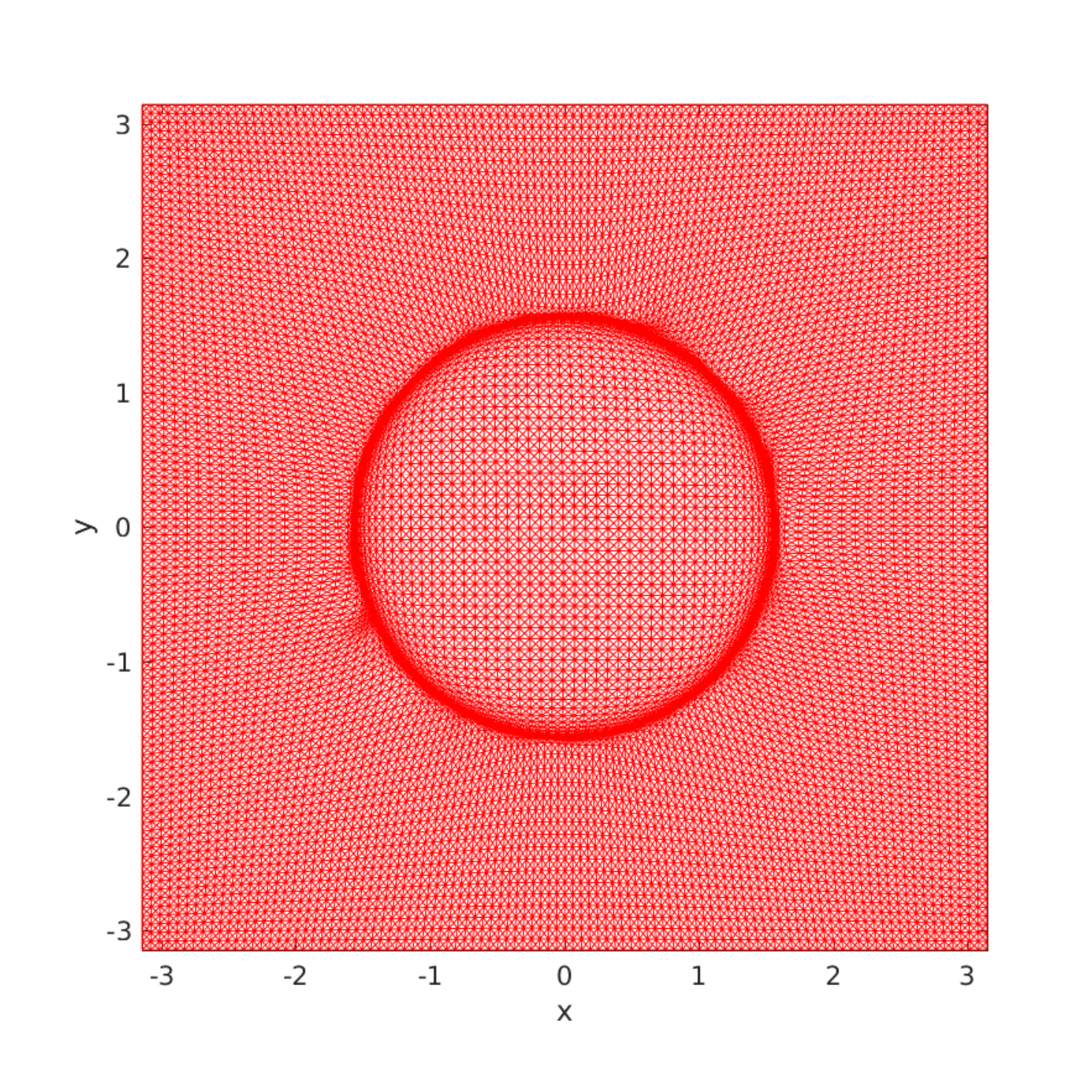}\caption{$t = 5$}\end{subfigure}\hspace{5mm}
        \begin{subfigure}[b]{0.35\linewidth}\includegraphics[scale=0.32]{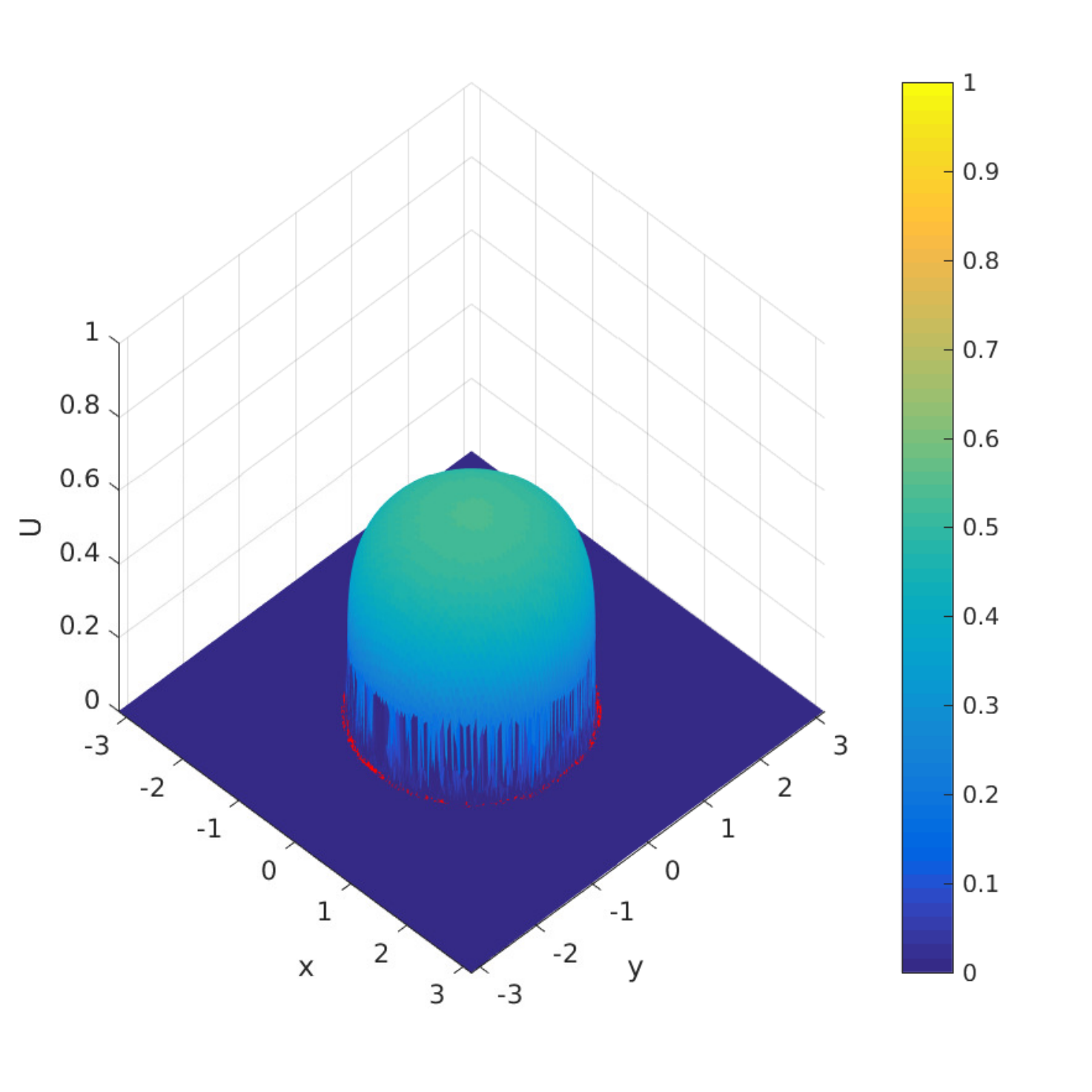}\caption{$t = 5$}\end{subfigure}\\%
         \vspace{-1mm}
        \begin{subfigure}[b]{0.35\linewidth}\includegraphics[scale=0.34]{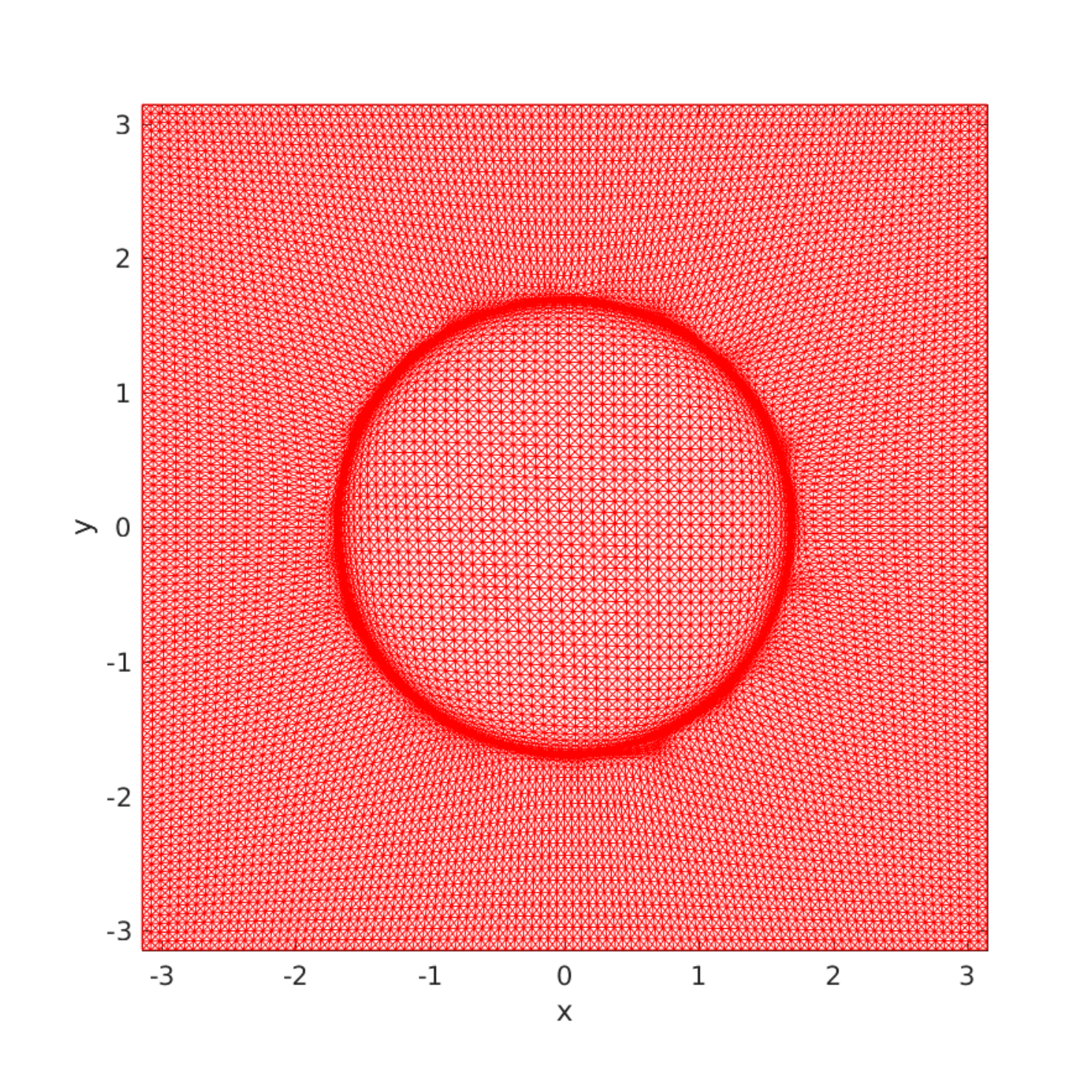}\caption{$t = 18.01$}\end{subfigure}\hspace{5mm}
        \begin{subfigure}[b]{0.35\linewidth}\includegraphics[scale=0.32]{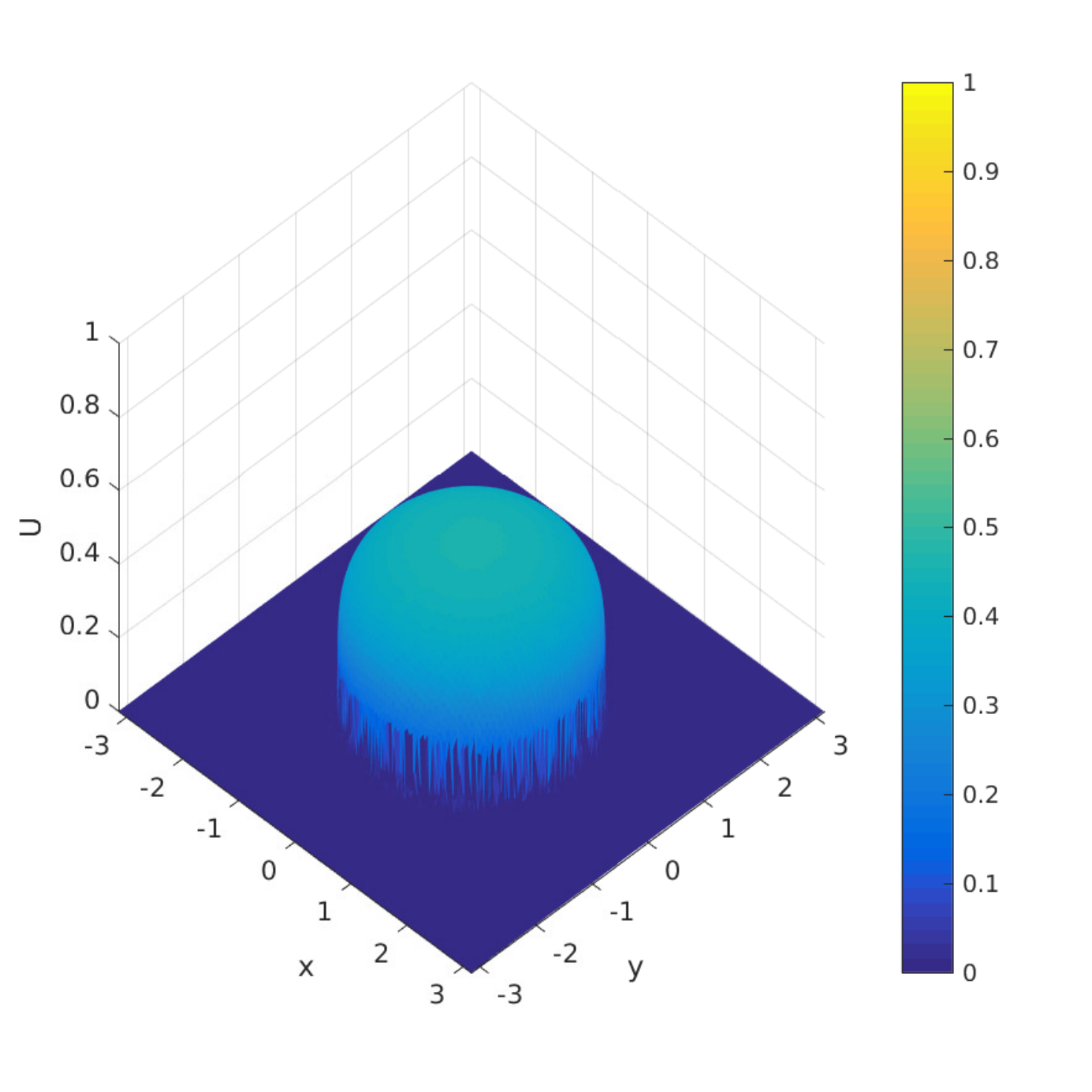}\caption{$t = 18.01$}\end{subfigure}%
        \caption{Example~\ref{exam4.4}. A computed solution is shown at various time instants ($N = 40000$).}
        \label{fig:waiting-2D}
    \end{figure}

\begin{exam}[Waiting-time phenomenon]
\label{exam4.4}
From Section~\ref{SEC:PME-theory} we recall that PME exhibits the waiting-time phenomenon for a certain type of
initial solutions. To see this, we consider
    \[
        m = 8, \quad \Omega = (-\pi, \pi) \times (-\pi, \pi) ,
    \]
    \begin{equation}
        u_0(x,y) = 
        \begin{cases}
            \cos(\sqrt{x^2 + y^2}) , & \quad \text{for} \quad \sqrt{x^2 + y^2} \leq \frac{\pi}{2} \\
            0 , & \quad \text{otherwise} .
        \end{cases}
    \label{waiting-time-1}
    \end{equation}
 We have
 \[
 \grad \cos^m(\sqrt{x^2+y^2}) = -\frac{m\, \cos^{m-1}(\sqrt{x^2+y^2})\, \sin(\sqrt{x^2+y^2})}{\sqrt{x^2+y^2}} 
 \begin{bmatrix} x \\ y \end{bmatrix} ,
 \]
which diminishes at $\sqrt{x^2 + y^2} = \frac{\pi}{2}$. From Darcy's law (\ref{Darcy-law}),
we do not anticipate that the free boundary moves initially.
 
In Figs.~\ref{fig:waiting-2D-cross-section} and Fig.~\ref{fig:waiting-2D} we show
the cross section at $y = 0$ of a computed solution and the solution itself.
The results show that the free boundary of the solution does not move until around $t = 10$.
Before this time, the solution is steepening. Interestingly, the steepening does not occur
on the whole initial support. Instead, it first occurs on a smaller region inside the support and
then this region is expanding until it fills the whole initial support. After that, the free boundary
waits until it becomes sufficiently steep and then moves.
\qed \end{exam}

%%%%%%%%%%%%%%%%%%%%%%%%%%%%%%%
%% Section: generalized PME
\section{Numerical experiment for PME with variable exponents and absorption}
\label{SEC:PME-numerics-2}
    
To demonstrate the robustness of the moving mesh finite element method described in Section~\ref{SEC:mmfem},
we consider its application to PME with absorption and/or variable exponents, 
\begin{equation}
u_t = \grad \cdot (|u|^{\gamma} \grad u) - \lambda u^{\sigma}, \quad \Omega \times (t_0, T]
\label{PME-2}
\end{equation}
subject to a homogeneous Dirichlet boundary condition and an initial condition.
Here, $\gamma = \gamma (\vect{x}, t)$ and $\sigma = \sigma (\vect{x}, t)$ are nonnegative bounded
functions and $\lambda$ is a constant.
PME in the form of (\ref{PME-2}) arises in continuum mechanics to model the motion of a barotropic gas
through a porous medium, where the pressure is considered to depend on the density and
temperature \cite{Antontsev2005}.
Like the standard PME, (\ref{PME-2}) with constant exponents (i.e., PME with absorption)
has been studied extensively; e.g., see \cite{Knerr1979,Shmarev2005}.
However, there are very few theoretical results for the case with variable exponents \cite{Antontsev2005,Lian2008}.
For example, there is no theoretical result on the movement of the free boundary (cf. (\ref{Darcy-law}))
although the solution to (\ref{PME-2}) is known to have the property of finite speed of propagation.
Neither is there much numerical work on this situation; see \cite{Duque2013,Duque2014,Duque2015}.

\begin{figure}[ht]
        \centering
        \begin{subfigure}[b]{0.35\linewidth}\includegraphics[scale=0.34]{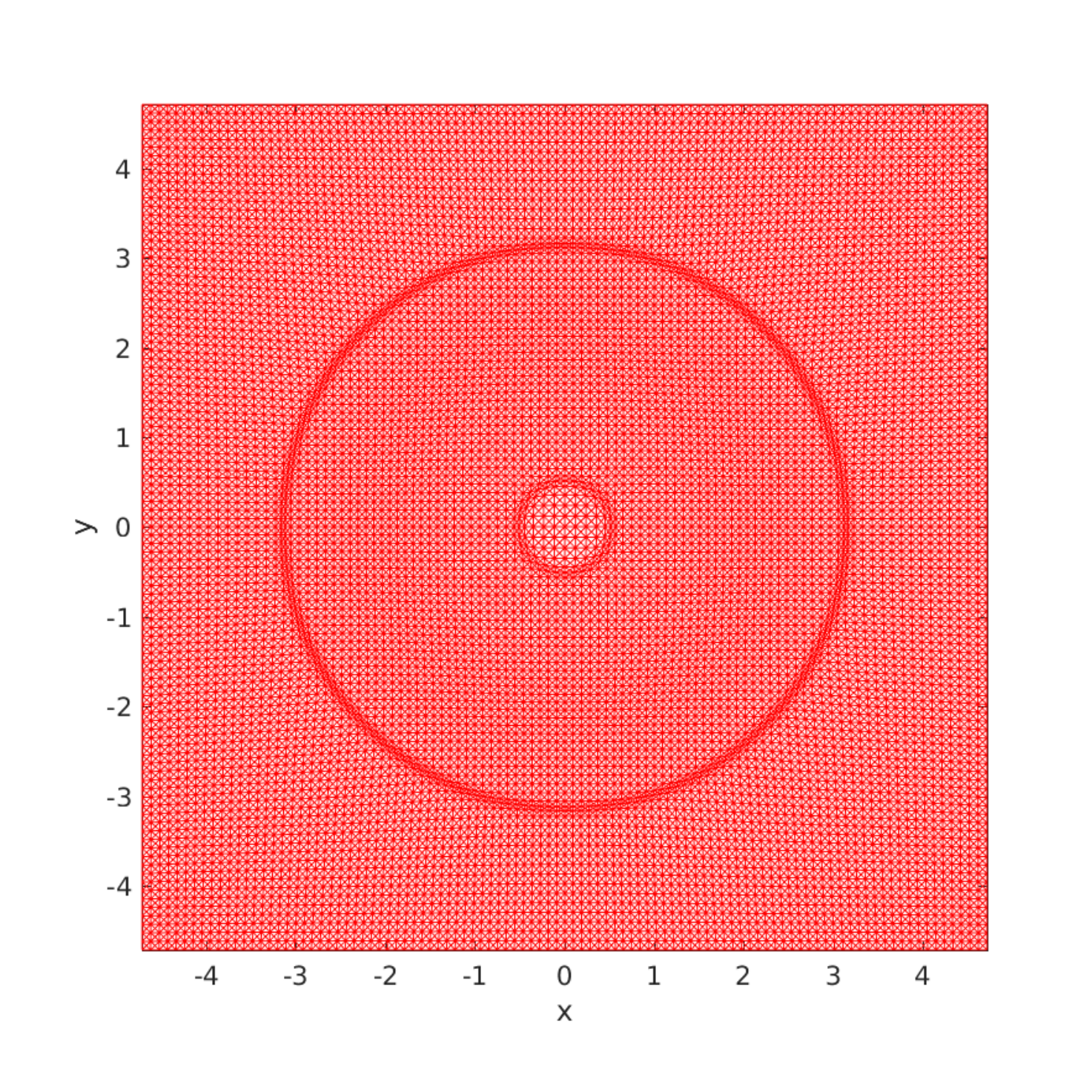}\caption{$t = 0$}\end{subfigure}\hspace{5mm}%
        \begin{subfigure}[b]{0.35\linewidth}\includegraphics[scale=0.32]{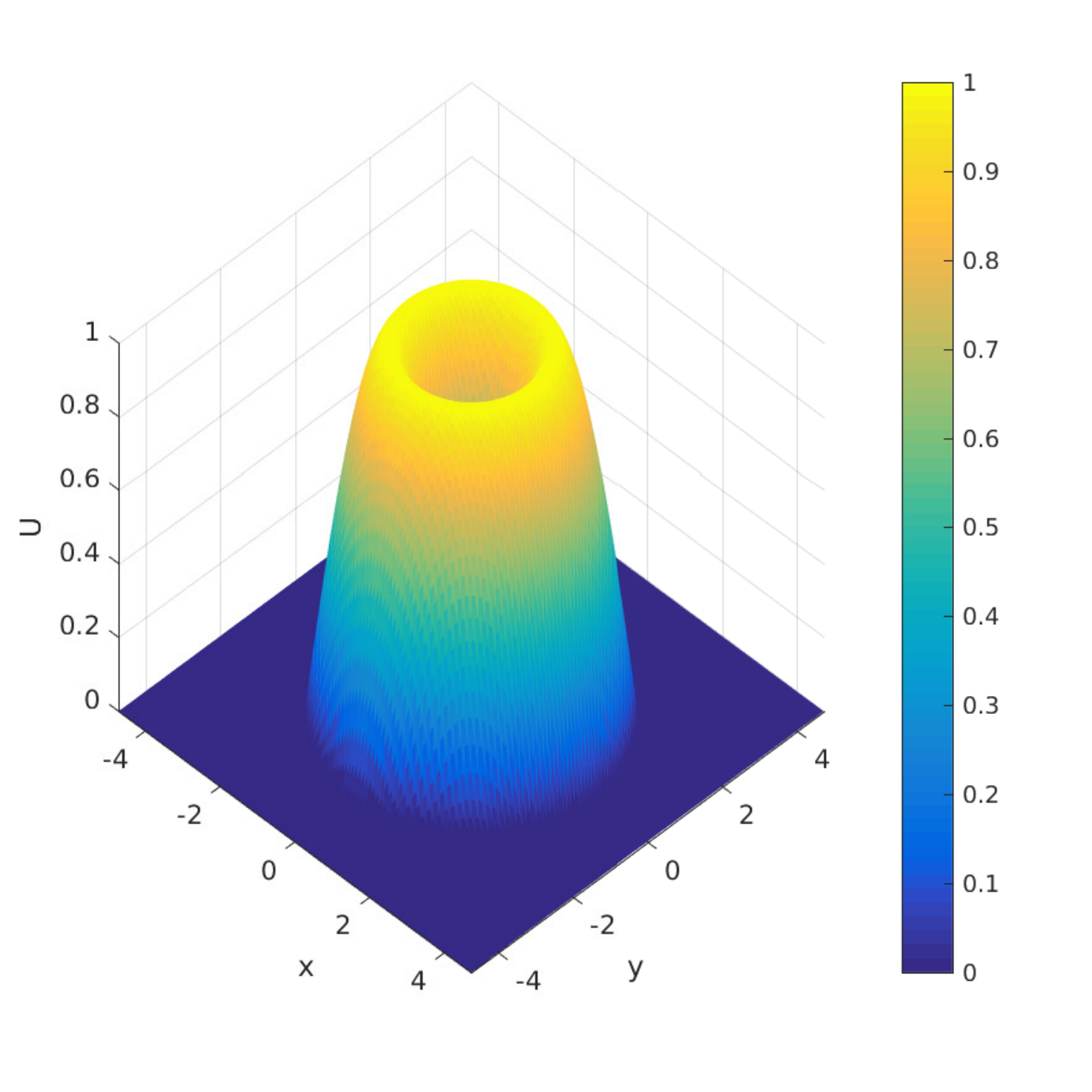}\caption{$t = 0$}\end{subfigure}\\%
        \vspace{-1mm}
        \begin{subfigure}[b]{0.35\linewidth}\includegraphics[scale=0.34]{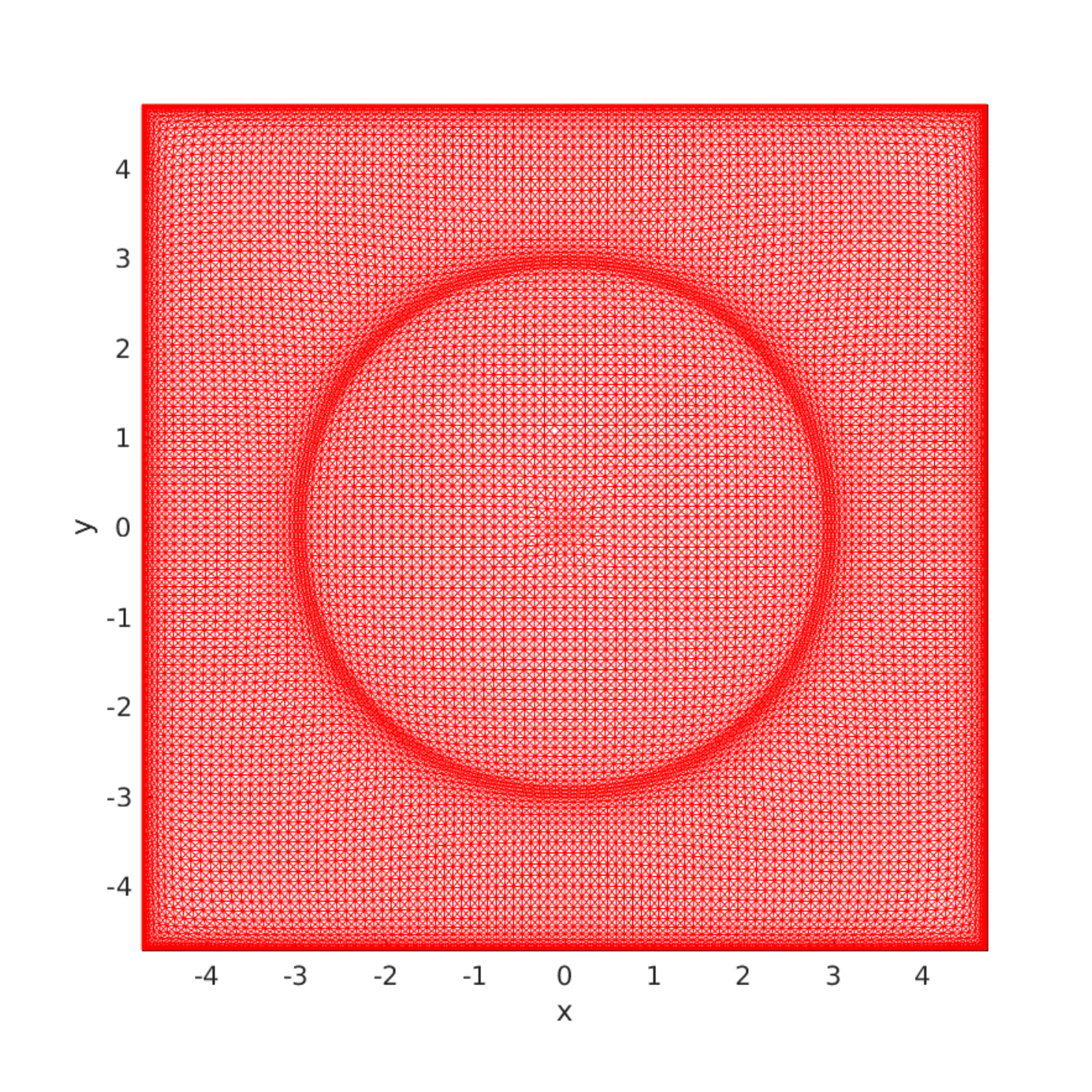}\caption{$t = 0.40$}\end{subfigure}\hspace{5mm}%
        \begin{subfigure}[b]{0.35\linewidth}\includegraphics[scale=0.32]{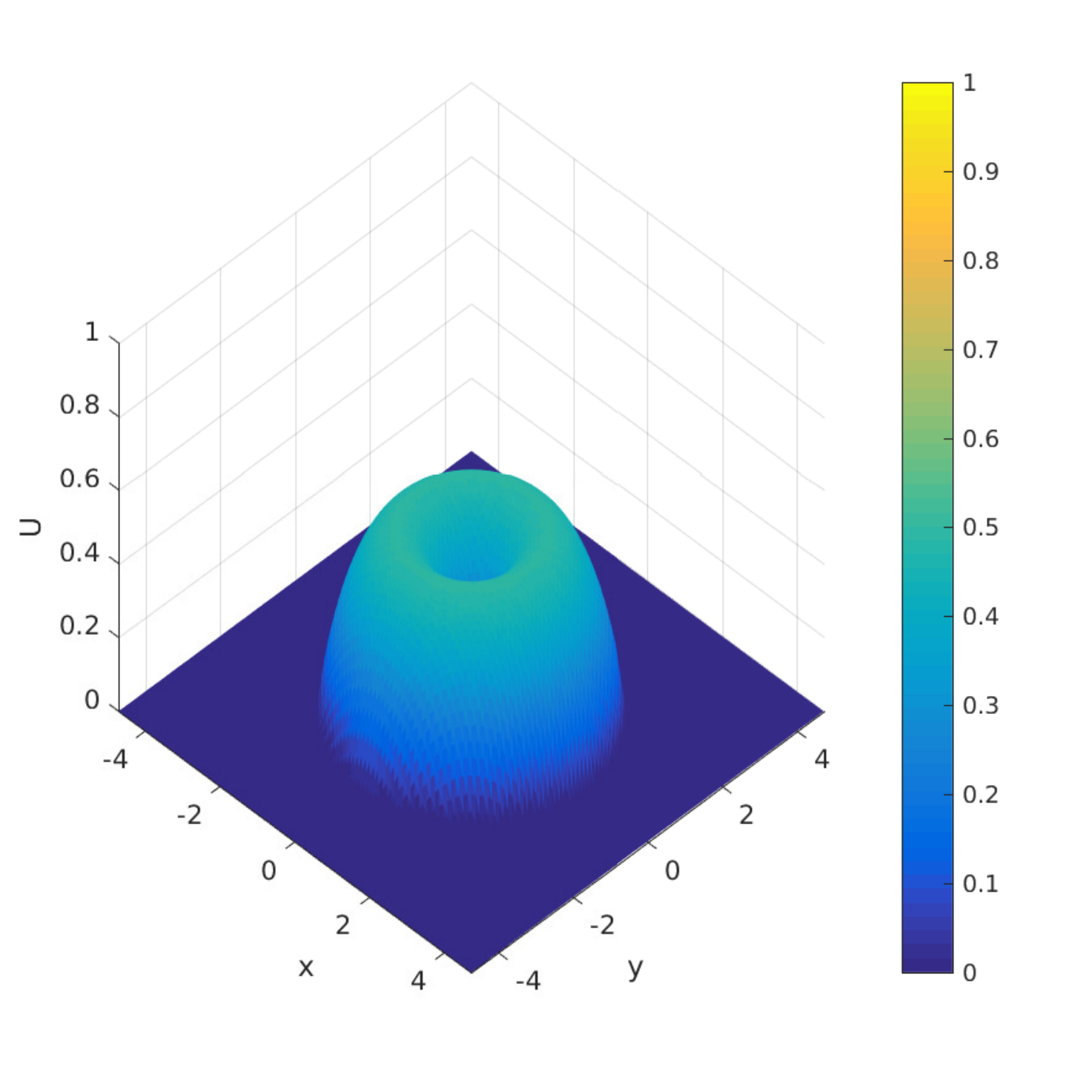}\caption{$t = 0.40$}\end{subfigure}\\%
        \vspace{-1mm}
        \begin{subfigure}[b]{0.35\linewidth}\includegraphics[scale=0.34]{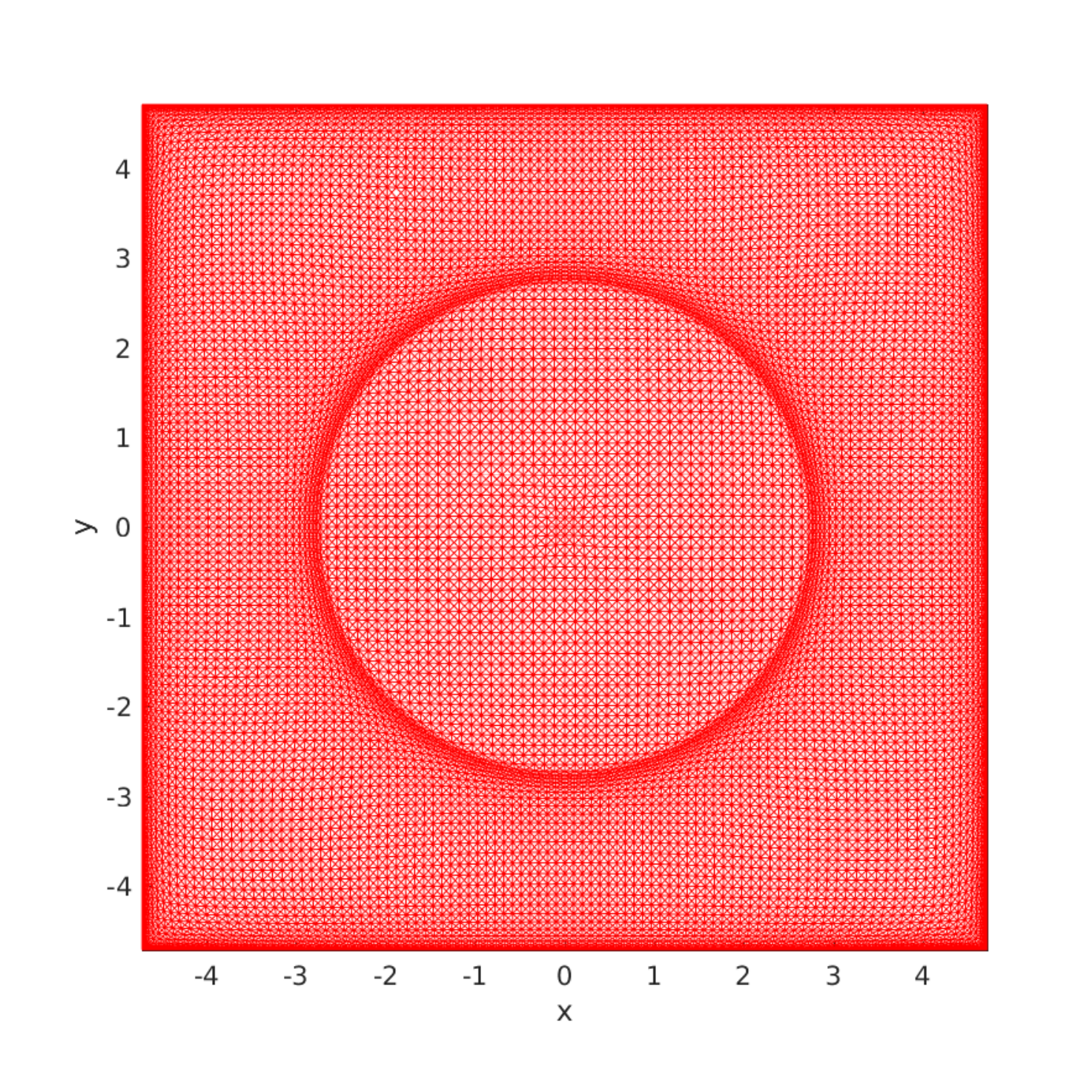}\caption{$t = 0.64$}\end{subfigure}\hspace{5mm}%
        \begin{subfigure}[b]{0.35\linewidth}\includegraphics[scale=0.32]{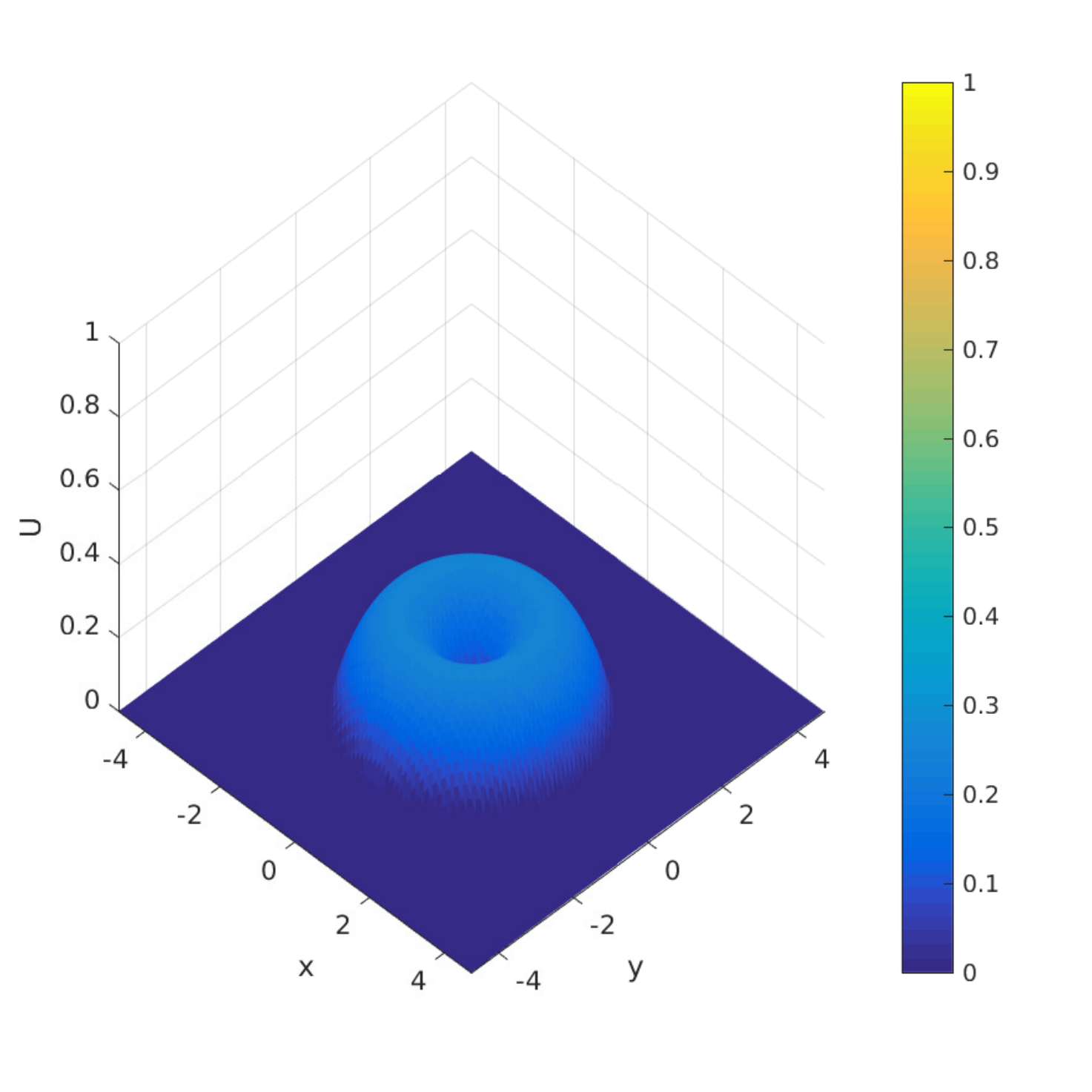}\caption{$t = 0.64$}\end{subfigure}\\%
        \vspace{-1mm}
        \begin{subfigure}[b]{0.35\linewidth}\includegraphics[scale=0.34]{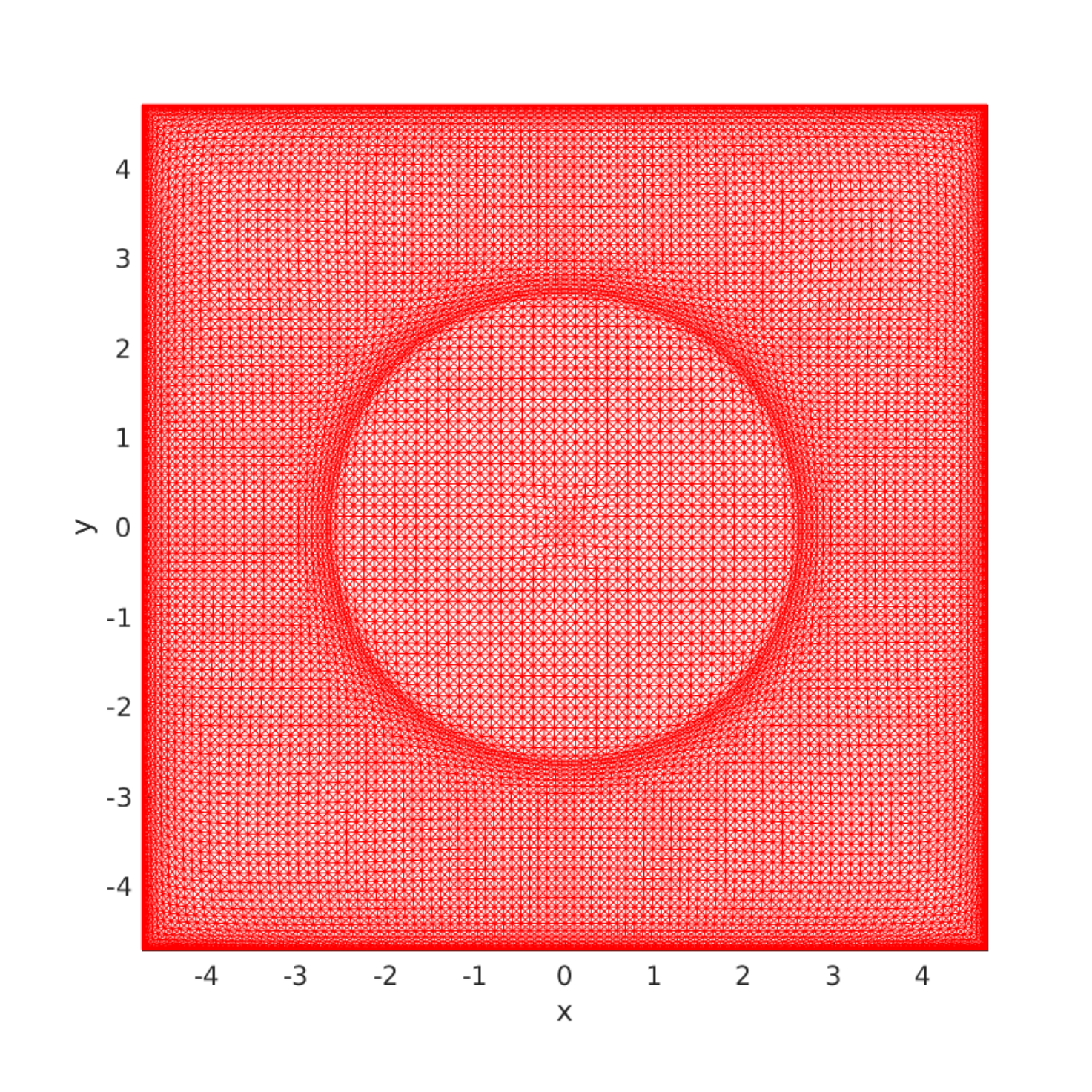}\caption{$t = 0.80$}\end{subfigure}\hspace{5mm}%
        \begin{subfigure}[b]{0.35\linewidth}\includegraphics[scale=0.32]{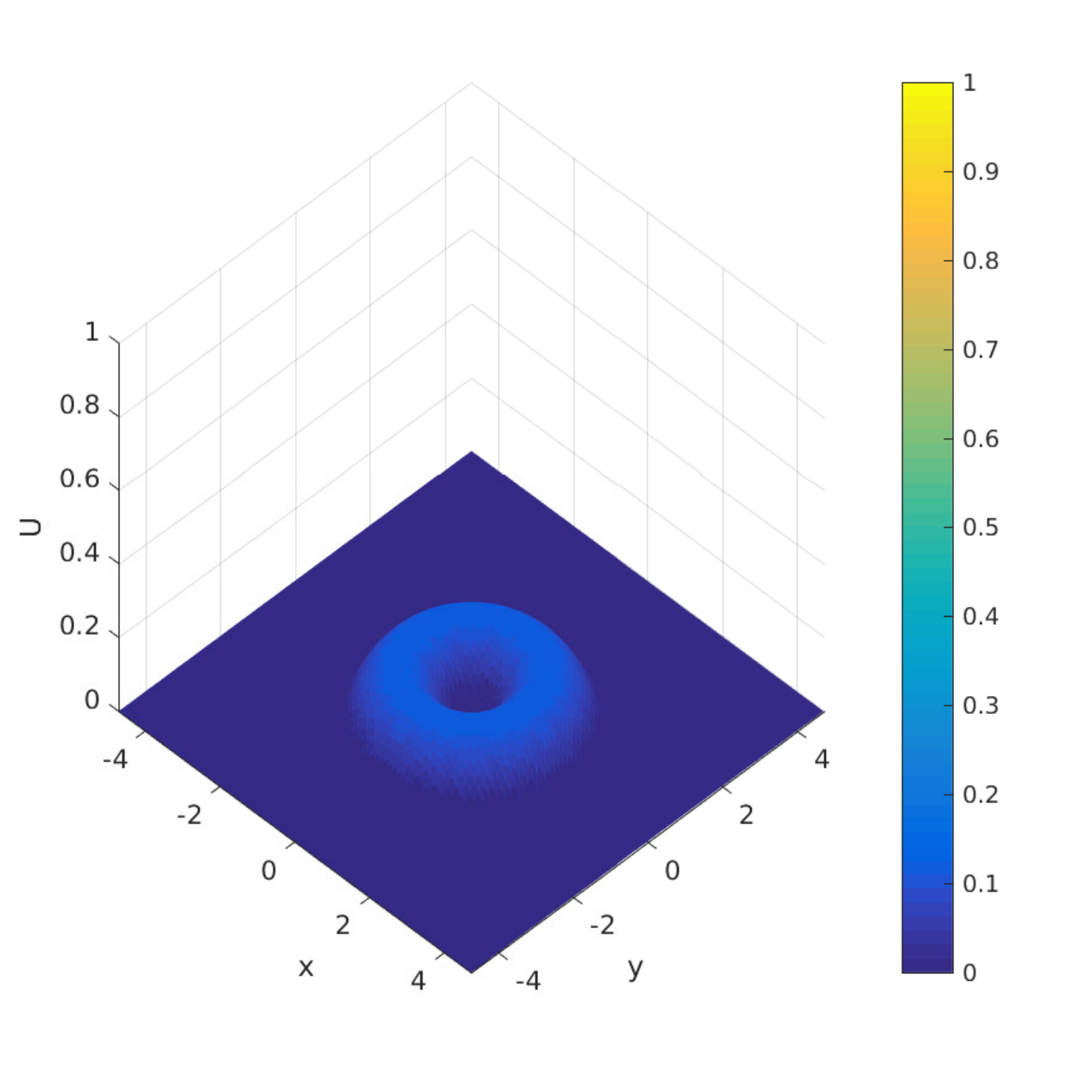}\caption{$t = 0.80$}\end{subfigure}%
        \caption{Example~\ref{exam5.1}. An adaptive mesh and the corresponding solution at various time instants ($N = 40000$).}
        \label{fig:Splitting-soln-mesh}
    \end{figure}

\begin{exam}[Constant exponents with absorption]
\label{exam5.1}
We first consider an example with an absorption term,
\begin{align*}
& \lambda = 1, \quad \gamma = 2,\quad \sigma = 0.1,
\quad \Omega =  (-1.5 \pi, 1.5 \pi) \times (-1.5 \pi, 1.5 \pi) , \\
&
        u_0 = 
        \begin{cases}
            | \sin(\sqrt{x^2 + y^2}) | , & \quad \text{for} \quad \sqrt{x^2 + y^2} \in  (\frac{\pi}{6}, \pi) \\
            0.5 , & \quad \text{for} \quad \sqrt{x^2 + y^2} \in [0, \frac{\pi}{6})  \\
            0 , & \quad \text{otherwise} .
        \end{cases}
\end{align*}
This example is the two-dimensional generalization of a one-dimensional example in \cite{Zhang-2009} that shows
a splitting phenomenon in the middle after a finite time. 

An adaptive mesh and the corresponding computed solution are shown in Fig.~\ref{fig:Splitting-soln-mesh}.
We can see that as time evolves, the solution is becoming lower and the support
is expanding on the outer boundary. Meanwhile, the solution is being ``punched through'' with a hole
at the middle of the support. This is an additional feature with the absorption term.
\qed \end{exam}

\begin{figure}[ht]
        \centering
        \begin{subfigure}[b]{0.35\linewidth}\includegraphics[scale=0.34]{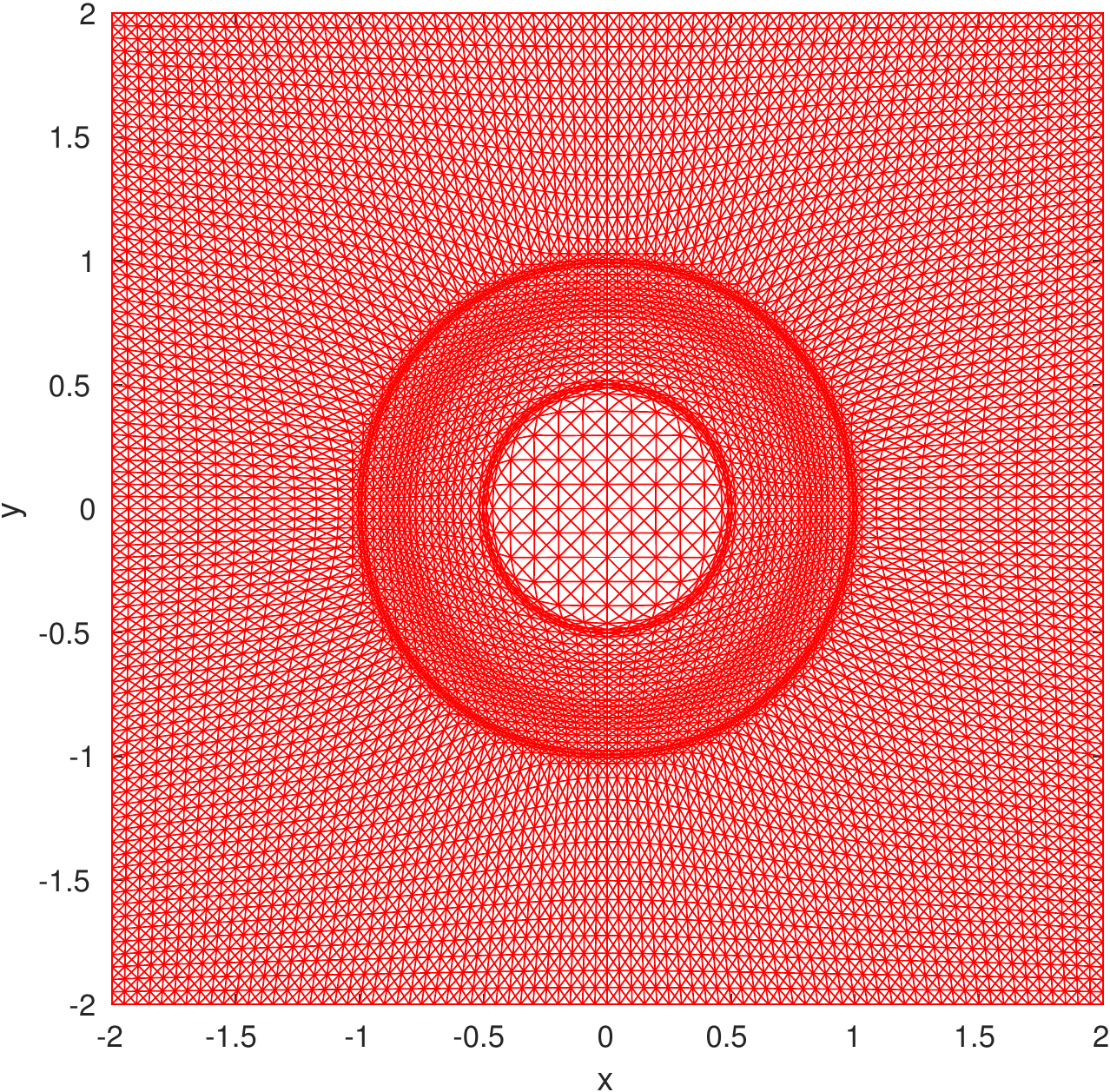}\caption{$t = 0$}\end{subfigure}\hspace{5mm}
        \begin{subfigure}[b]{0.35\linewidth}\includegraphics[scale=0.26]{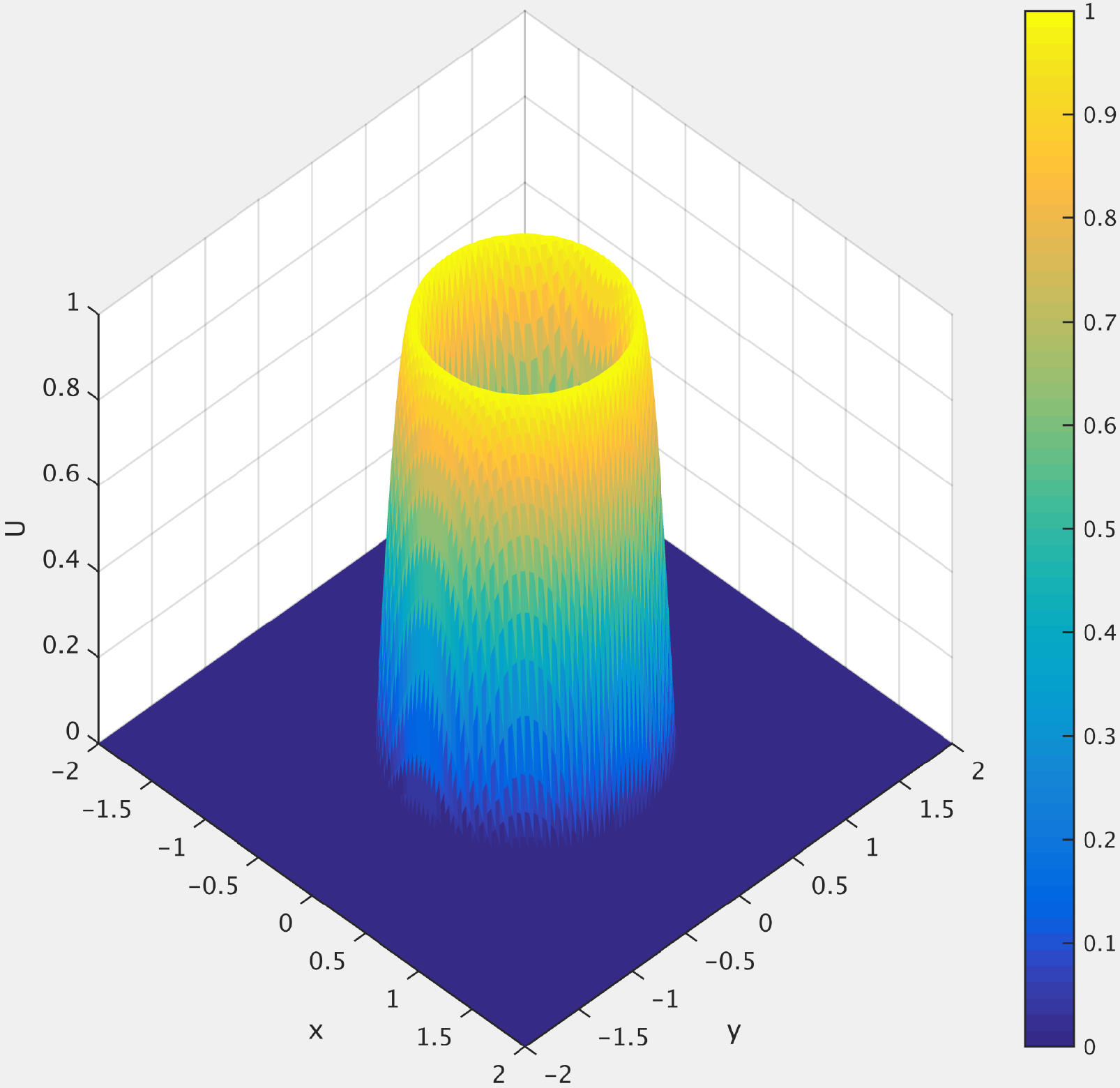}\caption{$t = 0$}\end{subfigure}\\%
        \vspace{2mm}
        \begin{subfigure}[b]{0.35\linewidth}\includegraphics[scale=0.34]{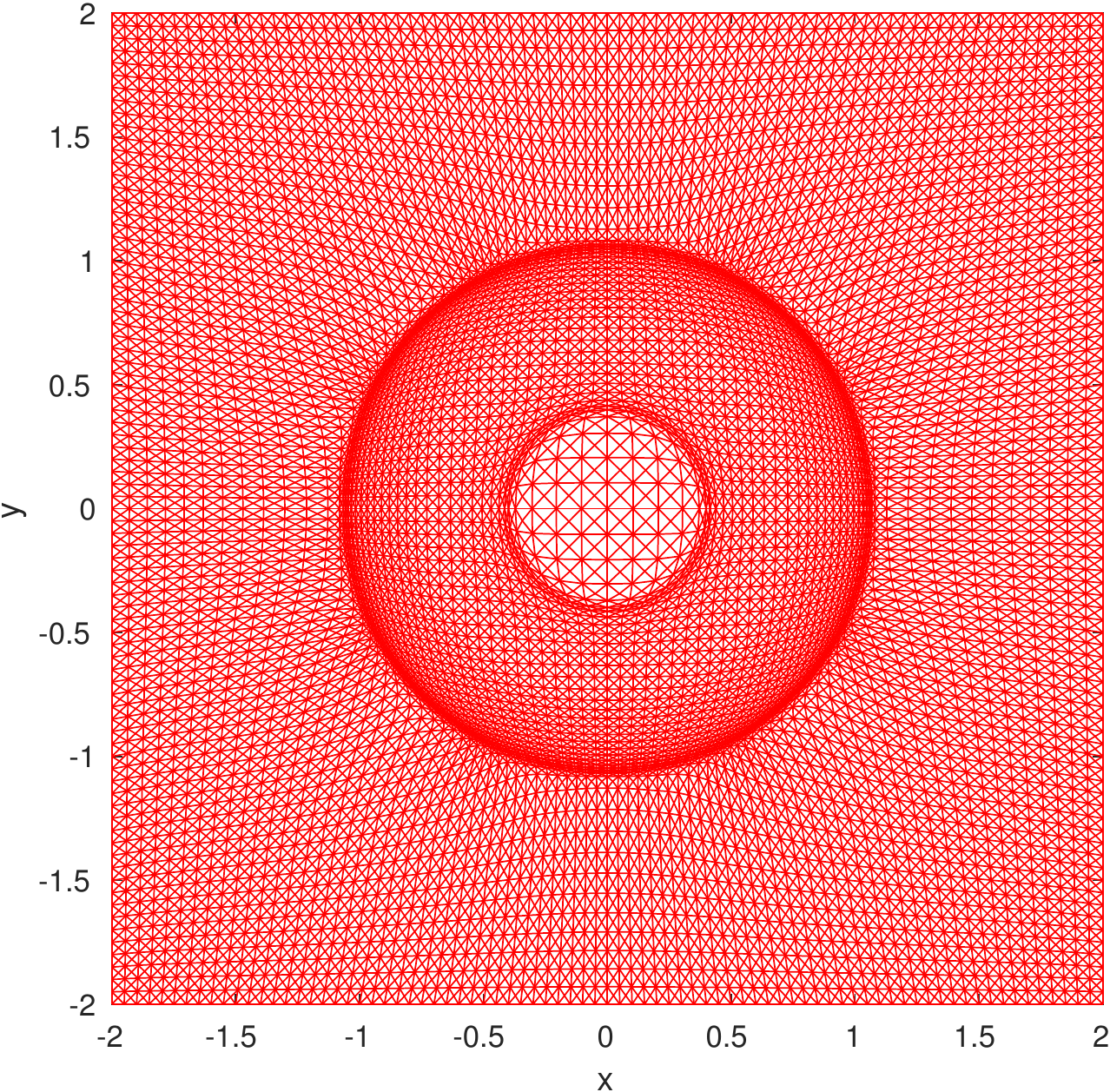}\caption{$t = 0.02$}\end{subfigure}\hspace{5mm}
        \begin{subfigure}[b]{0.35\linewidth}\includegraphics[scale=0.26]{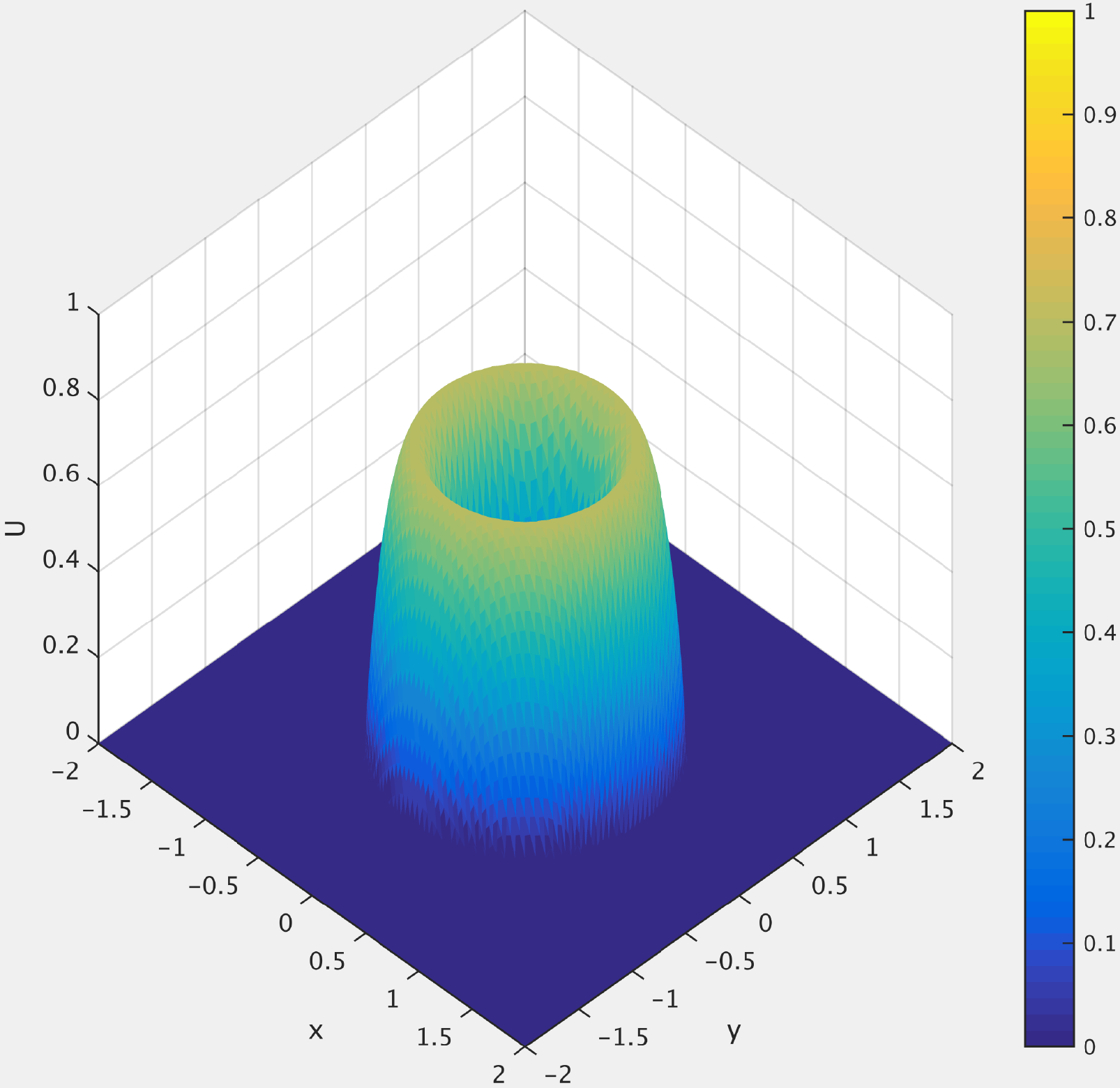}\caption{$t = 0.02$}\end{subfigure}\\%
        \vspace{2mm}
        \begin{subfigure}[b]{0.35\linewidth}\includegraphics[scale=0.34]{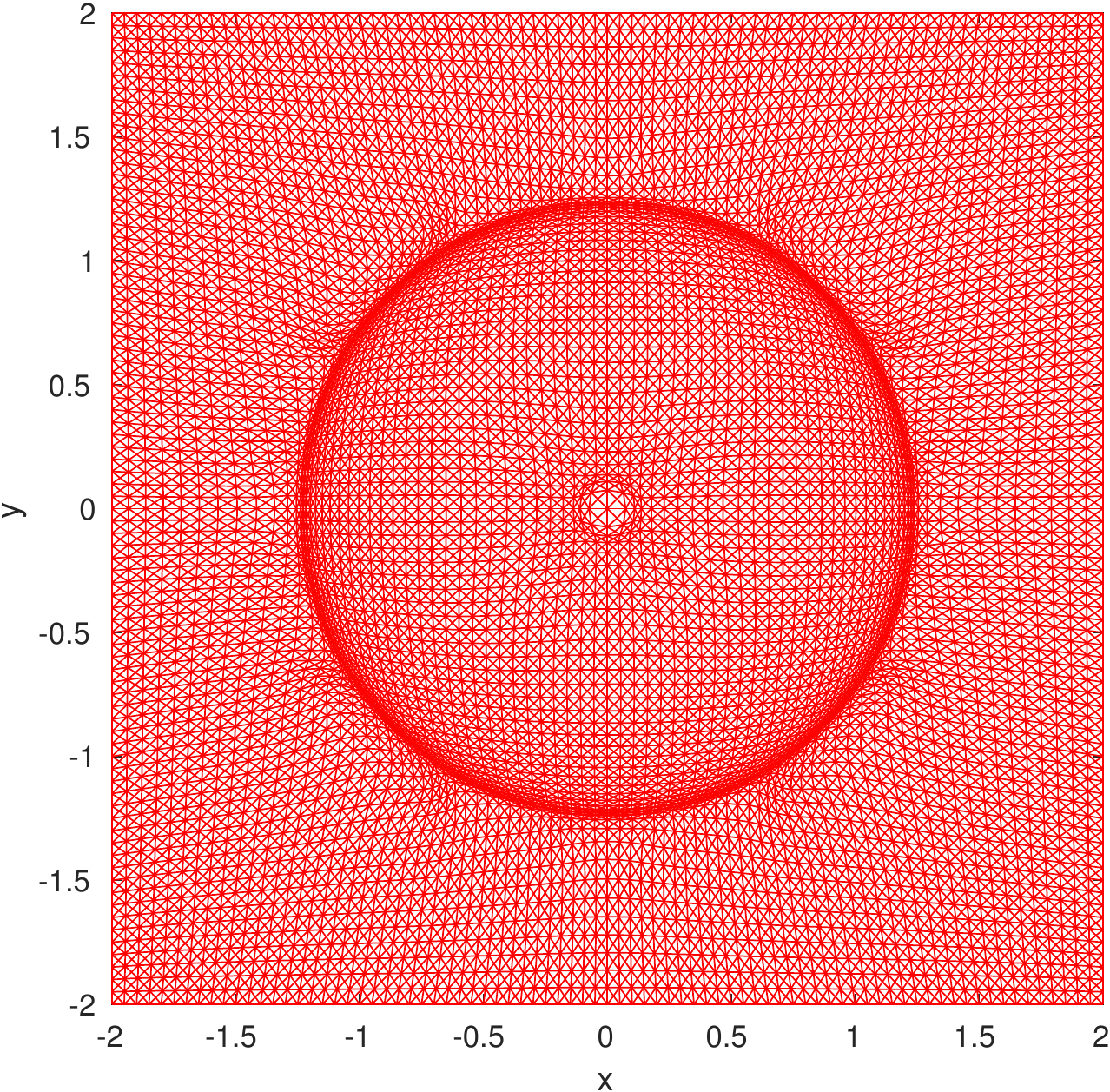}\caption{$t = 0.10$}\end{subfigure}\hspace{5mm}
        \begin{subfigure}[b]{0.35\linewidth}\includegraphics[scale=0.26]{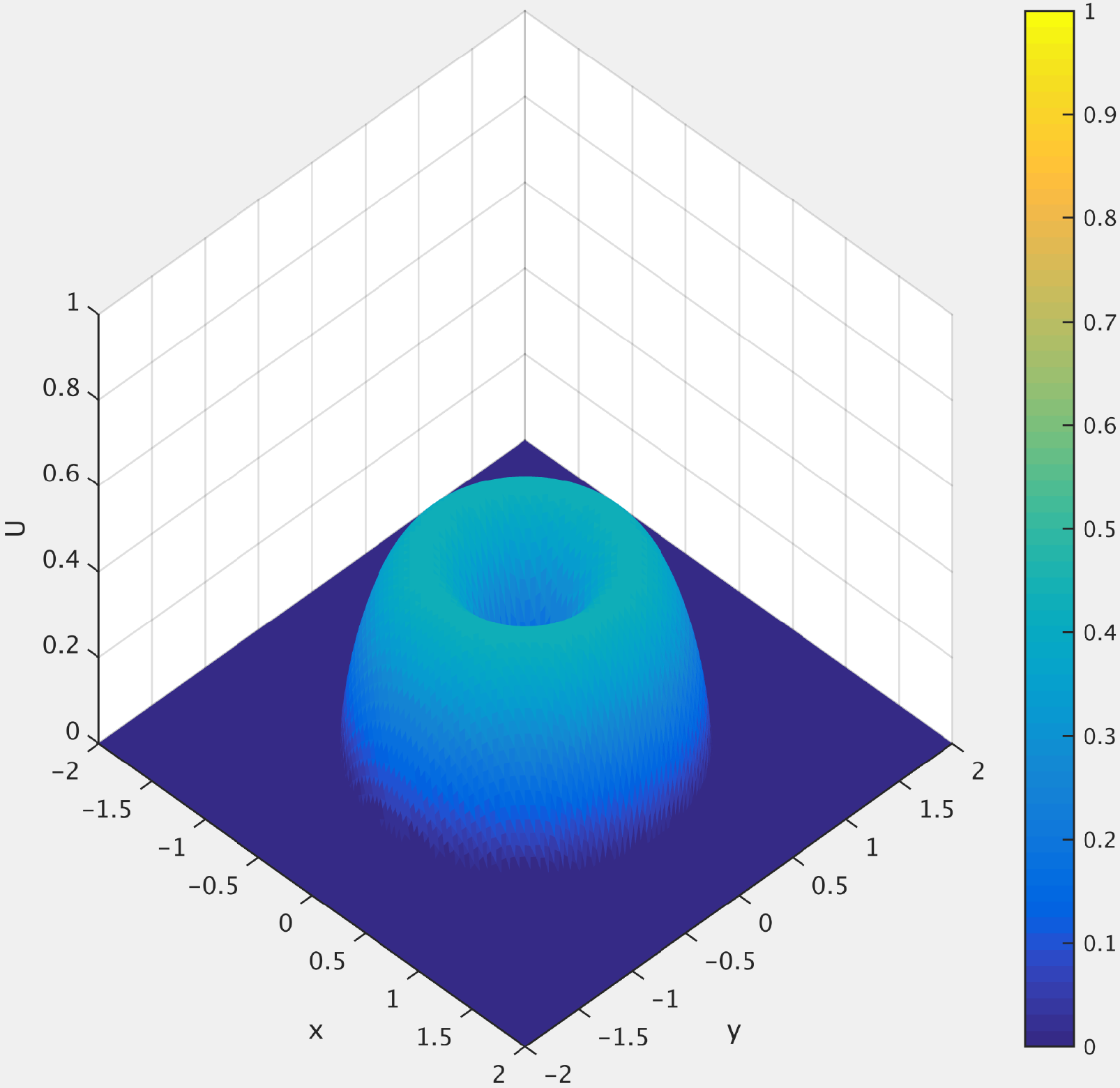}\caption{$t = 0.15$}\end{subfigure}\\%
        \vspace{2mm}
        \begin{subfigure}[b]{0.35\linewidth}\includegraphics[scale=0.34]{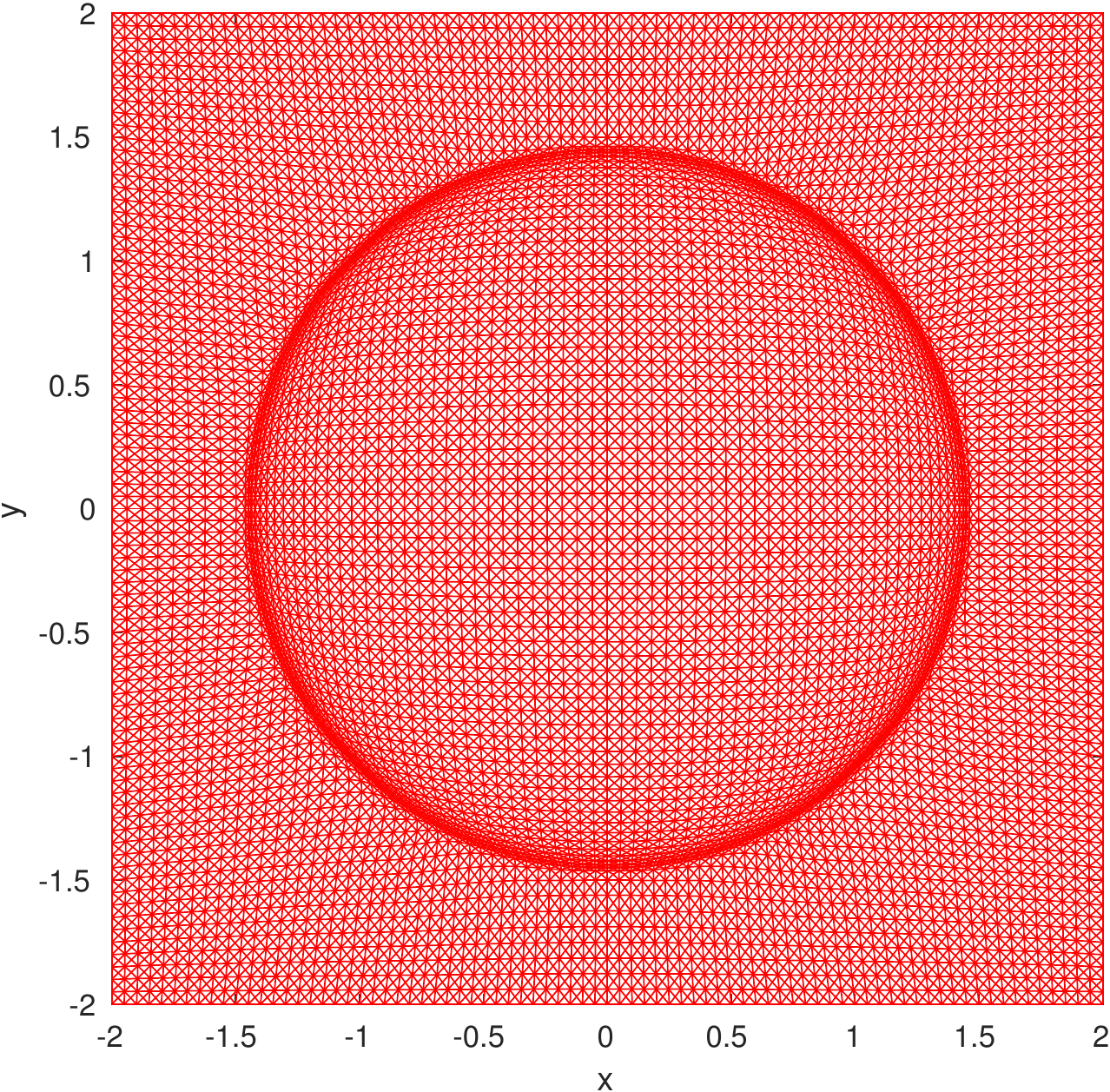}\caption{$t = 0.70$}\end{subfigure}\hspace{5mm}
        \begin{subfigure}[b]{0.35\linewidth}\includegraphics[scale=0.26]{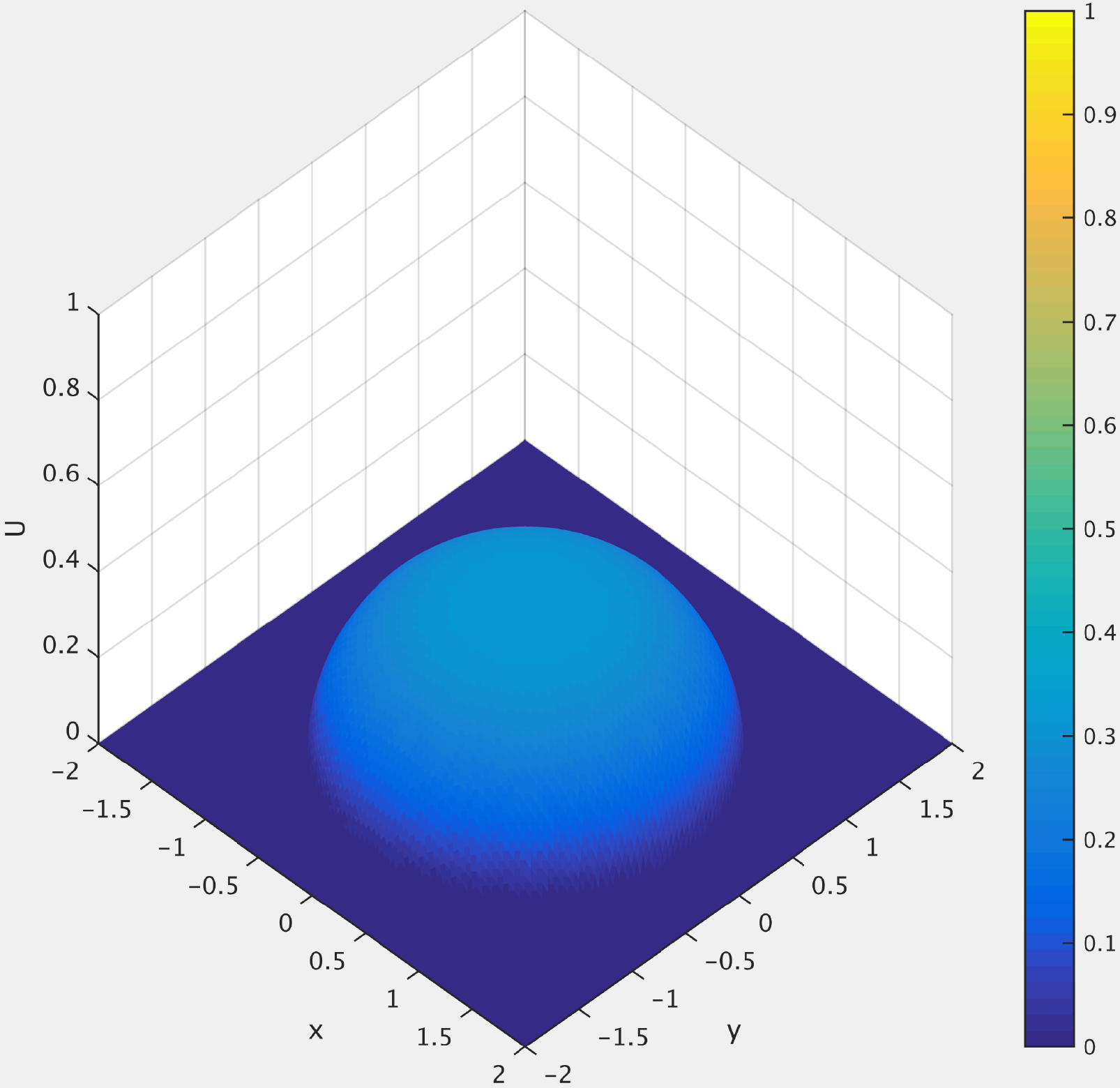}\caption{$t = 0.70$}\end{subfigure}%
        \caption{Example~\ref{exam5.2}. An adaptive mesh and the corresponding solution at various time instants ($N = 25600$).}
        \label{fig:VarExp1-example5.2}
    \end{figure}

\begin{exam}[Variable exponent without absorption]
\label{exam5.2}
For this example,
\begin{align*}
& \lambda = 0,\quad \gamma = \left (\frac{x}{2}\right )^2 + \left (\frac{y}{2}\right )^2 + 1.1,
\quad \Omega = (-2,2)\times (-2,2),
\\
& u_0 = \begin{cases} - \sin(2 \pi \sqrt{x^2 + y^2}),& \text{for} \quad 0.5 < \sqrt{x^2 + y^2} < 1 \\
	0,& \text{otherwise}.
\end{cases}
\end{align*}
This example has been studied in \cite{Duque2013,Duque2015}. The support of the solution
has a hole in the middle which disappears in a finite time.
We take $t \in [0, 0.2]$ in the computation.
An adaptive mesh and the corresponding numerical solution is shown in Fig. \ref{fig:VarExp1-example5.2}.
The result appears to have better resolution than that in \cite{Duque2013} where a uniform mesh has been used.
Moreover, our method works just fine through the closing of the inside hole (cf. Fig. \ref{fig:VarExp1-example5.2})
whereas the method in \cite{Duque2015} which explicitly traces the free boundary encounters
the mesh singularity problem near the time when the hole is closing.
\qed \end{exam}

\begin{figure}[ht]
        \centering
        \begin{subfigure}[b]{0.35\linewidth}\includegraphics[scale=0.34]{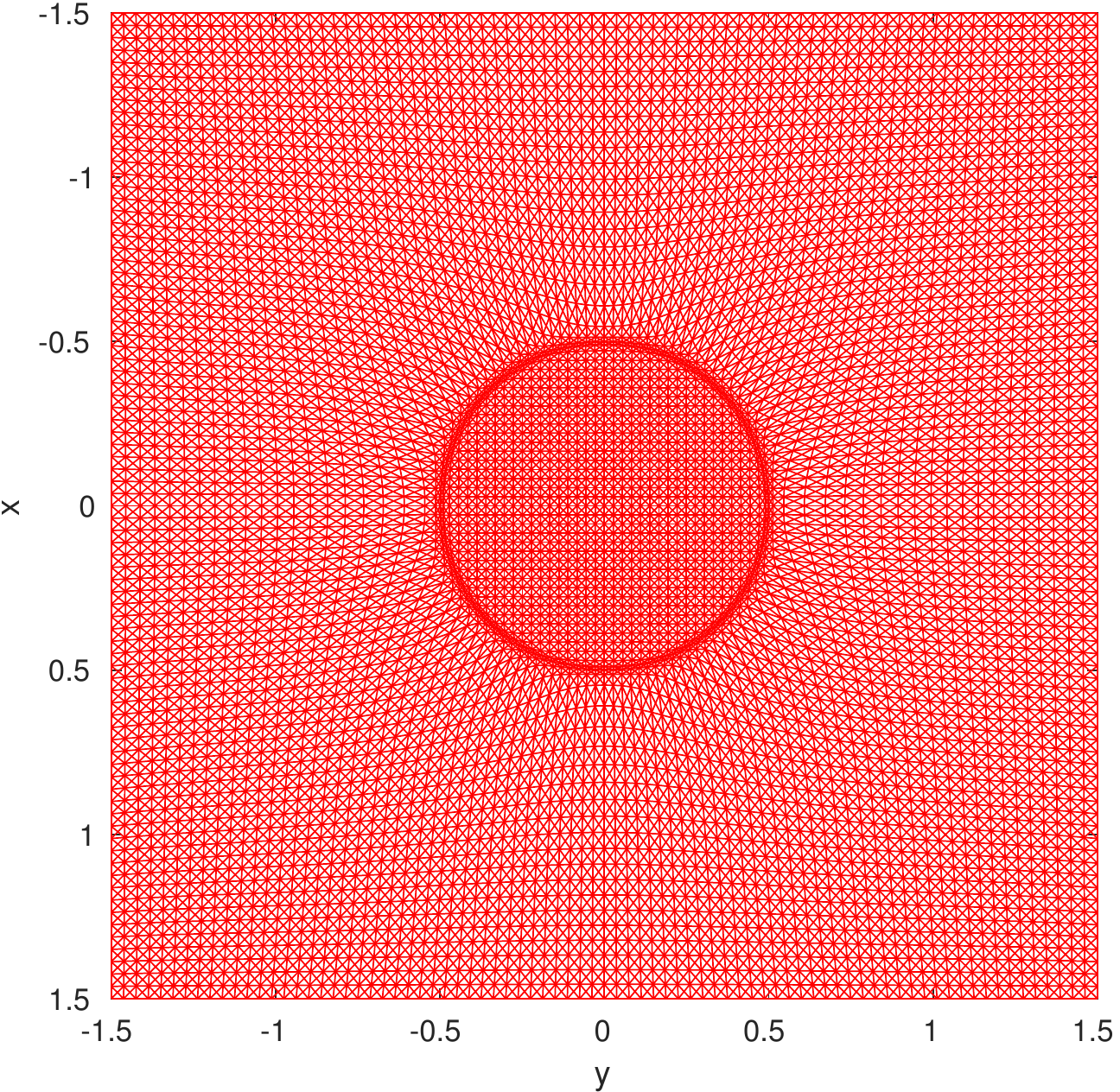}\caption{$t = 0$}\end{subfigure}\hspace{5mm}
        \begin{subfigure}[b]{0.35\linewidth}\includegraphics[scale=0.32]{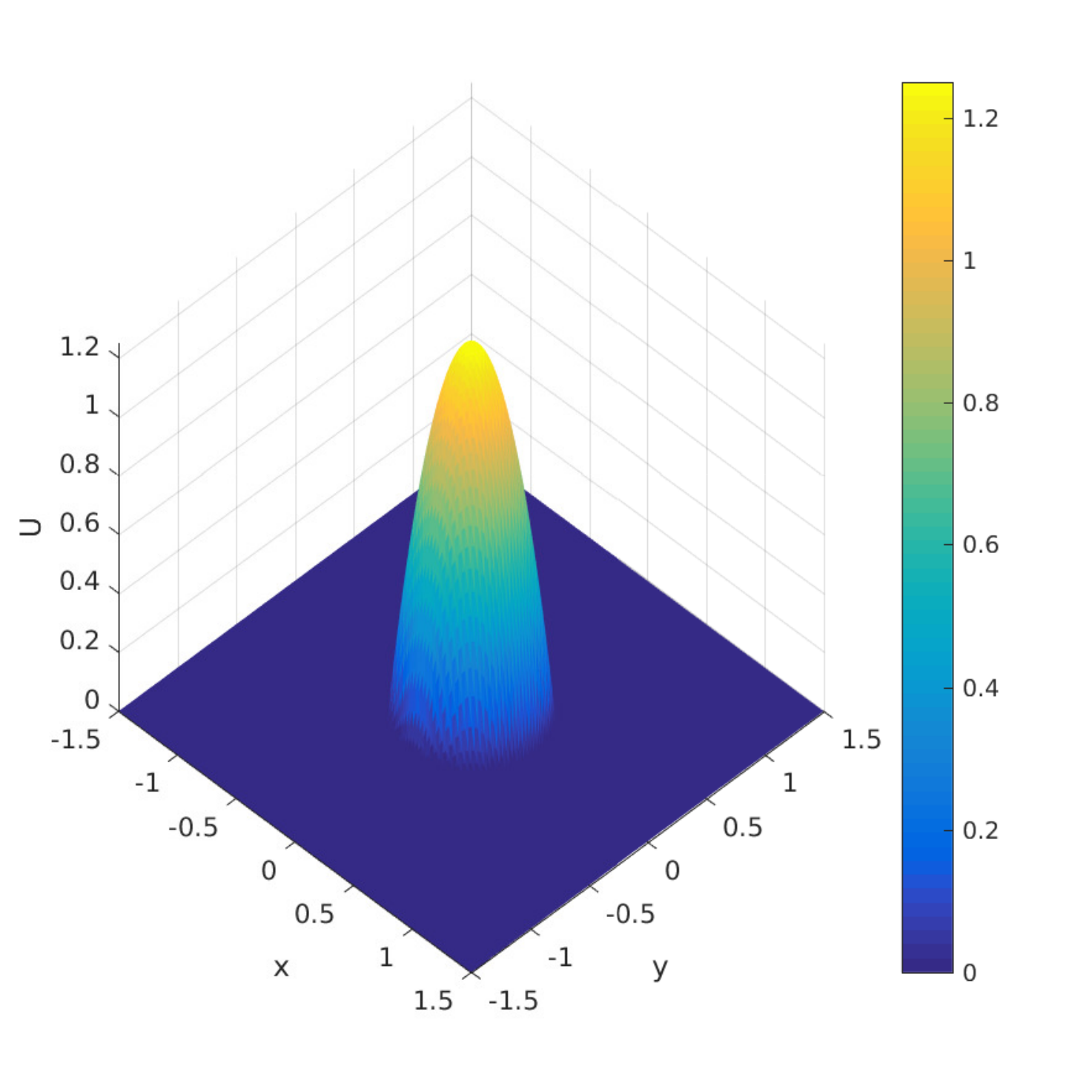}\caption{$t = 0$}\end{subfigure}\\
        \vspace{-1mm}
        \begin{subfigure}[b]{0.35\linewidth}\includegraphics[scale=0.34]{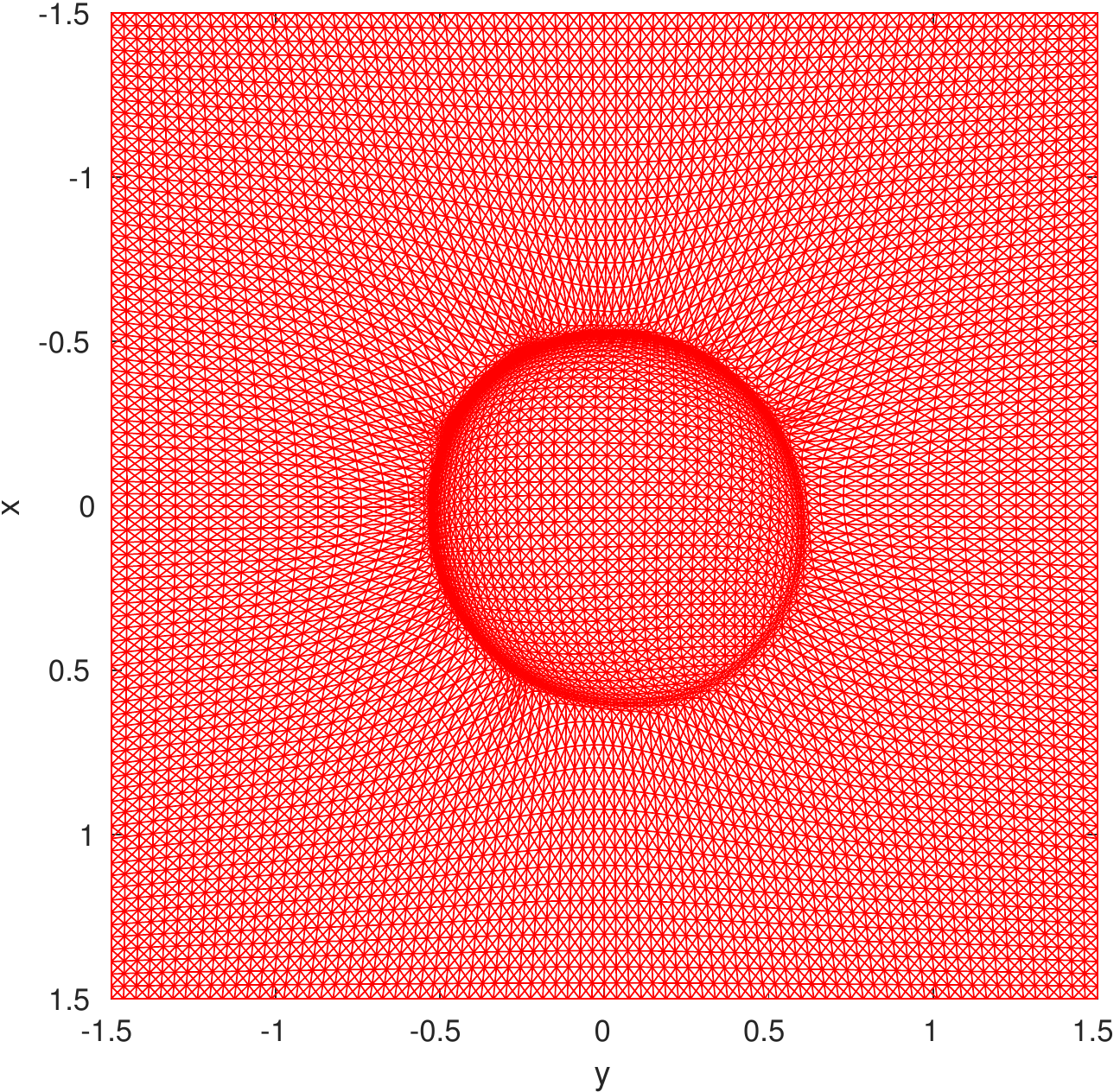}\caption{$t = 0.06$}\end{subfigure}\hspace{5mm}
        \begin{subfigure}[b]{0.35\linewidth}\includegraphics[scale=0.32]{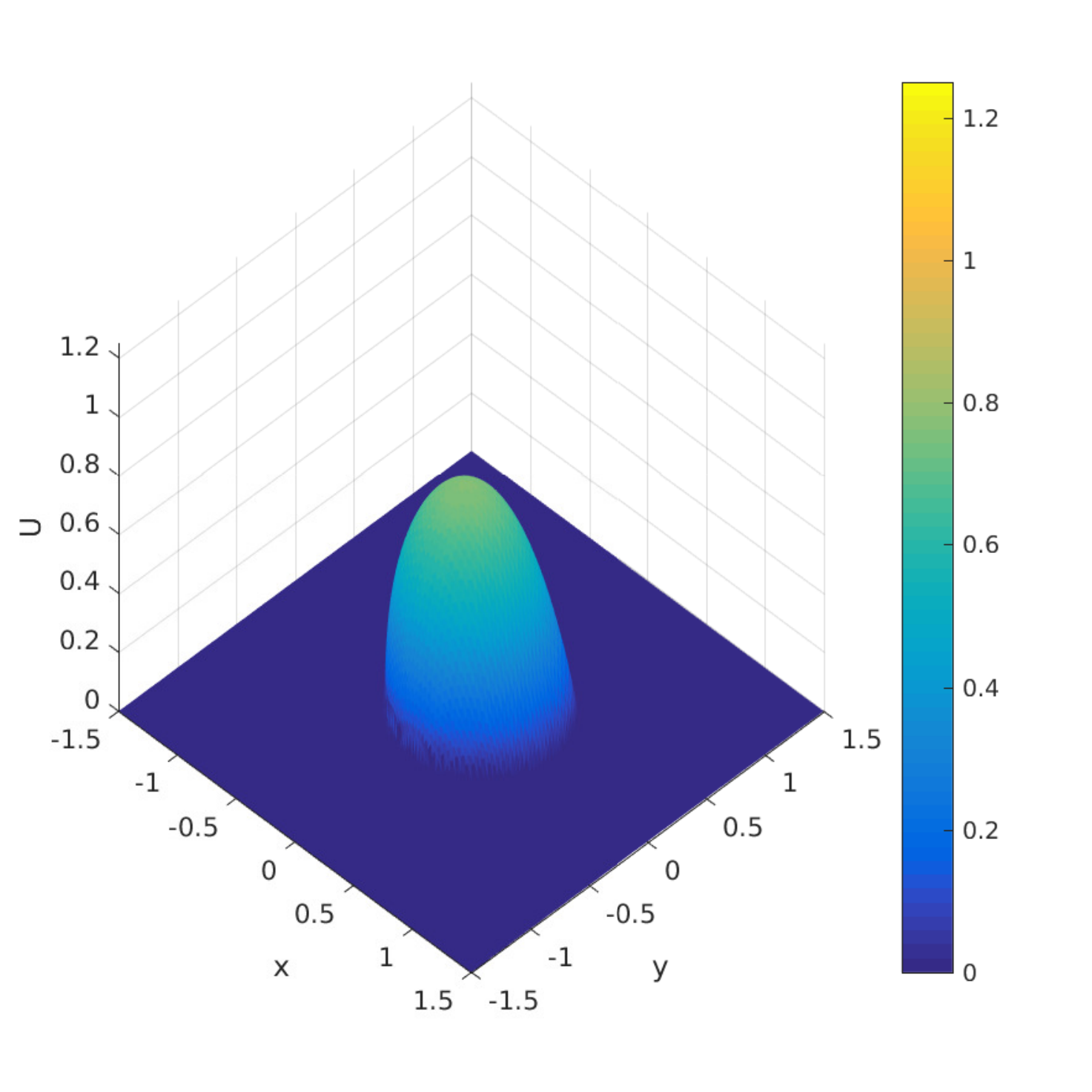}\caption{$t = 0.06$}\end{subfigure}\\%
        \vspace{-1mm}
        \begin{subfigure}[b]{0.35\linewidth}\includegraphics[scale=0.34]{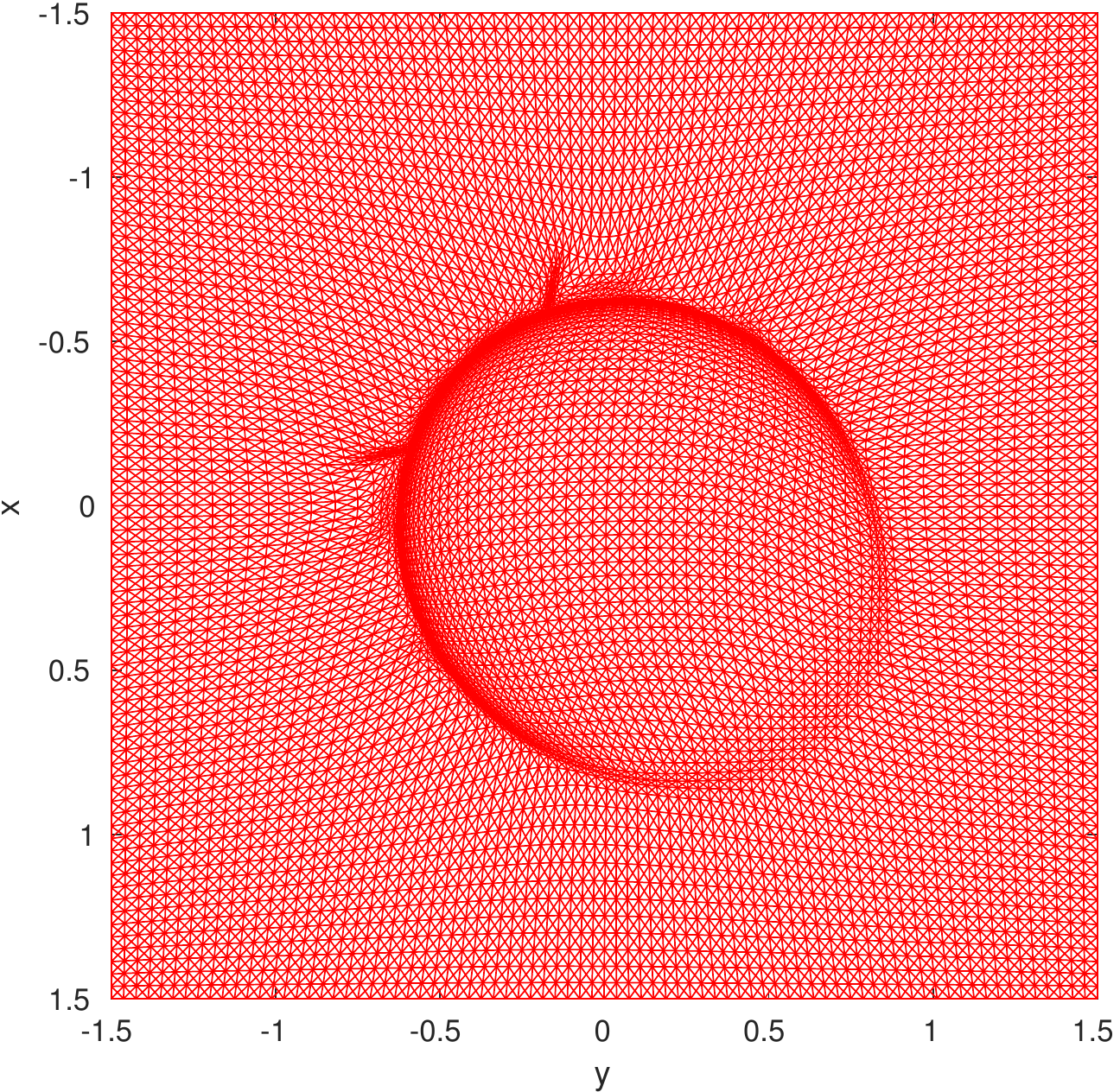}\caption{$t = 0.28$}\end{subfigure}\hspace{5mm}
        \begin{subfigure}[b]{0.35\linewidth}\includegraphics[scale=0.32]{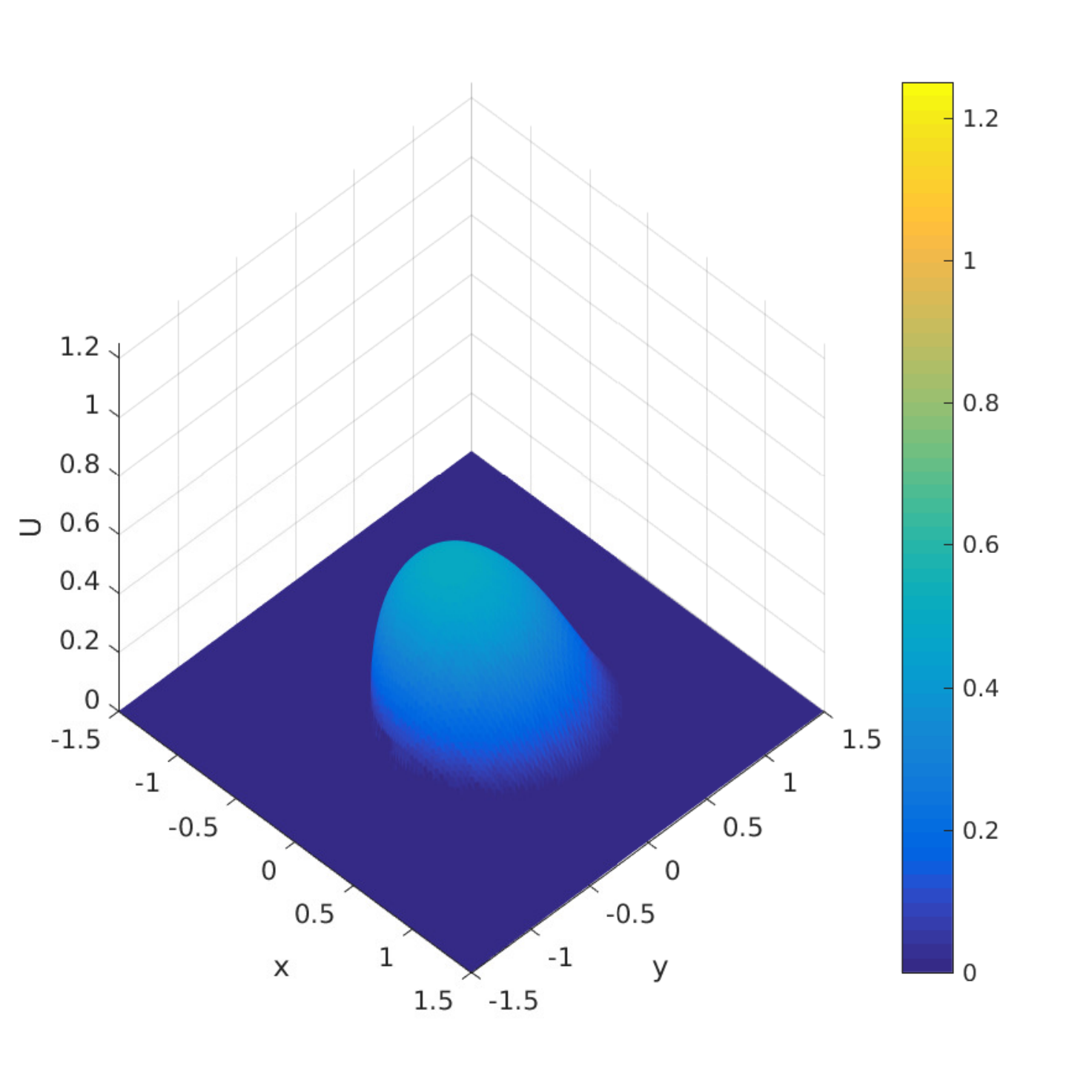}\caption{$t = 0.28$}\end{subfigure}\\
        \vspace{-1mm}
        \begin{subfigure}[b]{0.35\linewidth}\includegraphics[scale=0.34]{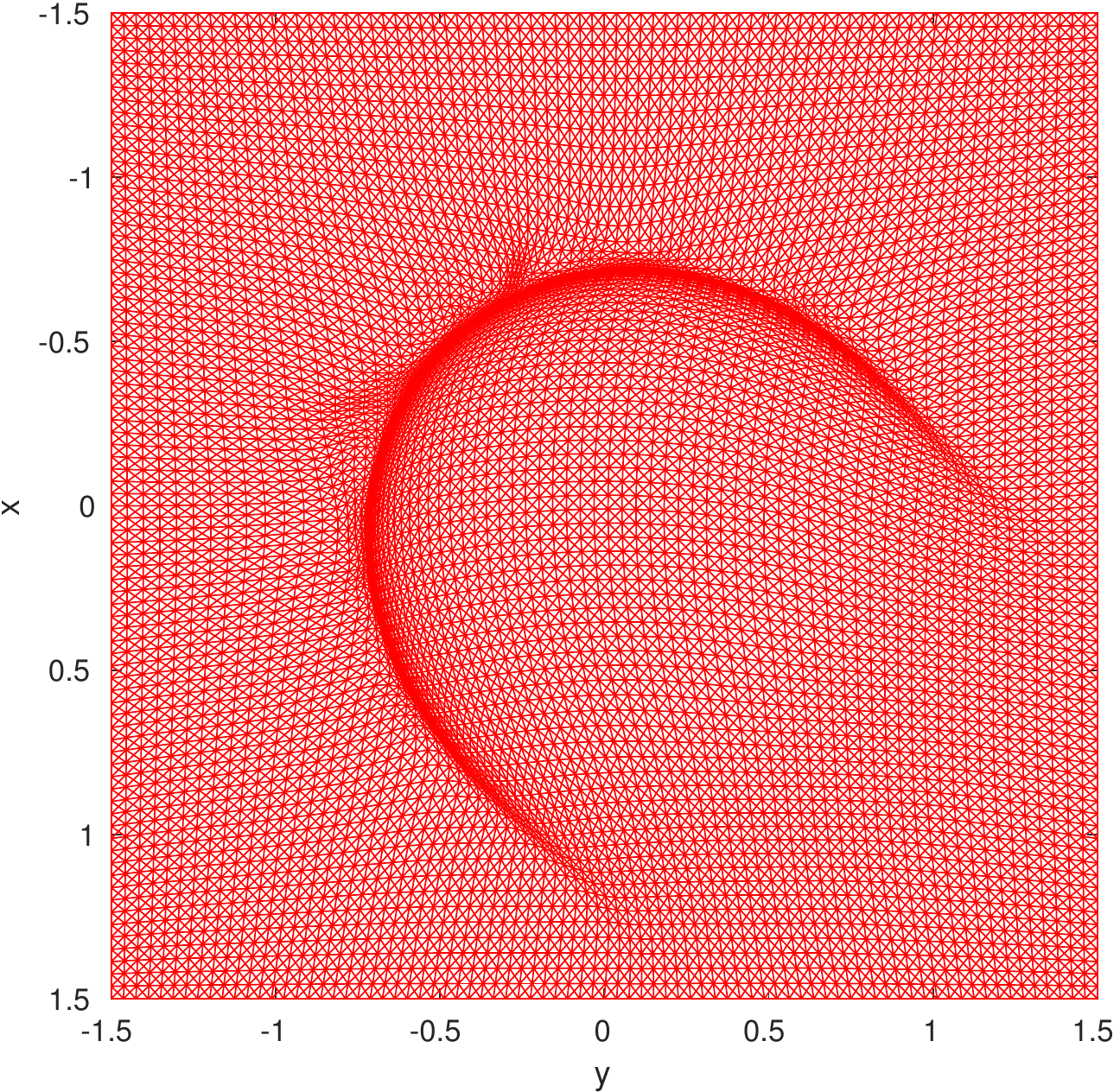}\caption{$t = 0.90$}\end{subfigure}\hspace{5mm}
        \begin{subfigure}[b]{0.35\linewidth}\includegraphics[scale=0.32]{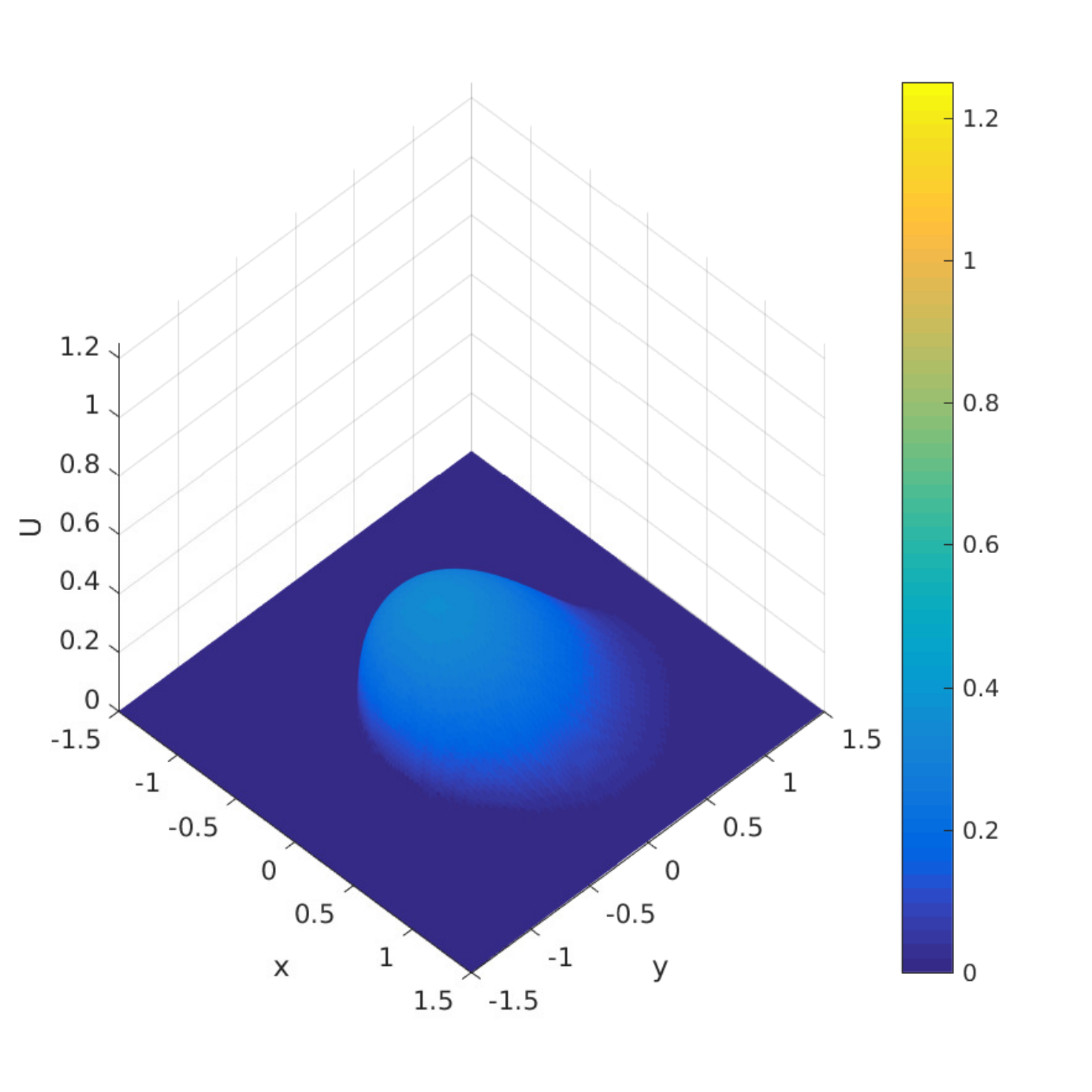}\caption{$t = 0.90$}\end{subfigure}%
        \caption{Example~\ref{exam5.3}. An adaptive mesh and the corresponding solution at various time instants ($N = 25600$).}
        \label{fig:VarExp2-example5.3}
    \end{figure}

\begin{figure}[ht]
        \centering
        \begin{subfigure}[b]{0.25\linewidth}\includegraphics[scale=0.20]{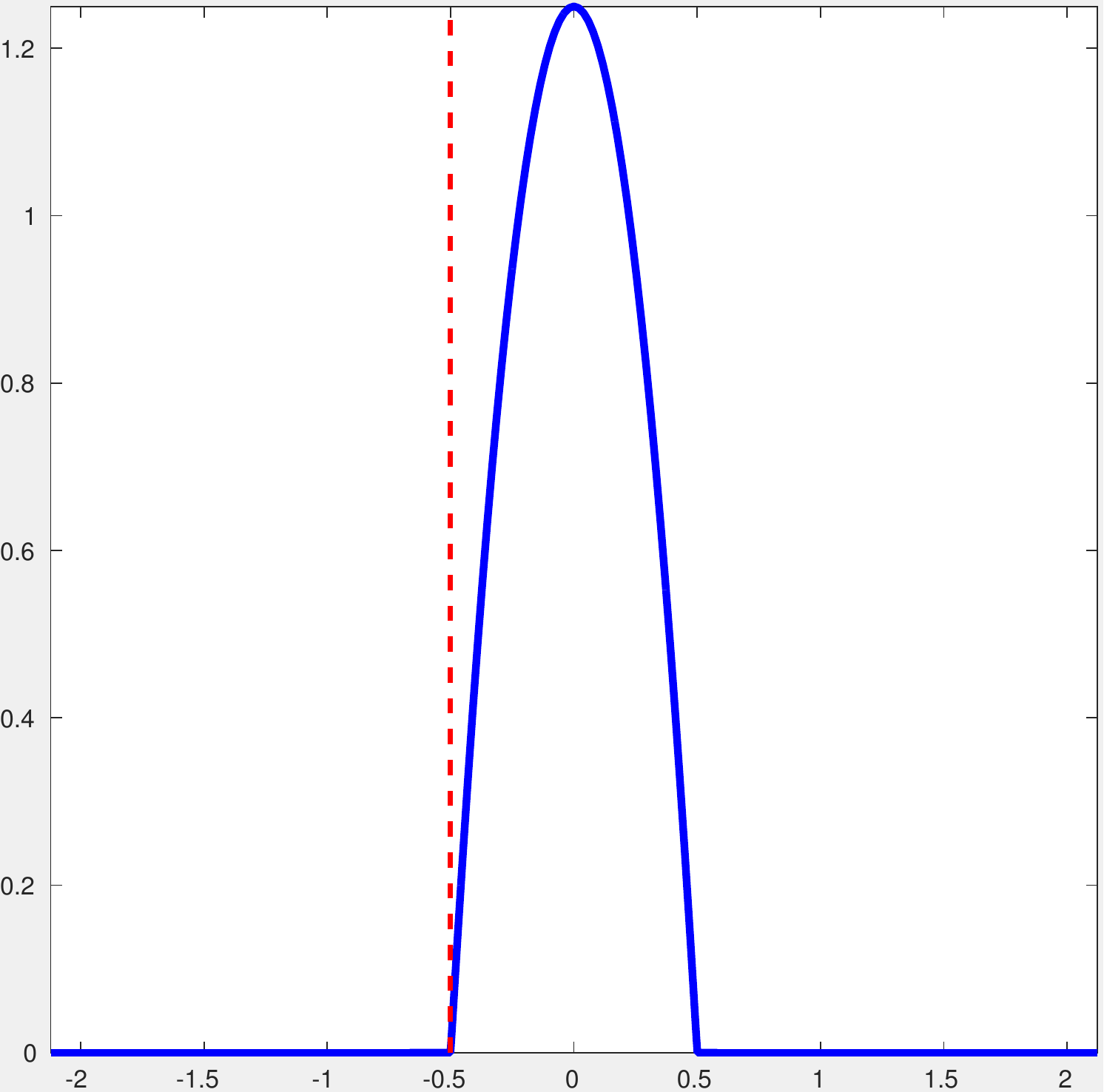}\caption{$t = 0$}\end{subfigure}%
        \begin{subfigure}[b]{0.25\linewidth}\includegraphics[scale=0.20]{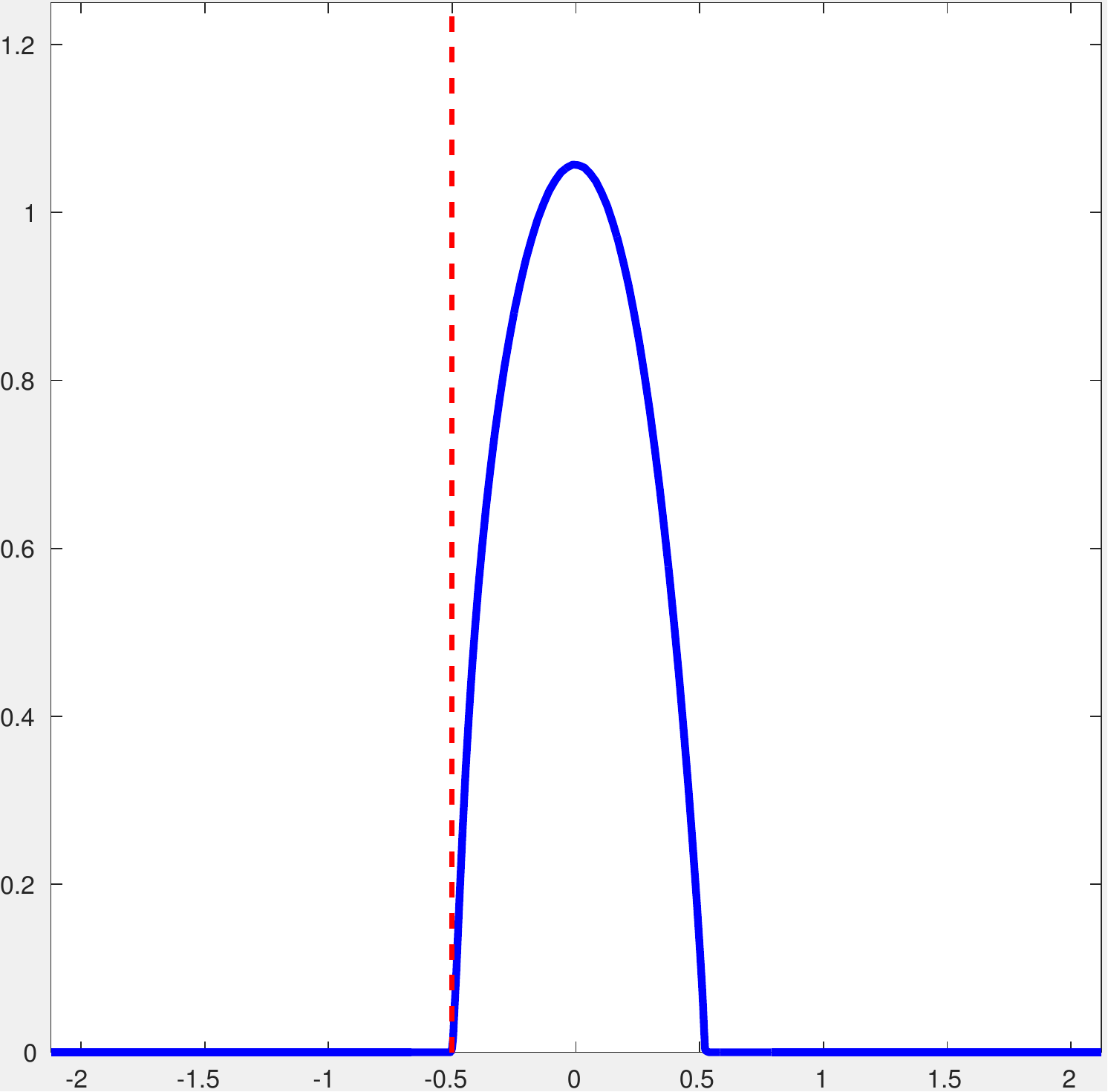}\caption{$t = 0.01$}\end{subfigure}%
        \begin{subfigure}[b]{0.25\linewidth}\includegraphics[scale=0.20]{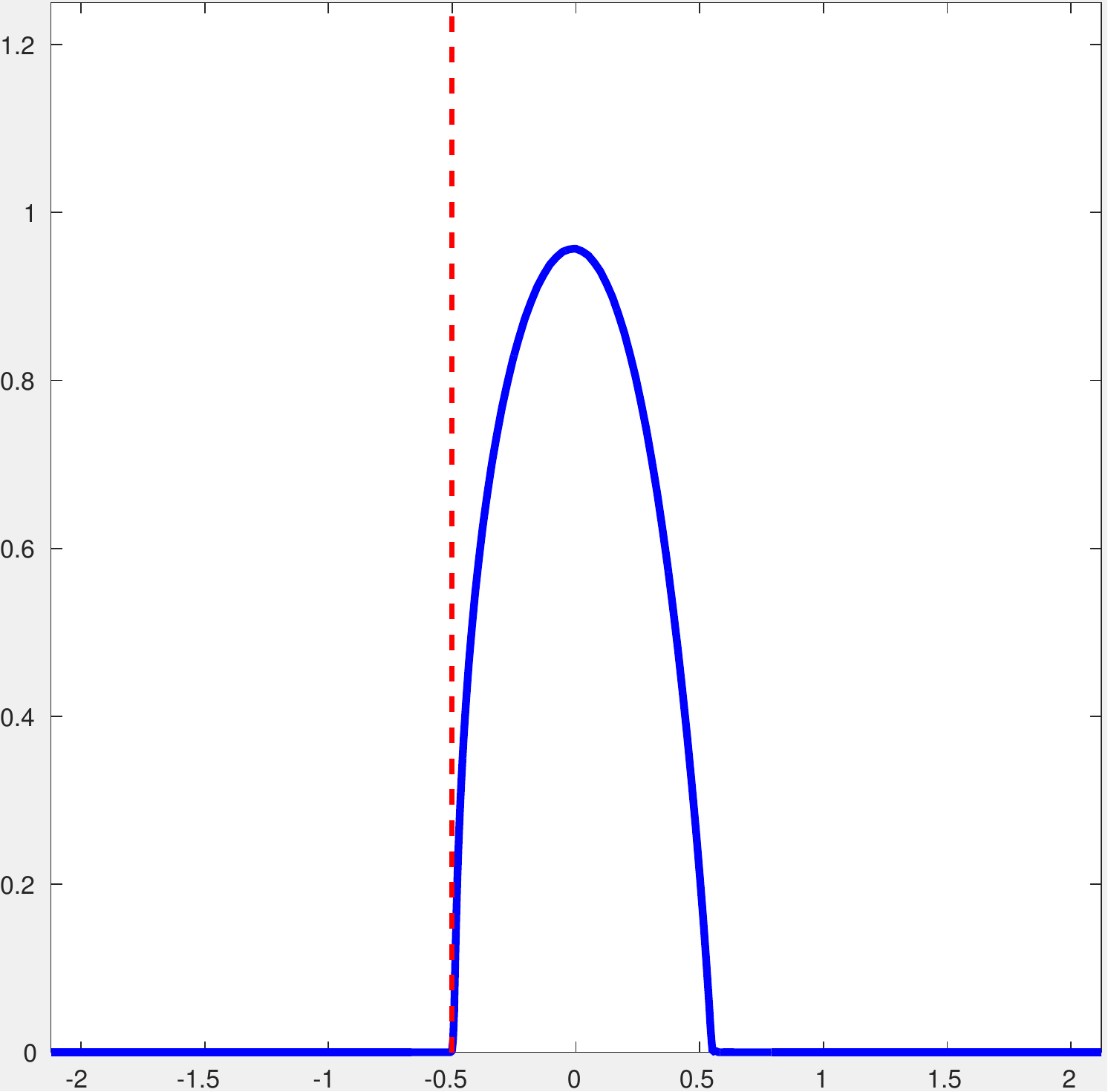}\caption{$t = 0.02$}\end{subfigure}%
        \begin{subfigure}[b]{0.25\linewidth}\includegraphics[scale=0.20]{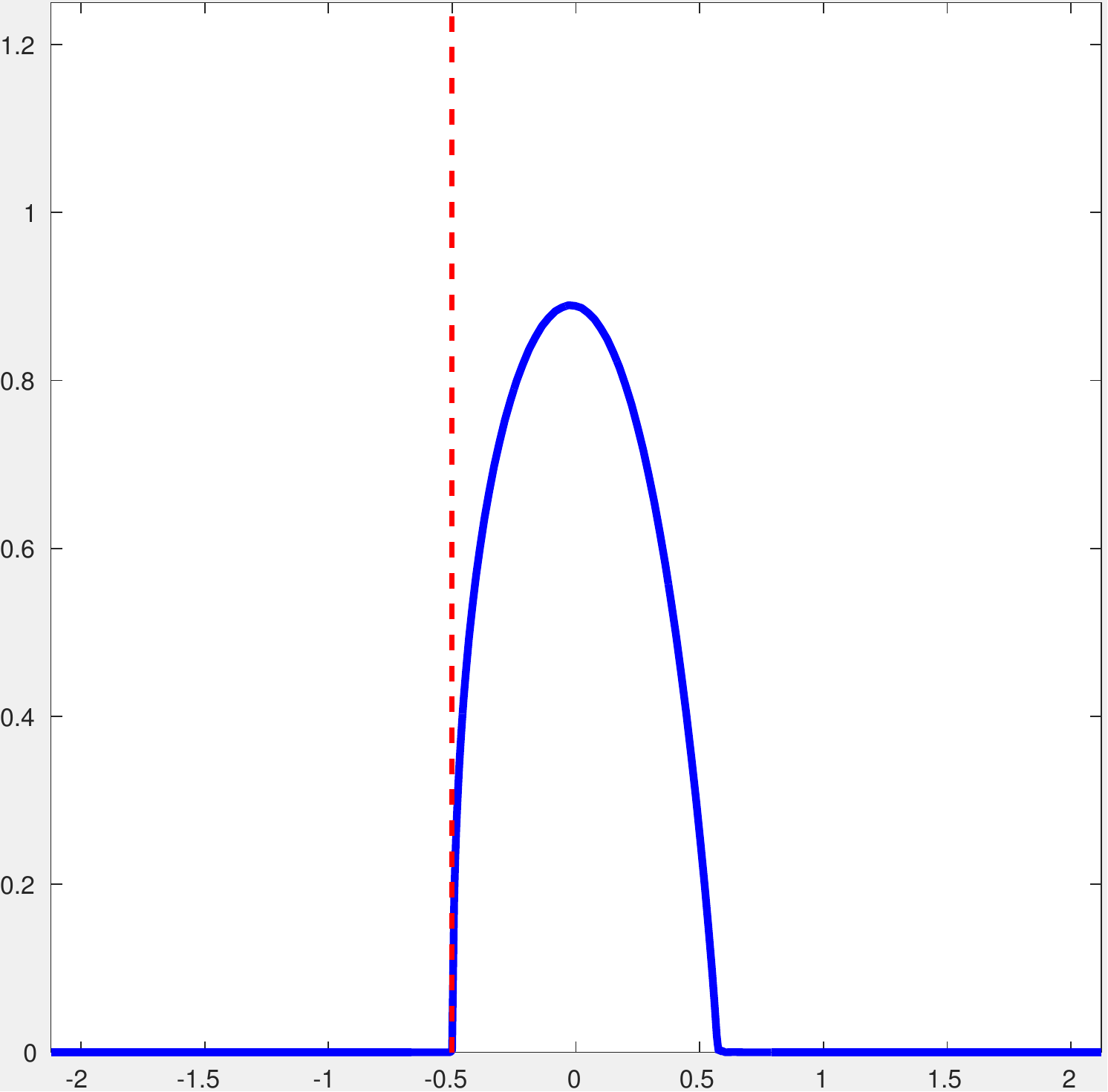}\caption{$t = 0.03$}\end{subfigure}\\%
        \begin{subfigure}[b]{0.25\linewidth}\includegraphics[scale=0.20]{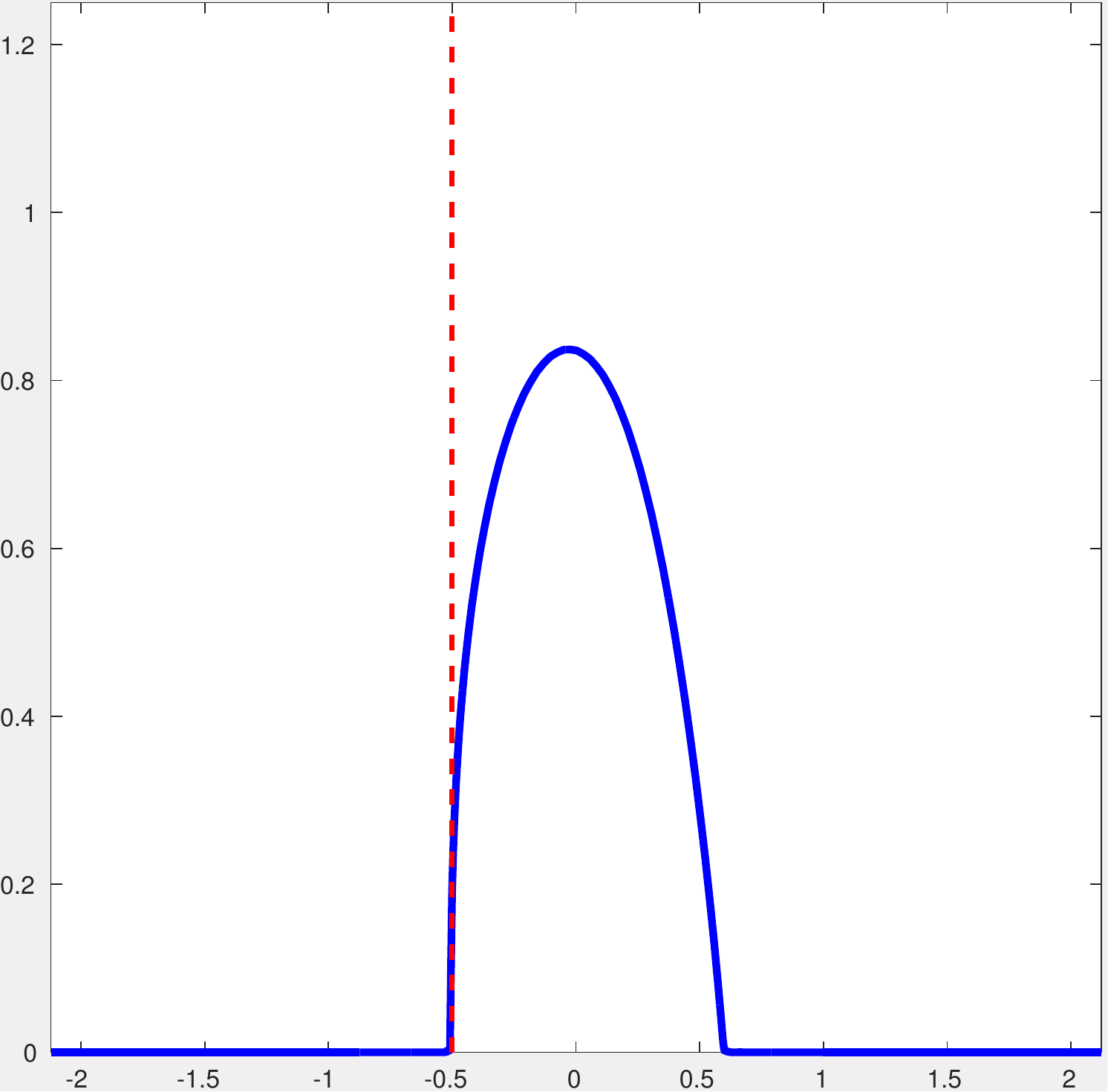}\caption{$t = 0.04$}\end{subfigure}%
        \begin{subfigure}[b]{0.25\linewidth}\includegraphics[scale=0.20]{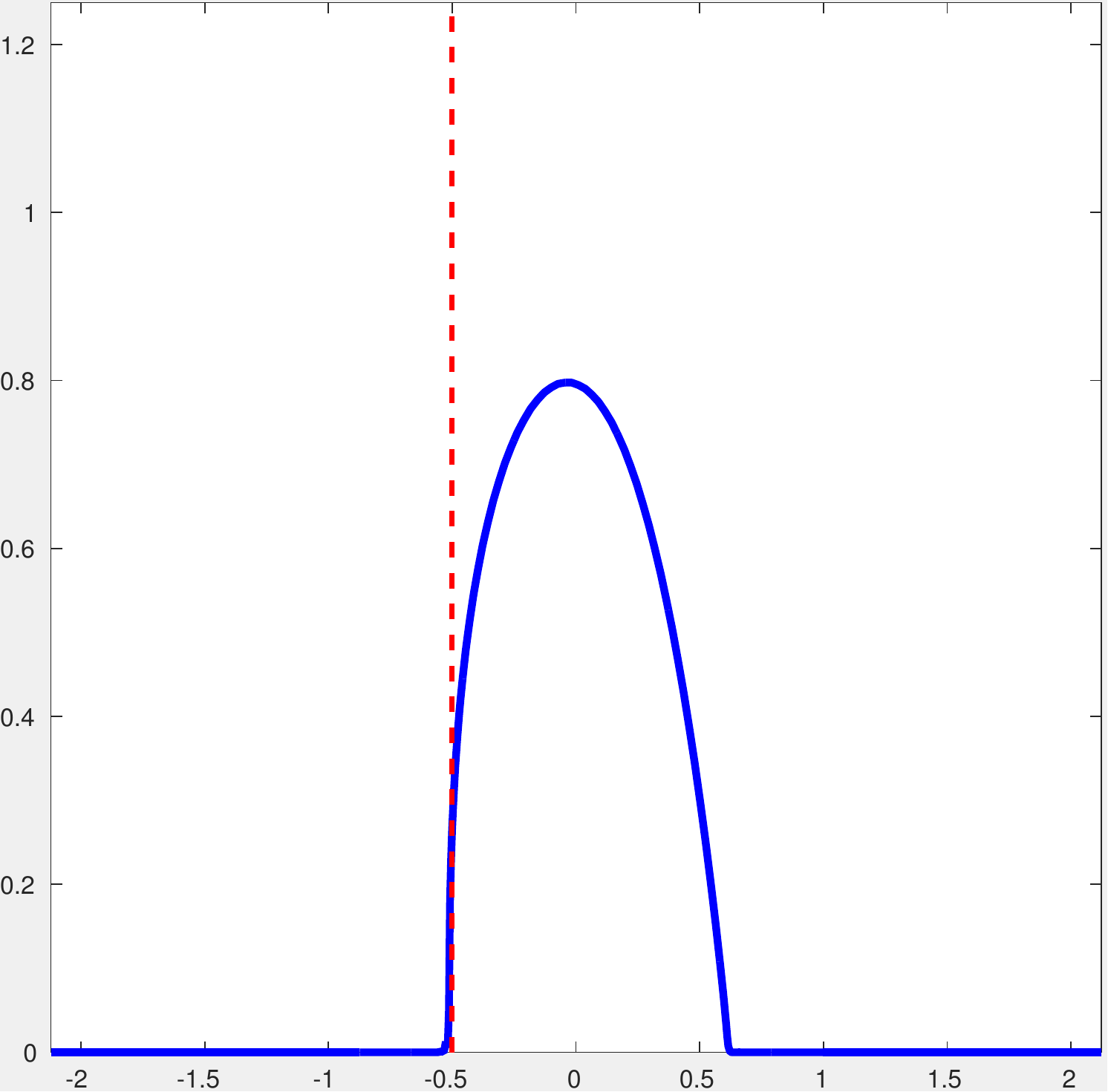}\caption{$t = 0.05$}\end{subfigure}%
        \begin{subfigure}[b]{0.25\linewidth}\includegraphics[scale=0.20]{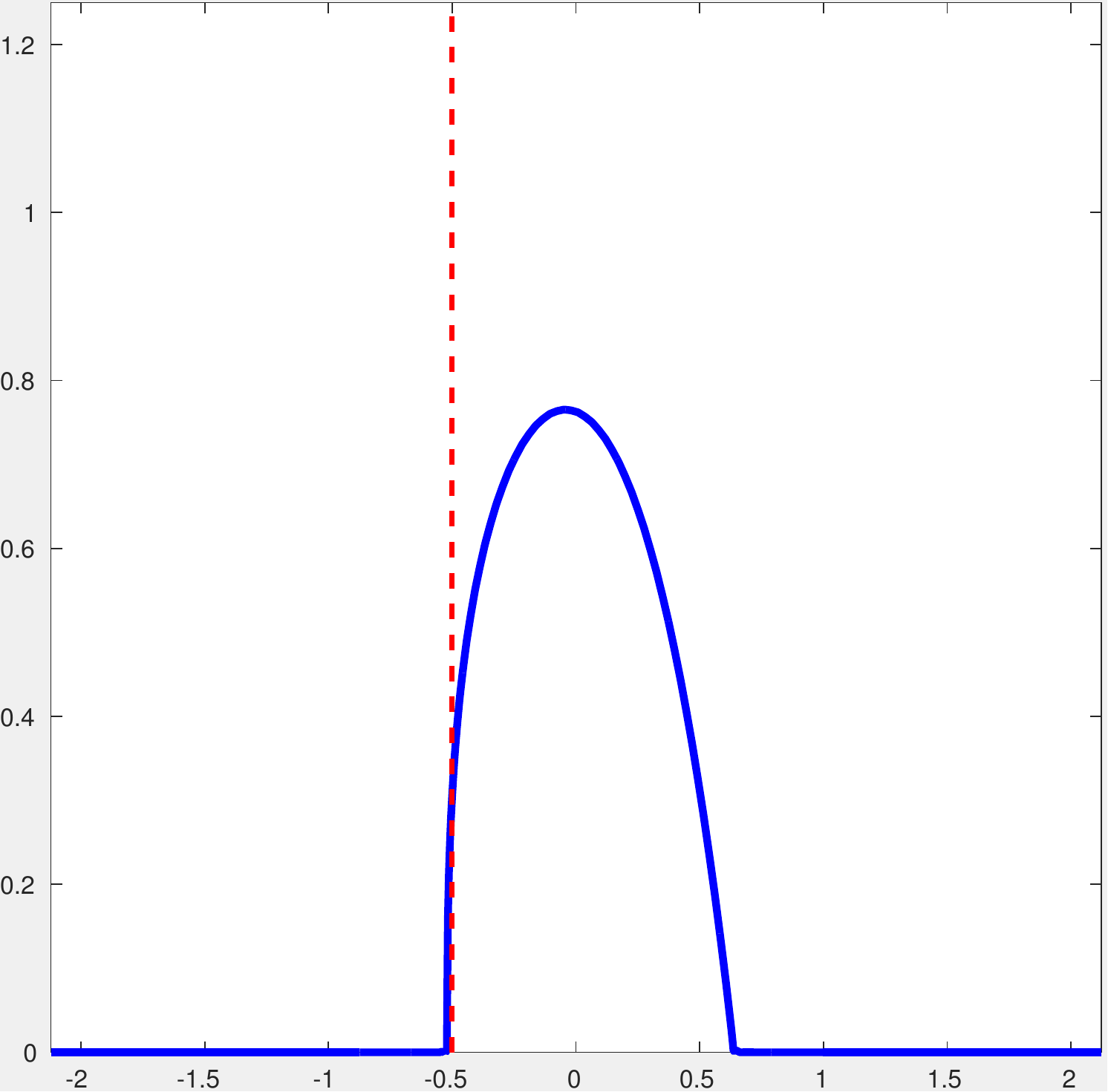}\caption{$t = 0.06$}\end{subfigure}%
        \begin{subfigure}[b]{0.25\linewidth}\includegraphics[scale=0.20]{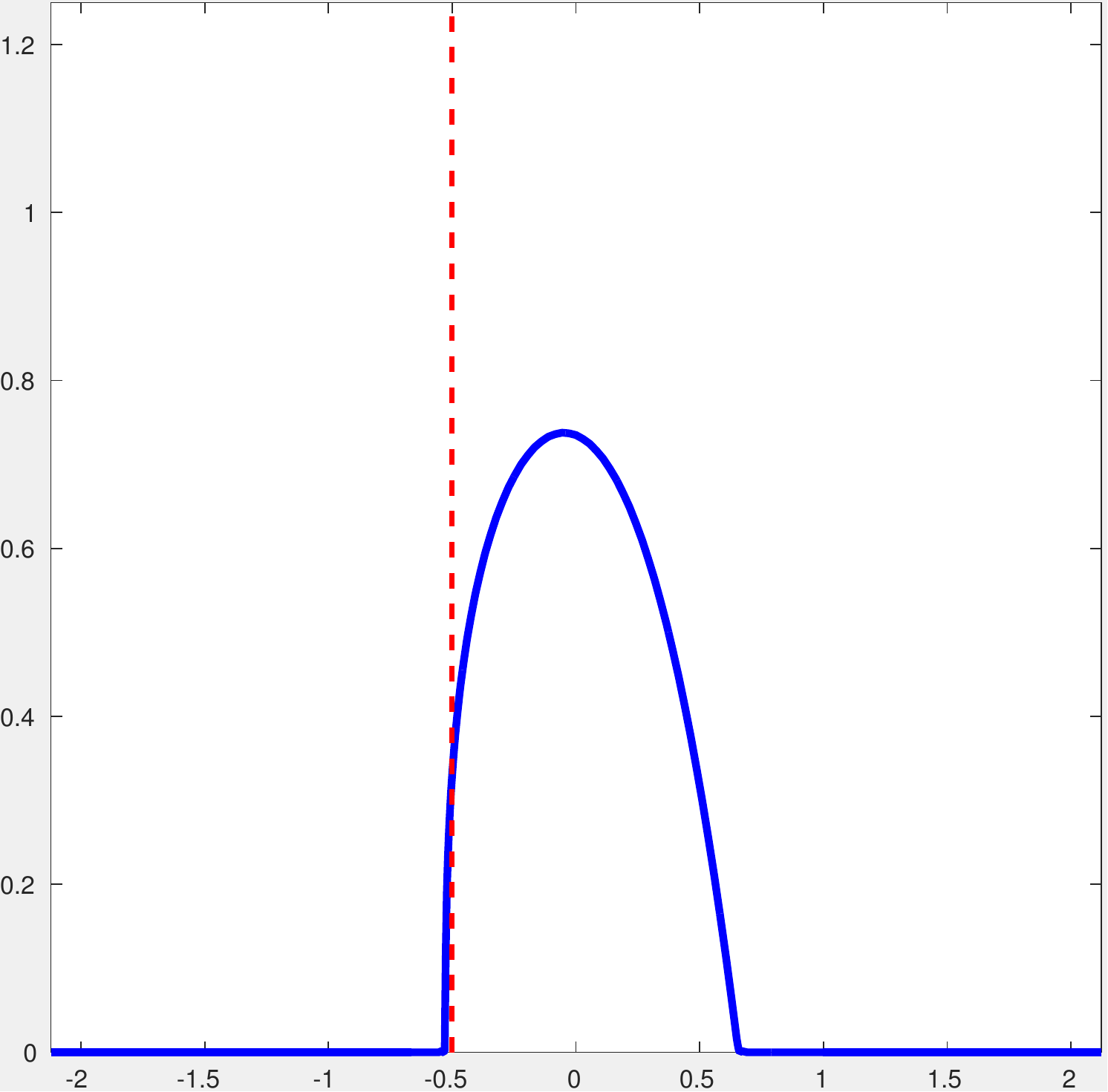}\caption{$t = 0.07$}\end{subfigure}\\%
        \begin{subfigure}[b]{0.25\linewidth}\includegraphics[scale=0.20]{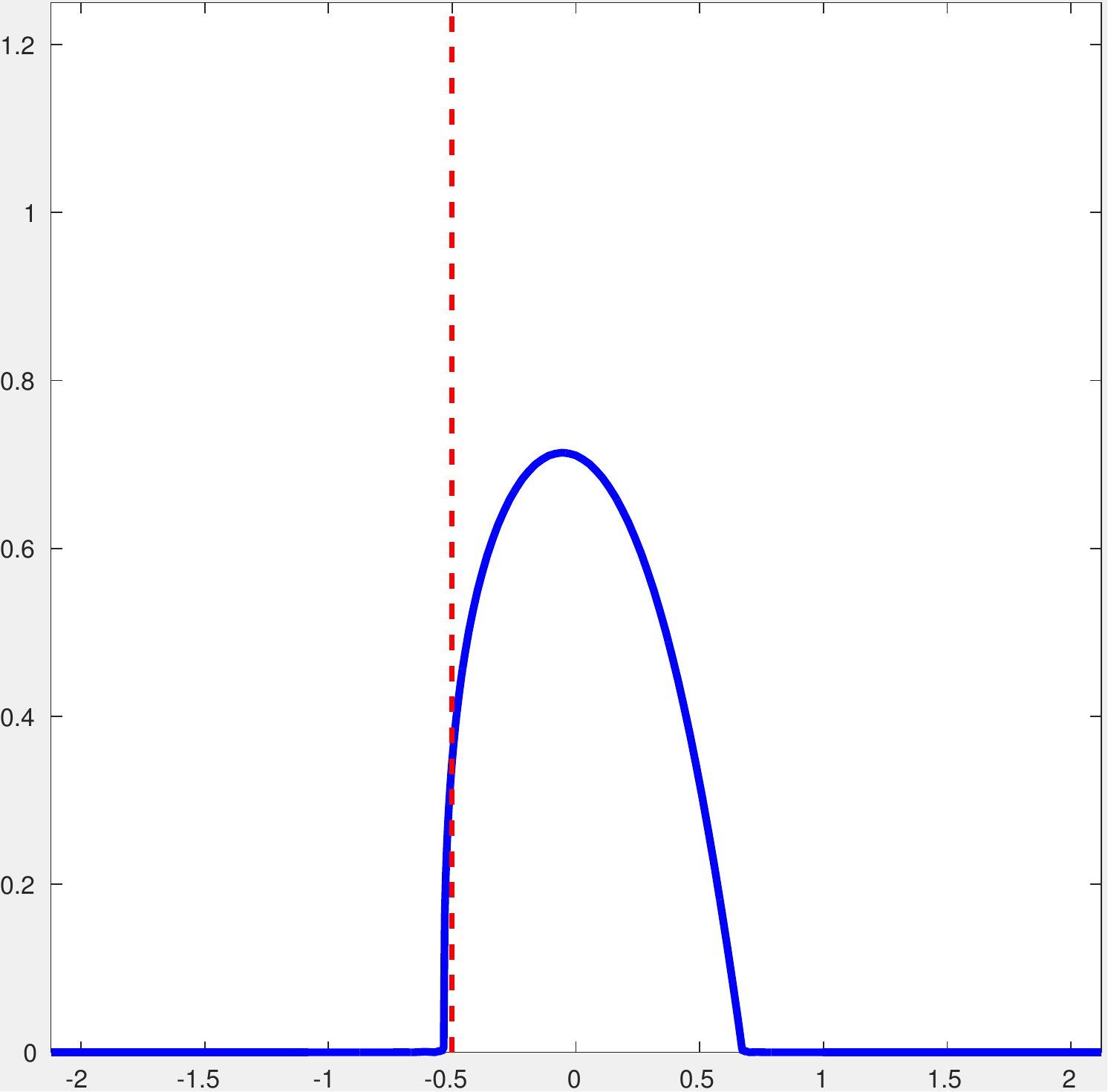}\caption{$t = 0.15$}\end{subfigure}%
        \begin{subfigure}[b]{0.25\linewidth}\includegraphics[scale=0.20]{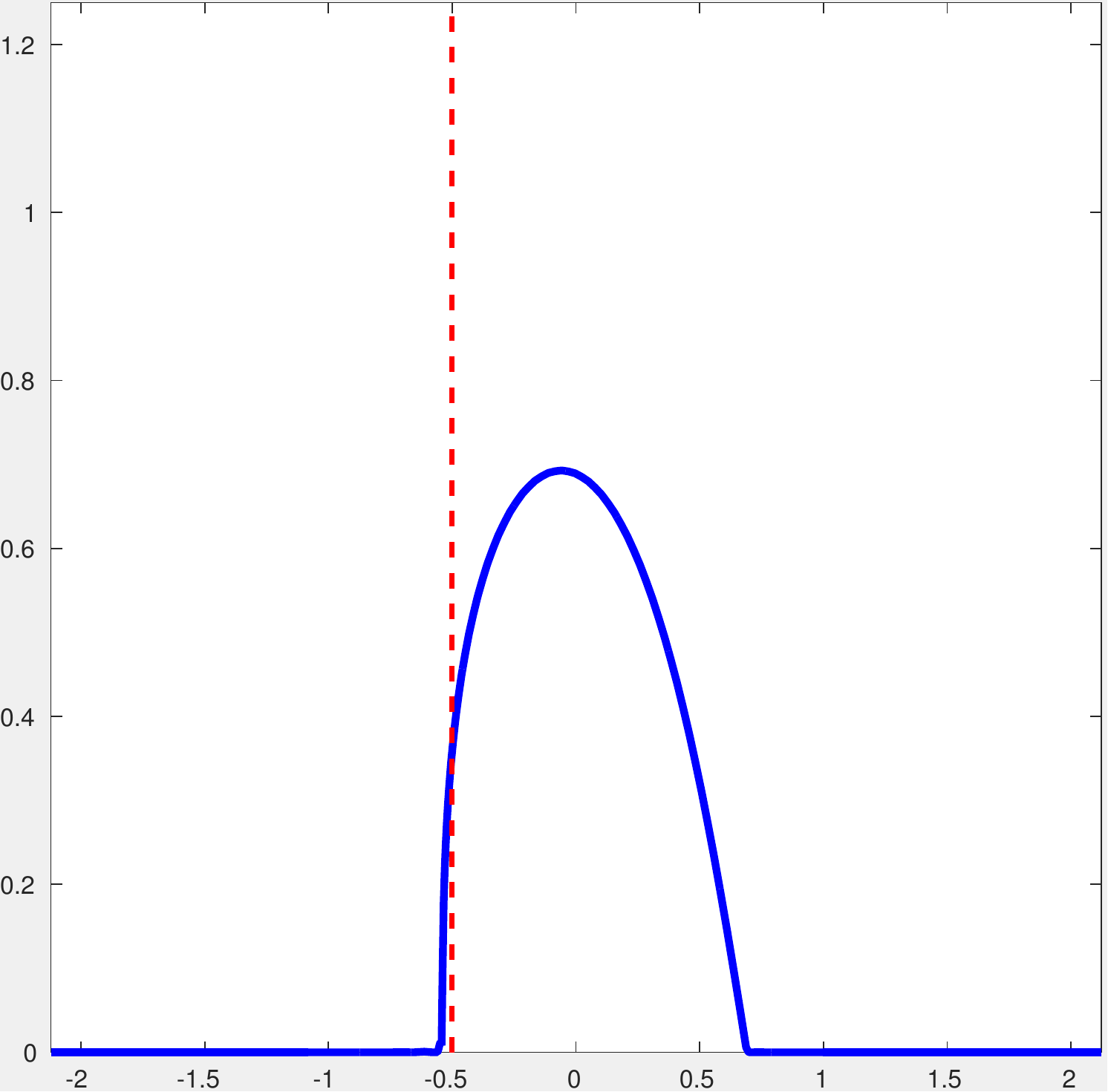}\caption{$t = 0.2$}\end{subfigure}%
        \begin{subfigure}[b]{0.25\linewidth}\includegraphics[scale=0.20]{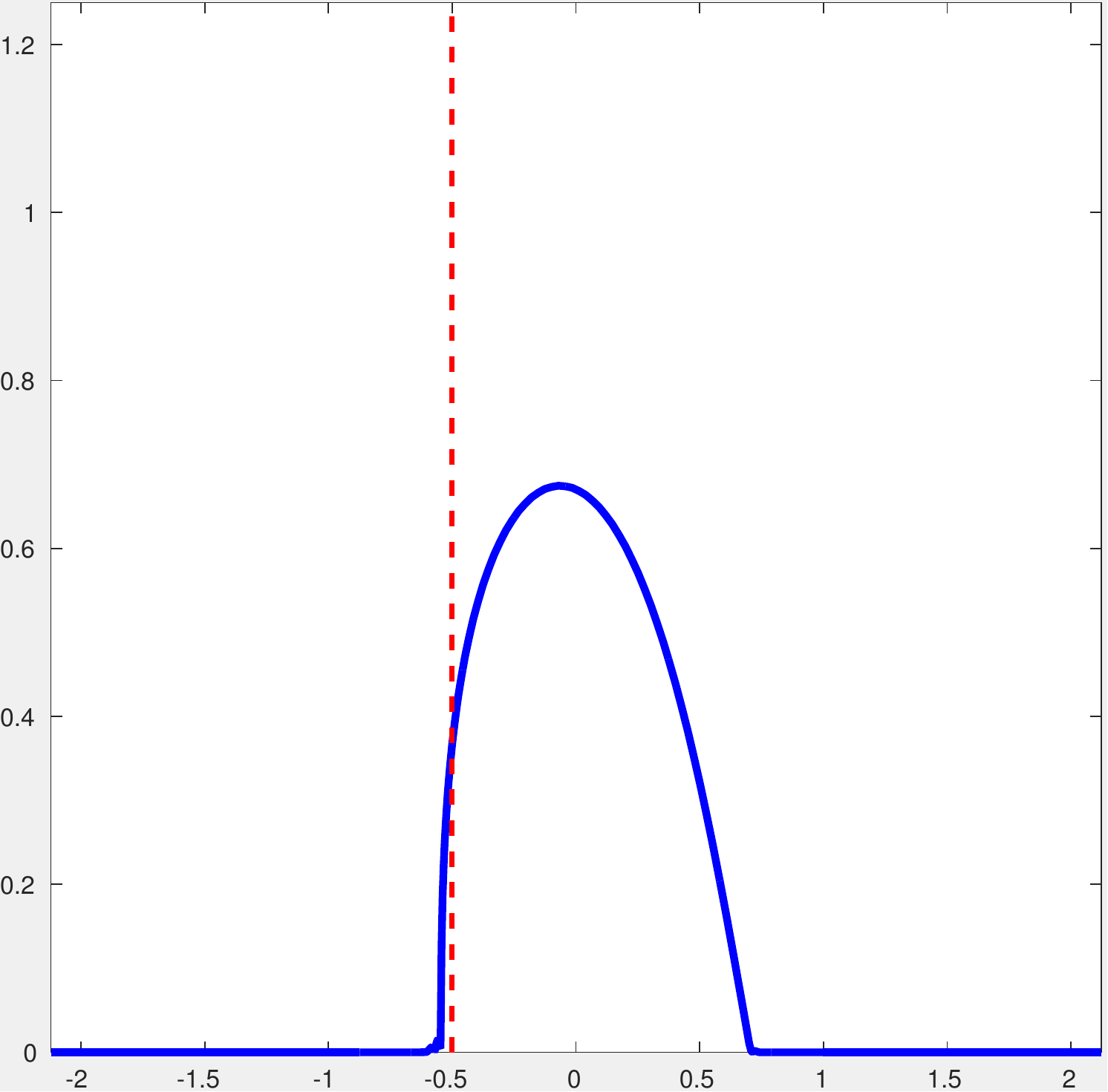}\caption{$t = 0.3$}\end{subfigure}%
        \begin{subfigure}[b]{0.25\linewidth}\includegraphics[scale=0.20]{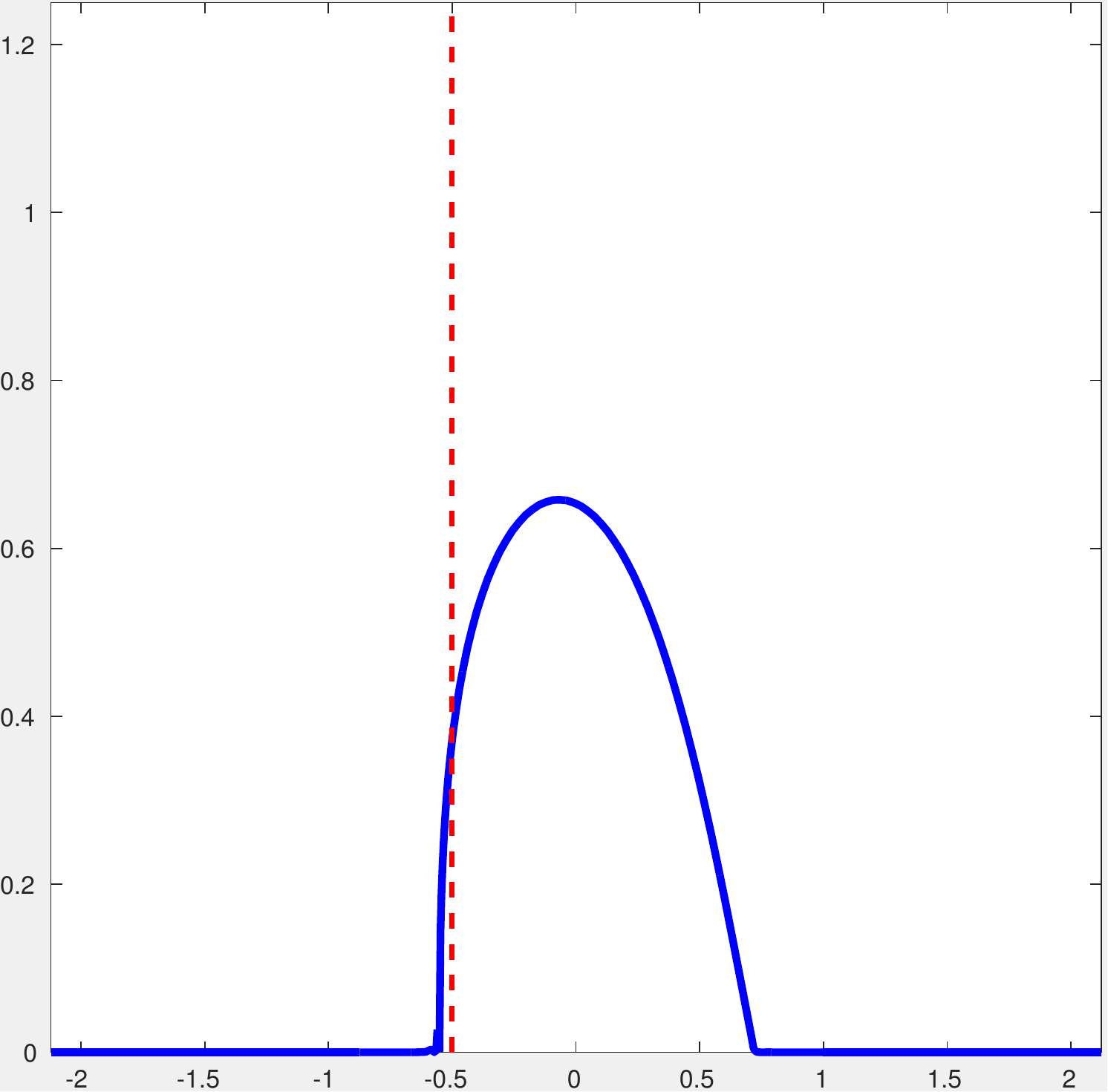}\caption{$t = 0.4$}\end{subfigure}\\%
        \begin{subfigure}[b]{0.25\linewidth}\includegraphics[scale=0.20]{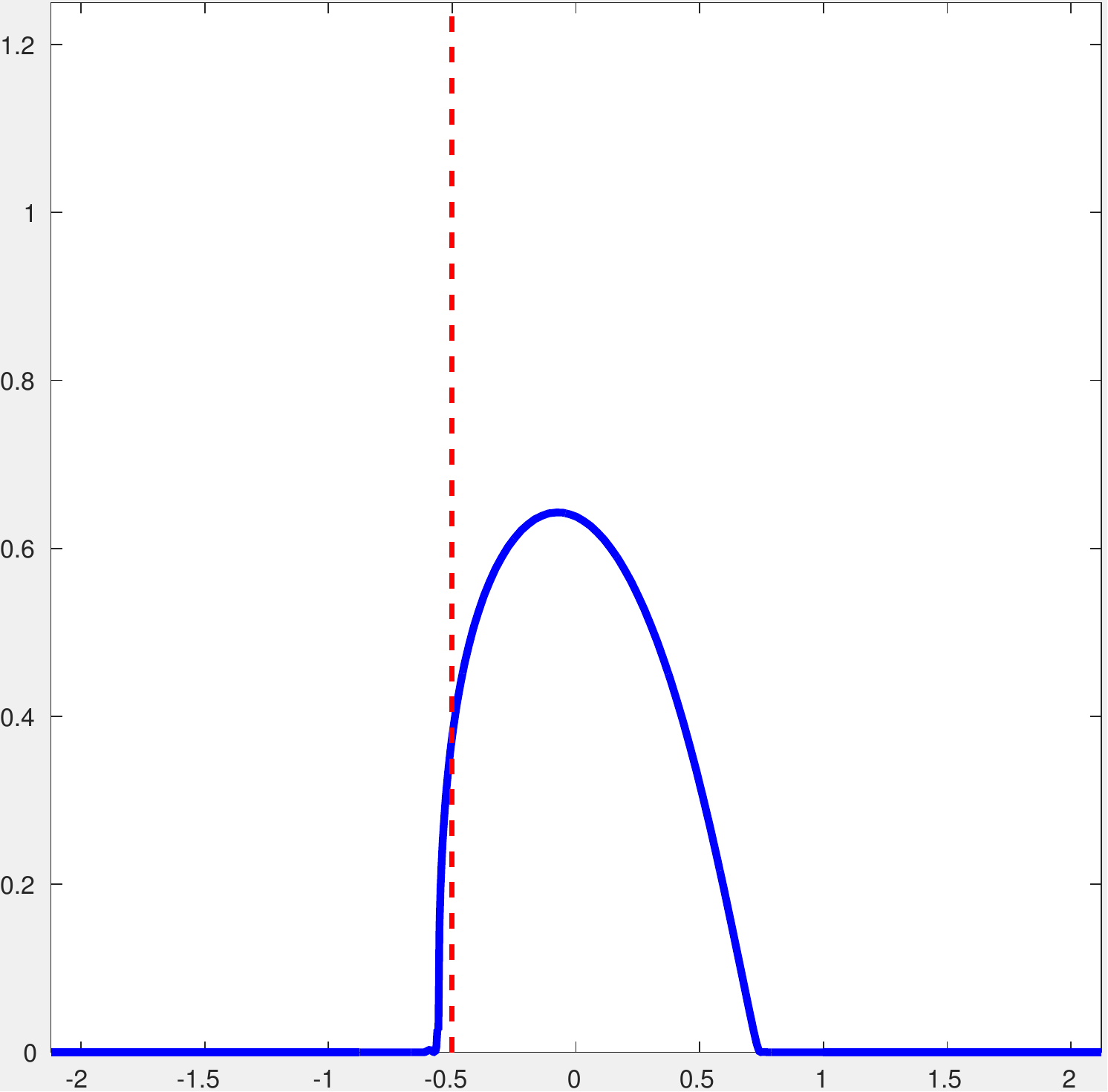}\caption{$t = 0.5$}\end{subfigure}%
        \begin{subfigure}[b]{0.25\linewidth}\includegraphics[scale=0.20]{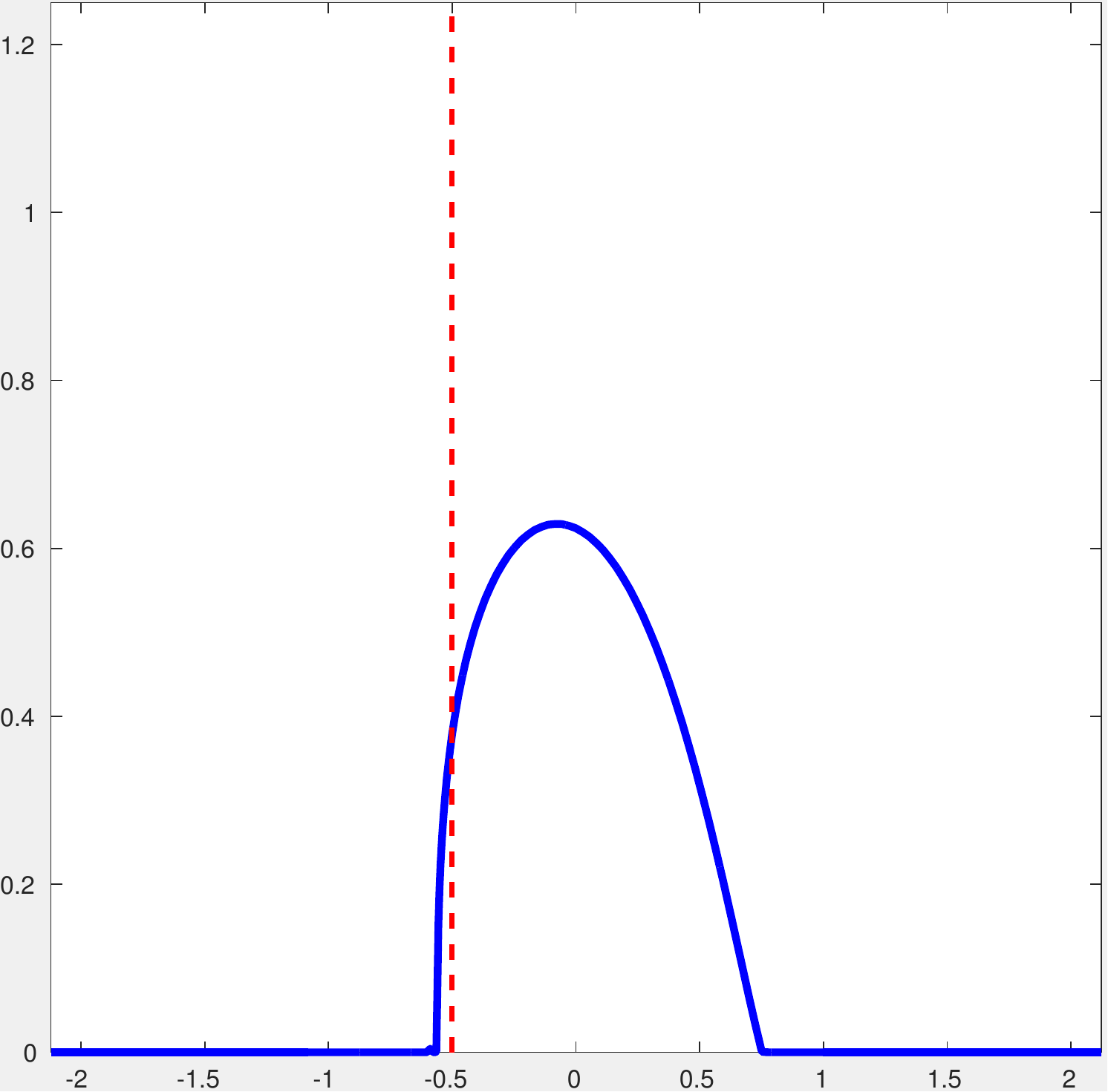}\caption{$t = 0.6$}\end{subfigure}%
        \begin{subfigure}[b]{0.25\linewidth}\includegraphics[scale=0.20]{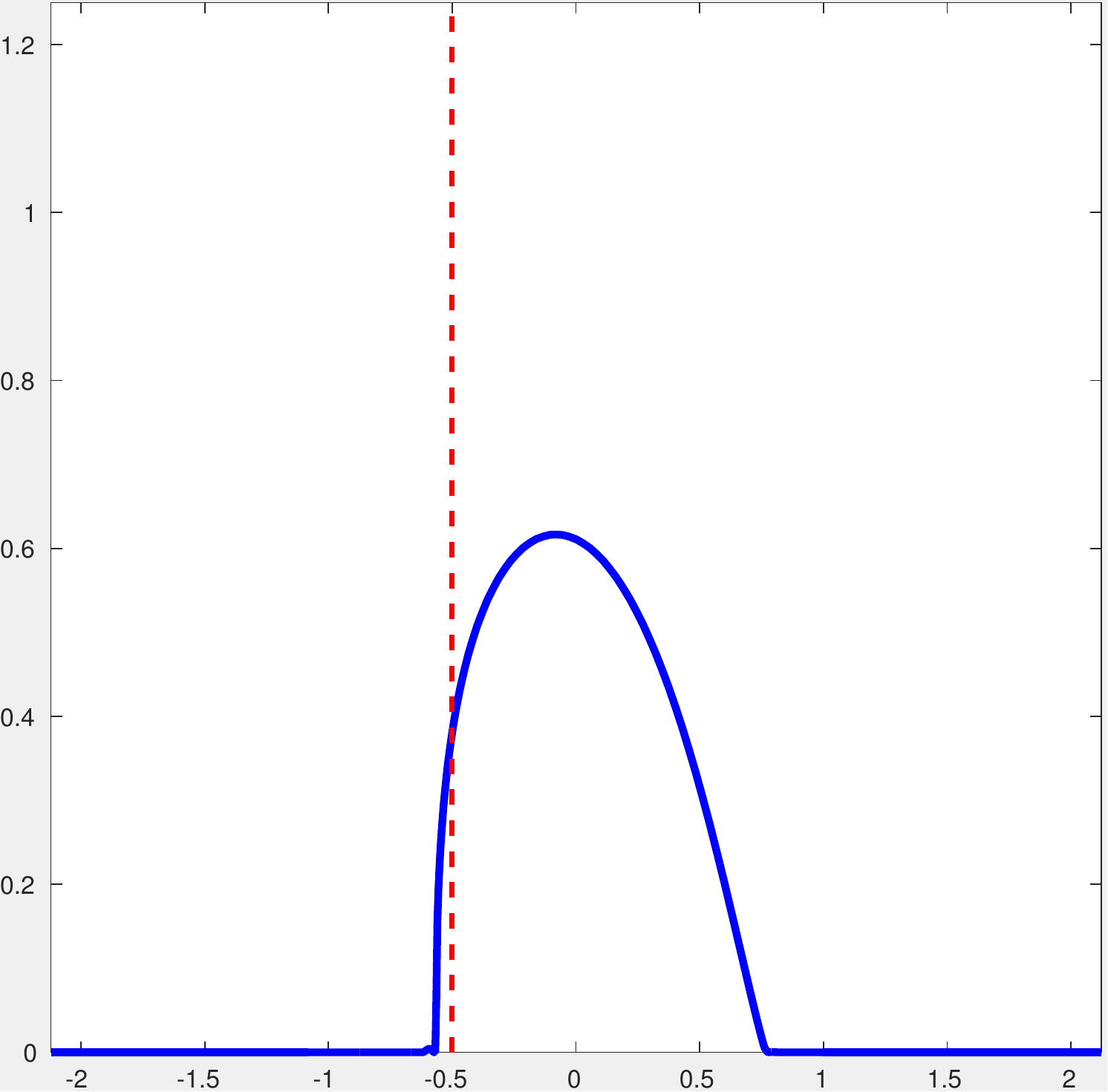}\caption{$t = 0.7$}\end{subfigure}%
        \begin{subfigure}[b]{0.25\linewidth}\includegraphics[scale=0.20]{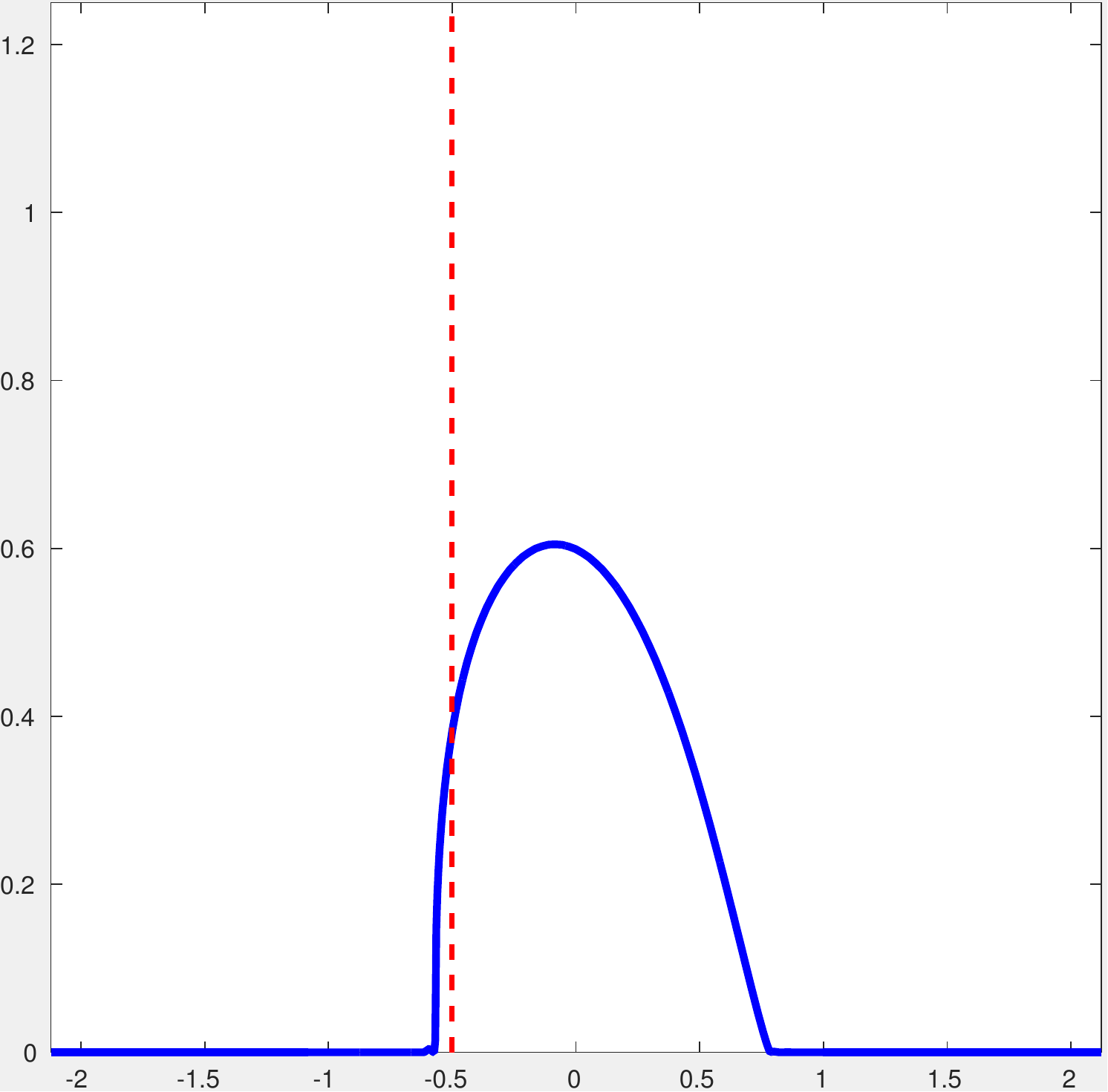}\caption{$t = 0.95$}\end{subfigure}\\%
        \caption{Example~\ref{exam5.3}. The cross section at $y = 0$ of a computed solution is shown at various time instants ($N = 40000$).}
        \label{fig:waiting-varExp2-2D-cross-section}
    \end{figure}

\begin{exam}[Waiting-time for variable exponent without absorption]
\label{exam5.3}
For this example, 
\begin{align*}
& \lambda = 0,\quad \gamma = 2-x-y,\quad \Omega = (-1.5,1.5)\times (-1.5,1.5),
\\
& u_0 = \begin{cases} 5 (0.25-x^2-y^2), & \text{for} \quad  \sqrt{x^2 + y^2} < 0.5 \\
	0,& \text{otherwise}.
\end{cases}
\end{align*}
This example has been studied in \cite{Duque2015}. 
We take $t \in [0, 0.05]$. The free boundary of the solution does not move  until $t \approx 0.02$.
A moving mesh and the corresponding computed solution are shown in Fig.~\ref{fig:VarExp2-example5.3}.
We can see that the variation of the exponent causes the free boundary to expand anisotropically and
the solution to have different steepness along the free boundary.
Moreover, a closer examination of the results confirms the waiting time phenomenon, where the interface
in the region $\{(x,y): x+y \leq 0 \}$ does not move until a finite time has elapsed. Fig. \ref{fig:waiting-varExp2-2D-cross-section} show the cross sections of the numerical solutions with the plane $y = x$ at various instants of time.
In the figure, the red dashed line refers to the position of the initial interface,
where the waiting time phenomenon subsequently occurs.
%Thus, our results agree with those in \cite{Duque2015}.
\qed \end{exam} 
    
\begin{figure}[ht]
        \centering
        \begin{subfigure}[b]{0.35\linewidth}\includegraphics[scale=0.34]{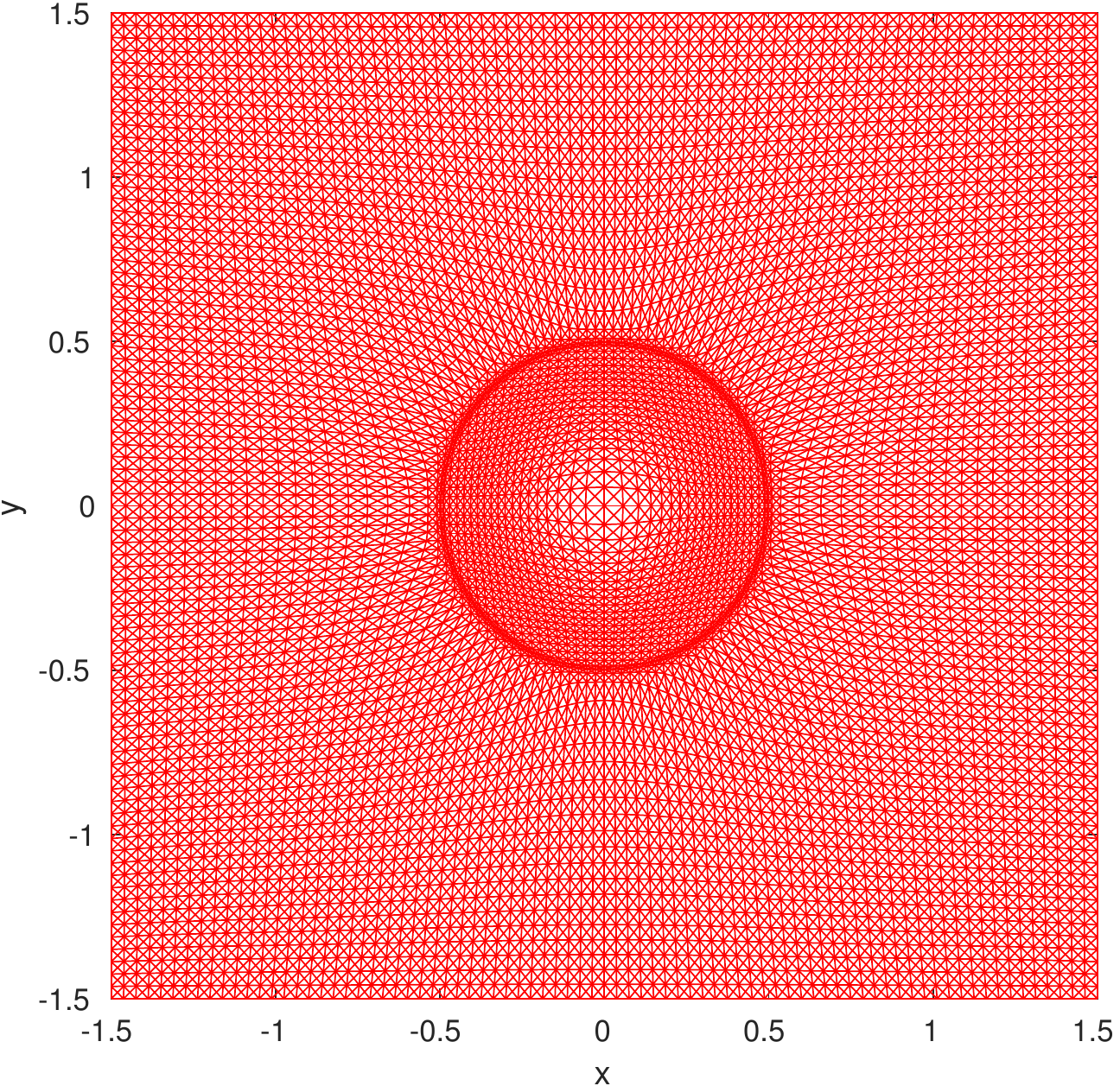}\caption{$t = 0.00$}\end{subfigure}\hspace{5mm}
        \begin{subfigure}[b]{0.35\linewidth}\includegraphics[scale=0.32]{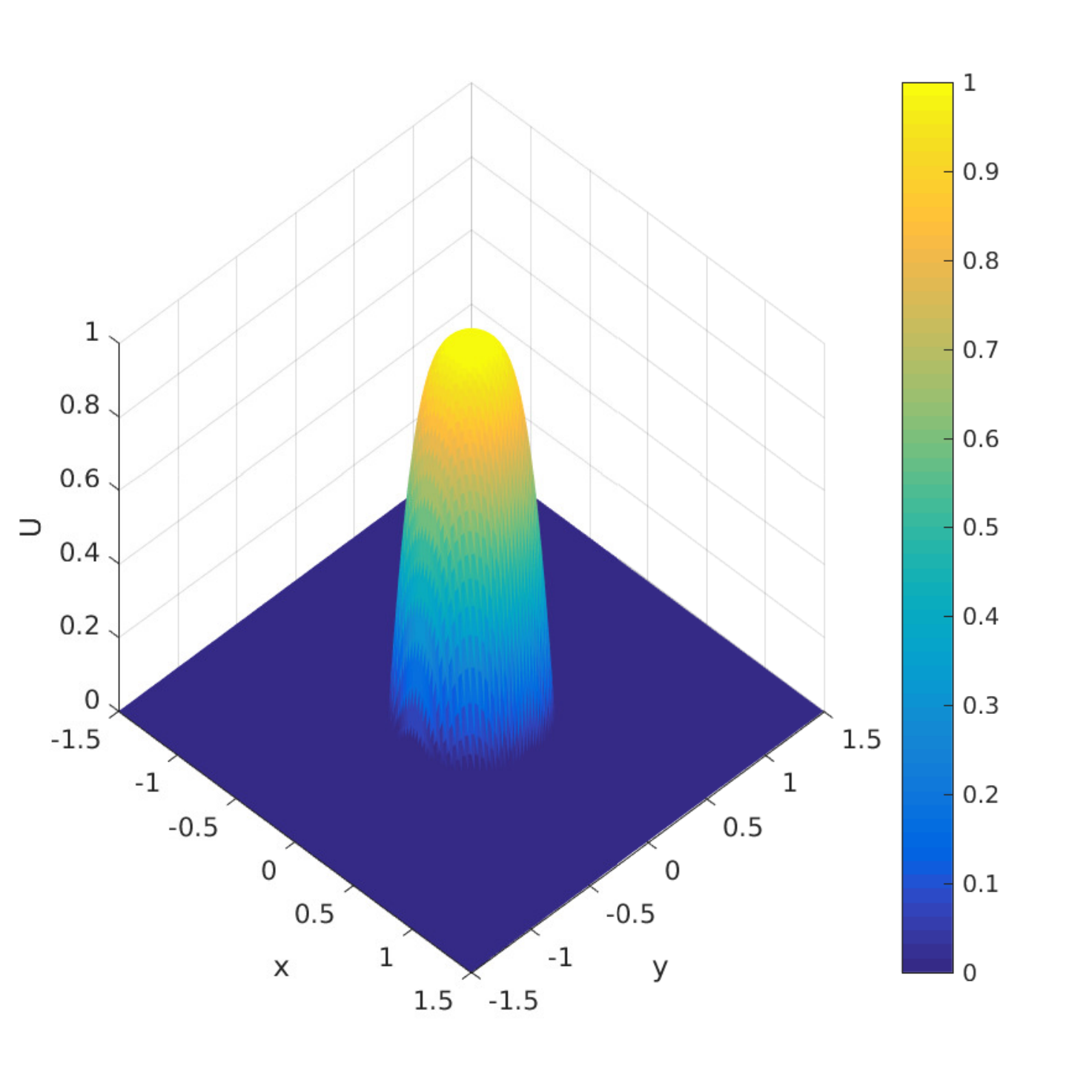}\caption{$t = 0$}\end{subfigure}\\
        \vspace{-1mm}
        \begin{subfigure}[b]{0.35\linewidth}\includegraphics[scale=0.34]{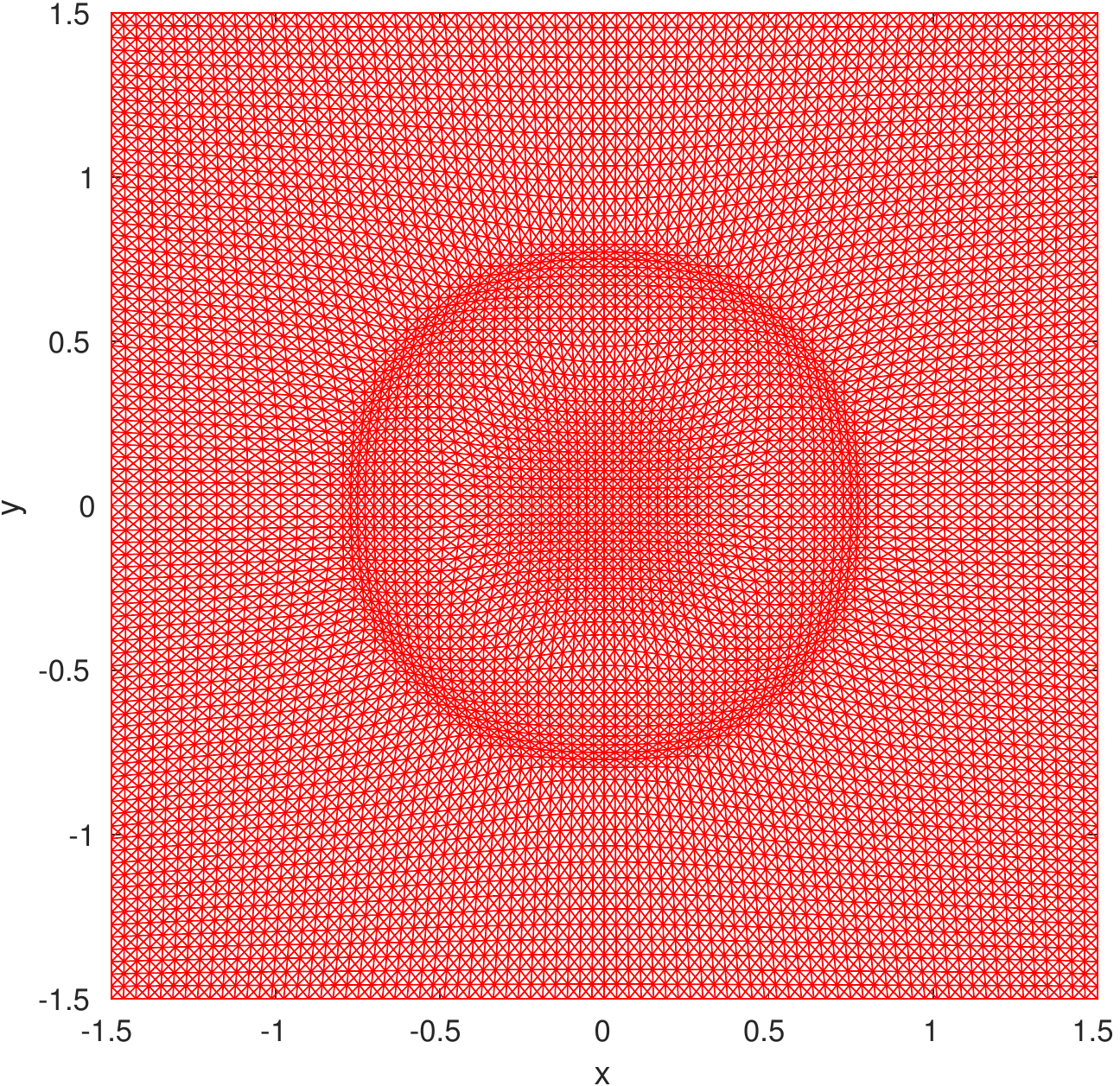}\caption{$t = 0.03$}\end{subfigure}\hspace{5mm}
        \begin{subfigure}[b]{0.35\linewidth}\includegraphics[scale=0.32]{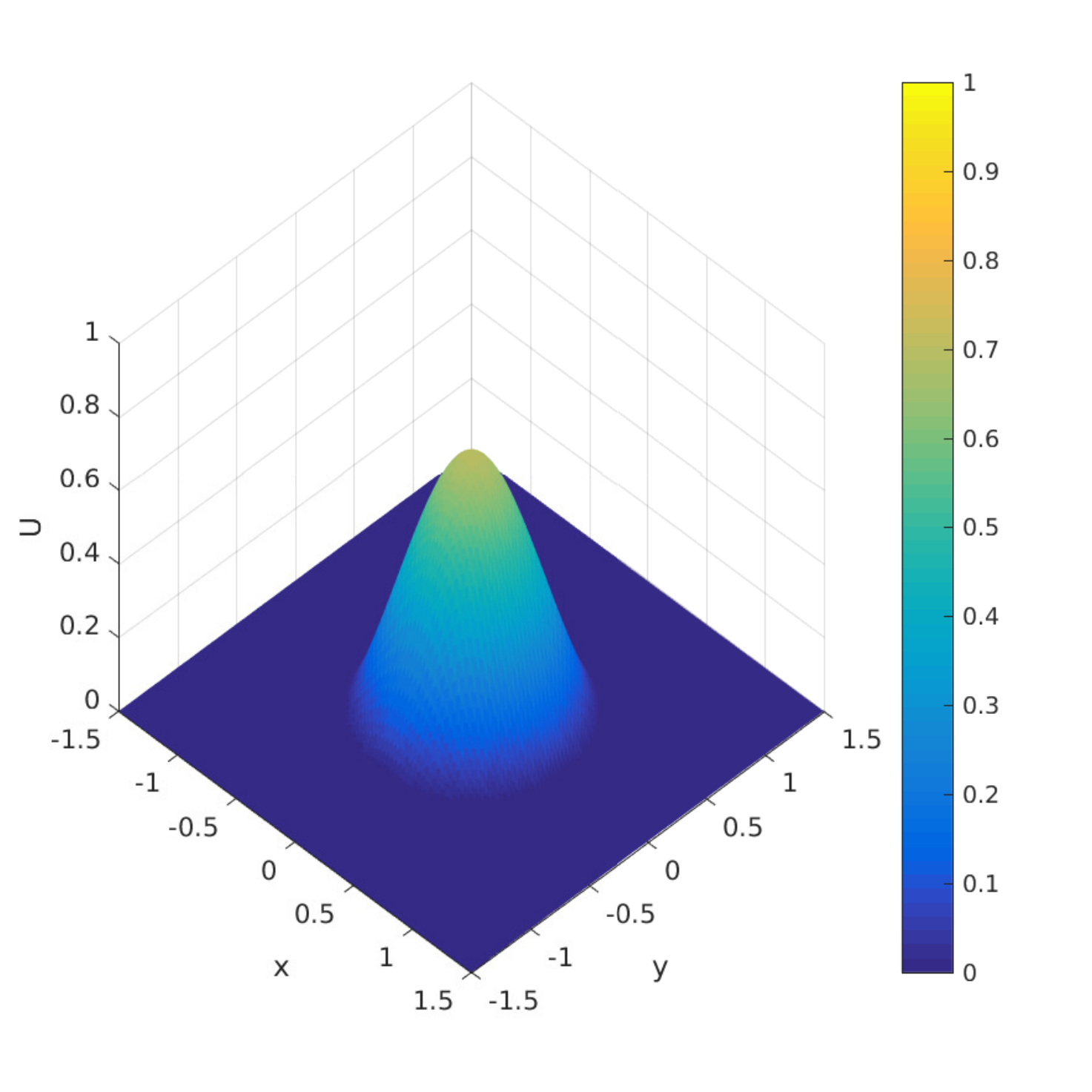}\caption{$t = 0.03$}\end{subfigure}\\%
        \vspace{-1mm}
        \begin{subfigure}[b]{0.35\linewidth}\includegraphics[scale=0.34]{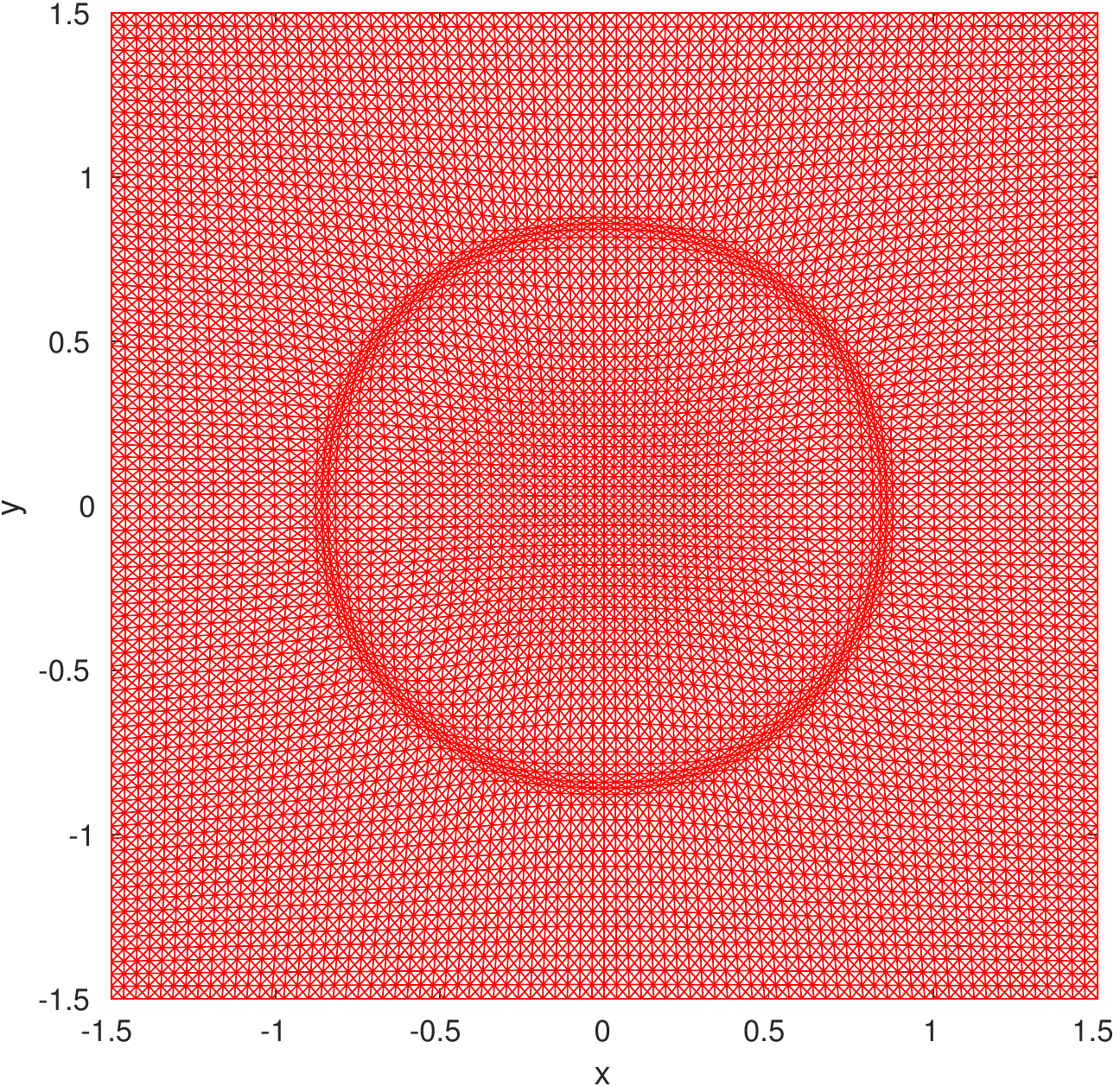}\caption{$t = 0.06$}\end{subfigure}\hspace{5mm}
        \begin{subfigure}[b]{0.35\linewidth}\includegraphics[scale=0.32]{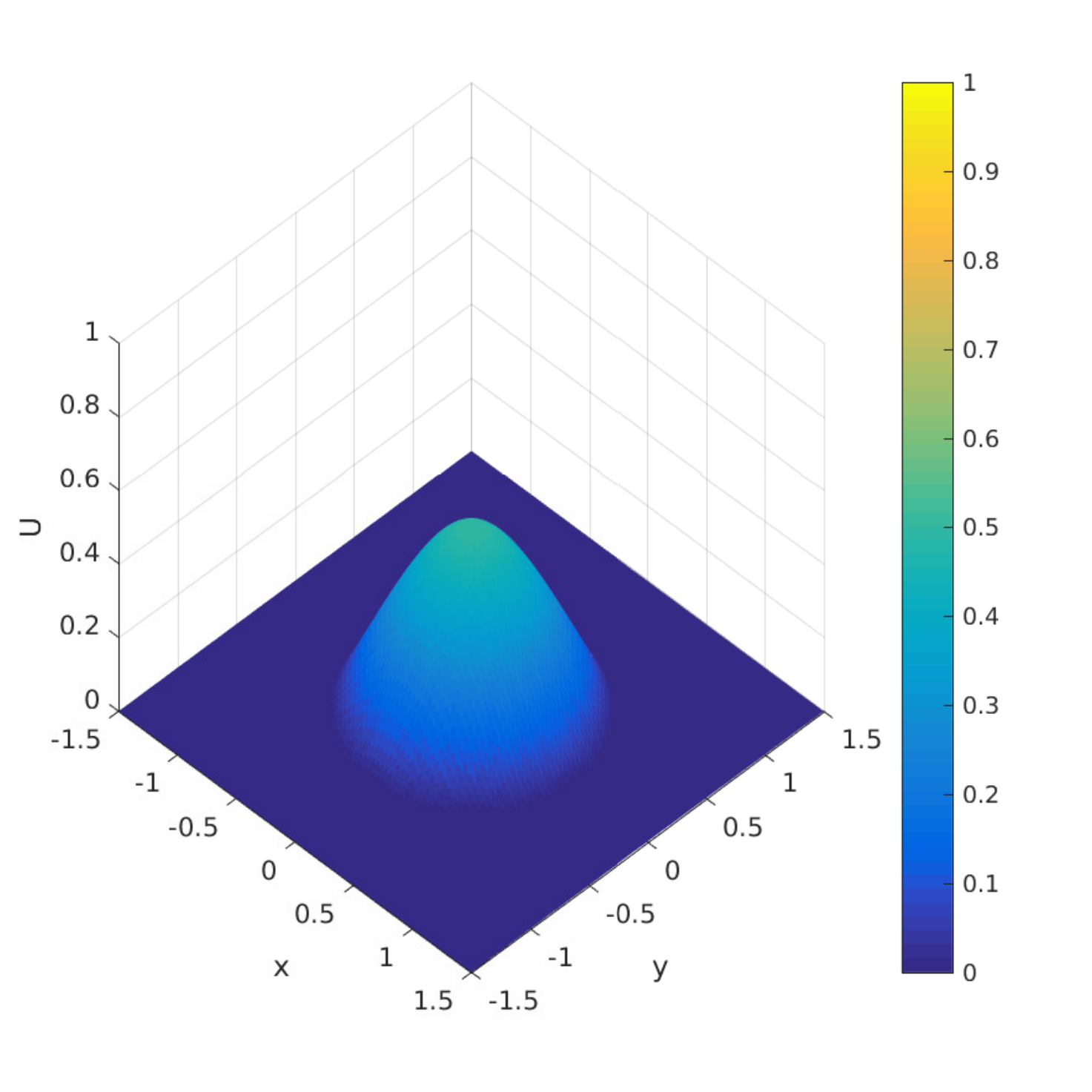}\caption{$t = 0.06$}\end{subfigure}\\
        \vspace{-1mm}
        \begin{subfigure}[b]{0.35\linewidth}\includegraphics[scale=0.34]{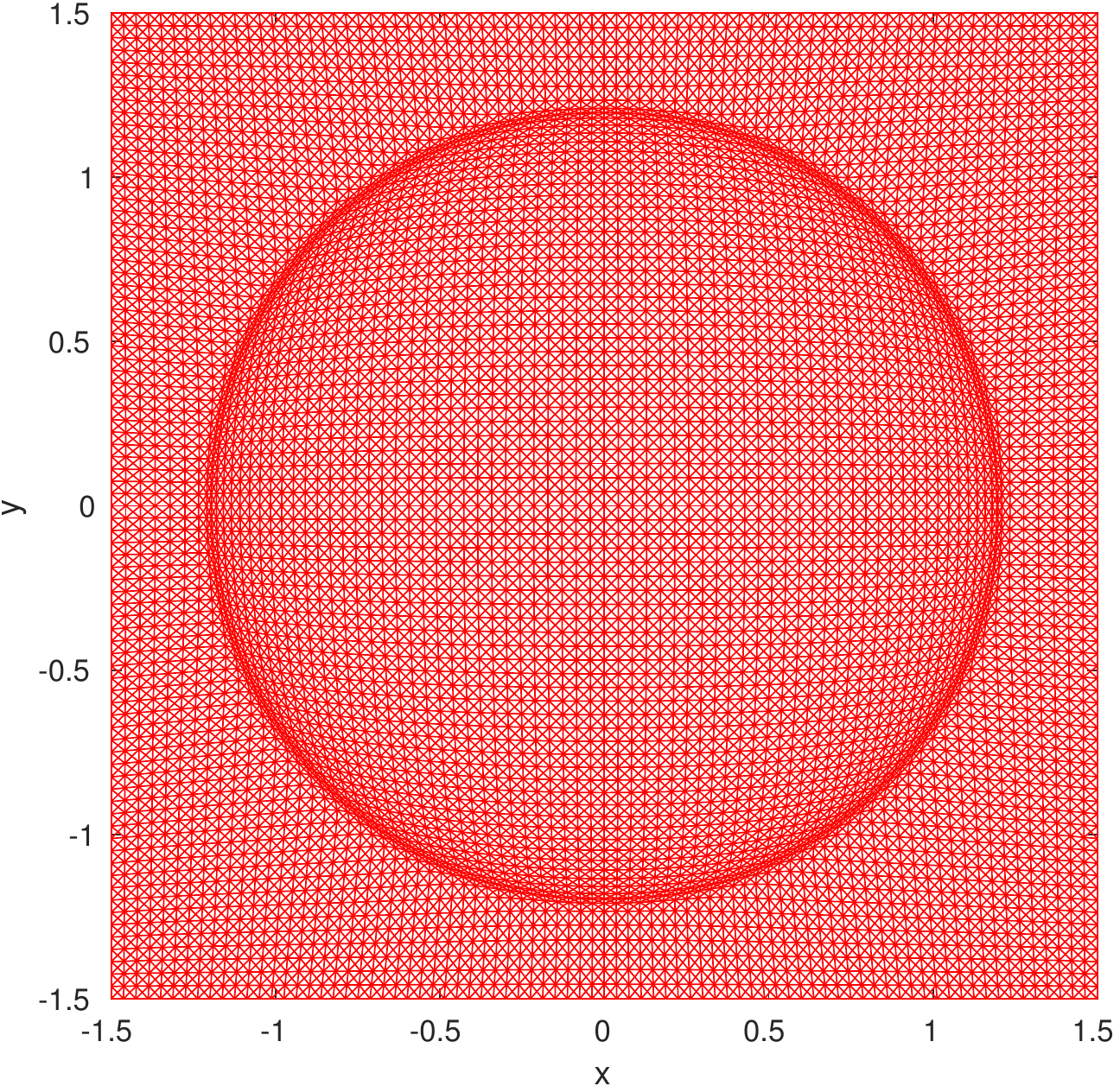}\caption{$t = 0.50$}\end{subfigure}\hspace{5mm}
        \begin{subfigure}[b]{0.35\linewidth}\includegraphics[scale=0.32]{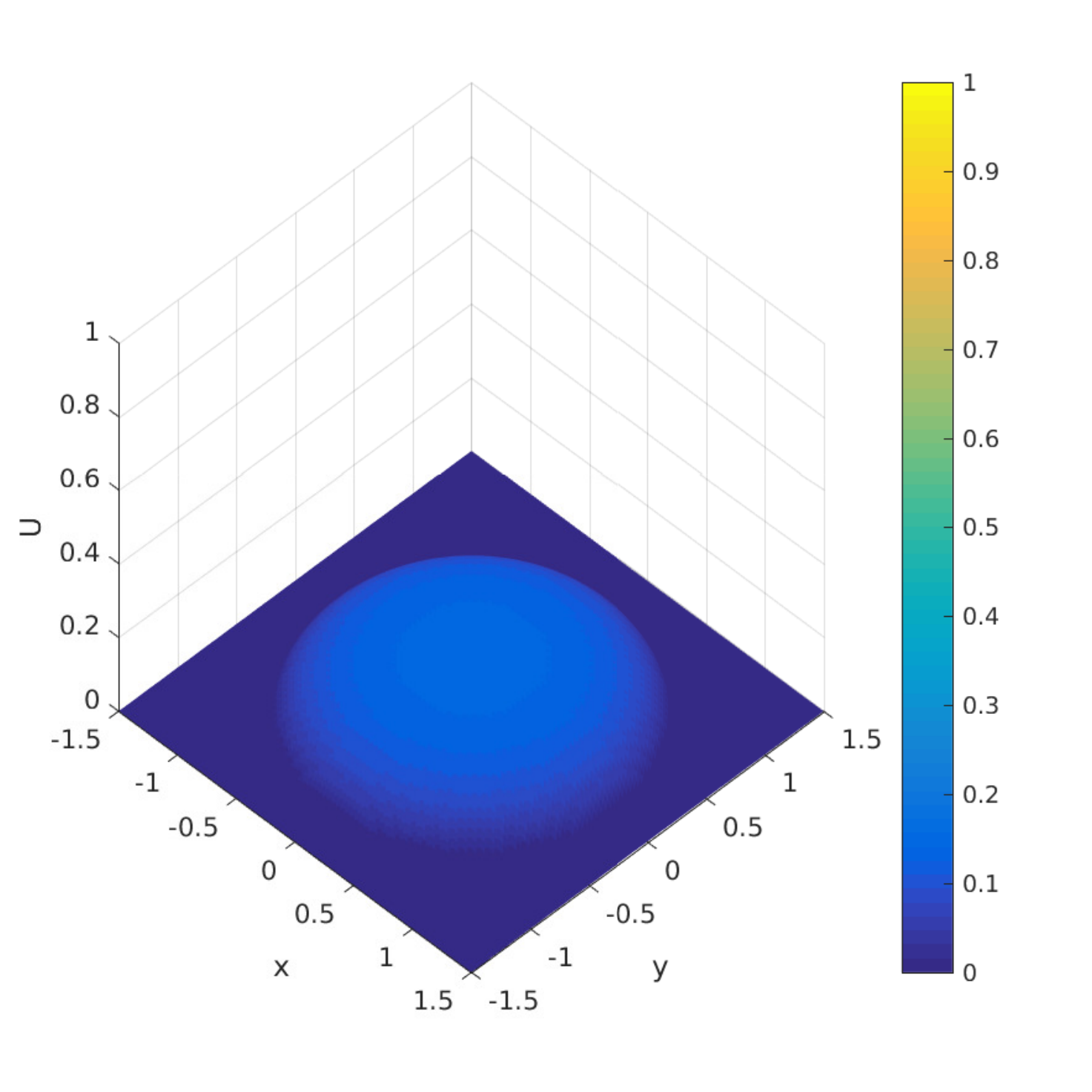}\caption{$t = 0.50$}\end{subfigure}%
        \caption{Example~\ref{exam5.4}. An adaptive mesh and the corresponding solution at various time instants ($N = 25600$).}
        \label{fig:VarExp3-example5.4-mesh}
    \end{figure}

\begin{exam}[Variable exponents with absorption]
\label{exam5.4}
The last example, taken from \cite{Duque2014}, has time dependent exponents, i.e., 
\begin{align*}
& \lambda = 1,\quad \gamma = \frac{x^2+y^2}{t^2+1},\quad \sigma = x^2+y^2 + 1 + e^{-t},
\quad \Omega = (-1.5,1.5)\times (-1.5,1.5),
\\
& u_0 = \begin{cases} \cos(2 \pi (x^2+y^2)), &  \text{for} \quad \sqrt{x^2 + y^2} < 0.5 \\
	0,& \text{otherwise}.
\end{cases}
\end{align*}
We take $t \in [0, 0.1]$. 
The numerical results are shown in Fig.~\ref{fig:VarExp3-example5.4-mesh}.
They are comparable with those in \cite{Duque2014}.
\qed \end{exam}  

%%%%%%%%%%%%%%%%%%%%%%%%%%%%%%%
%% Section: Conclusions
\section{Conclusions and further remarks}
\label{SEC:conclusion}

In the previous sections we have studied an adaptive moving mesh finite element method for the numerical
solution of PME. The method is based on
the MMPDE moving mesh strategy and its new implementation and uses a  linear finite element
method and the fifth-order Radau IIA scheme for the spatial and temporal discretization. 
Numerical results show that the method is able to produce correct mesh concentration around
the free boundary and deal with problems having complex solution support.
Three types of mesh have been considered, uniform and acrlength- and Hessian-based
adaptive meshes. The method shows a first-order convergence behavior as the mesh is refined
for uniform and arclength-based adaptive meshes and improves to a second-order convergence
when Hessian-based adaptive meshes are used. This indicates that mesh concentration around the free boundary
is important to the accuracy of the method. Moreover, the prompt response
of the mesh movement to the changes in the solution is also crucial, requiring that a small value
of the parameter $\tau$ in mesh movement (cf. (\ref{mmpde-1})) be used especially for the computation
with fine meshes.
We have also studied the application of the method to PME with variable exponents and absorption
for which there are very few theoretical results available. Numerical results demonstrate that the method
is robust and able to deal with PDEs having more complicated structures. 

%{\color{red}
It should be pointed out that there are small oscillations around the free boundary
in computed solutions; see the discussion in Sect.~\ref{SEC:PME-numerics}.
How to suppress these oscillations using a monotone or structure-preserving scheme
(e.g., see \cite{BaSo1991,LSSV07,LiHu2013,NH2015,Oberman2006,Le05,YuSh2008,Zhang-2009,ZZS2013})
and to combine them with adaptive mesh movement for PME are worth future investigations.
%}

\vspace{20pt}

%{\color{red}
{\bf Acknowledgment.} 
Support from US Army Research Office under grant W911-NF-1510377 is gratefully acknowledged.
The authors would also like to thank the anonymous referees for their valuable comments in improving the quality of the paper.
%}

%\bibliographystyle{amsplain}
%\bibliographystyle{abbrv}
%\bibliography{fdmref}

%\bibliographystyle{abbrv}
%\bibliography{/Users/huang/tex/bib/mmesh}

\end{document}